\documentclass[12pt]{report}
\usepackage{amsmath, amsthm, mathtools, amssymb, epsfig, sectsty, graphicx, float, url}
\usepackage{algorithm, algorithmicx, algpseudocode}
\usepackage[hypertexnames=false,hyperfootnotes=false]{hyperref}
\usepackage[hmargin=1in,vmargin=1in]{geometry}
\usepackage{float}
\usepackage{natbib}
\usepackage{eurosym}
\usepackage{wrapfig}
\usepackage{bm}

\usepackage[colorinlistoftodos,bordercolor=orange,backgroundcolor=orange!20,line
color=orange,textsize=scriptsize]{todonotes}
\catcode`\@=11
\catcode`\@=12

\usepackage{subcaption}

\numberwithin{equation}{chapter}%

\usepackage{booktabs}
\usepackage{tabularx}
\usepackage{array}
\usepackage{adjustbox}
\usepackage[noabbrev, capitalise, nameinlink]{cleveref}
\usepackage{multirow}

\usepackage[nottoc,notlof,notlot]{tocbibind}

\theoremstyle{plain}
\newtheorem{theorem}{Theorem}[chapter]
\newtheorem{proposition}{Proposition}[chapter]
\newtheorem{example}{Example}[chapter]
\newtheorem{lemma}{Lemma}[chapter]
\newtheorem{corollary}{Corollary}[chapter]
\theoremstyle{definition}
\newtheorem{definition}{Definition}[chapter]
\newtheorem{remark}{Remark}[chapter]

\numberwithin{equation}{chapter}

\newcommand{\continuanceref}{}

\theoremstyle{plain}
\newtheorem{continuancex}{Example}

\newenvironment{continuance}[1]
  {\renewcommand{\continuanceref}{\ref{#1}}\begin{continuancex}}
  {\end{continuancex}}

\newcommand{\dima}{{\sf dim aff}}

\newcommand{\conv}{{\sf conv}}

\newcommand{\ext}{{\sf ext}}
\newcommand{\opt}{^{\star}}
\newcommand{\tr}{^{\top}}
\newcommand{\kl}{{\sf KL}}

\newcommand{\bP}{\mathbb{P}}

\newcommand{\R}{\mathbb{R}}
\newcommand{\N}{\mathbb{N}}
\newcommand{\E}{\mathbb{E}}

\newcommand{\cS}{\mathcal{S}}
\newcommand{\cR}{\mathcal{R}}
\newcommand{\cA}{\mathcal{A}}

\newcommand{\cD}{\mathcal{D}}
\newcommand{\cI}{\mathcal{I}}
\newcommand{\cF}{\mathcal{F}}
\newcommand{\cJ}{\mathcal{J}}

\newcommand{\cH}{\mathcal{H}}
\newcommand{\cP}{\mathcal{P}}
\newcommand{\cL}{\mathcal{L}}

\newcommand{\cE}{\mathcal{E}}
\newcommand{\cU}{\mathcal{U}}
\newcommand{\cUf}{\mathcal{U}^{\sf f}}
\newcommand{\Tf}{T_{\sf f}}
\newcommand{\cW}{\mathcal{W}}

\newcommand{\cM}{\mathcal{M}}
\newcommand{\cX}{\mathcal{X}}
\newcommand{\cY}{\mathcal{Y}}
\newcommand{\cV}{\mathcal{V}}
\newcommand{\cVb}{\bar{\mathcal{V}}}
\newcommand{\cQ}{\mathcal{Q}}
\newcommand{\cG}{\mathcal{G}}

\newcommand{\PiS}{\Pi_{\sf S}}
\newcommand{\PiH}{\Pi_{\sf H}}
\newcommand{\PiHra}{\Pi_{\sf H}^{\sf ra}}
\newcommand{\PiSra}{\Pi_{\sf S}^{\sf ra}}
\newcommand{\PiM}{\Pi_{\sf M}}
\newcommand{\PiMD}{\Pi_{\sf MD}}
\newcommand{\PiSD}{\Pi_{\sf SD}}
\newcommand{\PiHD}{\Pi_{\sf HD}}

\newcommand{\hP}{\hat{P}}

\newcommand{\DSSA}{\Delta(\cS)^{\cS \times \cA}}
\newcommand{\DSA}{\Delta(\cS)^{\cA}}
\newcommand{\DS}{\Delta(\cS)}
\newcommand{\DA}{\Delta(\cA)}
\newcommand{\Dr}{\Delta([r])}

\newcommand{\RSAS}{\R^{\cS \times \cA \times \cS}}
\newcommand{\RS}{\R^{\cS}}
\newcommand{\RSA}{\R^{\cS \times \cA}}
\newcommand{\RAS}{\R^{\cA \times \cS}}

\newcommand{\cPSA}{\cP^{\sf SA}}
\newcommand{\cPS}{\cP^{\sf S}}
\newcommand{\cPG}{\cP^{\sf G}}
\newcommand{\cPC}{\cP^{\sf c}}

\newcommand{\cVS}{\cV_{\sf S}}
\newcommand{\cVNS}{\cV_{\sf NS}}

\newcommand{\tT}{\tilde{T}}
\newcommand{\tv}{\tilde{v}}
\newcommand{\bT}{\bar{T}}
\newcommand{\btheta}{\bar{\theta}}
\newcommand{\bv}{\bar{v}}
\newcommand{\bnu}{\bar{\nu}}
\newcommand{\bs}{\bar{s}}

\newcommand{\recs}{{\sf Rec}_{s}}
\newcommand{\recsa}{{\sf Rec}_{sa}}

\newcommand{\Ravg}{R_{\sf avg}}

\newcommand{\essinf}{{\sf essinf}}
\newcommand{\cvar}{{\sf CVaR}}
\newcommand{\evar}{{\sf EVaR}}
\newcommand{\erm}{{\sf ERM}}

\newcommand{\tX}{\tilde{X}}
\newcommand{\tR}{\tilde{R}}
\newcommand{\tY}{\tilde{Y}}
\newcommand{\ts}{\tilde{s}}
\newcommand{\ta}{\tilde{a}}

\title{Robust Markov Decision Processes: An Introduction using Convex Optimization and Game Theory}

\author{Julien Grand-Cl\'ement \\
HEC Paris \\ 
grand-clement@hec.fr\\}
\date{\today}

\begin{document}
\maketitle

\phantomsection
\addcontentsline{toc}{chapter}{Abstract}
\noindent\textbf{Abstract.}
The goal of these notes is to provide a unified and rigorous introduction to robust sequential decision-making. While there exist many introductory notes on reinforcement learning and classical Markov decision processes (MDPs), these notes focus specifically on {\em robust} MDPs, a field that has become an active area of research since the mid-2000s. We provide a self-contained introduction to the main results in the field, which are spread across several seminal papers. The notes are organized as follows:
\begin{itemize}
    \item \Cref{chap: lec 0 mdp}, MDPs: this chapter gives a concise introduction to nominal (classical) Markov Decision Processes and the operator approach to computing optimal policies.
    \item \Cref{chap: lec 1 robust MDP}, general model for robust MDPs: this chapter describes the robust MDP problem in the general setting for the case of uncertain transition probabilities and highlights the main departure from the tractable properties of nominal MDPs.
    \item \Cref{chap: lec 2 rectangular models}, rectangular models for robust MDPs: this chapter introduces the s-rectangularity and sa-rectangularity assumptions and shows how these assumptions enable the efficient computation of optimal policies by adapting classical fixed point algorithms.
    \item \Cref{chap: lec 3 connection}, connections with other models: this chapter shows the connections between robust MDPs and several other popular models for sequential decision-making: MDPs with uncertain rewards, Multi-Objective MDPs and Pareto optimality, MDPs with time-varying uncertain parameters, non-rectangular robust MDPs,  stochastic games, distributionally robust MDPs, and risk-averse MDPs.
\end{itemize}
The reader new to sequential decision-making should start from \Cref{chap: lec 0 mdp}, which provides a short introduction to MDPs and builds intuition useful in the rest of these notes.
The reader already versed in MDPs can refer to \Cref{chap: lec 0 mdp} only for notation and start from \Cref{chap: lec 1 robust MDP}. The reader already versed in robust MDPs will find interesting extensions at the end of \Cref{chap: lec 2 rectangular models} and in \Cref{chap: lec 3 connection}.
More advanced material is signaled with an asterisk (*), and exercises are given for the first three chapters. The vast majority of the results presented in these notes are already published in various venues, and the references are listed at the end of each chapter (where we also mention which results of these notes are new; note that we by no means attempt to provide an exhaustive literature review on MDPs, there already exist surveys for this purpose). Having all these results in one place, within a single unifying formalism, will hopefully be useful for researchers who want to start working on this topic. These notes were written during the summer 2026 and I hope to be able to add new chapters on some of the current research directions in robust MDPs in the near future.

\vspace{2mm}
\noindent
These notes are meant to be readable by a broad audience consisting of graduate students and faculty in optimization, operations research, game theory, computer science, and adjacent fields. A basic background in linear programming, convex analysis, and probability is required, although we usually reintroduce the notions where we need them. The appendices provide the statements and explanations of advanced technical results required for our purpose. These notes are entirely designed and written by the author; we acknowledge the use of artificial intelligence (more specifically, OpenAI's ChatGPT) for (a) finding the NP-hardness proofs for \Cref{prop: lec 3 hardness uncertain r p} and \Cref{prop: lec 3 hardness r rec}, and (b) reviewing the document for typos/inconsistencies and enhancing the rendering of tables.

\vspace{2mm}
\noindent
In order to keep technicalities to a minimum, we focus on ``the simplest setup where interesting things happen": sequential decision problems with finitely many states and actions, and most chapters deal with discounted return. In general we strive to give results in their most general forms, i.e., without unnecessary assumptions, and we reuse the same small number of examples/counterexamples throughout the notes. The exposition combines operator and fixed point approaches, convex optimization (more specifically, duality in saddle-point problems), and game-theoretic arguments. The latter identifies robust MDPs with a certain type of two-player game; this helps build useful intuition, especially when constructing counterexamples. My understanding of MDPs, robust MDPs, and games in general has greatly benefited from stimulating research discussions with many of my colleagues and coauthors, and particularly with Marek Petrik, St{\'e}phane Gaubert, and Nicolas Vieille.
\tableofcontents

\chapter{Introduction to Markov Decision Processes}\label{chap: lec 0 mdp}
This chapter introduces the model of Markov decision processes for sequential decision-making. We first define the notions of policies, the most classical objective function (discounted return) for infinite-horizon decision-making, and the value functions and state-action frequencies. We then provide the main properties of optimal policies, the optimality equations for the value functions of optimal policies and classical algorithms (value iteration, policy iteration, linear programming) that return uniformly optimal policies.
\section{Model, assumptions, and objective}\label{sec: lec 0 model}
\begin{figure}
    \centering
\includegraphics[width=0.5\linewidth]{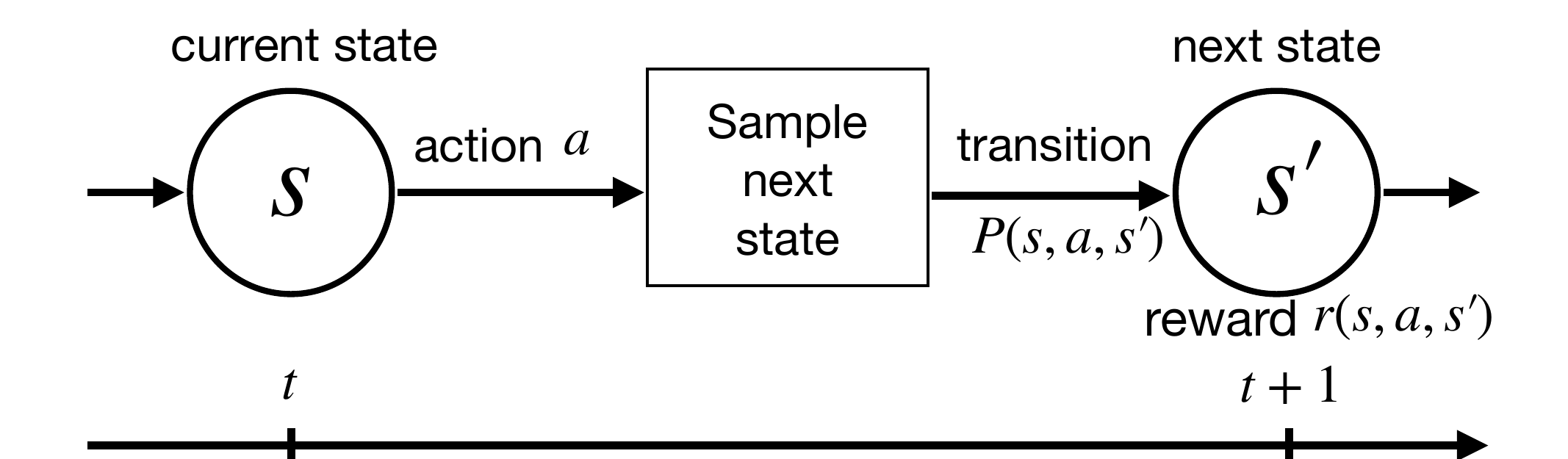}
    \caption{Dynamics for Markov decision processes. The next state is sampled with the probability distribution $P(s,a) \in \Delta(\cS)$. We choose action $a$ in state $s$ to balance the instantaneous reward $r(s,a,s')$ and the likelihood of obtaining good rewards in the future (from the next state $s'$).}
    \label{fig: lec 0 dynamics}
\end{figure}

Markov decision processes (MDPs) are used to model sequential decision-making problems.
An MDP $\cM$ is a tuple $\cM = (\cS,\cA,r,P,p_0)$ where:
\begin{itemize}
    \item $\cS$ and $\cA$ are the finite sets of states and actions,
    \item $r \in \R^{\cS \times \cA \times \cS}$ models the instantaneous rewards attached to each $(s,a,s')$ with $(s,a) \in \cS \times \cA$ the current state-action pair and $s'$ the next state,
    \item $P \in \Delta(\cS)^{\cS \times \cA}$, where $P(s,a,s')$ models the transition probabilities toward the next state $s'$ given a current state-action pair $(s,a)$,
    \item $p_0 \in \Delta(\cS)$ is the initial distribution over the set of states: $p_0(s)$ is the probability that we start in state $s \in \cS$.
    \end{itemize}
Recall that the simplex $\Delta(\cE)$ over a finite set $\cE$ is the set of all probability distributions over $\cE$ and can be identified with $\displaystyle \{x \in \R^{\cE} \; | \; \sum_{e \in \cE} x(e)=1, x(e) \geq 0, \forall \; e \in \cE\}$.
We usually represent the dynamics as in Figure~\ref{fig: lec 0 dynamics}. It is worth emphasizing the assumptions listed above:
\begin{enumerate}
    \item In these notes we focus on finite sets of states and actions, and discrete time periods. Already in this setting many interesting phenomena occur. We will later need compact action sets (see \Cref{sec: lec 2 adv}).
    \item The dynamics is assumed to be Markovian, in that a current state-action pair $(s_t,a_t)$ at period $t$ completely suffices to specify the distribution of the next state $s_{t+1}$ as a random variable with values in $\cS$. 
    \item For now, we assume that the same set of actions is available at each state. This is without loss of generality, and we do so to lighten notation for now, but we will need MDPs with state-dependent action sets when we study robust MDPs in \Cref{sec: lec 2 adv}.
\end{enumerate}
We will tend to emphasize the role of the transition probabilities since this will play a key role for robust MDPs.
\paragraph{Policies.} We now explain how we make decisions in an MDP. A policy $\pi$ should choose an action at each period $t \in \N$, based on all the information observed so far. At period $t \in \N$, we can observe $h_t = (s_0,a_0,s_1,a_1, ..., s_{t-1},a_{t-1},s_{t}) \in (\cS \times \cA)^{t} \times \cS$, i.e., we can observe the sequence of all states visited and actions chosen, we know the current state, and we need to decide what action to pick next. Thus in all generality a policy is a map from the set of all finite histories $\cH = \cup_{t \in \N} (\cS \times \cA)^{t} \times \cS$ to the simplex $\Delta(\cA)$ over the action set $\cA$.

We use the following notation for the sets of policies of interest:
\begin{itemize}
    \item $\PiH$ is the set of all policies, possibly history-dependent: $\PiH = \{\pi: \cH \to \Delta(\cA)\}$.
    \item $\PiM$ is the set of Markovian policies, i.e., policies where the action chosen at period $t$ may only depend on the current state $s_t \in \cS$ and period $t$, but not on the state-action pairs $(s_{0},a_{0},...,s_{t-1},a_{t-1})$ from the past.
    \item $\PiS$ is the set of stationary policies, i.e., policies where the action chosen at period $t$ may only depend on the current state $s_t \in \cS$.
    \item $\PiSD$ is the set of stationary deterministic policies, i.e., policies where the action chosen at period $t$ may only depend on the current state $s_t \in \cS$ and the same action is always chosen at a given state.
\end{itemize}
Of course we have $\PiSD \subset \PiS \subset \PiM \subset \PiH$.
Note that we can naturally identify the set $\PiM$ of Markovian policies with $\Delta(\cA)^{\cS \times \N}$, i.e., with maps returning a randomized action in $\Delta(\cA)$ for each pair $(s,t)$ of state-period in $\cS \times \N$. We can identify the set $\PiS$ of stationary policies with the compact set $\Delta(\cA)^{\cS}$ (returning a randomized action for each state) and the set $\PiSD$ of stationary deterministic policies with the finite set $\cA^{\cS}$.

Note that history-dependent policies require (in the worst case) keeping track of everything that happened in the past and may be difficult to store and deploy. Similarly, randomized stationary policies mean that every time we visit the same state $s$, we need to sample from $\pi(s) \in \Delta(\cA)$. This may be difficult to explain in some applications (e.g., in a healthcare application where the set of states is the set of health conditions for our patients and the set of actions is the set of feasible treatments), and in fact, randomization means delegating part of the decision to a generator of random numbers, which again may be problematic (how do you explain this to practitioners?). For all purposes, the set of stationary deterministic policies is the simplest, and its policies are the easiest to implement. 
\begin{example}\label{ex: lec 0 ex 1}
Consider the MDP instance from Figure~\ref{fig: lec 0 ex 1}. The set of states is $\cS = \{s_0, G,B\}$ and there are two actions: $\cA = \{a_1,a_2\}$. The rewards are $0$ starting from $s_0$ or $B$, and the reward from $G$ is always $1$. The transition probabilities are parametrized by some $p \in [0,1]$. Note that the only difference between $a_1$ and $a_2$ lies in the transitions out of $s_0$. For instance, if $p=0$, choosing action $a_1$ in $s_0$ leads to $B$ deterministically and choosing $a_2$ in $s_0$ leads to $G$ deterministically. If $p=1$ it is the reverse. The transition probabilities are given by $P(s_0,a_1) = (0,p,1-p),P(s_0,a_2) = (0,1-p,p), P(G,a_1) = P(G,a_2)= P(B,a_1) = P(B,a_2) = (1,0,0)$. The ``good" state $G$ yields a reward of $1$, while the ``bad" state yields a reward of $0$. Obviously we want to pick our actions so that we pass through state $G$ more often than state $B$. 
\begin{figure}[h!]
    \begin{center}
    \begin{subfigure}{0.45\textwidth}
    \begin{center}
        \includegraphics[width=\linewidth]{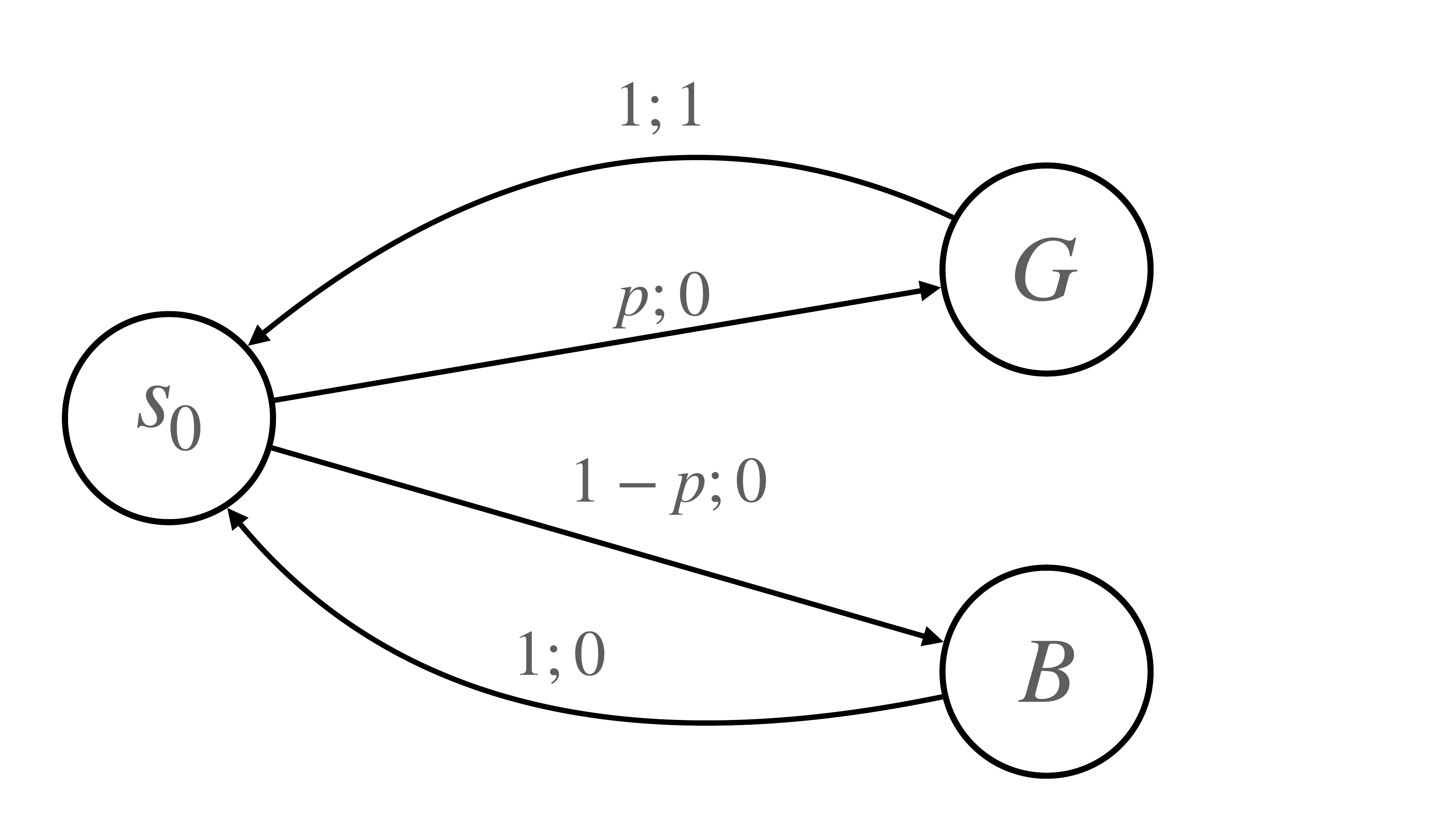}
    \end{center}
    \caption{Transitions and rewards for action $a_1$}
    \label{fig: lec 0 ex 1 a1}
    \end{subfigure}
    \begin{subfigure}{0.45\textwidth}
    \begin{center}
        \includegraphics[width=\linewidth]{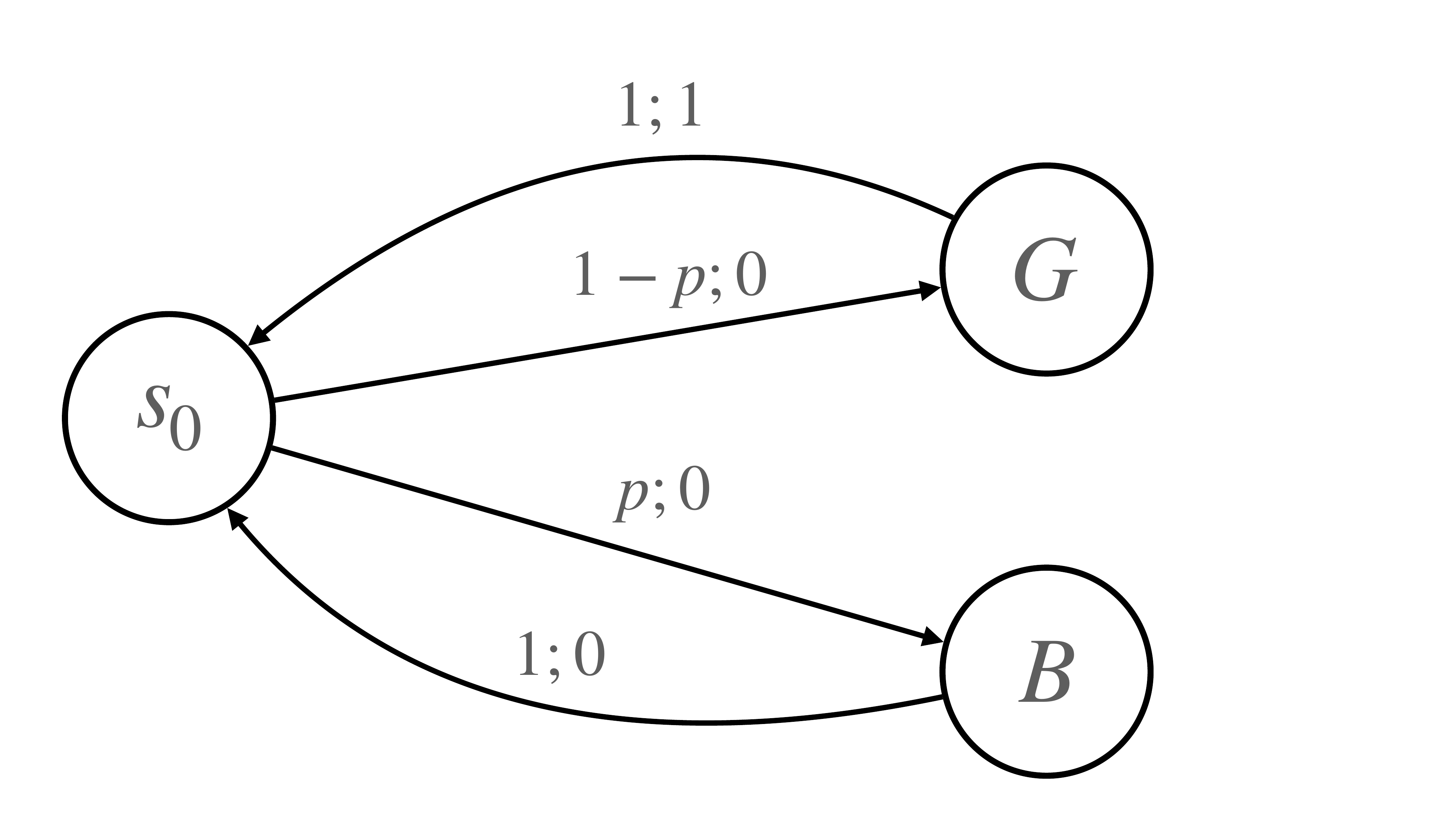}
    \end{center}
    \caption{Transitions and rewards for action $a_2$}
    \label{fig: lec 0 ex 1 a2}
    \end{subfigure}
    \end{center}
        \caption{MDP instance for \Cref{ex: lec 0 ex 1}. The transition probabilities and rewards are indicated beside each arc.}
        \label{fig: lec 0 ex 1}
    \end{figure}
\end{example}

\begin{example}\label{ex: lec 0 ex 2}
Consider the MDP instance from Figure~\ref{fig: lec 0 ex 2}, where there are five states $\{s_0,s_1,s_2,s_3,s_4\}$ and two actions $a_1$ and $a_2$. The transition probabilities are parametrized by $p,p' \in [0,1]$. Absorbing states are represented with squares. The rewards are only obtained in the absorbing states $s_2,s_3,s_4$. Actions $a_1$ and $a_2$ swap the probabilities towards the next states from $s_0$ and $s_1$. 
\begin{figure}[h!]
    \begin{center}
    \begin{subfigure}{0.45\textwidth}
    \begin{center}
        \includegraphics[width=\linewidth]{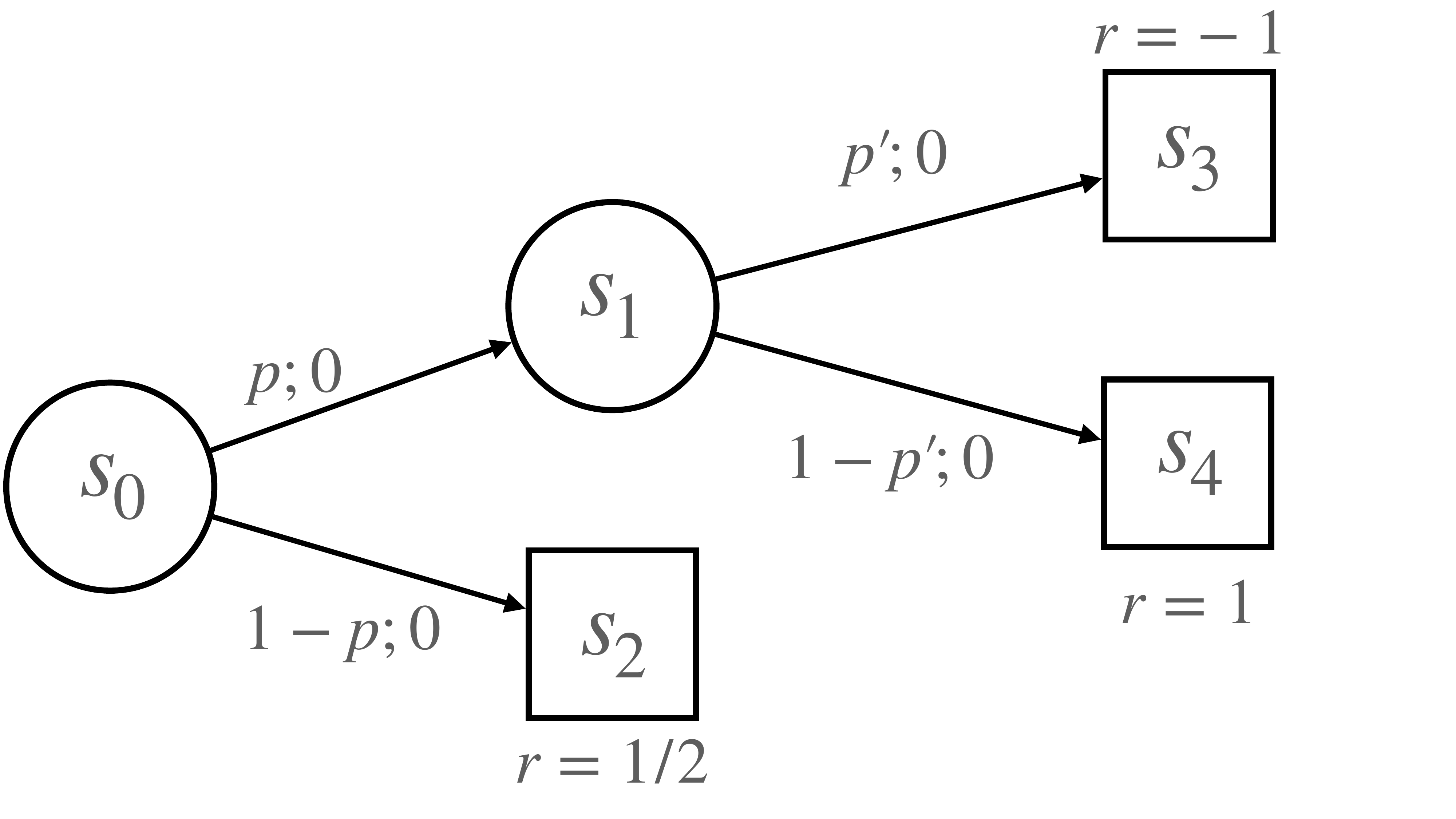}
    \end{center}
    \caption{Transitions and rewards for action $a_1$}
    \label{fig: lec 0 ex 2 a1}
    \end{subfigure}
    \begin{subfigure}{0.45\textwidth}
    \begin{center}
        \includegraphics[width=\linewidth]{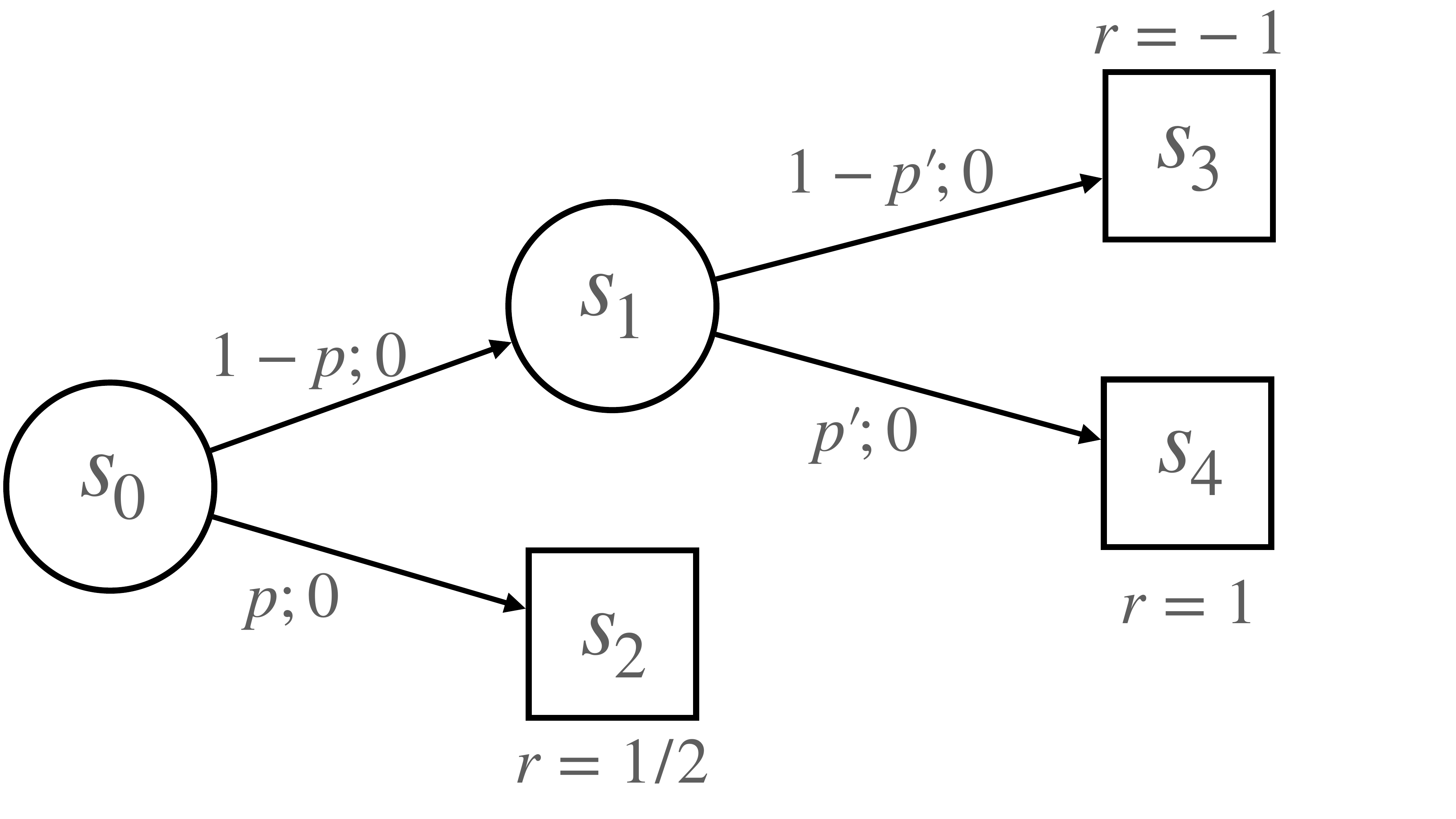}
    \end{center}
    \caption{Transitions and rewards for action $a_2$}
    \label{fig: lec 0 ex 2 a2}
    \end{subfigure}
    \end{center}
        \caption{MDP instance for \Cref{ex: lec 0 ex 2}. The transition probabilities and rewards are indicated beside each arc.}
        \label{fig: lec 0 ex 2}
    \end{figure}
\end{example}
\begin{remark}[Finite-horizon models]
    An additional property of \Cref{ex: lec 0 ex 2} (compared to \Cref{ex: lec 0 ex 1}) is that this MDP instance is absorbed after at most two periods. Hence, it actually represents a sequential decision problem {\em with a finite horizon}. There are two main models for finite-horizon problems:
    
    Model (A): we assume that we visit independent copies of the state set at each period (so that each state is visited only once), or
    
    Model (B): we assume that we can visit the same states multiple times. 

    Both models are equivalent, but stationary policies in case (A) correspond to Markovian policies in case (B). In these notes we'll focus on infinite-horizon models, and our results apply to finite-horizon models with model (A) (up to renormalizing the rewards at each state with the appropriate power of the discount factor).
\end{remark}
\paragraph{Dynamics and objective function.} Let's say we choose some policy $\pi \in \PiH$ and we ``let the dynamics run": with probability $p_{0}(s_0) \in [0,1]$ we start in state $s_0$, then we choose some action $a_0 \in \cA$ (possibly at random), transition to a next state $s_1$, pick another action, and so forth: $a_t \in \cA$ is sampled from the distribution $\pi(s_0,a_0,...,s_{t-1},a_{t-1},s_{t}) \in \Delta(\cA)$, and $s_{t+1} \in \cS$ is sampled from $P(s_t,a_t)\in \Delta(\cS)$. This results in an infinite trajectory of visited state-action pairs $(s_t,a_t)_{t \in \N}$ , an element of $(\cS \times \cA)^{\N}$, and given a policy $\pi \in \PiH$ and transition probabilities $P$, some trajectories are more likely than others. This can be made rigorous with the Ionescu-Tulcea theorem, which states that the pair $(\pi,P)$ induces a probability distribution over the set $(\cS \times \cA)^{\N}$ of infinite trajectories, compatible with the probabilities of any finite trajectories. We use the notation $\E^{\pi,P}$ for the expectation with respect to this probability distribution. 
A trajectory $(s_t,a_t)_{t \in \N}$ induces a stream of payoffs $(r(s_t,a_t,s_{t+1}))_{t \in \N}$, and given that $(\pi,P)$ induces a distribution over the possible $(s_t,a_t)_{t \in \N}$ then $(\pi,P)$ induces a distribution over the possible streams of payoffs $(r(s_t,a_t,s_{t+1}))_{t \in \N}$. We are now left with designing a scalar objective function. In these notes we focus on the following {\em discounted return}: for $\pi \in \PiH$ and $P \in \Delta(\cS)^{\cS \times \cA}$, 
\begin{equation}\label{eq:discount objective mdp}
    R(\pi,P) = \E^{\pi,P}\left[ \sum_{t=0}^{+\infty} \gamma^t r(s_t,a_t,s_{t+1}) \right]
\end{equation}
where $\gamma \in [0,1)$ is the discount factor and the expectation is over the streams of payoffs associated with (random) trajectories $(s_t,a_t)_{ t \in \N}$. The discount factor captures the diminishing importance of future rewards. In practice, this parameter governs the convergence rate and complexity of algorithms for discounted MDPs: most of them scale with $(1-\gamma)^{-1}$, so MDPs become increasingly hard to solve when $\gamma$ is close to $1$ and we want to be more farsighted.  
The main objective in an MDP is to solve the following maximization problem:
\begin{equation}\label{eq: lec 0 mdp}
    \sup_{\pi \in \PiH} R(\pi,P)
\end{equation}
A few modeling remarks:
\begin{itemize}
    \item We use the convention that MDPs are maximizing rewards here. Since $\E[X] = - \E[-X]$ we could have used the convention that MDPs minimize cost, and all results presented here remain the same. In particular we will use ``minimizing MDPs" to model the policy evaluation problem in \Cref{sec: lec 2 adv}.
    \item By linearity of the expectation, scaling by a positive factor and translating the rewards (i.e. $r \gets \alpha r + \beta$ with $\alpha>0, \beta \in \R$) does not change the set of optimal solutions in~\eqref{eq: lec 0 mdp}. 
    \item It is good to keep in mind that $\sum_{t=0}^{+\infty} \gamma^t r(s_t,a_t,s_{t+1})$ is a random variable, given that the trajectories $(s_t,a_t)_{ t \in \N}$ are themselves random.  
    \end{itemize}
The following definition will simplify the statements of our theorems.
\begin{definition}
\begin{itemize}
    \item     A policy $\pi \in \PiH$ is {\em optimal} if it is an optimal solution of~\eqref{eq: lec 0 mdp}.
    \item A policy $\pi \in \PiH$ is {\em uniformly optimal}  if it is an optimal solution of~\eqref{eq: lec 0 mdp} for all initial distributions $p_0 \in \Delta(\cS)$ at the same time.
\end{itemize}
\end{definition}    
When $p_0$ has full support, optimal policies are also uniformly optimal, but in general uniform optimality is a much stronger condition than optimality.
While, in principle, we are interested in computing optimal policies (and at this point it may not even be obvious that uniformly optimal policies exist), we will see that most algorithms actually return uniformly optimal policies. For some $\epsilon>0$ we will also use the notion of $\epsilon$-optimal policies, which are policies $\pi \in \PiH$ such that $R(\pi,P) \geq \sup_{\pi' \in \PiH} R(\pi',P) - \epsilon$, and uniformly $\epsilon$-optimal policies, which are $\epsilon$-optimal for all initial distributions.
\begin{remark}
    It is worth comparing the notion of uniform optimality with the notion of subgame perfect equilibrium for Nash equilibria in extensive-form games, which are equilibria that play optimally even in parts of the game trees that are not visited. Optimal policies (for a given initial distribution $p_0$) need not take ``good actions" in parts of the set of states that are not visited under $p_0$, in contrast to stationary uniformly optimal policies.
\end{remark}
The main difficulty in computing optimal policies lies in the multistage aspect of~\eqref{eq: lec 0 mdp}. Indeed when visiting a state $s$ we have to balance (a) achieving good instantaneous rewards by choosing actions yielding large values of $r(s,a,s')$, and (b) the probability to transition to a next state $s'$ where we may be cornered and unable to derive good rewards from there on. This sort of tradeoff is best expressed in the framework of value functions.
\paragraph{Value functions.} The objective function $\pi \mapsto R(\pi,P)$ can be decomposed as follows by conditioning on the realization of the first state:
\begin{align*}
    \E^{\pi,P}\left[ \sum_{t=0}^{+\infty} \gamma^t r(s_t,a_t,s_{t+1}) \right] & = \sum_{s \in \cS} p_{0}(s) \E^{\pi,P}\left[ \sum_{t=0}^{+\infty} \gamma^t r(s_t,a_t,s_{t+1}) \;\middle| \; s_0 = s \right] \\
    & = \sum_{s \in \cS} p_{0}(s) v^{\pi,P}(s) \\
    & = p_{0}\tr v^{\pi,P}
\end{align*}
with $v^{\pi,P} \in \R^{\cS}$ the {\em value function} associated with a policy $\pi \in \PiH$ and a transition kernel $P$:
\begin{equation}\label{eq:value funciton}
    v^{\pi,P}(s) = \E^{\pi,P}\left[ \sum_{t=0}^{+\infty} \gamma^t r(s_t,a_t,s_{t+1}) \; \middle| \; s_0 = s \right], \quad \forall \; s \in \cS.
\end{equation}
Note that in principle the trajectory laws depend on $p_0$ so that conditioning on $s_0=s$ can be ill-defined if $p_0(s)=0$; in the expression above with the conditioning on $s_0=s$ the expectation over the trajectory law is taken with initial state $s$, irrespective of the initial distribution.
The scalar $v^{\pi,P}(s)$ represents the ``value" that the decision maker would obtain if starting in state $s$ (with probability $1$) and playing according to $\pi$ along the trajectories. For stationary policies $\pi \in \PiS$ this has a closed-form. Indeed, for $\pi \in \PiS$ the pair $(\pi,P)$ induces a Markov chain over the finite set $\cS$, with transition matrix $L^{\pi,P} \in \R^{\cS \times \cS}$ such that
\[L^{\pi,P}_{ss'} = \sum_{a \in \cA} \pi(s,a) P(s,a,s'), \quad \forall \; s,s' \in \cS.\]
Indeed if we are in state $s$, then with probability $\pi(s,a)$ we pick action $a$ and then with probability $P(s,a,s')$ we transition to $s'$. Associated with this transition $s \to s'$ is the function $r^{\pi,P} \in \R^{\cS}$ defined as 
\begin{equation}\label{eq: lec 0 def r pi P}
    r^{\pi,P}(s) = \sum_{a \in \cA} \pi(s,a) P(s,a)\tr r(s,a),
\end{equation}
which captures the expected reward that we obtain at (only) one period in each state, with $r(s,a) = (r(s,a,s'))_{s' \in \cS} \in \R^{\cS}$. Finally, recall that the $t$-th power $(L^{\pi,P})^{t}$ captures the distribution over the state at period $t$. With this notation, we get that 
\begin{equation}\label{eq: lec 0 value function formula}
    v^{\pi,P} = \sum_{t=0}^{+\infty} \gamma^t (L^{\pi,P})^t r^{\pi,P} = (I - \gamma L^{\pi,P})^{-1}r^{\pi,P}.
\end{equation}
The first equality is by definition of value functions as in~\eqref{eq:value funciton} and interchanging the summation and the expectation (we can because of the discount factor). The second equality is by $(I-A)^{-1} = \sum_{t=0}^{+\infty} A^t$ for $A$ a matrix with spectral radius $\sigma(A)$ strictly less than $1$ (again, $\sigma(L^{\pi,P})=1$ because $L^{\pi,P}$ is a transition matrix, and $\sigma(\gamma A) = \gamma \sigma(A)$ for $\gamma>0$).
This shows that we can compute $v^{\pi,P}$, at least for stationary policies, by inverting a matrix. This can become expensive for large state sets.

Note that uniform optimality can be interpreted in terms of value functions. 
In particular, the policy $\pi^{P}  \in \PiH$ is uniformly optimal if and only if
\[\forall \; p_0 \in \Delta(\cS),  \forall \; \pi \in \PiH, R(\pi,P) \leq R(\pi^{P} ,P).\]
Hence, choosing $p_0$ that puts probability $1$ on a single state $s$, for each $s \in \cS$, we get that $\pi^{P} $ is uniformly optimal if and only if its value function dominates componentwise the value function of any other policy:
\[\forall \; \pi \in \PiH, v^{\pi,P} \leq v^{\pi^{P} ,P}.\]
In particular this implies that the value functions of all uniformly optimal policies must coincide.
In some lecture notes (and in the literature on stochastic games) the optimality criteria are usually presented as finding a stationary uniformly optimal policy.
In robust MDPs we will use optimality criteria based on $R(\pi,P)$ rather than on the componentwise comparisons of value functions, mostly because we are interested in situations where the transition probabilities are not independent across states (so-called {\em non-rectangular} instances) so that componentwise comparison of (worst-case) value functions may be meaningless (or rather, they only provide lower bounds on the worst-case returns). We will be able to recover uniformly optimal policies for robust MDPs only under some specific assumptions (like rectangularity).
\begin{continuance}{ex: lec 0 ex 1}
Consider the trajectories in the MDP instance from Figure~\ref{fig: lec 0 ex 1} starting in $s_0$. At even periods we are in state $s_0$ and at odd periods we are in state $G$ or $B$; we pick an action in $\{a_1,a_2\}$ at each period. A stream of payoffs is therefore equal to a sequence of $0$ and $1$, where even elements are $0$'s, and odd elements are $1$ (if we visit $G$ at this period) or $0$ (if we visit $B$ at this period). 
A stationary policy maps each element of $\{s_0,G,B\}$ to a probability distribution over $\{a_1,a_2\}$. Note that the choice of actions is irrelevant at $G,B$ (the rewards and transitions are the same for $a_1$ and $a_2$) so we can parametrize stationary policies by $x \in [0,1]$ modeling $\pi(s_0,a_1)$ (we can then recover $\pi(s_0,a_2)$ as $1-\pi(s_0,a_1)$).

We can compute the value function of a stationary policy as follows: 
$v(G) = 1+\gamma v(s_0), v(B) = 0+\gamma v(s_0)$, and 
\[v(s_0) = x \left( p \cdot \gamma v(G) + (1-p) \cdot \gamma v(B) \right) + (1-x) \left( (1-p) \cdot \gamma v(G) + p \cdot \gamma v(B)\right).\]
Replacing $v(G),v(B)$ by their expressions in terms of $v(s_0)$ in the equation above and rearranging the terms, we get
\[v(s_0) = \gamma \frac{xp + (1-p)(1-x) }{1-\gamma^2}\]
and we observe that $xp + (1-p)(1-x)$ is the probability to reach $G$ from $s_0$. Maximizing over $x \in [0,1]$ yields an optimal policy among stationary policies. Hence if we maximize over $x \in [0,1]$ we end up choosing $x=0$ if $p \in [0,1/2)$, $x=1$ if $p \in (1/2,1]$, and at $p=1/2$ we are indifferent and we can pick any $x \in [0,1]$ as an optimal solution.
\end{continuance}
\begin{continuance}{ex: lec 0 ex 2}
To compute value functions one can use backward induction, that is, we first compute the value functions of the terminal states then propagate the value backward to their predecessors. Let $\pi \in \PiS$ be a policy, parametrized by $x \in [0,1]$ the probability to play $a_1$ in $s_0$ and $y \in [0,1]$ the probability to play $a_1$ in $s_1$. We get $v(s_4) = \frac{1}{1-\gamma},v(s_3) = -\frac{1}{1-\gamma},v(s_2) = \frac{1}{2}\frac{1}{1-\gamma}.$  This yields
\[v(s_1) = 0+ y \gamma ( p' v(s_3) + (1-p')v(s_4)) + (1-y)\gamma ((1-p') v(s_3) + p' v(s_4)) = \frac{\gamma}{1-\gamma}(2y-1)(1-2p').\]
Now $v(s_0) = x \gamma(p v(s_1) + (1-p)v(s_2)) + (1-x) \gamma ((1-p) v(s_1) +p v(s_2))$.
This yields the following closed-form expression for $v(s_0)$:
\begin{equation}\label{eq: lec 0 ex 2 vs0}
    v(s_0) = \frac{\gamma}{1-\gamma}\left( \gamma (1-2p') \left( 2(2p-1) xy - (2p-1)x + 2 (1-p)y - (1-p) \right) + \frac{1-2p}{2} x + \frac{p}{2}\right).
\end{equation}
\end{continuance}
A very important fact about the function $(\pi,P) \mapsto R(\pi,P)$ is that it is somewhat ``ill-behaved" from a classical optimization standpoint. In particular, the function $\pi \mapsto R(\pi,P)$ need not be concave (we want to maximize it) and the function $P \mapsto R(\pi,P)$ need not be convex (we will want to minimize later on). 
This is exemplified by the computations of the value functions from \Cref{ex: lec 0 ex 2}, where we see that starting in $s_0$ yields a value function $v(s_0)$ that depends on some products of the probabilities chosen by the policy (terms in $x \cdot y$ in Equation~\eqref{eq: lec 0 ex 2 vs0}) and some products of the transition probabilities (terms in $p \cdot p'$ above). In fact if we had more periods the return $R(\pi,P)$ would involve higher-order products of the entries of $\pi \in \PiS$. Thus in general maximizing $\pi \mapsto R(\pi,P)$ over stationary policies does not fall in the class of ``usual easy problems" that we are used to in convex optimization (let alone optimizing over the entire class of history-dependent policies). The next section shows how to solve MDPs.
\paragraph{State-action frequencies.}
The previous reformulation relating returns and value functions is the cornerstone of the analysis of robust MDPs. We derive here another very useful equality that relates the return and the {\em state-action frequencies} (also called occupancy frequencies). Let $\pi \in \PiH$. Then
\begin{align*}
    R(\pi,P) & = \E^{\pi,P}\left[ \sum_{t=0}^{+\infty} \gamma^t r(s_t,a_t,s_{t+1}) \right] \\
    & =  \sum_{t=0}^{+\infty} \gamma^t \E^{\pi,P}\left[ r(s_t,a_t,s_{t+1}) \right] \\
    & =  \sum_{t=0}^{+\infty} \gamma^t  \sum_{s,a} P(s,a) \tr r(s,a) \cdot \bP^{\pi,P}((s_t,a_t) = (s,a))  \\
    & =   \sum_{s,a} P(s,a) \tr r(s,a) \sum_{t=0}^{+\infty}  \gamma^t     \bP^{\pi,P}((s_t,a_t) = (s,a))  \\
    & = \sum_{s,a} \mu^{\pi}(s,a) P(s,a) \tr r(s,a)
   \end{align*}
where the first equality is by definition, the second equality is to be justified in \Cref{sec: lec 0 fdt}, and the other equalities are from reformulating and rearranging the terms (note $\bP^{\pi,P}$ is the probability with respect to the distribution over the infinite trajectories induced by the pair $(\pi,P)$), and where we define $\mu^{\pi} \in \R^{\cS \times \cA}_{+}$ as, for any pair $(s,a) \in \cS \times \cA$,
\begin{equation}\label{eq: lec 0 state-action occupancy}
    \mu^{\pi}(s,a) = \sum_{t=0}^{+\infty} \gamma^t \bP^{\pi,P}((s_t,a_t) = (s,a))
\end{equation}
The scalar $\mu^{\pi}(s,a)$ is the discounted sum of probabilities of visiting the pair $(s,a)$ across the entire infinite horizon of play, and $\mu^{\pi}$ is referred to as the {\em state-action frequencies.} Overall we obtained that for $\pi \in \PiH$ we have
    \begin{equation}\label{eq: lec 0 return frequencies}
    R(\pi,P) = \sum_{s \in \cS} \sum_{a \in \cA} \mu^{\pi}(s,a) P(s,a) \tr r(s,a)
\end{equation}
and in fact the above expression is also linear in the rewards (which will play a role in \Cref{sec: lec 3 uncertain rewards} for the case of uncertain rewards).
Note that \eqref{eq: lec 0 return frequencies} has a natural interpretation: the state-action frequencies $\mu^{\pi}(s,a)$ count in a discounted manner how much time we spend at a state-action pair, and $P(s,a)\tr r(s,a)$ is the expected instantaneous reward we get from $(s,a)$. We can also derive \eqref{eq: lec 0 return frequencies} for stationary policies by noting that $R(\pi,P) = p_0 \tr v^{\pi,P}$, using \eqref{eq: lec 0 value function formula} to get $R(\pi,P) = p_0 \tr (I - \gamma L^{\pi,P})^{-1} r^{\pi,P} = r^{\pi,P \; \top} (I - \gamma L^{\pi,P \; \tr})^{-1}p_0$ and then reformulating this last expression by isolating the factor in front of $P(s,a) \tr r(s,a)$ as $\mu^{\pi}(s,a) = \pi(s,a)\left((I - \gamma L^{\pi,P \; \tr})^{-1}p_0\right)(s)$ for each $(s,a)$. The matrix-vector product $(I - \gamma L^{\pi,P \; \tr})^{-1}p_0$ is precisely the term giving the next flow constraints~\eqref{eq: lec 0 flow constraint mu}.
Note that for $s' \in \cS$ we have
\begin{align*}
    \sum_{a \in \cA} \mu^{\pi}(s',a) & = \sum_{t=0}^{+\infty} \gamma^t \bP^{\pi,P}((s_t) = s') \\
    & = p_0(s') +  \sum_{t=1}^{+\infty} \gamma^t \bP^{\pi,P}((s_t) = s') \\
    & = p_0(s') +  \gamma \sum_{t=0}^{+\infty} \gamma^t \bP^{\pi,P}((s_{t+1}) = s') \\
    & = p_0(s') +  \gamma \sum_{s,a} P(s,a,s') \sum_{t=0}^{+\infty} \gamma^t \bP^{\pi,P}((s_{t},a_t) = (s,a)) \\
    & =  p_0(s') + \sum_{s \in \cS} \sum_{a \in \cA} \gamma P(s,a,s')\mu^{\pi}(s,a)
\end{align*}
i.e., we have shown that the following natural flow constraints hold:
\begin{equation}\label{eq: lec 0 flow constraint mu}
    p_0(s') + \sum_{s \in \cS} \sum_{a \in \cA} \gamma P(s,a,s')\mu^{\pi}(s,a) = \sum_{a \in \cA} \mu^{\pi}(s',a), \quad \forall \; s' \in \cS
\end{equation}
where the right-hand side is the (discounted) amount of time spent in state $s'$ and the left-hand side is the incoming flow at $s'$. Therefore the set $\cQ \coloneqq \{ \mu^{\pi} \in \R^{\cS \times \cA}_{+} \; | \; \pi \in \PiH\}$ can be written as the following polytope:
\begin{equation}\label{eq: lec 0 frequency polytope}
    \cQ = \left\{ \mu \in \RSA \; | \; \mu \geq 0, p_0(s') + \sum_{s \in \cS} \sum_{a \in \cA} \gamma P(s,a,s')\mu(s,a) = \sum_{a \in \cA} \mu(s',a), \quad \forall \; s' \in \cS \right\}
\end{equation}
and given $\mu \in \cQ$ we can recover a feasible {\em stationary} policy $\pi \in \PiS$ as \[\pi(s,a) = \mu(s,a) / \sum_{b \in \cA} \mu(s,b)\] if $\sum_{b \in \cA} \mu(s,b) \neq 0$, and $\pi(s,a) = 1/|\cA|$ otherwise (we won't visit state $s$ anyway so we can play anything there). Conversely, given $\pi \in \PiS$, \eqref{eq: lec 0 flow constraint mu} implies that we can compute $\mu^{\pi}$ by solving a linear system of equations, which can be done efficiently. In turn, the extreme points of $\cQ$ correspond to the state-action frequencies induced by stationary deterministic policies.

Overall we have seen three crucial properties of state-action frequencies:
\begin{itemize}
    \item we can map policies $\pi \in \PiH$ to elements $\mu^{\pi}$ of the polytope $\cQ$,
    \item we can map each element $\mu \in \cQ$ to a {\em stationary policy} $\pi$,
    \item the return $R(\pi,P)$ of a policy $\pi \in \PiH$ is a linear function of the state-action frequency $\mu^{\pi}$.
\end{itemize}
In fact the three points above already suggest that stationary policies will be sufficient for optimality in MDPs, as we will prove next.
The state-action frequencies are powerful tools in the MDP toolkit and have a natural linear programming interpretation (see \Cref{sec: lec 0 algo}).
The state-action frequencies will not play a role for robust MDPs with uncertain transitions. This is because the polytope $\cQ$ actually depends on $P$ for the expressions of the flow constraints (in a ``non-rectangular manner": every $s'$-th constraint depends on $P(s,a,s')$ across all $(s,a) \in \cS \times \cA$). Additionally $\cQ$ depends on $p_0$ so uniform optimality is harder to prove based on frequencies than based on value functions (which are defined component-wise and independently of $p_0$). Still the polytope $\cQ$ will play a crucial role for robust MDPs with uncertain rewards since the return is linear in the rewards (see \Cref{sec: lec 3 uncertain rewards}).
\section{Main results for MDPs}\label{sec: lec 0 main}
\subsection{Operators and optimal policies.}
The first crucial result relates the value function of a stationary policy to the fixed point of an operator. For $\pi \in \PiS$, we write $T^{\pi,P}:\R^{\cS} \to \R^{\cS}$ for the operator defined as
\begin{equation}\label{eq: lec 0 T pi P}
    T^{\pi,P}v(s) = \sum_{a \in \cA} \pi(s,a) P(s,a)\tr (r(s,a) + \gamma v), \quad \forall \; v \in \R^{\cS}, s\in \cS.
\end{equation}
We can interpret $T^{\pi,P}v(s)$ as the expected return that one would obtain if starting in state $s \in \cS$, playing according to $\pi \in \PiS$ for the first period, transitioning to a next state $s'$ with distribution $P(s,a) \in \Delta(\cS)$, and quitting the dynamics with the termination value $v(s') \in \R$, i.e., we have
\[T^{\pi,P}v(s) = \E_{a \sim \pi(s),s ' \sim P(s,a)} \left[ r(s,a,s') + \gamma  v(s')\right].\]
Remark that the $s$-th component of $T^{\pi,P}v$ only depends on the $s$-th component of $\pi$ (i.e., how we choose actions at state $s$) and on the $s$-th component of $P$ (i.e., the transition probabilities out of state $s$), and not on the other components of $\pi$ and $P$ related to other states (this will play a crucial role in the later chapters).

From there we get that if we define the sequence $v^{t+1} := T^{\pi,P}v^t, v^0 = 0 \in \R^{\cS}$ then $v^{t}$ is the total discounted return obtained up to period $t-1$ (and terminating with a reward of $0$ in the last state): 
\[    v^{t}(s) = \E^{\pi,P}\left[ \sum_{t'=0}^{t-1} \gamma^{t'} r(s_{t'},a_{t'},s_{t'+1}) \; \middle| \; s_0 = s \right], \quad \forall \; s \in \cS.\]
This suggests that we can interpret $v^{\pi,P}$ as the limit of a sequence $\left(T^{\pi,P}\right)^{t}0_{\cS}$ for $t \to + \infty$ (with $0_{\cS}=(0,...,0) \in \R^{\cS}$). In fact, if we consider a sequence of policies $\pi^0,..., \pi^t \in \PiS$ then $T^{\pi^0,P} \circ ... \circ T^{\pi^t,P} 0_{\cS}$ is the reward  obtained after $t+1$ periods (recall period indices start at $0$) by playing according to policy $\pi^0$ at period $0$, then according to policy $\pi^1$ at period $1$, etc.

To make this rigorous, we note that $T^{\pi,P}$ satisfies the following two key properties, defined for general operators $F:\R^{\cS} \to \R^{\cS}$.
\begin{definition}\label{def: lec 0 monotone contraction}
Let $F:\R^{\cS} \to \R^{\cS}$ be an operator.
\begin{itemize}
    \item $F$ is {\em monotone} if
    \begin{equation}\label{eq:monotone}
        u \leq v \Rightarrow Fu \leq Fv, \forall \; u,v \in \R^{\cS}.
    \end{equation}
    
    \item $F$ is a {\em contraction} if
    \begin{equation}\label{eq:contraction}
        \exists \; \lambda \in [0,1), \forall \; u,v \in \R^{\cS}, \|Fu-Fv\|_{\infty} \leq \lambda \|u-v\|_{\infty}.
    \end{equation}
\end{itemize}
\end{definition}
As highlighted from the calculations of the value functions for our simple examples above, we can compute $v^{\pi,P}$ as the fixed point of a monotone contraction. This yields the next proposition.
\begin{proposition}\label{prop: lec 0 operators contraction}
    Let $\pi \in \PiS$.
    \begin{enumerate}
        \item The operator $T^{\pi,P}$  defined in~\eqref{eq: lec 0 T pi P} is a monotone contraction.
        \item The value function $v^{\pi,P} \in \R^{\cS}$ is the unique fixed point of the operator $T^{\pi,P}$.
    \end{enumerate}
\end{proposition}
\begin{proof}
\begin{enumerate}
    \item Monotonicity follows from $P(s,a) \in \R^{\cS}_{+}$ for any $(s,a) \in \cS \times \cA$:
    \[u \leq v \Rightarrow P(s,a) \tr u \leq P(s,a) \tr v \Rightarrow T^{\pi,P}u \leq T^{\pi,P}v.\]
    The contraction property follows from $P(s,a) \in \R^{\cS}_{+}$ and $\sum_{s' \in \cS} P(s,a,s') = 1$, i.e., from $P(s,a) \in \Delta(\cS)$. To see this, note that 
    \[T^{\pi,P}u = T^{\pi,P}(v + (u-v)) \leq T^{\pi,P}v + \gamma \| u-v\|_{\infty} 1_{\cS}\]
    where the inequality is from the componentwise inequality $u-v \leq \| u-v\|_{\infty} 1_{\cS}$. Rearranging the terms we obtain $T^{\pi,P}u - T^{\pi,P}v \leq  \gamma \| u-v\|_{\infty} 1_{\cS}.$
    Reversing the roles of $u$ and $v$ yields the contraction property for $T^{\pi,P}$.
    \item We have $v^{\pi,P} = (I-\gamma L^{\pi,P})^{-1}r^{\pi,P}\iff(I-\gamma L^{\pi,P})v^{\pi,P}=r^{\pi,P} \iff v^{\pi,P} = r^{\pi,P} + \gamma L^{\pi,P}v^{\pi,P} \iff v^{\pi,P} = T^{\pi,P}v^{\pi,P}$ where the first equivalence is by multiplying on both sides by $(I-\gamma L^{\pi,P})$, the second equivalence is by reordering the terms on both sides, and the last equivalence uses the expression of $T^{\pi,P}$. 
\end{enumerate}
\end{proof}
The second important result for MDPs pertains to the existence of stationary deterministic optimal policies. We will use the following {\em Bellman optimality operator} $T^{P}:\R^{\cS} \to \R^{\cS}$ defined as
\begin{equation}\label{eq: lec 0 T P}
    T^{P}v(s) = \max_{\pi \in \PiS} T^{\pi,P}v(s)
\end{equation}
For $s \in \cS$, note that the $s$-th component of $T^{\pi,P}v$ only depends on $\pi(s) \in \Delta(\cA)$ and $(P(s,a))_{a \in \cA} \in \Delta(\cS)^{\cA}$, i.e., on the $s$-th component of $P \in \Delta(\cS)^{\cS \times \cA}$. Hence we actually have
\begin{align*}
    T^{P}v(s)  & =  \max_{\pi(s) \in \Delta(\cA)}\sum_{a \in \cA} \pi(s,a) P(s,a)\tr (r(s,a) + \gamma v) \\ 
    & = \max_{a \in \cA} P(s,a)\tr (r(s,a) + \gamma v)
\end{align*}
for each $v \in \RS$ and $s \in \cS$,
where the second equality uses the fact that the maximum of a linear form over a polytope is attained (among other solutions) by at least one extreme point, and we can identify the extreme points of $\Delta(\cA)$ with the deterministic choice of actions. In short this shows that
\[T^{P}v(s) = \max_{\pi \in \PiSD} T^{\pi,P}v(s).\]
Since the maximum over a finite set is always attained, we get the following {\em attainability property} for the family $\{T^{\pi,P} \; | \; \pi \in \PiSD\}$ with respect to $T^{P}$: for any $v\in \R^{\cS}$, $T^P v$ is actually $T^{\pi,P}v$ for some stationary deterministic $\pi \in \PiSD$, i.e.,
\begin{equation}\label{eq:attainability T P}
    \forall \; v \in \R^{\cS}, \exists \; \pi \in \PiSD, T^{P}v=T^{\pi,P}v.
\end{equation}
As a side remark, note that we should not hope to be able to reverse the quantifiers in the above equations except in some specific cases. We can now state the ``fundamental theorem of MDPs".
\begin{theorem}\label{th: lec 0 mdp main}
    Let $\cM=(\cS,\cA,r,P,p_0)$ be an MDP with finite sets of states and actions and consider computing an optimal policy for the discounted return as in~\eqref{eq: lec 0 mdp}.
    \begin{enumerate}
        \item {\em Bellman equation.} The Bellman optimality operator defined in~\eqref{eq: lec 0 T P} is a monotone contraction with unique fixed point $v^{P}  \in \R^{\cS}$, satisfying the following Bellman equations in the variable $v \in \R^{\cS}$:
        \begin{equation}\label{eq: lec 0 Bellman mdp}
            v(s) = \max_{a \in \cA} P(s,a)\tr (r(s,a)  + \gamma v), \quad \forall \; s \in \cS.
        \end{equation}
        \item {\em Uniform optimality.} A uniformly optimal policy $\pi^{P} $ exists and can be chosen stationary and deterministic, and the value function $v^{\pi^{P},P}$ of a stationary uniformly optimal policy $\pi^{P} $ coincides with $v^{P} \in \R^{\cS}$, the unique fixed point of $T^{P}$.
        \item {\em Optimality conditions.} We can compute a stationary uniformly optimal policy by choosing any policy that jointly attains the $\arg \max$ in $T^{P}v^{P} (s)$ across all $s \in \cS$, i.e., we have:
        \begin{align*}
            \pi^{P}  \in \PiS \text{ is uniformly optimal } & \iff  T^{P}v^{P} =T^{\pi^{P} ,P}v^{P}   \\
            & \iff \text{the support of } \pi^{P} (s) \in \DA \text{ is contained in} \\
            & \; \arg \max_{a \in \cA} P(s,a)\tr(r(s,a) + \gamma v^{P} ), \forall\; s \in \cS \\
            & \iff \pi^{P}  \in \cap_{s \in \cS} \arg \max_{\pi \in \PiS} T^{\pi,P}v^{P} (s)
        \end{align*}
    \end{enumerate}
\end{theorem}
\begin{proof}
We actually prove all statements at the same time. 

The proof that $T^{P}$ is a monotone contraction is identical to that of the same properties of $T^{\pi,P}$ for $\pi \in \PiS$ and is omitted. Let $v^{P}  \in \R^{\cS}$ be its unique fixed point. Let $\pi \in \PiS$ be any policy such that $T^{P}v^{P}  = T^{\pi,P}v^{P} $ (the existence of such $\pi$ is from the attainability property). Note that $T^{P}v^{P}  = T^{\pi,P}v^{P} $ implies in particular that $v^{P}  = T^{\pi,P}v^{P} $ since $v^{P} $ is the fixed point of $T^{P}$. Hence by uniqueness of the fixed point of $T^{\pi,P}$ we actually get $v^{P}  = v^{\pi,P}$. We now prove that $\pi$ is a uniformly optimal policy. The fact that $v^{\pi,P} \geq v^{\pi',P}$ for any stationary $\pi' \in \PiS$ is easy to prove: indeed we have 
\[v^{\pi,P}=v^{P}  =T^Pv^{P}  \geq T^{\pi',P}v^{P} \]
which implies, using monotonicity, that $v^{P} \geq \left(T^{\pi',P}\right)^t v^{P}$ for each $t \in \N$, and taking the limit as $t \to + \infty$ and using the contraction property, that $v^{\pi,P} \geq v^{\pi',P}$ for any $\pi' \in \PiS$. We must also show that $v^{\pi,P} \geq v^{\pi',P}$ for $\pi' \in \PiH$, i.e. we must show that $\pi$ is better than history-dependent policies.

Let $\pi' \in \PiH$ and $h_t$ be a history ending at period $t$ at $s_t$. Consider playing for one period then quitting the game with the value given by $v^{P} $. Then from $s_t$ the policy $\pi'$ would obtain
\begin{align*}
    \E^{\pi',P} \left[ r(s_t,a_t,s_{t+1}) + \gamma v^{P} (s_{t+1}) \middle| h_t \right] 
    & = \sum_{a \in \cA} \pi'(h_t,a)P(s_t,a)\tr\left( r(s_t,a) + \gamma v^{P}  \right) \\
    & \leq \max_{a \in \cA} P(s_t,a)\tr\left( r(s_t,a) + \gamma v^{P}  \right) = T^{P}v^{P} (s_t) = v^{P} (s_t)
\end{align*}
where the expectation in the first line is over the chosen action $a_t$ and the next state $s_{t+1}$, and where $\pi'(h_t,a)$ is the probability that $\pi'$ chooses action $a$ for $h_t$. Repeating this argument over $H+1$ periods we get 
\[\E^{\pi',P} \left[ \sum_{t=0}^{H} \gamma^tr(s_t,a_t,s_{t+1}) + \gamma^{H+1} v^{P} (s_{H+1}) \; \middle| s_0 = s \right] \leq v^{P} (s).\]
Since $\gamma^H \to 0$ and the rewards are bounded, we can conclude that 
\[v^{\pi',P}(s)=\E^{\pi',P} \left[ \sum_{t=0}^{+\infty}   \gamma^t r(s_t,a_t,s_{t+1}) \; \middle| s_0 = s \right] \leq v^{P} (s).\]
Since $v^{P}  = v^{\pi,P}$ we conclude that $\pi$ is a uniformly optimal policy. 

We now prove that any stationary uniformly optimal policy $\bar{\pi} \in \PiS$ is such that $T^{P}v^{P}=T^{\bar{\pi},P}v^{P}$. For $\bar{\pi} \in \PiS$ stationary uniformly optimal we must have $v^{\bar{\pi},P} = v^{P}$ (since $v^{P}$ is attained by some policies). This implies that $T^{P}v^{P}=v^{P} = v^{\bar{\pi},P} = T^{\bar{\pi},P}v^{\bar{\pi},P}=T^{\bar{\pi},P}v^{P}$, i.e., $T^{P}v^{P} = T^{\bar{\pi},P}v^{P}$, as claimed.

Finally, we can find a stationary {\em deterministic} uniformly optimal policy since each component of the Bellman operator~\eqref{eq: lec 0 T P} optimizes a linear form over the set $\PiS$, and the extreme points of the set of stationary policies $\PiS$ coincide with the set of stationary deterministic policies $\PiSD$.
\end{proof}
From the above theorem we get that the equation $v = T^{P}v$ in the variable $v \in \R^{\cS}$ has a unique solution which is equal to $v^{P} \in \R^{\cS}$; this equation is referred to as the ``Bellman optimality equations".
There are some misconceptions that we should avoid when thinking about the above theorem:
\begin{enumerate}
    \item {\em Non-stationary optimal policies.} While an optimal policy may be chosen stationary and deterministic, there may exist stationary randomized optimal policies, and even history-dependent optimal policies. This is because {\em any} policy that chooses (randomly) an action in $\arg \max_{a \in \cA} P(s_t,a)\tr (r(s_t,a) + \gamma v^{P})$ at each state $s_t \in \cS$ {\em visited} at period $t \in \N$ is optimal across all initial distributions $p_0$. Hence in general it is improper to talk of ``the" solution of \eqref{eq: lec 0 mdp} (e.g. ``the optimal policy ..."), one should talk of
``a" solution.

I emphasized the word ``visited" above because, obviously, we don't need to play optimally at states that are never visited! This highlights that while we have defined optimal policies as optimal solutions in~\eqref{eq: lec 0 mdp} (which in principle depends on the choice of $p_0$), in practice we will most of the time be interested in computing a uniformly optimal solution $\pi^{P}  \in \PiSD$ that is optimal in a stronger sense (optimal for all $p_0$ at the same time).

\item {\em Optimality vs. uniform optimality.} Note that the second statement is false if we replace {\em ``... of a stationary uniformly optimal policy $\pi^{P} $ ..."} by {\em ``... of an optimal stationary policy $\pi^{P} $ ..."}, i.e., the value function of an optimal policy need not be equal to the fixed point of the Bellman operator!

Indeed an optimal policy that is non-uniformly optimal must choose suboptimal actions in some states, which from the above item means it doesn't play an action in the $\arg \max$ of the Bellman operator $T^P$ at some states.

A more mathematical reformulation of the last paragraph is as follows: given $v^{P} $ the fixed point of $T^{P}$, any policy $\pi \in \PiS$ such that $T^{\pi,P}v^{P}  = T^{P}v^{P} $ is uniformly optimal. However, there may be optimal (non-uniformly optimal) policies $\pi' \in \PiS$ that do not satisfy $T^{\pi',P}v^{P}  = T^{P}v^{P} $. This is because 
\[R(\pi,P) = \sum_{s \in \cS} p_0(s) v^{\pi,P}(s)\]
hence if a policy chooses suboptimal actions at a state that is never visited, this does not impact the discounted return $R(\pi,P)$. I give an example of this type of issue below (see continued \Cref{ex: lec 0 ex 2}). We also note that the existence of uniformly optimal policies makes the initial distribution $p_0$ irrelevant in the description of the MDP instance. This will not be true in general for robust MDPs, but this will be true for robust MDPs with rectangular models of uncertainty (see \Cref{chap: lec 1 robust MDP} and \Cref{chap: lec 2 rectangular models}).

\item {\em Beyond stationarity.} It is worth emphasizing that the value functions of non-stationary (history-dependent or Markovian) policies $\pi \in \PiH$ are well-defined objects. However we can't define the Bellman operator $T^{\pi,P}$ as in~\eqref{eq: lec 0 T pi P} because~\eqref{eq: lec 0 T pi P} involves the notation $\pi(s,a)$, which is ill-defined for non-stationary policies: by definition for non-stationary policies we need more information than just the current state to pick some actions (for Markovian policies we need additionally the current period, and for history-dependent policies we need the full history up to the current state).
    \end{enumerate}
We compute the optimal policies and uniformly optimal policies in the next example.
\begin{continuance}{ex: lec 0 ex 1}
We saw that the value function of a policy choosing $a_1$ with probability $x \in [0,1]$ in state $s_0$ is $v(s_0) = \frac{\gamma}{1-\gamma^2} (xp + (1-p)(1-x))$. Note that $xp + (1-p)(1-x) = (2p-1)\cdot x + (1-p)$. Maximizing this last expression in $x$ allows us to compute optimal policies: $if p < 1/2$ then $2p-1 < 0$ and we pick action $a_2$ with probability $1$ (i.e., we choose $x\opt=0$), otherwise if $p > 1/2$ we pick action $a_1$ (i.e., $x\opt=1$). At $p=1/2$ we have $2p-1=0$, we are indifferent and we always get a payoff proportional to $1-p=1/2$. Note that here we also computed a uniformly optimal policy (since from $G,B$ we transition deterministically to $s_0$).
\end{continuance}
\begin{continuance}{ex: lec 0 ex 2}
The optimal policies depend on the scalar values $p,p' \in [0,1]$ and it is a bit tedious to enumerate all possibilities. Assume for now that we start in $s_0$ ($p_0(s_0)=1$ and $p_0(s)=0$ for $s \in \{s_1,s_2,s_3,s_4\}$) and that $p=0,p'=1/4$. The value functions that we already computed simplify as follows. Choosing action $a_1$ at $s_1$ ($y=1$) yields a value from $s_1$ of $\frac{\gamma}{2(1-\gamma)}$ while action $a_2$ ($y=0$) yields a value from $s_1$ of $-\frac{\gamma}{2(1-\gamma)}$. Hence we want to pick action $a_1$ at $s_1$, i.e., $y=1$. This makes sense just looking at the probabilities: since $p' = 1/4 < 1/2$ then $a_1$ gives us more chance than $a_2$ to end up in the ``good" state $s_4$ with reward $+1$ instead of ending up in the ``bad" state $s_3$ with reward $-1$.

Now let us compare the two actions at state $s_0$. If we choose $a_1$ at $s_0$ we transition deterministically to $s_2$ and we get a return of $\gamma v(s_2) = \frac{\gamma}{2(1-\gamma)}$. Otherwise if we pick $a_2$ at $s_0$ we transition deterministically to $s_1$ and we get a return of $\gamma v(s_1) =  \frac{\gamma^2}{2(1-\gamma)}$. Since $\gamma^2 < \gamma$ for $\gamma \in (0,1)$, we see that we prefer to transition to $s_2$ from $s_0$, i.e., we pick action $a_1$ in $s_0$. 

This shows that the action that we choose at $s_1$ is irrelevant: any policy that chooses $a_1$ with probability $1$ at $s_0$ (i.e. $x=1$) is optimal, whatever the choice of actions at $s_1$. However if we want to choose a {\em uniformly} optimal policy we must also play optimally at $s_1$ (i.e., choose action $a_1$ at $s_1$).
\end{continuance}
We conclude by highlighting another important difference between the set of stationary uniformly optimal policies and the set of stationary optimal policies: the former set is always convex, while the latter may not be convex. This is because optimal policies need not take optimal actions in parts of the set of states that they don't visit, whereas uniform optimality imposes componentwise dominance of the value function, i.e., a uniformly optimal policy plays optimally starting from any state. This convexity property is why it is often much simpler to work with uniform optimality than with vanilla optimality. 
\begin{lemma}\label{lem: lec 0 convex sets unf policy}
    \begin{enumerate}
        \item The set of stationary uniformly optimal policies is always a convex closed set.
        \item The set of stationary optimal policies is always closed but may not be convex.
    \end{enumerate}
\end{lemma}
\begin{proof}
    \begin{enumerate}
        \item From the third statement in \Cref{th: lec 0 mdp main} we know that the set of stationary uniformly optimal policies is the finite intersection over $s \in \cS$ of the sets $\arg \max_{\pi \in \PiS} T^{\pi,P}v^{P} (s)$. Each of these optimization programs optimizes a linear form $\pi \mapsto T^{\pi,P}v^{P} (s)$ over the polytope $\PiS$, hence these sets are convex, and the intersection of these sets is convex and closed hence it is a convex closed set.
        \item The set of stationary optimal policies is closed by continuity of $\pi \mapsto R(\pi,P)$. 
        For the non-convexity of the set of optimal policies, consider the MDP instance depicted in \Cref{fig: lec 0 ex 4}.
        \begin{figure}[h!]
    \begin{center}
    \begin{subfigure}{0.45\textwidth}
    \begin{center}
        \includegraphics[width=\linewidth]{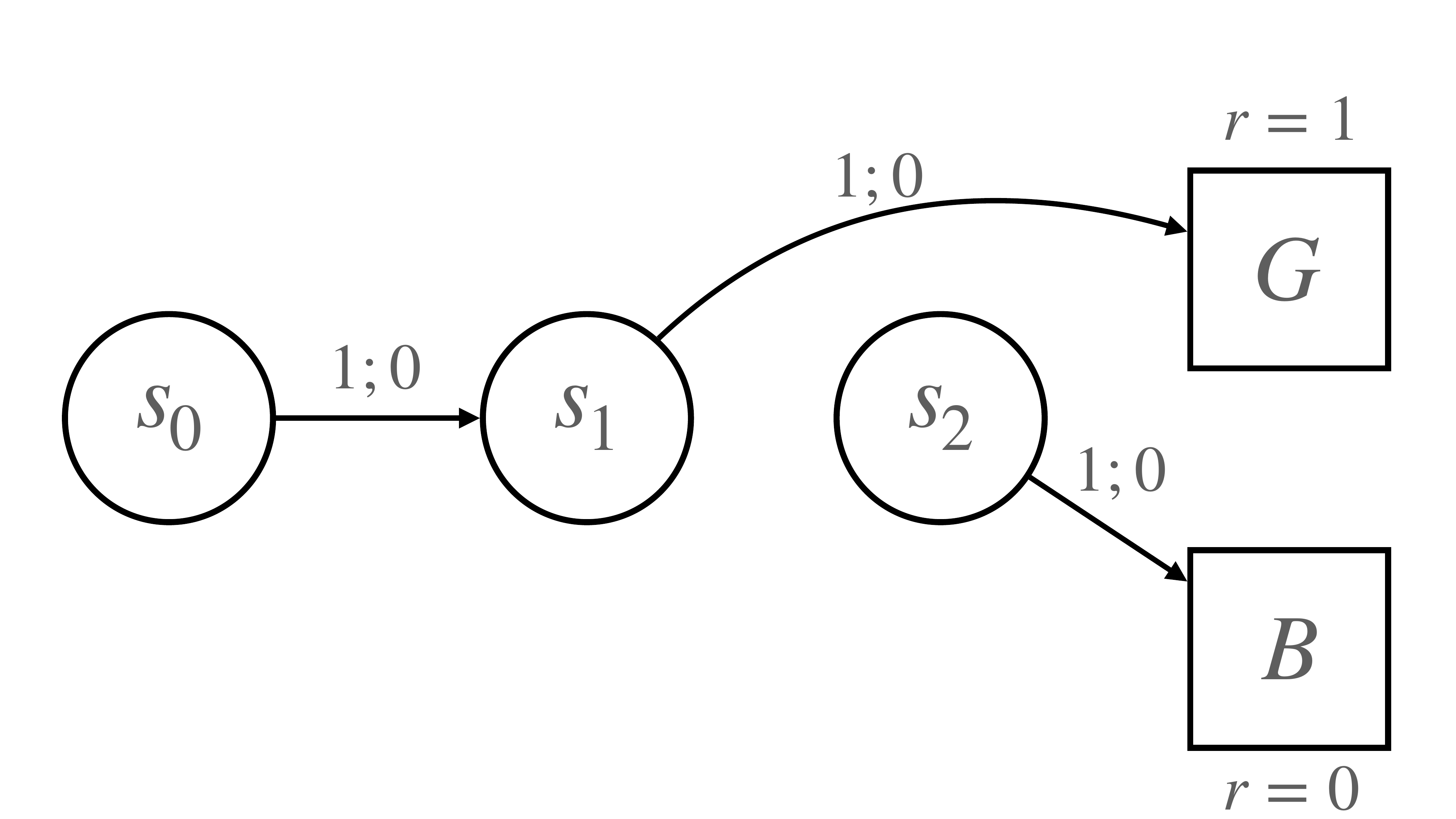}
    \end{center}
    \caption{Transitions and rewards for action $a_1$}
    \label{fig: lec 0 ex 4 a1}
    \end{subfigure}
    \begin{subfigure}{0.45\textwidth}
    \begin{center}
        \includegraphics[width=\linewidth]{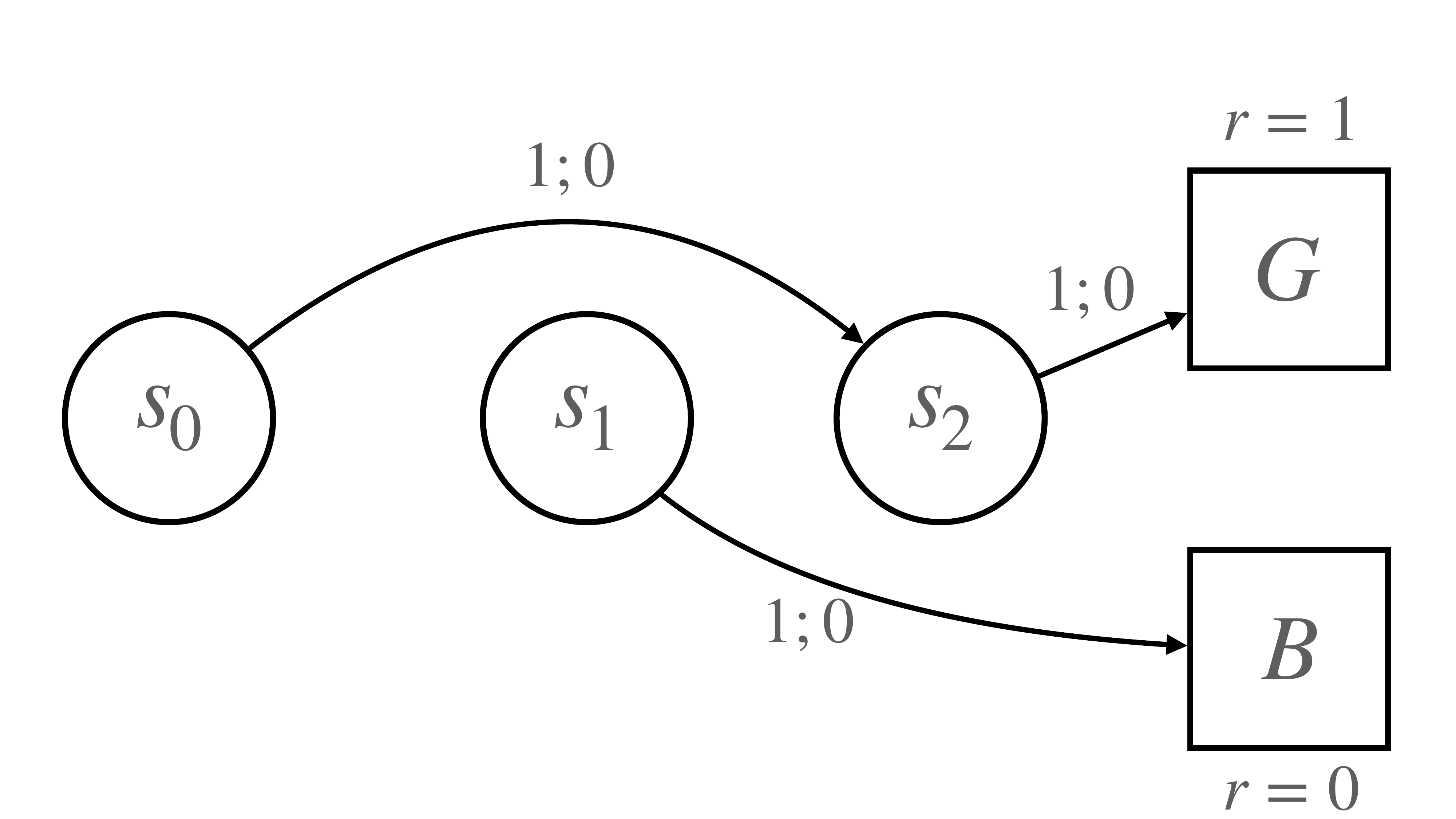}
    \end{center}
    \caption{Transitions and rewards for action $a_2$}
    \label{fig: lec 0 ex 4 a2}
    \end{subfigure}
    \end{center}
        \caption{The transition probabilities and rewards are indicated beside each arc. Absorbing states are represented by squares.}
        \label{fig: lec 0 ex 4}
    \end{figure}
        There are five states $\{s_0,s_1,s_2,G,B\}$ and two actions $a_1,a_2$. The states $G,B$ are absorbing and the MDP starts in state $s_0$ with probability $1$. The rewards are $0$ everywhere except at $G$ where we obtain a reward of $+1$. The transitions are deterministic as follows: under $a_1$ the transitions are $s_0 \to s_1, s_1 \to G, s_2 \to B$, and under $a_2$ the transitions are $s_0 \to s_2, s_1 \to B, s_2 \to G.$ 

        A policy is optimal if and only if with probability $1$ we end up in state $G$. Let $\pi^1 \in \PiS$ be defined as $\pi^1(s_0,a_1)=\pi^1(s_1,a_1) =\pi^1(s_2,a_1) = 1$. The path followed under $\pi^1$ is $s_0 \to s_1 \to G$ and $s_2$ is never visited (even though we take a suboptimal action at $s_2$). Hence $\pi^1$ is optimal. Let $\pi^2 \in \PiS$ be defined as $\pi^2(s_0,a_2)=\pi^2(s_1,a_2) = \pi^2(s_2,a_2) = 1$. The path followed under $\pi^2$ is $s_0 \to s_2 \to G$ and $s_1$ is never visited (even though we take a suboptimal action at $s_1$). However it is clear that $\frac{1}{2}\pi^1 + \frac{1}{2}\pi^2$, the policy that plays $a_1$ and $a_2$ with probability $1/2$ at $s_0,s_1$ and $s_2$, is strictly suboptimal, because there is a positive probability (here, $1/2$) that we end up in state $B$ under this policy.
    \end{enumerate}
\end{proof}
Recall that a function $f: \R^{n} \to \R$ is said to be {\em quasi-convex} if its sublevel set $\{x \in \R^{n} \; | \; f(x) \leq \alpha \}$ is convex for each $\alpha \in \R$ (this is equivalent to saying that the maximum of $f$ over any segment is attained at an endpoint of this segment). The second statement in \Cref{lem: lec 0 convex sets unf policy} implies that $\pi \mapsto R(\pi,P)$ is not a quasi-convex function in general.
\begin{corollary}\label{cor: lec 0 return not quasi cvx}
    The function $\pi \mapsto R(\pi,P)$ may not be quasi-convex over $\PiS$.
\end{corollary}
\begin{proof}
    Take $\alpha = -\max_{\pi \in \PiS} R(\pi,P)$ and consider the MDP from the second statement in \Cref{lem: lec 0 convex sets unf policy} but with a minus sign for the rewards. Then the sublevel set associated with $\alpha$ for the return in this ``flipped" MDP instance is the set of optimal policies, which is not convex, hence the return is not quasi-convex.
\end{proof}
\subsection{Algorithms and complexity.}\label{sec: lec 0 algo}
\Cref{th: lec 0 mdp main} shows how to compute an optimal policy as soon as we know $v^{P} $, the fixed point of $T^{P}$. An interesting observation is that computing $v^{P} $ requires (a priori) solving a non-linear program (computing the fixed point of the non-linear operator $T^P$). However the same theorem tells us that once we know $v^{P} $ we can compute a uniformly optimal policy $\pi \in \PiSD$ by choosing, for each $s \in \cS$, an action in
\[\arg \max_{a \in \cA} P(s,a)\tr \left( r(s,a) + \gamma v^P \right).\]
The classical algorithms for solving MDPs are variants of the following three methods for computing the fixed point of the Bellman optimality operator $T^P$: value iteration, policy iteration, and linear programming. We do not cover (yet?) policy mirror descent in these notes.
\paragraph{Value iteration.}
Value iteration is based on the contraction property of $T^P$. 
\begin{algorithm}
\caption{Value iteration for nominal MDPs}\label{alg: lec 0 VI}
\begin{algorithmic}[1]
\State Initialize $v^{0} \in \R^{\cS}$.
\For{$t \in \N$}
    \State $v^{t+1} = T^{P} v^{t}$
\EndFor
\end{algorithmic}
\end{algorithm}
The convergence rate of \Cref{alg: lec 0 VI} is given next. It is a straightforward corollary of a more general version of Value Iteration, given in Appendix~\ref{app: algorithms}.
\begin{proposition}\label{prop: lec 0 VI rate}
    Let $(v^t)_{t \in \N}$ be the sequence of vectors generated by \Cref{alg: lec 0 VI}.
    \begin{enumerate}
        \item The sequence $(v^t)_{t \in \N}$ converges to $v^{P} $ at a geometric rate of $\gamma$:
        \[\|v^{P}  - v^{t}\|_{\infty} \leq \gamma^t \|v^{P}  - v^{0}\|_{\infty}\]
    \item Let $\epsilon>0$. Assume that $\|v^{P}  - v^{t}\|_{\infty} \leq \epsilon$ for some $t \in \N$ and let $\pi \in \PiS$ such that $T^{\pi,P}v^t = T^{P}v^t$. Then $\pi$ is a uniformly $2\epsilon \gamma (1-\gamma)^{-1}$-optimal policy: $\|v^{P}-v^{\pi,P}  \|_{\infty} \leq 2\epsilon \gamma (1-\gamma)^{-1}$.

    Additionally, we can check that $\|v^P - v^{t}\|_{\infty} \leq \epsilon$ by checking $\|T^P v^{t}-v^{t}\|_{\infty} \leq \epsilon (1-\gamma)$, and that last inequality holds after at most $\lfloor \frac{1}{1-\gamma} \log\left(\frac{\|v^1-v^0\|_{\infty}}{\epsilon (1-\gamma)}\right) \rfloor + 1$ iterations of \Cref{alg: lec 0 VI}.
    \end{enumerate}
\end{proposition}
Note that the second item of the above proposition provides an implementable stopping criterion for \Cref{alg: lec 0 VI} to ensure that we obtain a uniformly $\epsilon$-optimal policy.
\paragraph{Policy iteration.}
{\em Policy iteration} is an exact method that terminates with an optimal policy (whereas value iteration needs to be stopped and returns an $\epsilon$-optimal policy, for some $\epsilon>0$ that we can control). In policy iteration we generate a sequence of stationary deterministic policies $(\pi^t)_{t \in \N}$, where the current policy $\pi^t$ greedily improves upon the previous one $\pi^{t-1}$. A natural stopping criterion is to stop as soon as we can choose $\pi^{t} = \pi^{t-1}$.
\begin{algorithm}
\caption{Policy iteration for nominal MDPs}\label{alg: lec 0 PI}
\begin{algorithmic}[1]
\State Initialize $\pi^0 \in \PiSD$, $v^{0} = v^{\pi^{0},P}$
\For{$t \in \N$}
    \State Choose $\pi^{t+1} \in \PiSD$ such that $T^{P}v^{t}(s) = T^{\pi^{t+1},P}v^{t}(s)$ for all $s \in \cS$ (and equal to $\pi^t$ if we can)
    \State Choose $v^{t+1} = v^{\pi^{t+1},P}$
    \If{$\pi^{t+1}=\pi^{t}$}{ return $\pi^{t+1}$}
    \EndIf
\EndFor
\end{algorithmic}
\end{algorithm}
Policy iteration can be interpreted as a ``block simplex" algorithm on the linear programming formulation of MDPs (see next algorithm). The proof of the next proposition is given (in more generality) in Appendix~\ref{app: algorithms}.
The convergence of \Cref{alg: lec 0 PI} is a bit longer to prove than for \Cref{alg: lec 0 VI} but uses the same ideas; again we only use the monotonicity, contractivity, and attainability properties of the family of operators $\{T^{\pi,P} \; | \; \pi \in \PiSD\}$. More refined convergence analysis exists, based on analyzing the connection between policy iteration and (block) simplex algorithms for the next linear programming formulation, see \Cref{sec: lec 0 ref}.
\begin{proposition}\label{prop: lec 0 PI rate}
 Let $(\pi^t)_{t \in \N}$ be the sequence of policies generated by~\Cref{alg: lec 0 PI}. Then 
         $\|v^{P}-v^{\pi^t,P}\|_{\infty} \leq \gamma^t \|v^{P} -v^{\pi^0,P}\|_{\infty}.$
\end{proposition}
\paragraph{Linear programming and the contraction lemma.}
The linear programming formulation of MDPs is based on the following {\em contraction lemma}, which shows that the fixed point of a monotone contraction is the unique solution to minimizing a componentwise increasing function over the set of super-solutions.
\begin{lemma}[Contraction lemma]\label{lem: lec 0 contraction lemma}
    Let $F:\R^{\cS} \to \R^{\cS}$ be a monotone contraction with fixed point $v^{P}  \in \R^{\cS}$. Then 
    \[\{v^{P} \} = \arg \min \left\{ \sum_{s \in \cS} v(s) \; \middle| \; v \geq Fv, v \in \R^{\cS}\right\}.\]
\end{lemma}
\begin{proof}
    By definition $v^{P} $ is a feasible solution of the optimization problem on the right-hand side of the above equation. Additionally for $v \in \R^{\cS}$ feasible we get, using monotonicity, that $v \geq Fv \geq F^t v $ and since $\lim_{t \to + \infty} F^t v = v^{P} $ we get that $v \geq v^{P} $.
\end{proof}
 Applying the contraction lemma to the Bellman operator $T^P$ yields the following proposition.
\begin{proposition}\label{prop: lec 0 LP MDP}
    The optimal value function $v^{P} $ of an MDP is the unique optimal solution to the following linear program:
    \begin{equation}\label{eq: lec 0 LP MDP}
        \min \left\{ \sum_{s \in \cS} v(s) \; \middle| \; v(s) \geq P(s,a) \tr (r(s,a) + \gamma v), \forall \; (s,a) \in \cS \times \cA, v \in \R^{\cS}\right\}
    \end{equation}
\end{proposition}
\begin{proof}
    Immediate from the contraction lemma and the fact that $ v(s) \geq \max_{a \in \cA} P(s,a) \tr (r(s,a) + \gamma v), \forall\; s \in \cS \iff v(s) \geq P(s,a) \tr (r(s,a) + \gamma v), \forall\; (s,a) \in \cS \times \cA $.
\end{proof}
Note that we use the objective function $v \mapsto \sum_{s \in \cS} v(s)$. It would be more natural to use $v \mapsto p_0 \tr v$ but if some coefficients of $p_0$ are equal to $0$ then there may be other optimal solutions than $v^{P} $, from which we wouldn't necessarily recover a uniformly optimal policy. In fact we may not recover an optimal policy at all if $p_0$ does not have full support and we just consider an arbitrary optimal solution to $\min \{ p_0 \tr v \;  | \;  v \geq T^{P}v \}$, as the next example shows.
\begin{continuance}{ex: lec 0 ex 2}
    Consider \Cref{ex: lec 0 ex 2} with $p=0,p'=1/4,\gamma=1/2,p_0(s_0)=1$ (the same parameters that we used right before \Cref{lem: lec 0 convex sets unf policy}). The linear program $\min \{ p_0 \tr v \;  | \;  v \geq T^{P}v \}$ becomes
    \begin{equation}\label{eq: lec 0 LP example 1.2}
        \begin{aligned}
        \min \; & v(s_0) \\
        & v(s_0) \geq \frac{1}{2}v(s_2)\\
        & v(s_0) \geq \frac{1}{2}v(s_1)\\
        & v(s_1) \geq \frac{1}{8}v(s_3) + \frac{3}{8}v(s_4) \\
        & v(s_1) \geq \frac{3}{8}v(s_3) + \frac{1}{8}v(s_4) \\
        & v(s_2) \geq 1, v(s_3) \geq -2, v(s_4) \geq 2, \\
        & v \in \RS.
\end{aligned}
    \end{equation}
We have already computed an optimal policy and we know that $v^P = \left( \frac{1}{2}, \frac{1}{2}, 1, -2, 2 \right)$ with optimal value $v^{P}(s_0)=\frac{1}{2}$. We already argued that $v^{P}$ is an optimal solution for the linear program \eqref{eq: lec 0 LP example 1.2}. However note that for $v^{P}$ there is a positive slack in the second constraint, involving $v(s_1)$ and $v(s_0)$. Hence we can increase $v^{P}(s_1)$ and still get an optimal solution to \eqref{eq: lec 0 LP example 1.2} (the objective only involves $v(s_0)$). In particular, $\hat{v} = \left( \frac{1}{2}, 1, 1, -2, 2 \right)$ is feasible and optimal in \eqref{eq: lec 0 LP example 1.2}, however, based on this vector, at $s_0$ action $a_1$ and $a_2$ both lead to a continuation value of $1/2$, so assume that a greedy policy for this vector chooses action $a_2$ in $s_0$; at $s_1$, $a_1$ yields a continuation value of $1/2$ and $a_2$ yields a continuation value of $-1/2$, so we would take $a_1$. But the resulting policy has return $1/4$, which is strictly suboptimal (the optimal objective is $1/2$). The issue here comes from the fact that $\hat{v}$ is not a fixed point of $T^P$ so that in particular a greedy policy for $\hat{v}$ need not have value function $\hat{v}$.
\end{continuance}
Finally, it is worth mentioning that the dual of \eqref{eq: lec 0 LP MDP} (if we use $v \mapsto p_0 \tr v$ in the objective, instead of the sum of the components of $v$) yields a linear program whose variables are the state-action frequencies $\mu  \in \cQ \subset \R_{+}^{\cS \times \cA}$, defined in \eqref{eq: lec 0 frequency polytope}:
\begin{equation}\label{eq: lec 0 LP dual mu}
    \max \{ \sum_{s \in \cS} \sum_{a \in \cA} \mu(s,a) P(s,a) \tr r(s,a) \; | \; \mu \in \cQ\}.
\end{equation}
So we can obtain two linear programs to solve MDPs: one with variables $v \in \R^{\cS}$ (given in~\eqref{eq: lec 0 LP MDP}) and one with variables $\mu  \in \R_{+}^{\cS \times \cA}$. While we can recover optimal policies from any optimal solutions of these linear programs, to the best of our knowledge, there is no linear programming formulation with variables $\pi \in \PiS$ directly (which would be the most natural).
\paragraph{Complexity.}
Both value iteration and policy iteration require evaluating the Bellman operator $T^{P}v$ at some point $v \in \R^{\cS}$. This can be done by computing, for each $s \in \cS$, the maximum over actions in $\cA$, of values involving the scalar product of vectors in $\R^{\cS}$. This can be done in $O(|\cS|^2 |\cA|)$ arithmetic operations. Policy Iteration additionally requires computing the value function for some stationary deterministic policies at each iteration, which may be expensive, as it involves inverting a matrix of dimension $|\cS| \times |\cS|$. To avoid doing so it is common to resort to {\em modified} policy iteration, and instead of computing $v^{\pi,P}$, to just compute a few consecutive images of $T^{\pi,P}$. The linear program for MDPs can be solved by any existing methods for linear programs, such as the simplex algorithm or interior point methods. 
\section{Additional properties*}\label{sec: lec 0 additional prop}
We present here some additional properties of MDPs, which are not necessary for a first reading but provide interesting insights.
\paragraph{Behavioral policies.}
(A variant of) Kuhn's theorem asserts that to describe a policy in our sequential decision-making process it is equivalent to (A) assign a (possibly randomized) action to each finite history $(s_0,a_0,...,s_{t-1},a_{t-1},s_t)$ for each period $t$ (so-called ``behavioral form"), or (B) to directly choose a probability distribution over trajectories (compatible with the transition probabilities, e.g., we can't assign positive probabilities to sets of trajectories that contain transitions $s,a \to s'$ that are forbidden, i.e. for which $P(s,a,s')=0$). 

For instance consider \Cref{ex: lec 0 ex 2} starting from $s_0$ with probability $1$. The set of {\em all} trajectories is:
\begin{align*}
    & (s_0,a_1,s_2^{\infty}), (s_0,a_2,s_2^{\infty}) \\
    & (s_0,a_1,s_1,a_1,s_3^{\infty}), (s_0,a_1,s_1,a_2,s_3^{\infty}) \\
    & (s_0,a_2,s_1,a_1,s_3^{\infty}), (s_0,a_2,s_1,a_2,s_3^{\infty}) \\
    & (s_0,a_1,s_1,a_1,s_4^{\infty}), (s_0,a_1,s_1,a_2,s_4^{\infty}) \\
    & (s_0,a_2,s_1,a_1,s_4^{\infty}), (s_0,a_2,s_1,a_2,s_4^{\infty})
\end{align*}
where we use $s^{\infty}$ to indicate that we have reached an absorbing state (hence from there on the state $s$ is repeated in the trajectory and the actions are irrelevant).
For some choices of $p,p' \in [0,1]$ some of these trajectories must have probability $0$. For instance if $p=0$ then choosing action $a_2$ cannot lead to $s_2$ from $s_0$ hence $(s_0,a_2,s_2^{\infty})$ must have probability $0$. For the choice $p=p'=0$ the set of feasible trajectories is given by
\begin{align*}
    & h_1 \coloneqq (s_0,a_1,s_2^{\infty}) \\
    & h_2 \coloneqq (s_0,a_2,s_1,a_1,s_4^{\infty}) \\
    & h_3 \coloneqq (s_0,a_2,s_1,a_2,s_3^{\infty})
\end{align*}
and specifying a probability distribution over these elements yields a policy. For instance choosing $\bP(h_1)=1/2,\bP(h_2)=\bP(h_3)=1/4$ is equivalent to playing $a_1$ with probability $1/2$ at $s_0$ and $s_1$. Note that we must have $1-\bP(h_1) = \bP(h_2) + \bP(h_3)$ given we play either $a_1$ {\em or} $a_2$ in $s_0$. Later in these notes, we will touch upon contexts where it is easier to reason in terms of distributions over trajectories than in terms of behavioral policies, typically when the objective function is not continuous, such as for mean-payoff objectives.
\paragraph{Adjacent policies.} When a policy $\pi \in \PiS$ is not uniformly optimal, we can always find a policy $\pi' \in \PiS$ that coincides with $\pi$ at all states but one state $\bar{s} \in \cS$, is deterministic at $\bar{s}$, and such that $v^{\pi',P} \geq v^{\pi,P}$ with $v^{\pi',P}(\bar{s}) > v^{\pi,P}(\bar{s})$. Therefore, we can ``strictly improve" upon $\pi$ by changing the action chosen at only one state; this is the basis of a ``lazy" version of policy iteration, which updates the running policy at only one state (instead of potentially changing the running policy entirely as in \Cref{alg: lec 0 PI}). We will prove a stronger (robust) version of this result in \Cref{lem: lec 2 adjacent policies}. It is worth emphasizing that the same statement for {\em return} (instead of value functions) is false! That is, if a policy $\pi \in \PiS$ is such that $R(\pi,P) < \max_{\pi' \in \PiS} R(\pi',P)$, to improve $R(\pi,P)$ we may need to change the actions in at least two states, as shown in the next example.
\begin{example}\label{ex: lec 0 ex 3}
    Consider the MDP instance described in \Cref{fig: lec 0 ex 3}, and assume that we start at $s_0$ and that $\pi$ chooses action $a_2$ at $s_0$ and $s_1$. Then if we change only the action chosen at $s_0$ or only the action chosen at $s_1$ we would still get a return of $0$. Here to get a positive return we need to take $a_1$ in {\em both} $s_0$ and $s_1$.
\begin{figure}[h!]
    \begin{center}
    \begin{subfigure}{0.45\textwidth}
    \begin{center}
        \includegraphics[width=\linewidth]{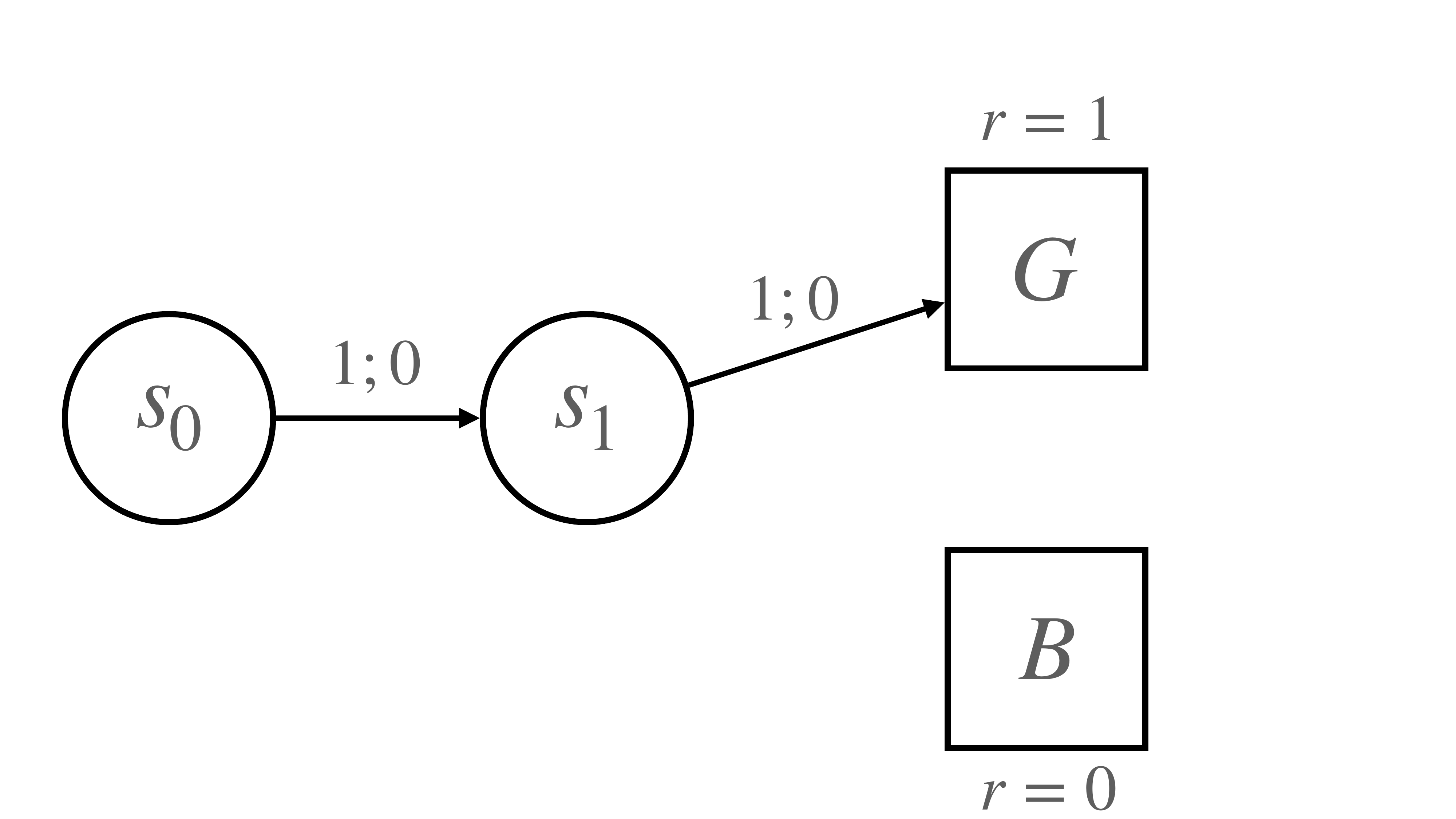}
    \end{center}
    \caption{Transitions and rewards for action $a_1$}
    \label{fig: lec 0 ex 3 a1}
    \end{subfigure}
    \begin{subfigure}{0.45\textwidth}
    \begin{center}
        \includegraphics[width=\linewidth]{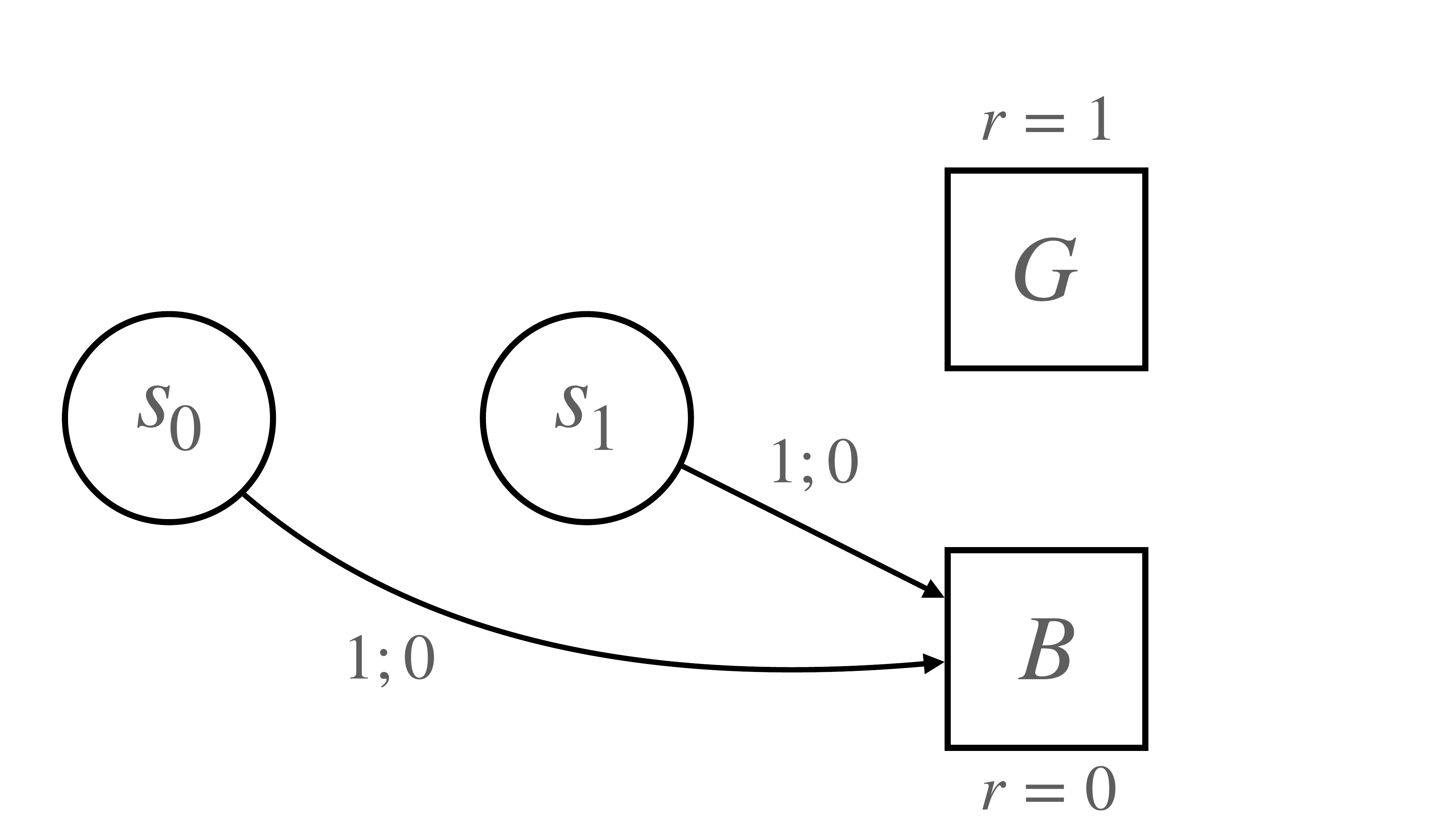}
    \end{center}
    \caption{Transitions and rewards for action $a_2$}
    \label{fig: lec 0 ex 3 a2}
    \end{subfigure}
    \end{center}
        \caption{MDP instance for \Cref{ex: lec 0 ex 3}. The transition probabilities and rewards are indicated beside each arc. Absorbing states are represented by squares.}
        \label{fig: lec 0 ex 3}
    \end{figure}
\end{example}
\paragraph{Choice of rewards.} We have introduced MDPs with instantaneous rewards $r \in \R^{\cS \times \cA \times \cS}$ indexed by the current state-action pair $(s,a)$ and the next state $s'$. For {\em nominal} MDPs, it is without loss of generality to assume that the rewards only depend on the current state-action pair $(s,a)$ and {\em not} the next state $s'$, simply by replacing the rewards by their expectation over the next states. There are two ways to prove this fact. The first way is to notice that
\begin{align*}
    T^{P}v(s) & = \max_{a \in \cA} P(s,a) \tr (r(s,a) + \gamma v) \\
    & = \max_{a \in \cA}  P(s,a) \tr r(s,a) + \gamma P(s,a) \tr v \\
    & = \max_{a \in \cA}  \tilde{r}(s,a) + \gamma P(s,a) \tr v 
\end{align*}
where in the last line we introduce $\tilde{r} \in \R^{\cS \times \cA}$ defined as $\tilde{r}(s,a) \coloneqq P(s,a) \tr r(s,a)$. Starting from an MDP with instantaneous reward $r \in \R^{\cS \times \cA \times \cS}$, we can replace it by an MDP with the same optimal policies but instantaneous rewards $\tilde{r}$ that are independent of the next state. 

A second proof is by augmenting the set of states with intermediate ``fictitious states" $(s,a,s') \in \cS \times \cA \times \cS$; rewards are obtained only at the ``state" $(s,a,s')$, so that in this model rewards only depend on the current state (and not even on the action chosen), and transitions are chosen to make this new model consistent with the original MDP; that is, $s \to (s,a,s')$ with probability $P(s,a,s')$ when choosing action $a \in \cA$ in state $s$. One also needs to replace the discount factor $\gamma$ by $\sqrt{\gamma}$ (since rewards are obtained every other period). It is worth noting here that for robust MDPs we {\em cannot} make such reductions in our choices of rewards in all generality, as will be made clear in \Cref{sec: lec 3 discussion rewards}. We also note that this reduction may not be very practical as it introduces many dummy states, which can slow down the algorithms to compute optimal policies.
\section{Main takeaways}\label{sec: lec 0 main takeaways}
We summarize the main takeaways for this chapter below:
\begin{itemize}
    \item {\em Return reformulation.} The discounted return can be modeled via two complementary approaches: one approach is based on value functions \eqref{eq:value funciton} and one approach is based on state-action frequencies \eqref{eq: lec 0 return frequencies}. The first approach leads to studying the fixed point of monotone contraction operations defined componentwise, known as Bellman operators (\eqref{eq: lec 0 T pi P} and \eqref{eq: lec 0 T P}); the second approach makes the return linear in the frequencies.
    \item {Optimality equations and computation.} An optimal policy can be chosen stationary, deterministic, and {\em uniformly} optimal, that is, optimal for all choices of initial distributions; it is characterized by the so-called Bellman equations, related to the fixed point of the Bellman operator (\Cref{th: lec 0 mdp main}). One can then use Value Iteration, Policy Iteration, and linear programming (\Cref{alg: lec 0 VI}, \Cref{alg: lec 0 PI}, \eqref{eq: lec 0 LP MDP}) to compute (approximate) fixed points of the Bellman operators.
\end{itemize}
\section{Exercises}\label{sec: lec 0 fdt}
The following is a list of observations/facts that help understand MDPs beyond what is presented above. They are presented in the same order as the objects introduced above. 
\begin{enumerate}
    \item {\em Definitions of the discounted return.} Recall that we have defined the discounted return as~\eqref{eq:discount objective mdp}. 
    \begin{enumerate}
        \item Justify that the limit exists in $ \lim_{T \to + \infty}  \E^{\pi,P}\left[ \sum_{t=0}^{T} \gamma^t r(s_t,a_t,s_{t+1})\right]$.
        \item Justify that \[\E^{\pi,P}\left[ \sum_{t=0}^{+\infty} \gamma^t r(s_t,a_t,s_{t+1}) \right] =  \sum_{t=0}^{+\infty} \gamma^t \E^{\pi,P}\left[r(s_t,a_t,s_{t+1}) \right]\]
and explain the different interpretations for the terms on each side of the above equation. 
        \item What would change if we introduced state-dependent discount factors $\gamma(s) \in [0,1)$ for each state $s \in \cS$?
    \end{enumerate}  
    \item {\em Supremum is attained.} Give the topological arguments to claim that the supremum in~\eqref{eq: lec 0 mdp} is always attained (over $\PiH$) (of course, do not use \Cref{th: lec 0 mdp main}).
    \item {\em Compact heterogeneous action sets.} Go through the main results in these notes and justify that \Cref{th: lec 0 mdp main} remains true when considering compact heterogeneous action sets, i.e., when the action set at each state $s \in \cS$ is $\cA_s$, a compact action set (not necessarily finite), and that rewards and transition probabilities are continuous in the actions.
    \item {\em Operator properties.}
Let $F:\R^{\cS} \to \R^{\cS}$ be an operator.
\begin{enumerate}
    \item {\em Relaxed value iteration / Krasnoselskii-Mann.} Let $F_{\alpha} = \alpha I + (1-\alpha) F$ for $\alpha \in \R$. Assume that $F$ is a monotone contraction. Show that $F_{\alpha}$ is a monotone contraction for $\alpha \in [0,1)$. What happens for $\alpha$ outside this interval? Compare the fixed point of $F_{\alpha}$ with that of $F$. Can we obtain better convergence rates for \Cref{alg: lec 0 VI} by considering $F_{\alpha}$ instead of $F$?
    \item {\em Crandall-Tartar theorem.} We say that $F$ is additively homogeneous if $F(x + \lambda 1_{\cS}) = Fx + \lambda 1_{\cS}, \forall \; \lambda \in \R, \forall \; x \in \R^{\cS}$, with $1_{\cS}= (1,...,1) \in \R^{\cS}$. We say $F$ is non-expansive if $\|Fx-Fy\|_{\infty} \leq \|x-y\|_{\infty}$. Note $F$ non-expansive implies that $v \mapsto F(\gamma v)$ is a contraction (hence non-expansiveness is a natural notion to study in our area).

Show the following theorem: assume that $F$ is additively homogeneous. Then $F$ is monotone if and only if $F$ is non-expansive.

Hint: use the following characterization of the $\ell_{\infty}$-norm for $x \in \R^{\cS}$: \[\|x \|_{\infty} = \inf \{ \lambda \in \R \; | \; - \lambda  \leq x(s) \leq \lambda, \forall \; s \in \cS\}.\]
\end{enumerate}
    \item {\em About discount factors.} Justify that we can interpret the quantity $1-\gamma$ as a probability of ``quitting the system with reward $0$ onward" (this is the original model for MDPs back to the 1950s). Can the results in this section be adapted to state-dependent discount factors?     
    \item {\em About the contraction lemma.}
    \begin{enumerate}
        \item Let $F$ be a monotone contraction with fixed point $v^{P} $. Show that $v^{P} $ is the maximum of the sub-solutions:
\[\{v^{P} \} = \arg \max \{ \sum_{s \in \cS} v(s) \; | \; v \leq Fv, v \in \R^{\cS}\}.\]
Is this a better formulation than the minimization program already given in the contraction lemma? 

Show that $\min \{ p_0 \tr v \; | \; v \geq Fv, v \in \R^{\cS}\} = p_0 \tr v^P$ for any $p_0 \in \DS$  and that $v^P$ is the unique minimizer if $p_0$ has full support. What are the other minimizers otherwise?
\item Let $\pi \in \PiS$. Write the contraction lemma applied to the monotone contraction $T^{\pi,P}$. Does this yield a ``fast" algorithm for computing the fixed point $v^{\pi,P}$?
    \end{enumerate}

\end{enumerate}
\section{References}\label{sec: lec 0 ref}
MDPs date back to the late 1940s-1950s where they were used to model operations research problems (typically how to manage inventories over time; this particular example comes back very often in papers from that period). 
There are so many seminal references from the 1950s to early 1980s that I do not attempt to review all of them here. For historical background on MDPs I refer to Chapter 1 and Chapter 2 in~\cite{feinberg2012handbook}, and the ``Bibliographic Remarks" at the end of each chapter in~\cite{puterman2014markov}. For a good modern introduction to MDPs I recommend the lecture notes by Joon Kwon~\cite{kwon2026introduction}, which are available on GitHub. The lecture notes in \cite{kallenberg2011markov} provide a comprehensive overview of MDPs.
Bellman's book~\cite{Bellman1957} and Howard's book~\cite{howard1960dynamic} played an important role and Bellman is credited with coining the term ``dynamic programming". Shapley's work on stochastic games~\cite{shapley1953stochastic} introduced a two-player version of MDPs.
\cite{d1960probleme} introduced the linear programming formulation for MDPs with discounted objectives (both the primal and the dual formulations; at the same time \cite{de1960problemes,manne1960linear} introduced a linear programming formulation for irreducible MDPs with a mean-payoff objective).
We refer to~\cite{ye2011simplex,scherrer2013improved} for a refined complexity analysis of policy iteration. The contraction lemma we give here is from~\cite{nilim2005robust} but it is already present in~\cite{d1960probleme}. The Crandall-Tartar theorem can be found in \cite{crandall1980some}. Kuhn's theorem is from~\cite{kuhn1953extensive}.

\newpage
\chapter{Robust MDPs with general uncertainty sets}\label{chap: lec 1 robust MDP}
This chapter introduces the concept of uncertainty set, interpreted as a set of all plausible transition probabilities (e.g., consistent with some historical dataset) and describes the robust MDP problem. The tractability properties of worst-case kernels and optimal policies are then derived for general uncertainty sets. 
\section{Main model and assumptions}
In many real-world setups, we may face some uncertainty on the ``true" value of the transition parameters. This is the case for instance if we only have access to a finite dataset of historical transitions. What should we do in this case? A natural approach would be to compute $\hat{P}$ the empirical average of the observed transitions and compute an optimal policy for the MDP with transition probabilities $\hat{P}$. However, we know that the value functions may depend on the entries of $\hat{P}$ in a rather complicated manner (see the formula in~\eqref{eq: lec 0 value function formula}), and it is unclear if a policy that performs well for $\hat{P}$ would perform well for other transition probabilities, even close to $\hat{P}$.

Hence in our search for a ``good" policy we may want to incorporate additional information, for instance we may want to take into account the confidence intervals that one can construct around the empirical average $\hat{P}$, and search for a policy with ``good" performance for any realization of the transition probabilities in these confidence intervals. For instance we may want to consider the set of transition probabilities that are $\epsilon$-close to $\hat{P}$:
\[\{P \in \Delta(\cS)^{\cS \times \cA} \; | \; |P(s,a,s') - \hat{P}(s,a,s') | \leq \epsilon, \forall \; (s,a,s') \in \cS \times \cA \times \cS\}\]
and we could adjust the radius $\epsilon$ to each triplet $(s,a,s')$ depending on the number of observations for the transition $(s,a) \to s'$.

More generally, we may want to construct a set $\cP \subset \Delta(\cS)^{\cS \times \cA}$ of feasible transition probabilities, called the {\em uncertainty set}, which captures all the plausible realizations of the transitions (plausible given the dataset we have). The {\em robust MDP problem} is then to solve the following optimization problem:
\begin{equation}\label{eq: lec 1 rmdp}
    \sup_{\pi \in \PiH} \inf_{P \in \cP} R(\pi,P)
\end{equation}
An optimal solution of~\eqref{eq: lec 1 rmdp} (provided it exists) guarantees the highest worst-case return possible across all realizations of $P \in \cP$; a uniformly optimal policy is optimal for all initial distributions. This robust model is also related to percentile guarantees (see~\Cref{sec: lec 1 fdt}). We emphasize that there are {\em many} other models that can be used to cope with parameter uncertainty (e.g. multi-objective models, partially observable MDPs, percentile optimization, distributionally robust MDPs, etc.), many of these problems having a deep connection to robust MDPs as we will show in~\Cref{chap: lec 3 connection}.

For $\pi \in \PiH$, the {\em policy evaluation problem} asks to evaluate the worst-case return of $\pi$ across all transition probabilities in $\cP$ and is defined as follows:
\begin{equation}\label{eq: lec 1 policy evaluation}
    \inf_{P \in \cP} R(\pi,P)
\end{equation}
It is common to interpret~\eqref{eq: lec 1 rmdp} as a ``game" (loosely defined for now) where, before the dynamics start, we first choose a policy $\pi \in \PiH$ and an ``adversary" chooses the transition probabilities $P \in \cP$; the adversary is also sometimes called ``Nature". The game is zero-sum, in that we are maximizing the objective and the adversary is minimizing it. Note that there is an asymmetry in the decisions of the first player (choosing $\pi \in \PiH$, i.e. eventually changing their decision rules at each period, depending on the realization of the current trajectory) and the second player (choosing $P \in \cP$ which is independent of time). We will treat the case of time-varying transition probabilities in \Cref{sec: lec 3 time var P}.

The rest of this chapter is devoted to understanding the fundamental properties (existence and characterization) of the optimal solutions in the robust MDP problem~\eqref{eq: lec 1 rmdp} and the policy evaluation problem~\eqref{eq: lec 1 policy evaluation}. In this chapter, we limit the assumptions to what is strictly necessary for stating our results (in particular, many of the existence results for optimal policies don't even need compactness of $\cP$); in the next chapter, we'll make additional assumptions to simplify the proofs. The topological results needed below (e.g. Weierstrass' theorem or Tykhonoff's theorem) are given and explained in Appendix~\ref{app:sec analysis}.
\begin{remark}
    In these notes I mostly focus on the case of uncertain transition probabilities and known rewards. The case of uncertain rewards is somewhat easier (at least for the main results) and is treated in~\Cref{chap: lec 3 connection}.
\end{remark}
\section{Existence of optimal solutions}
In this section we show that under mild technical conditions, worst-case kernels and optimal policies exist for the robust MDP problem. Here we'll always assume that $\cP \neq \emptyset$. We start by studying the continuity properties of the return.
\begin{proposition}\label{prop: lec 1 continuity}
    Let $\cP \subset \Delta(\cS)^{\cS \times \cA}$.
    \begin{enumerate}
        \item The function $(\pi,P) \mapsto R(\pi,P)$ is continuous over $\PiH \times \cP$.
        \item The function $\pi \mapsto \inf_{P \in \cP} R(\pi,P)$ is upper semicontinuous over $\PiH$.
    \end{enumerate}
\end{proposition}
\begin{proof}
\begin{enumerate}
    \item    Let $r_{\infty} = \max_{s,a,s'} |r(s,a,s')|$. We have $r_{\infty} < + \infty$ since $\cS$ and $\cA$ are finite sets. For $H \in \N$, consider the truncated return $R_{H}(\pi,P) = \E^{\pi,P}\left[ \sum_{t=0}^{H-1} \gamma^t r(s_t,a_t,s_{t+1})\right]$. The function $(\pi,P) \mapsto R_{H}(\pi,P)$ is a polynomial in finitely many entries of $\pi \in \PiH$ and $P \in \cP$ since it depends only on the distributions induced by $\pi,P$ over the first $H+1$ states visited and first $H$ actions chosen. Moreover we have
    \[|R(\pi,P) - R_{H}(\pi,P)| \leq r_{\infty} \frac{\gamma^H}{1-\gamma}, \forall \; (\pi,P) \in \PiH \times \cP.\]
    Hence $(\pi,P) \mapsto R(\pi,P)$ is the uniform limit, as $H \to + \infty$, of the continuous functions $(\pi,P) \mapsto R_{H}(\pi,P)$, and therefore it is continuous.
    \item This follows because the pointwise infimum of some continuous functions is always upper semicontinuous.
\end{enumerate}
\end{proof}
We first show that the policy evaluation problem~\eqref{eq: lec 1 policy evaluation} always has a solution when $\cP$ is compact, i.e., a worst-case transition kernel always exists. 
\begin{proposition}\label{prop: lec 1 existence P worst-case}
    Let $\cP \subset \Delta(\cS)^{\cS \times \cA}$ be a compact set and $\pi \in \PiH$. Then $\inf_{P \in \cP} R(\pi,P)$ is attained.
\end{proposition}
\begin{proof}
Let $\pi \in \PiH$. By assumption $\cP$ is compact, and from the first statement in \Cref{prop: lec 1 continuity} the function $P \mapsto R(\pi,P)$ is continuous. Weierstrass' theorem states that a continuous function on a (nonempty) compact set attains its minimum, which yields the proposition.
\end{proof}
This enables us to talk later of ``worst-case transition probabilities" given a fixed policy $\pi \in \PiH$. It is worth noting that in general there may be several different worst-case transition probabilities, leading to different value functions (see \Cref{prop: lec 1 multiple value functions}).
Regarding the existence of optimal policies, we can prove the following two positive results.
First, we note that an optimal robust policy always exists, i.e. the supremum in~\eqref{eq: lec 1 rmdp} is always attained.
\begin{proposition}\label{prop: lec 1 existence piopt}
    Let $\cP \subset \Delta(\cS)^{\cS \times \cA}$. An optimal policy (eventually in $\PiH$) always exists for the robust MDP problem~\eqref{eq: lec 1 rmdp}.
\end{proposition}
\begin{proof}
The set $\PiH$ is the product over all finite trajectories $h \in \cH$ of the compact set $\Delta(\cA)$. Note that $\cH$ is countable since $\cS,\cA$ are finite sets. From Tykhonoff's theorem, $\PiH$ is a compact set for the product topology (see \Cref{rmk: lec 0 topology}). Additionally, from the second statement in \Cref{prop: lec 1 continuity} we have that $\pi \mapsto \inf_{P \in \cP} R(\pi,P)$ is an upper semicontinuous function. Weierstrass' theorem gives that it attains its supremum over a compact set.
\end{proof}
Note that the above proposition does not provide any characterization of optimal policies (e.g., we do not know if they can be chosen stationary) and does not require $\cP$ to be compact (we optimize over $\PiH$ in~\eqref{eq: lec 1 rmdp}).
The next proposition shows the existence of a policy that is optimal across all choices of {\em deterministic} initial distributions (note the difference with uniform optimality, which requires to be optimal across all choices of initial distributions).
\begin{proposition}\label{prop: lec 1 exist uniform opt policies Dirac}
Let $\cP \subset \Delta(\cS)^{\cS \times \cA}$. There exists a policy $\pi \in \PiH$ that is optimal in~\eqref{eq: lec 1 rmdp} for all choices of {\em deterministic} initial distribution.
\end{proposition}
\begin{proof}
The proof follows by patching together optimal policies for each Dirac distribution. For $s \in \cS$ let us consider $\pi^{\star,s} \in \PiH$ an optimal policy (possibly history-dependent) when starting in state $s$; the existence of such a policy follows from \Cref{prop: lec 1 existence piopt}. Now define a single policy $\bar{\pi} \in \PiH$ that coincides with $\pi^{\star,s}$ along trajectories starting in $s \in \cS$. By definition $\bar{\pi}$ is optimal across all deterministic initial distributions.
\end{proof}
\begin{remark}\label{rmk: lec 0 topology}
    Note that $\PiH$ is a set contained in an infinite-dimensional space; the rigorous definition given in \Cref{sec: lec 0 model} is $\PiH = \times_{h \in \cH} \Delta(\cA)$ where $\cH$ is the (infinite, countable) set of all histories of finite length. Hence in principle one needs to be careful when choosing a topology for defining the notion of compactness and continuity; the product topology suffices for our purpose (and Tykhonoff's theorem gives compactness of $\PiH$ as a Cartesian product of compact sets).
\end{remark}
\section{Properties of optimal solutions}\label{sec: lec 1 properties optim sol}
We now analyze the properties of optimal policies for the robust MDP problem (the existence of which we proved in~\Cref{prop: lec 1 existence piopt}).
In all generality, we first note that we may need non-stationary optimal policies. This is in contrast to the case of nominal MDPs studied in~\Cref{chap: lec 0 mdp}.
\begin{proposition}\label{prop: lec 1 nonstationary opt policies}
    For a general uncertainty set $\cP$, an optimal policy for the robust MDP problem~\eqref{eq: lec 1 rmdp} may need to be history-dependent (non-stationary and non-Markovian).
\end{proposition}
The proof of~\Cref{prop: lec 1 nonstationary opt policies} follows from analyzing the optimal policies in the MDP instance from \Cref{ex: lec 0 ex 1}.
\begin{continuance}{ex: lec 0 ex 1}
Consider the MDP instance from Figure~\ref{fig: lec 0 ex 1} (starting in $s_0$) and consider the uncertainty set of cardinality $2$, corresponding to the transition matrices induced by the choice $p=0$ or by the choice $p=1$. We first compute the payoff of the best stationary policies. A stationary policy $\pi \in \PiS$ is determined by $x \in [0,1]$ the probability that $\pi$ chooses action $a_1$ in state $s_0$, and we computed the discounted return of such a stationary policy as (for fixed $p$) $\frac{\gamma}{1-\gamma^2} (xp + (1-x)(1-p))$. Hence, the best stationary policy solves
\[\max_{x \in [0,1]} \min_{p \in \{0,1\}} xp + (1-x)(1-p) = \max_{x \in [0,1]} \min 
\{x,1-x\}.\] 
From this, we get that the best stationary policy chooses $x=1/2$, for a discounted return of $\frac{1}{2}\frac{\gamma}{1-\gamma^2}$. Now consider the following choice of non-stationary policy. We choose $a_1$ with probability $1/2$ at period $0$ (to match what the best stationary policy does). Then we get back to state $s_0$ at period $2$, and what we do depends on the observed trajectory so far (below we omit the actions chosen at $G$ and $B$ since they are irrelevant):
\begin{itemize}
    \item If we observe $(s_0,a_1,G,s_0)$: we play $a_1$ forever onward.
    \item If we observe $(s_0,a_1,B,s_0)$: we play $a_2$ forever onward.
    \item If we observe $(s_0,a_2,G,s_0)$: we play $a_2$ forever onward.
    \item If we observe $(s_0,a_2,B,s_0)$: we play $a_1$ forever onward.
\end{itemize}
In short, the decisions outlined above check whether we played optimally during the first period (i.e. whether we reached $G$ instead of $B$ at period $t=1$) and adjust the actions chosen next accordingly.
This ensures that we get the same reward as the best stationary policy in the first period ($t=0$), and then we'll always visit the good state $G$ after $s_0$, so we get the maximum reward possible starting from period $t=2$, i.e. we get a discounted return of $\gamma(\frac{1}{2} +  \frac{\gamma^2}{1-\gamma^2})$. Therefore this non-stationary policy strictly outperforms the best stationary policy. Note that the worst-case transition probabilities here are attained by both choices of $p=0$ and $p=1$, since at period $t=0$ the choice of $p$ doesn't matter (we play uniformly) and for later periods we always visit $G$ instead of $B$ no matter the choice of $p$.
\end{continuance}
Intuitively, non-stationary policies beat stationary policies in the above example because they can exploit information revealed over time (in this example, we can determine whether the adversary chose $p=0$ or $p=1$ by observing the first transition, i.e., after observing the realized state at period $1$).
Note that in the example above, we need genuinely history-dependent policies (the optimal decisions at any period $t \geq 1$ depend on the realized transitions at time $t=0$ so optimal policies are non-Markovian). It is also worth emphasizing that if we are only interested in $\epsilon$-optimal policies, we can play {\em eventually stationary} policies, i.e., policies that become stationary for large enough periods. 
\begin{proposition}\label{prop: lec 1 eventual stationarity}
    Let $\epsilon>0$. We can find an $\epsilon$-optimal policy that is eventually stationary.
\end{proposition}
\begin{proof}
    Let $\epsilon>0$. We can just truncate the return up to period $H-1$ such that the remaining discounted return obtained after the first $H$ periods is always smaller than $\epsilon$. We can find such $H \in \N$ because the rewards are bounded and discounted by some $\gamma \in [0,1)$. Now consider the policy that plays like an optimal (potentially history-dependent) policy for the first $H$ periods, then always play action $a_1$ afterwards (or any other stationary policy). By construction this is an $\epsilon$-optimal policy.
\end{proof}
Note that from the above construction, as $\epsilon$ goes to $0$, the period at which we can start to play stationary while guaranteeing $\epsilon$-optimality increases to $+\infty$, which is not very practical. 
Additionally, there need not exist a {\em uniformly} optimal policy.
\begin{proposition}\label{prop: lec 1 non exist uniform opt policies}
    For a general uncertainty set $\cP$, there may not exist any uniformly optimal policy that is optimal for the robust MDP problem~\eqref{eq: lec 1 rmdp}, even among history-dependent policies.
\end{proposition}
The proof is based on an MDP close to \Cref{ex: lec 0 ex 2}.
\begin{continuance}{ex: lec 0 ex 2}
    Consider an MDP instance with the same transition probabilities as in Figure~\ref{fig: lec 0 ex 2}, and take $p=p',p \in [0,1]$. For the rewards, choose the following intermediate rewards: for $s' \in \{s_1,s_2\}$, $r(s_0,a_1,s_2) = 1, r(s_0,a_2,s')=-1, r(s_1,a_2,s_3)=-1, r(s_1,a_2,s_4)=1$, and all other rewards are $0$. We choose $\gamma=0$ so that only the action chosen at period $t=0$ impacts the return (and in particular when starting from $s_0$ or $s_1$ the rewards obtained in the absorbing states $s_2,s_3,s_4$ are irrelevant). Similarly as in \Cref{ex: lec 0 ex 2}, stationary policies $\pi \in \PiS$ can be parametrized by $x \in [0,1]$ the probability to play $a_1$ in $s_0$ and $y \in [0,1]$ the probability to play $a_1$ in $s_1$. Because $\gamma=0$, stationary policies and history-dependent policies induce the same set of returns, and we can focus on searching for optimal policies among stationary ones.

    Let us compute the set of all optimal policies when starting in $s_1$ with probability $1$. Action $a_1$ gives a return of $0$ while action $a_2$ gives a return of $(1-p)\cdot (-1) + p \cdot 1 = 2p-1.$ Hence the worst-case return of $\pi$ is $\inf_{p \in [0,1]} (1-y)(2p-1) = y-1$ attained at $p=0$ so that optimality requires $y=1$. Similarly it is easy to see that starting in $s_0$ the decision maker would choose $a_1$ and the adversary would choose $p=1$.

    We now compute the set of all optimal policies when starting in $s_0$ or in $s_1$ with probability $1/2$. A few lines of algebra show that 
    \[R(\pi,P) = \frac{1}{2}\left( -1 + x(2-p) + (1-y)(2p-1) \right) \]
    This expression is affine in $p \in [0,1]$ hence the worst-case over $[0,1]$ is always attained in $\{0,1\}$. Therefore 
    \[\inf_{P \in \cP} R(\pi,P) = \min \lbrace x + \frac{y}{2}-1,\frac{x-y}{2} \rbrace\]
    This expression always increases in $x \in [0,1]$ hence a maximizer must take $x=1$ for which the objective becomes $\frac{1}{2} \min \{y,1-y\}$, maximized at $y=1/2$. Thus the sets of optimal policies when starting in $s_1$ with probability $1$ or in $\{s_0,s_1\}$ uniformly are disjoint, and a uniformly optimal policy cannot exist in this example.
\end{continuance}
The previous proposition should be contrasted with \Cref{prop: lec 1 exist uniform opt policies Dirac}, which shows that we can find policies that are optimal across all {\em deterministic} initial distributions; it is the fact that initial distributions may be randomized that prevents the existence of uniformly optimal policies in all generality.

The key in the above example is that the same parameter $p \in [0,1]$ appears in both the transitions out of the two states $s_0$ and $s_1$. From $s_0$ we want to choose $a_1$ and the adversary wants to choose $p=1$; for $s_1$ we want to choose $a_1$ and the adversary wants to choose $p=0$. But if we start in $s_1$ or $s_0$, each with  probability $1/2$, then the adversary has ``conflicting" worst-case kernels from these two states ($p=0$ from $s_1$ and $p=1$ from $s_0$), which we can exploit with the choice $y=1/2$. 
It is perhaps also worth noting that the previous example shows, in passing, that we may need randomized (non-deterministic) policies in general (this is the case when starting in $s_0$ or $s_1$ with probability $1/2$). 

An additional potential hiccup is that optimal solutions of both~\eqref{eq: lec 1 rmdp} and~\eqref{eq: lec 1 policy evaluation} may have irrational entries, even when all the inputs (rewards, uncertainty set, discount factor) have rational (or even integer) entries. 
\begin{proposition}\label{prop: lec 1 irrational entries}
  There exists a robust MDP with two states and two actions, binary rewards, such that for a choice of discount factor $\gamma=1/2$ and $\cP$ a convex compact polytope with constraints described by integers, every stationary optimal policy has irrational entries.
\end{proposition}
The proof is based on the following example.
\begin{example}\label{ex: lec 1 irrational}
    Consider an MDP instance with states $s_0,s_1$ and actions $a_1,a_2$, represented in Figure~\ref{fig: lec 1 ex 3}. 
    \begin{figure}[h!]
    \begin{center}
    \begin{subfigure}{0.45\textwidth}
    \begin{center}
        \includegraphics[width=\linewidth]{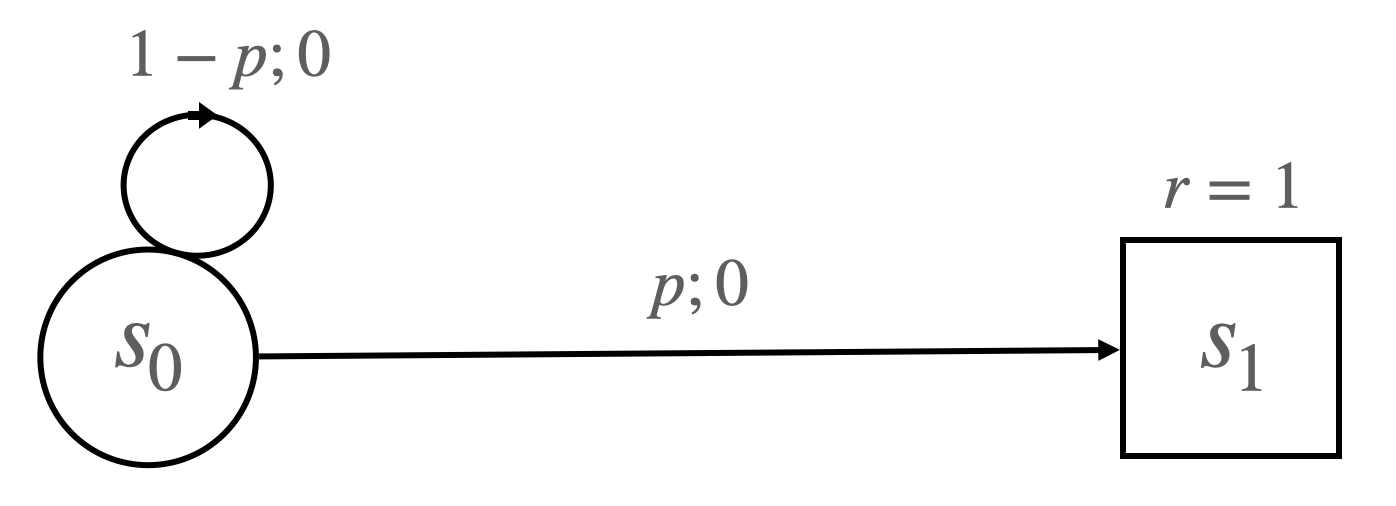}
    \end{center}
    \caption{Transitions and rewards for action $a_1$}
    \label{fig: lec 1 ex 3 a1}
    \end{subfigure}
    \begin{subfigure}{0.45\textwidth}
    \begin{center}
        \includegraphics[width=\linewidth]{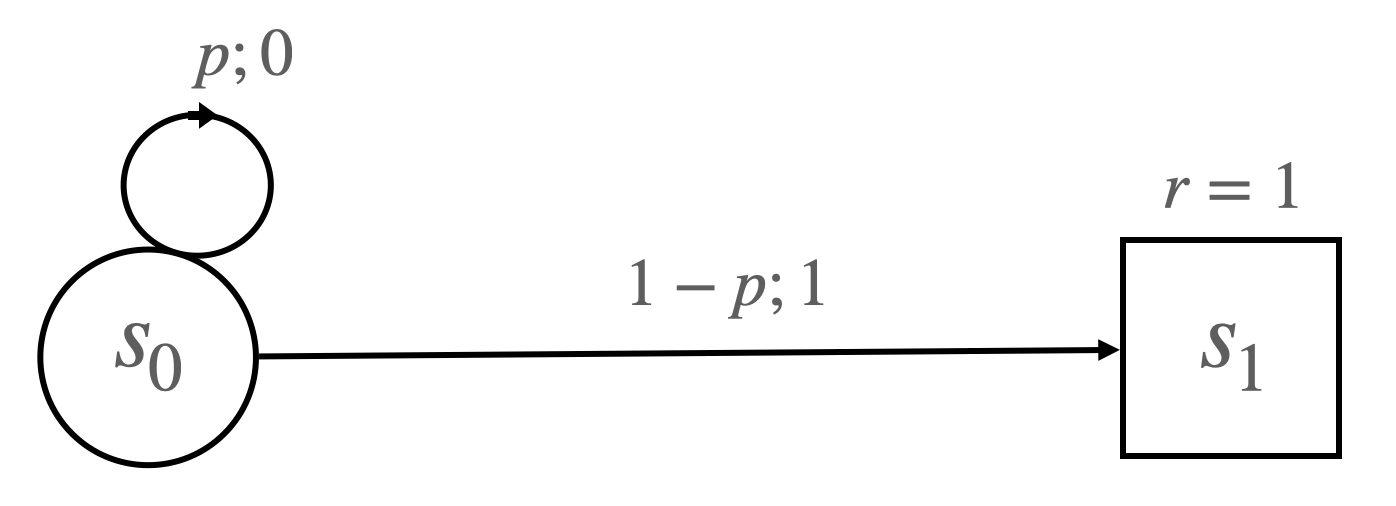}
    \end{center}
    \caption{Transitions and rewards for action $a_2$}
    \label{fig: lec 1 ex 3 a2}
    \end{subfigure}
    \end{center}
        \caption{MDP instance for \Cref{ex: lec 1 irrational}. The transition probabilities and rewards are indicated beside each arc.}
        \label{fig: lec 1 ex 3}
    \end{figure}
    For now we'll assume that an optimal policy can be chosen stationary (this will be justified in the next chapter in \Cref{th: lec 2 main rmdps}). Let $x = \pi(s_0,a_1) \in [0,1]$ and $p \in [0,1]$. Choose $\gamma=1/2$. The value function $v^{x,p}(s_0)$ is 
    \[v^{x,p}(s_0) = 2 \frac{xp + 2(1-x)(1-p)}{1+xp + (1-x)(1-p)}.\]
    This is linear-fractional in $p$ so the worst-case over $p \in [0,1]$ is attained at some $p \in \{0,1\}$. Consequently,
    \[\max_{x \in [0,1]} \min_{p \in [0,1]} v^{x,p}(s_0) = \max_{x \in [0,1]} \min \lbrace \frac{4(1-x)}{2-x}, \frac{2x}{1+x} \rbrace.\] 
    This yields a maximizer $x\opt \in [0,1]$ satisfying $4(1-x\opt)(1+x\opt) = 2x\opt (2-x\opt)$ so that $x^{\star 2} + 2x\opt - 2 =0$ i.e. $x\opt = \sqrt{3}-1 \in (0,1)$, which has irrational entries.
\end{example}
This situation is in line with the case of Nash equilibria of finite games with more than two players, as shown by Nash in his original 1951 paper, and in contrast with the case of linear programs with rational inputs (where there exist some optimal solutions with rational entries) and MDPs: they can be solved by linear programming, as we saw in the previous chapter, hence we can find optimal policies with rational entries. It suggests that even representing optimal solutions for robust MDPs may be non-trivial.
An additional difficulty comes from the fact that, for a fixed policy $\pi \in \PiS$, worst-case transition probabilities need not be unique; in the subparts of the set of states that are not visited (or whose visits are irrelevant), the adversary may choose any transition probabilities. However value functions are vector-valued functions of the policies and transition probabilities, indexed by the {\em entire} set of states. Consequently, {\em worst-case value functions need not be unique} (although of course the worst-case return is unique)! This is described in the next proposition.
\begin{proposition}\label{prop: lec 1 multiple value functions}
Let $\pi \in \PiS$. In general there may be multiple solutions in $\arg \min_{P \in \cP} R(\pi,P)$, and for $P_1,P_2 \in \arg \min_{P \in \cP} R(\pi,P)$ we may have $v^{\pi,P_1} \neq v^{\pi,P_2}$.
\end{proposition}
\begin{proof}
    It is easy to observe this phenomenon for the robust MDP instance from \Cref{fig: lec 0 ex 2}, with an uncertainty set described by $p=0,p' \in [0,1]$. Say we start in $s_0$ with probability $1$ and $\gamma \in (0,1)$. Take for instance $\pi \in \PiS$ that chooses action $a_1$ with probability $1$ in $s_0$ and with probability $3/4$ in $s_1$. With this choice of policy every kernel attains the worst-case return: $R(\pi,P) = \frac{\gamma}{2(1-\gamma)}$ for any $P \in \cP$ but we get $v^{\pi,P}(s_1) = \frac{\gamma (1-2p')}{2(1-\gamma)}$.
\end{proof}
Hence in general it is more proper to talk of the {\em set} of worst-case value functions of a policy: $\displaystyle \{v^{\pi,P} \; | \; P \in \arg \min_{P \in \cP} R(\pi,P)\}$, rather than ``the" worst-case value function of a policy (except if we proved either (a) that there is a unique worst-case kernel, or (b) that all worst-case kernels yield the same value function). One may wonder whether there is a worst-case value function of particular interest, for instance, whether we can define a notion of ``smallest" (i.e., worst) worst-case value function. We will see in the next chapter that for rectangular uncertainty sets, there is indeed a worst-case value function which is componentwise smaller than any other worst-case value function, and how to compute this particular vector.

We conclude this section by showing that even evaluating the worst-case return of a policy may be a difficult task\footnote{The concept of {\em NP-hardness} describes a large class of decision problems such that every problem in NP (intuitively: a decision problem is in NP if, if we are given a candidate solution to the decision problem, we can verify it in polynomial time) admits a polynomial-time reduction to this problem. The class of NP-hard problems is widely regarded as a class of decision problems that are ``difficult to solve".}.
\begin{proposition}\label{prop: lec 1 np hard}
Let $\pi \in \PiSD$.
    For a general convex compact uncertainty set $\cP$ and $c \in \R$, it may be NP-hard to decide if $\inf_{P \in \cP} R(\pi,P) \leq c$.
\end{proposition}
There exist several proofs of this statement (see references at the end of this chapter). I adapt here a proof  by Joel Goh and coauthors, based on the Hamiltonian cycle problem.
\begin{proof}[Proof of \Cref{prop: lec 1 np hard}]
Given a directed graph $G=(V,E)$ with non-empty edge set $E \subset V \times V$, the Hamiltonian cycle problem asks whether there exists a Hamiltonian cycle in $G$, i.e. a finite sequence of  vertices $(s_0,s_1,...,s_H,s_0)$ where $(s_i,s_{i+1}) \in E$ for $i \in \{0,...,H-1\}$, $(s_{H},s_0) \in E$, such that this cycle visits all states at most once: $|\{ s_0,s_1,...,s_H\}|=|V|$.

We can reduce this problem to comparing the optimal objective of a robust policy evaluation problem as follows. Choose $\cS = V$ and $\cA=\{a_1\}$ (i.e. there is no action to choose). Let $\pi:\cS \to \{a_1\}$ be the constant policy. Let us write $\cS = \{s_0,...,s_{H}\}$ for $H = |V|-1$. The uncertainty set is the set of bistochastic matrices compatible with the set of edges $E$:
    \[\cP = \{ P \in \R^{\cS \times \cS} \; | \; P(s,s') \in \R_{+}, P 1_{\cS} = 1_{\cS}, P\tr 1_{\cS} = 1_{\cS}, P(s,s') = 0 \text{ if } (s,s') \notin E\}\]
    with $1_{\cS} = (1,...,1) \in \R^{\cS}$.
Note that we can omit the dimension associated with $\cA$ since $|\cA|=1$. Choose the rewards as follows: there is a reward of $1$ when we are at state $s_0$ and a reward of $0$ everywhere else. Since the adversary solves $\inf_{P \in \cP} R(\pi,P)$, the adversary will choose a matrix in $\cP$ to spend as little time as possible in state $s_0$ (the only state giving us a positive reward). The initial distribution starts deterministically in $s_0$. The return $R(\pi,P)$ for $P \in \cP$ is
\[R(\pi,P) = \E^{P}\left[ \sum_{t=0}^{+\infty} \gamma^t 1_{\{ s_t = s_0\}} \right].\]
We first claim that $\inf_{P \in \cP} R(\pi,P) \geq \frac{1}{1-\gamma^{|V|}}$.
Indeed for $\tau$ the random variable modeling the return time to $s_0$ we get that 
\[R(\pi,P) = \E^{P}\left[ \sum_{t=0}^{+\infty} \gamma^{t \tau} \right] =  \sum_{t=0}^{+\infty} \E^{P}\left[ \gamma^{t \tau} \right] \geq \sum_{t=0}^{+\infty} \gamma^{\E^{P}\left[ t \tau \right]} = \frac{1}{1-\gamma^{\E^{P}\left[\tau\right]}}\]
where the second equality is from the Monotone Convergence Theorem and the inequality is from Jensen's inequality. Finally we note that we must have $\E^{P}\left[\tau\right] \leq |V|$, because the expected return time $\E^{P}\left[\tau\right]$ is the inverse of the $s_0$-component in the invariant distribution of the Markov chain induced by $P \in \cP$, and since $P$ is bistochastic (i.e. $P$ is a stochastic matrix and $P\tr$ too), we must have that this invariant distribution is uniform over the closed communicating class containing $s_0$. 
Now we claim that there exists a Hamiltonian cycle in $G$ if and only if $\inf_{P \in \cP} R(\pi,P) = \frac{1}{1-\gamma^{|V|}}$.

We first prove the easy implication.
If there is a Hamiltonian cycle in $G$, then there is a corresponding matrix $P \in \cP$ such that the dynamics is deterministic under $P$ and follows the Hamiltonian cycle. By definition for this choice of $P$ we will pass by $s_0$ every $|V|$ iterations and therefore the discounted return is $R(\pi,P) = \frac{1}{1-\gamma^{|V|}}$.

Now assume that $\inf_{P \in \cP} R(\pi,P) = \frac{1}{1-\gamma^{|V|}}$. By \Cref{prop: lec 1 existence P worst-case}, the infimum is attained at some $P\opt \in \cP$. The proof of the lower bound on $\inf_{P \in \cP} R(\pi,P)$ shows that we must have $\E^{P\opt}[\tau]=|V|$ as well as equality in Jensen's inequality, which implies the stronger property that $\tau = |V|$ almost surely. This means that starting from $s_0$ we come back to $s_0$ exactly at period $|V|$. It is not straightforward to construct a Hamiltonian cycle from $P\opt$, because $P\opt$ may include randomized transitions (the matrix is bistochastic but the entries need not be $0$ or $1$, we can have entries strictly between $0$ and $1$). However a theorem by Birkhoff says that the polytope $\cP$ admits the permutation matrices (compatible with $E$) as extreme points, i.e. we can write $P\opt = \sum_{\ell=1}^{m} \lambda_{\ell} Q^{\ell}$ with $Q^{\ell} \in \cP$ a permutation matrix for each $\ell \in \{1,...,m\}$: $Q^{\ell}_{ss'} = 1_{\{ s=\sigma^{\ell}(s')\}}$ for some permutation $\sigma^{\ell}$ of $\cS$, and $\lambda \geq 0, \sum_{\ell=1}^{m} \lambda_{\ell}=1$. Each permutation matrix decomposes into disjoint directed cycles, and for some $\ell \in \{1,...,m\}$ consider the cycle of $Q^{\ell}$ containing $s_0$. If this cycle has length strictly smaller than $|V|$ then we would have that $\bP_{P\opt}(\tau < |V|)>0$, which is impossible since $\tau = |V|$ almost surely (which exactly means $\bP_{P\opt}(\tau = |V|)=1$). Hence for each $Q^{\ell}$ the cycle containing $s_0$ induced by $Q^{\ell}$ has length $|V|$, i.e., there exists a Hamiltonian cycle.

Therefore we can verify if there exists a Hamiltonian cycle by verifying $\inf_{P \in \cP} R(\pi,P) \leq  1 / (1-\gamma^{|V|})$ (in which case we must have $\inf_{P \in \cP} R(\pi,P) =  1 / (1-\gamma^{|V|})$). This shows that solving the policy evaluation problem is NP-hard.
\end{proof}
\section{Main takeaways}
Let us summarize the main results in this section. For a general uncertainty set $\cP$:
\begin{itemize}
    \item {\em Properties of optimal policies.} Optimal policies for~\eqref{eq: lec 1 rmdp} always exist but may require being non-stationary, even when there are only two possible transition kernels (\Cref{prop: lec 1 existence piopt} and \Cref{prop: lec 1 nonstationary opt policies}). Approximately optimal policies may be chosen ``eventually stationary" but we may need longer and longer horizons to better approximate optimal policies (\Cref{prop: lec 1 eventual stationarity}). Even for instances described by simple (rational) entries, we may need irrational real values to describe stationary optimal policies (\Cref{prop: lec 1 irrational entries}).
    \item {\em Uniformly optimal policies.} Uniformly optimal policies for~\eqref{eq: lec 1 rmdp} may not exist (\Cref{prop: lec 1 non exist uniform opt policies}), unless we restrict ourselves to deterministic initial probability distributions and we play history-dependent policies (\Cref{prop: lec 1 exist uniform opt policies Dirac}).
    \item {\em Hardness and multiple worst-case kernels.} Evaluating the worst-case return of a stationary policy may be computationally challenging, even when the uncertainty set is a convex compact polytope, but worst-case transition probabilities always exist for compact uncertainty sets (\Cref{prop: lec 1 existence P worst-case} and \Cref{prop: lec 1 np hard}). Finally, there may exist multiple worst-case value functions, since worst-case transition kernels need not be unique (\Cref{prop: lec 1 multiple value functions}).
\end{itemize}
Let's take a step back. Overall it appears that solving the robust MDP problem in general is a difficult task under the various settings considered in this chapter. This is in stark contrast to the numerous favorable tractability properties of the nominal MDP setting studied in~\Cref{chap: lec 0 mdp}. At this point this discrepancy may seem surprising. Indeed, the roles of $\pi \in \PiS$ and $P \in \cP$ in the expression of the value function $v^{\pi,P}$ as in~\eqref{eq: lec 0 value function formula} are similar: $v^{\pi,P}$ depends on the transition matrix $L^{\pi,P} \in \R^{\cS \times \cS}$ of the Markov chain over $\cS$ induced by $\pi,P$, and on $r^{\pi,P} \in \R^{\cS}$ the expected instantaneous rewards. Both $L^{\pi,P}$ and $r^{\pi,P}$ are bilinear functions of $(\pi,P)$. So why is optimizing $R(\pi,P)$ over $P \in \cP$ difficult (for fixed $\pi$), while optimizing $R(\pi,P)$ over $\pi \in \PiS$ is tractable (for fixed $P$)? The main issue must come from the properties of the domains $\PiS$ and $\cP$ of each variable $\pi$ and $P$ (given the ``symmetry" of $(\pi,P) \mapsto v^{\pi,P}$ in $\pi$ and $P$). What is the fundamental structure of $\PiS$ that makes nominal MDPs tractable, and that is lacking in general $\cP$? Answering this question is the focus of the next chapter.
\section{Exercises}\label{sec: lec 1 fdt}
\begin{enumerate}
    \item {\em Worst-case transitions in the relative interior of the uncertainty sets.}
    Consider the robust MDP instance from \Cref{fig: lec 1 exo 0}, with an uncertainty set $\cP$ parametrized by $p \in [0,1]$. Use this instance to show that (optimal) worst-case transition probabilities may belong to the relative interior of $\cP$.
    \begin{figure}[h!]
    \begin{center}
    \begin{subfigure}{0.45\textwidth}
    \begin{center}
        \includegraphics[width=\linewidth]{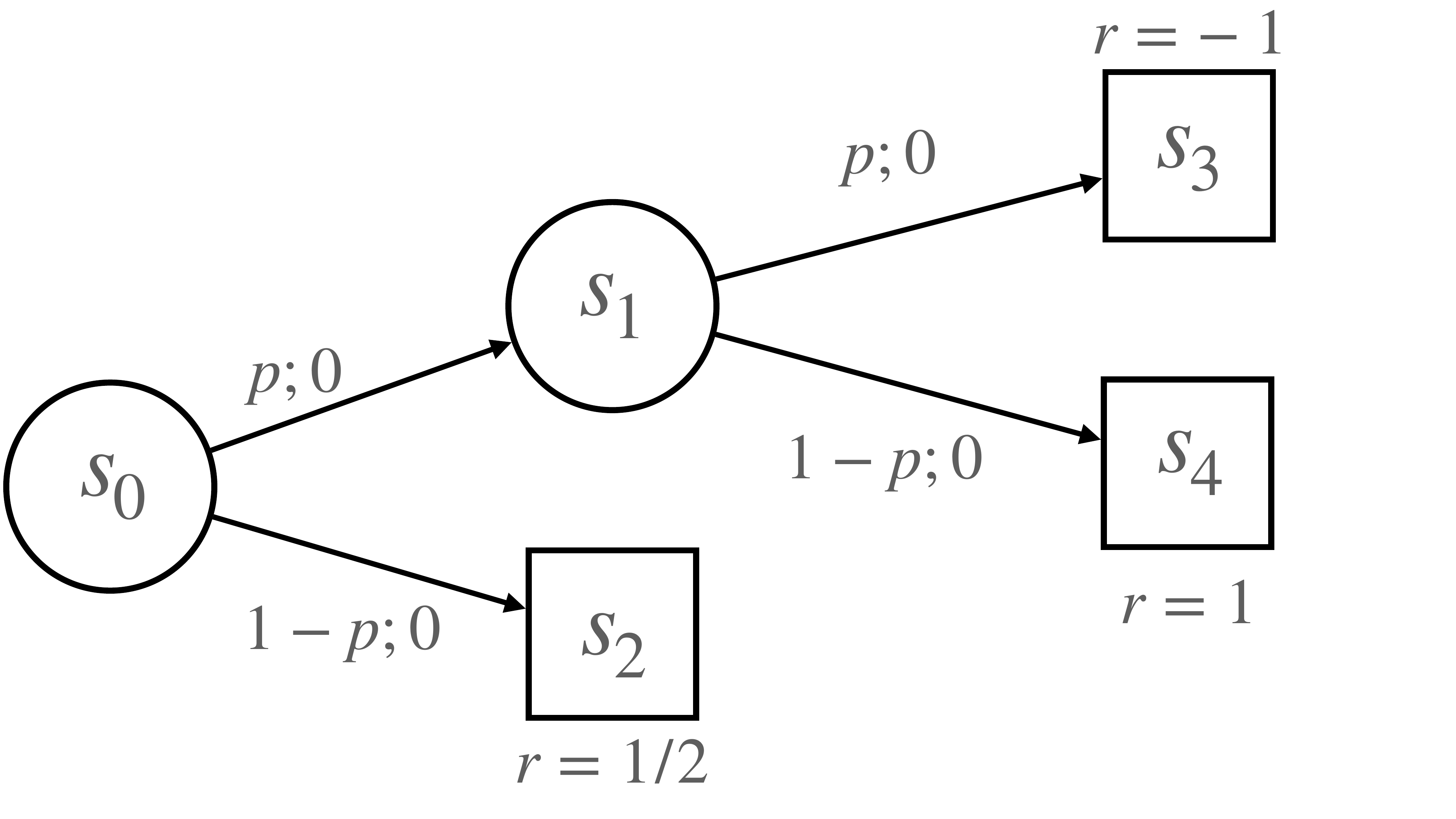}
    \end{center}
    \caption{Transitions and rewards for action $a_1$}
    \label{fig: lec 1 exo 0 a1}
    \end{subfigure}
    \begin{subfigure}{0.45\textwidth}
    \begin{center}
        \includegraphics[width=\linewidth]{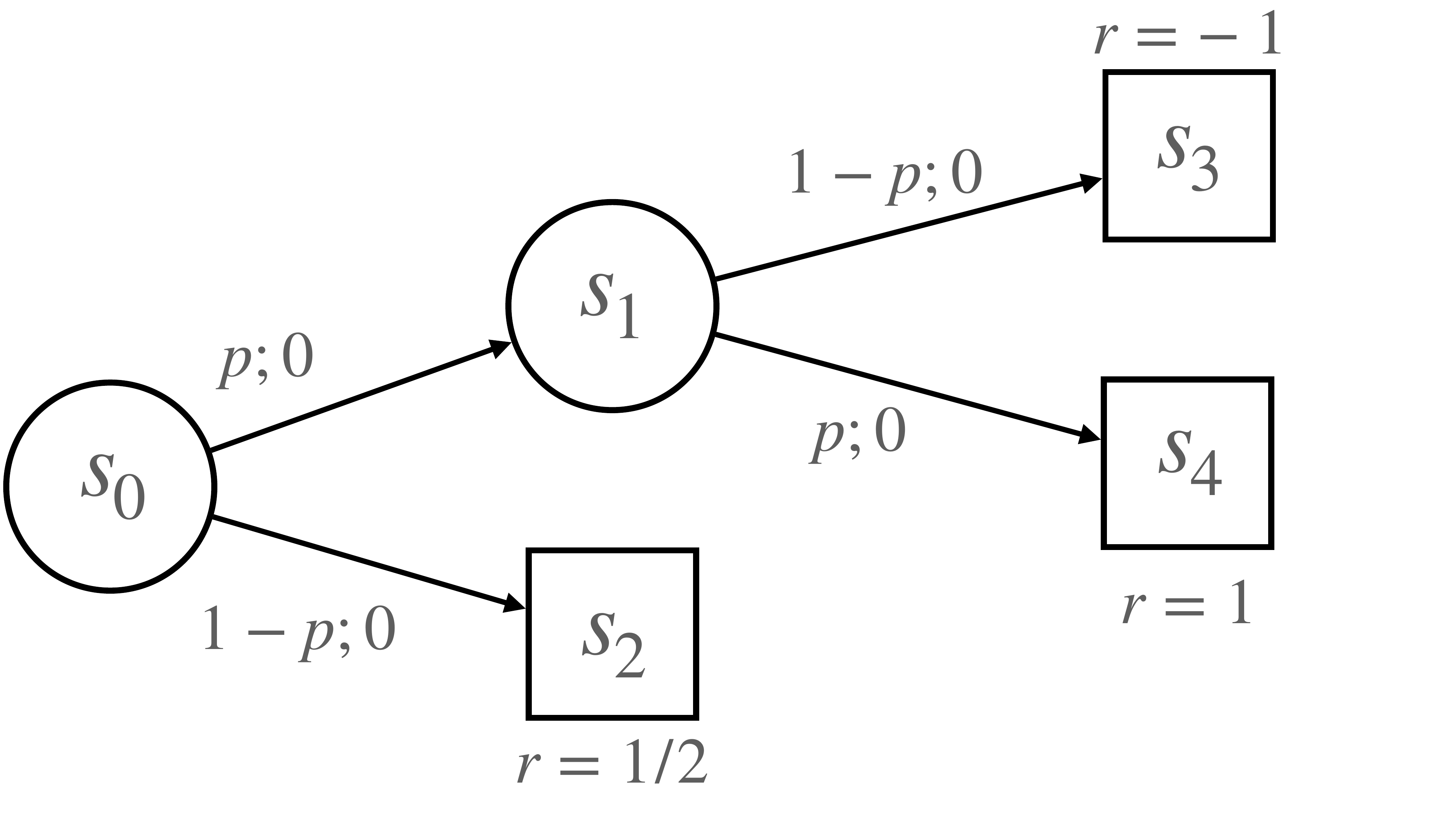}
    \end{center}
    \caption{Transitions and rewards for action $a_2$}
    \label{fig: lec 1 exo 0 a2}
    \end{subfigure}
    \end{center}
        \caption{MDP instance for \Cref{sec: lec 1 fdt}. The transition probabilities and rewards are indicated beside each arc.}
        \label{fig: lec 1 exo 0}
    \end{figure}
    
    \item {\em Percentile guarantees.} Assume that $\cP$ is constructed from some random observations sampled from a fixed true constant transition kernel $P^0$, so that $\cP$ is a random variable. Let $z\opt$ be the optimal objective value and $\pi^{\cP} $ be an optimal solution in~\eqref{eq: lec 1 rmdp}; note that $z\opt,\pi^{\cP} $ are random variables (they depend on the random variable $\cP$). Show that $\pi^{\cP} $ guarantees its worst-case return over $\cP$ with a probability at least as large as $\bP(P^0 \in \cP)$, i.e., show that $\bP\left(R(\pi^{\cP}, P^0) \geq z\opt\right) \geq \bP(P^0 \in \cP)$, where the probability $\bP$ is over the random variables $\pi^{\cP} ,z\opt,\cP$.
    \item {\em Entries of an optimal policy.} 
    Show that for any $d\in \N$, there exist simple robust MDP instances described by rational entries such that the entries of any stationary optimal policy may need to be algebraic numbers of degree $d$.
\end{enumerate}
\section{References}
Robust MDPs date as far back as~\cite{satia1973markov}. The fact that history-dependent policies may be necessary to solve robust MDPs with general uncertainty sets is proved in~\cite{wiesemann2013robust}; the proof presented above is adapted from an example in~\cite{wang2023foundation}. The potential non-existence of uniformly optimal policies follows from more general arguments given in \cite{benyamine2026existence}.
The irrationality of the entries of optimal policies/transition probabilities mirrors existing results for concurrent stochastic games (e.g. Example 3.2.1 in~\cite{filar2012competitive}) and finite games with more than two players~\cite{nash1951non}. Finally, there are many NP-hardness results for solving robust MDPs or policy evaluation for general uncertainty sets; I use here the elegant proof from~\cite{goh2018data} and \cite{borkar2009hamiltonicity} based on Hamiltonian cycles.
I also refer to the thesis of Le Tallec~\cite{le2007robust} which contains a flurry of counterexamples and hardness results. The instance from \Cref{fig: lec 1 exo 0} is from \cite{grand2024tractable}. A recent survey on robust MDPs can be found in \cite{suilen2024robust}.
\chapter{Robust MDPs with rectangular uncertainty sets}\label{chap: lec 2 rectangular models}
This chapter introduces the s-rectangular and sa-rectangular models of uncertainty and proves several of the favorable tractability properties implied by these rectangularity assumptions (strong duality, optimality of stationary policies, existence of uniformly optimal policies, fixed point optimality equations). It ends with a more advanced section showing that in some sense, rectangularity is the ``weakest" assumption one needs to make to ensure these favorable properties. We assume $\cP \neq \emptyset$ throughout.
\section{Modeling assumptions}
In the previous chapter, we have seen several examples of uncertainty sets leading to difficult optimization problems. For instance, the uncertainty set for~\Cref{ex: lec 0 ex 2} with $p=p',p\in [0,1]$ is convex but the same parameter $p$ appears for the transitions out of the two different states $s_0$ and $s_1$.
The class of {\em s-rectangular} uncertainty sets precisely remedies this issue.
\begin{definition}\label{def: lec 2 srec}
    An uncertainty set $\cP \subset \DSSA$ is s-rectangular if we can write $\cP = \times_{s \in \cS} \cP_{s}$ for some sets $\cP_{s} \subseteq \DSA$ for each $s \in \cS$.
\end{definition}
Let us parse this definition. If $\cP$ is s-rectangular then for each $s \in \cS$ the set $\cP_s \subseteq \DSA$ must be equal to the set $\{(P(s,a))_{a \in \cA} \; | \; P \in \cP\}$, the set of feasible transitions out of state $s$. The Cartesian product in the definition then means that the transitions out of each state $s$ can be chosen (by the adversary) independently of the transitions out of the other states. We can revisit the examples studied so far and inspect their rectangularity or non-rectangularity properties:
\begin{itemize}
    \item \Cref{ex: lec 0 ex 1} with $p \in \{0,1\}$ is s-rectangular. Let us index the columns by $ (a_1,a_2)$ and the rows by $(s_0,G,B)$. Then $P(G) = P(B) = \begin{pmatrix}
        1 & 1 \\
        0 & 0 \\
        0 & 0
    \end{pmatrix}$ since there is no uncertainty for the transitions out of the states $G,B$. Additionally for $p=0$ the transitions out of $s_0$ are $\begin{pmatrix}
        0 & 0 \\
        0 & 1 \\
        1 & 0
    \end{pmatrix}$, and for $p=1$ the transitions out of $s_0$ are $\begin{pmatrix}
        0 & 0 \\
        1 & 0 \\
        0 & 1
    \end{pmatrix}$. Hence $\cP_{s_0} = \left\{ \begin{pmatrix}
        0 & 0 \\
        0 & 1 \\
        1 & 0
    \end{pmatrix}, \begin{pmatrix}
        0 & 0 \\
        1 & 0 \\
        0 & 1
    \end{pmatrix} \right\}$. 
    Hence $\cP = \cP_{s_0} \times \{P(G)\} \times \{P(B)\}$ is s-rectangular. Note that $\cP$ is not convex (it is finite).
    \item \Cref{ex: lec 0 ex 2} with $p,p ' \in [0,1]$ is s-rectangular and convex: we can choose the transitions out of state $s_0$ (parameter $p$) and out of state $s_1$ (parameter $p'$) independently.
    \item \Cref{ex: lec 0 ex 2} with $p=p'$ is {\em not} s-rectangular: if we choose a value for $p$ we fix the transitions out of both $s_0$ and $s_1$. Still we would get a convex set even for $p=p'$.
    \item \Cref{ex: lec 1 irrational} with $p \in [0,1]$ is s-rectangular and convex: there is no uncertainty for the transitions out of the absorbing state $s_1$, so that $p$ only influences the transitions out of $s_0$. Indexing columns by $(a_1,a_2)$ and rows by $(s_0,s_1)$, the set of feasible transitions out of $s_0$ is $\cP_{s_0} = \left\{ \begin{pmatrix}
        1-p & p \\
        p & 1- p \\
    \end{pmatrix} \; | \; p \in [0,1] \right\}$.
\end{itemize}
In applications of robust MDPs, one typically considers sets based on distance from an empirical estimation $\hat{P} \in \DSSA$, of the form
\[\{P \in \Delta(\cS)^{\cS \times \cA} \; | D(P_s,\hat{P}_s) \leq \epsilon_s, \forall \; s \in \cS\}\]
with $D$ that can be the $\ell_{\infty}$-distance, the $\ell_{1}$-distance, the $\ell_2$-distance, or based on Kullback-Leibler or $\phi$-divergences.

One way to ``visually check" the s-rectangularity of a robust MDP instance represented on a graph (as in most examples in these notes) is to check that there is no parameter (such as some $p \in [0,1]$) involved in the transitions out of two different states, and no constraints across the transition parameters out of different states.
A  more abstract way to define s-rectangularity is to say that $\cP$ is s-rectangular if and only if we can mix the entries (indexed by states) of any two kernels in $\cP$ and still recover a feasible kernel in $\cP$ (see \Cref{sec: lec 2 fdt} for a rigorous statement). While we have defined s-rectangularity for uncertainty sets, this property is merely saying that $\cP$ is a Cartesian product, and we can wonder if other objects studied so far have this property. Crucially, this is indeed the case for the set of policies:
\begin{itemize}
    \item The set of stationary policies $\PiS$ can be identified with $\Delta(\cA)^{\cS}$, which is a Cartesian product: $\Delta(\cA)^{\cS} = \times_{s \in \cS} \Delta(\cA)$.
    \item The set of history-dependent policies $\PiH$ is also a Cartesian product of copies of the simplex $\Delta(\cA)$, taken over all possible finite trajectories in $\cH$.
\end{itemize}
Since s-rectangularity makes the uncertainty set $\cP$ (which we can interpret as the ``decision set" of the adversary choosing a worst-case transition kernel) much like the set of policies for nominal MDPs, we can expect that we can now derive more favorable tractability properties than what we've seen in the previous chapter. Regarding assumptions, we will assume compactness since this simplifies the exposition of the results; we only assume convexity where we need it (in particular we don't need convexity in~\Cref{sec: lec 2 adv} but we need it in \Cref{sec: lec 2 rmdp}).
\begin{remark}
Rectangularity is in fact encountered in many other research areas. For instance:
\begin{itemize}
    \item The two-dimensional rectangle $\{(x,y) \in \R^{2} \; | \; a \leq x \leq b, c \leq y \leq d\}$ is rectangular, and more generally $n$-dimensional boxes are rectangular (this is the prototypical example).
    \item The ``independent row uncertainty" and ``independent column uncertainty" assumptions in the literature on computing joint spectral radii exactly say the set of admissible matrices is a Cartesian product (over rows or columns).
    \item A multivariate random variable with independent components (i.e., its distribution is a product distribution) has a rectangular support.
\end{itemize}
\end{remark}
\section{Policy evaluation}\label{sec: lec 2 adv}
For now we focus on solving the policy evaluation problem~\eqref{eq: lec 1 policy evaluation} for a given stationary policy. Under the s-rectangularity assumption, we will see that this problem can be reformulated as an {\em adversarial MDP}, defined next.
\begin{definition}
    Let $\cM = (\cS,\cA,r,P,p_0)$ be an MDP and let $\cP \subseteq \DSSA$ be a compact s-rectangular uncertainty set. Let $\pi \in \PiS$ be a stationary policy.
    The adversarial MDP is defined as the following MDP instance $\cM'=(\cS',\cA',r',P',p'_0)$:
    \begin{itemize}
        \item The set of states $\cS'$ is equal to $\cS$.
        \item The set of actions $\cA'_s$ available at state $s\in \cS'$ is the compact set $\cP_s \subseteq \DSA$.
        \item The instantaneous reward for choosing action $P_s \in \cP_s$ at state $s \in \cS'$ is $r'(s,P_s) = - \sum_{a \in \cA} \pi(s,a) P(s,a)\tr r(s,a)$, which is independent of the next state. We will write $r'(s,a')1_{\cS'}$ for the vector $(r'(s,a'),...,r'(s,a')) \in \R^{\cS'}$, given some ``action" $a' \in \cP_s$ of the adversary.
        \item The transition probability toward $s' \in \cS$ after choosing action $P_s \in \cP_s$ at state $s \in \cS'$ is $\sum_{a \in \cA} \pi(s,a)P(s,a,s') \in \R$.
        \item The initial distribution $p'_0$ is equal to $p_0$.
    \end{itemize}
\end{definition}
In the adversarial MDP, a stationary deterministic policy is a choice of $P_s \in \cP_s$ for each state $s \in \cS$, i.e., it is a choice of $P \in \cP = \times_{s \in \cS} \cP_s$. As we discussed at the end of the previous chapter, the main missing ingredient to make~\eqref{eq: lec 1 policy evaluation} tractable for general $\cP$ was the (potential) lack of rectangularity. Now that we have assumed $\cP$ to be s-rectangular, we can indeed interpret the set $\cP$ as a set of stationary {\em deterministic} policies in the adversarial MDP; the set of stationary randomized policies in the adversarial MDP is $\times_{s \in \cS}\Delta(\cP_s)$.
We need the minus sign to turn the minimization problem \eqref{eq: lec 1 policy evaluation} into the adversarial MDP maximization problem \eqref{eq: lec 0 mdp}. The following theorem is the main result of this section.
\begin{theorem}\label{th: lec 2 adversarial mdp}
Let $\pi \in \PiS$ and $\cP \subset \DSSA$ be a compact s-rectangular uncertainty set.
    \begin{enumerate}
        \item {\em Value function correspondence.} Let $v^{P,P'} \in \R^{\cS}$ be the value function of the stationary deterministic adversarial policy $P \in \cP$ in the adversarial MDP. Then $v^{P,P'} = - v^{\pi,P}$ for $v^{\pi,P} \in \R^{\cS}$ the value function of the policy $\pi \in \PiS$ in the nominal MDP with transitions $P$.
        \item {\em Optimal objective correspondence.} We have $\arg \min_{P \in \cP} R(\pi,P)  = \arg \max_{P \in \cP} p_0 \tr v^{P,P'}$.
    \end{enumerate}
\end{theorem}
\begin{proof}
\begin{enumerate}
    \item By definition, the value function $v^{P,P'} \in \R^{\cS}$ of a stationary deterministic adversarial policy $P \in \cP$ for the transition probability $P'$ in the adversarial MDP satisfies the following fixed point equation:
\begin{align*}
    v^{P,P'}(s) & = \E_{a' = P(s),s ' \sim P'(s,a')} \left[ r'(s,a') + \gamma  v^{P,P'}(s') \right] \\
    & = P'(s,a') \tr (r'(s,a')1_{\cS'} + \gamma  v^{P,P'}) \\
    & = \left(\sum_{a \in \cA} \pi(s,a)P(s,a) \right) \tr \left( - \sum_{a \in \cA} \pi(s,a) P(s,a) \tr r(s,a) + \gamma v^{P,P'}(s')\right)_{s' \in \cS} \\
    & = - \sum_{a \in \cA} \pi(s,a)P(s,a)\tr (r(s,a) + \gamma (-v^{P,P'})) \\
    & = - T^{\pi,P}(-v^{P,P'})(s).
\end{align*}
The last equality shows that $v^{P,P'}$ is a fixed point of the operator $v \mapsto -T^{\pi,P}(-v)$, with $T^{\pi,P}$  defined in \eqref{eq: lec 0 T pi P}. Hence $v^{P,P'}=-v^{\pi,P}$ (since $v^{\pi,P}$ is the unique fixed point of the contraction $T^{\pi,P}$). 
    \item We have 
    \[  \min_{P \in \cP} R(\pi,P)  = \min_{P \in \cP} p_0 \tr v^{\pi,P} = - \max_{P \in \cP} p_0 \tr (-v^{\pi,P})  = - \max_{P \in \cP} p_0 \tr v^{P,P'}\]
where the first equality is by definition of the return in the nominal MDP, the second equality is by $\min f = - \max - f$, and the third equality is from the value function correspondence (first statement in this theorem). Hence the sets of optimal solutions for $\min_{P \in \cP} R(\pi,P) $ and $\max_{P \in \cP} p_0 \tr v^{P,P'}$ coincide.
\end{enumerate}
\end{proof}
Note that $\max_{P \in \cP} p_0 \tr v^{P,P'}$ computes the best {\em stationary deterministic} adversarial policies in the adversarial MDP (the transition probabilities are constant over time; defining history-dependent adversaries requires a bit more care and will be done in Chapter~\ref{chap: lec 3 connection}). Since in the adversarial MDP the action sets are compact and the rewards and transitions are continuous in the actions chosen, we can apply the results from Chapter~\ref{chap: lec 0 mdp} (see Section~\ref{sec: lec 0 fdt} for the case of nominal MDPs with state-dependent compact action sets) and we know that in the adversarial MDP a uniformly optimal policy $P$ exists and can be chosen stationary deterministic (i.e. $P \in \cP$).
We define the adversarial Bellman operator $T^{\pi}:\R^{\cS} \to \R^{\cS}$ as
\begin{equation}\label{eq: lec 2 T pi}
    T^{\pi}v(s) = \min_{P \in \cP} T^{\pi,P}v(s), \forall \; v \in \R^{\cS}, s \in \cS.
\end{equation}
Note that the family $\{T^{\pi,P} \; | \; P \in \cP\}$ satisfies the attainability property with respect to $T^{\pi}$, as we state formally next.
\begin{lemma}\label{lem: lec 2 attainability}
    Let $\cP \subset \DSSA$ be a compact s-rectangular uncertainty set and let $\pi \in \PiS$. For any $v \in \R^{\cS}$, there exists $P \in \cP$ such that $T^{\pi}v=T^{\pi,P}v$.
\end{lemma}
\begin{proof}
    Let $s \in \cS$ and $v \in \R^{\cS}$.  Note that the minimization program in the definition of $T^{\pi}$ as in~\eqref{eq: lec 2 T pi} only depends on $P(s)$ the $s$-th component of $\cP$.
    Since $\cP$ is s-rectangular compact we get that we can solve each of these minimization programs independently and recover a joint solution that is optimal simultaneously across all $s \in \cS$.
\end{proof}
In fact one can easily verify that the set $\{T^{\pi,P} \; | \; P \in \cP\}$ is a Cartesian product when $\cP$ is s-rectangular; for a non-rectangular set $\cP$, this set is still well-defined, as well as $T^{\pi}$ as in~\eqref{eq: lec 2 T pi}, but the attainability property may not hold (see~\Cref{sec: lec 2 minimality rec} for more on this topic).

It is worth emphasizing one more time the role of rectangularity. For nominal MDPs, the Bellman operator $T^P$ imposes {\em componentwise} maximization over the policies, i.e., it imposes that we can choose actions independently across each state. Rectangularity of the uncertainty set ensures the same property, that is, the adversary can independently minimize the values obtained at each state. Since we have already noted right before Section~\ref{sec: lec 1 fdt} how similar the two functions $\pi \mapsto v^{\pi,P}$ and $P \mapsto v^{\pi,P}$ are, the interpretation of the policy evaluation as an adversarial MDP is natural (once we allow for rectangular sets of decisions for the adversary, i.e., once we have s-rectangular uncertainty). \Cref{tab: lec 2 comparisons mdp adv} summarizes the similarities between nominal MDPs and adversarial MDPs. Note that in general it is much easier to work with value functions than directly with the return function; in some sense, non-rectangular models typically lead to difficult optimization problems because one cannot, in general, rely on value functions  (which are defined componentwise, i.e., in a rectangular manner).

The following direct corollary is a reformulation of~\Cref{prop: lec 0 operators contraction} and~\Cref{th: lec 0 mdp main} for the adversarial MDP; note that we do not need any convexity assumption on $\cP$ for now (indeed in \Cref{chap: lec 0 mdp} we worked with $\cA$ finite and then extended the results to $\cA$ compact in \Cref{sec: lec 0 fdt}).
\begin{corollary}\label{cor: lec 2 adversarial mdp - alg}
Let $\cP \subset \DSSA$ be a compact s-rectangular uncertainty set and let $\pi \in \PiS$.
\begin{enumerate}
    \item {\em Bellman equations.} The operator $T^{\pi}$ defined in \eqref{eq: lec 2 T pi} is a monotone contraction with unique fixed point $v^{\pi} \in \R^{\cS}$, satisfying the following worst-case Bellman equations in the variable $v \in \R^{\cS}$:\begin{equation}\label{eq: lec 2 adv bellman}
        v(s) = \min_{P \in \cP} \sum_{a \in \cA} \pi(s,a) P(s,a) \tr (r(s,a) + \gamma v), \quad \forall \; s \in \cS.
    \end{equation}
    \item {\em Uniform optimality.} There exists a {\em uniform} worst-case transition probabilities $P^{\pi} \in \cP$, that solves the policy evaluation problem \eqref{eq: lec 1 policy evaluation} simultaneously for all initial distributions $p_0 \in \DS$, and the value function $v^{\pi,P^{\pi}} \in \R^{\cS}$ coincides with $v^{\pi}$.
    \item {\em Optimality conditions.} We can compute a uniform worst-case transition kernel $P^{\pi} \in \cP$ by choosing any transition kernel that jointly attains that $\arg \min$ in $T^{\pi}v^{\pi}(s)$ across all $s \in \cS$, i.e., we have that $P^{\pi}$ is a uniform worst-case transition kernel if and only if $T^{\pi}v^{\pi} = T^{\pi,P^{\pi}}v^{\pi}$.
\end{enumerate}
\end{corollary}
Note that the first and second statements in \Cref{cor: lec 2 adversarial mdp - alg} show that $v^{\pi} \in \R^{\cS}$ is the ``smallest"  worst-case value function of $\pi \in \PiS$, in that
\[v^{\pi}(s) = \min_{P \in \cP} v^{\pi,P}(s), \quad \forall \; s \in \cS\]
and for any $p_0 \in \Delta(\cS)$ we must have
\[\min_{P \in \cP} R(\pi,P) = p_0 \tr v^{\pi}.\]
\begin{remark}
    \Cref{prop: lec 1 multiple value functions} shows that even for s-rectangular uncertainty sets, there may be multiple worst-case value functions (for a single policy $\pi \in \PiS$), since worst-case kernels need not be unique. As such, the vector $v^{\pi} \in \R^{\cS}$ is in some sense the component-wise smallest such worst-case value function (as shown by the last-but-one equation above).
\end{remark}
\paragraph{Computing worst-case kernels.}
Since we interpret the policy evaluation problem as an adversarial MDP, we use the whole machinery developed in \Cref{chap: lec 0 mdp} for computing a uniform worst-case optimal kernel. 
The value iteration and policy iteration algorithms and their convergence rates are given next.
\begin{algorithm}
\caption{Value iteration for adversarial MDPs}\label{alg: lec 1 VI}
\begin{algorithmic}[1]
\State {\bf Input:} $\pi \in \PiS$ and $\cP$ compact s-rectangular
\State Initialize $v^{0} \in \R^{\cS}$.
\For{$t \in \N$}
    \State $v^{t+1} = T^{\pi} v^{t}$
\EndFor
\end{algorithmic}
\end{algorithm}
\begin{proposition}\label{prop: lec 1 VI rate}
Let $\pi \in \PiS$ and let $\cP$ be a compact s-rectangular uncertainty set.
   Let $(v^t)_{t \in \N}$ be the sequence of vectors generated by \Cref{alg: lec 1 VI}.
    \begin{enumerate}
        \item The sequence $(v^t)_{t \in \N}$ converges to $v^{\pi}$ at a geometric rate of $\gamma$:
        \[\|v^{\pi} - v^{t}\|_{\infty} \leq \gamma^t \|v^{\pi} - v^{0}\|_{\infty}\]
    \item Assume that $\|v^{\pi} - v^{t}\|_{\infty} \leq \epsilon$ for some $t \in \N$ and let $P \in \cP$ such that $T^{\pi,P}v^t = T^{\pi}v^t$. Then $P$ is a uniform $2\epsilon \gamma (1-\gamma)^{-1}$ worst-case transition probabilities for $\pi$: $\|v^{\pi}-v^{\pi,P}\|_{\infty} \leq 2\epsilon \gamma (1-\gamma)^{-1}$.
    \end{enumerate}
\end{proposition}
\begin{algorithm}
\caption{Policy iteration for adversarial MDPs}\label{alg: lec 1 PI}
\begin{algorithmic}[1]
\State {\bf Input:} $\pi \in \PiS$ and $\cP$ compact s-rectangular
\State Initialize $P^0 \in \cP$, $v^0 = v^{\pi,P^0}$.
\For{$t \in \N$}
    \State Choose $P^{t+1} \in \cP$ such that $T^{\pi}v^{t}(s) = T^{\pi,P^{t+1}}v^{t}(s)$ for all $s \in \cS$
    \State Choose $v^{t+1} = v^{\pi,P^{t+1}}$
\EndFor
\end{algorithmic}
\end{algorithm}
\begin{proposition}\label{prop: lec 1 PI rate}
Let $\pi \in \PiS$ and let $\cP$ be a compact s-rectangular uncertainty set.
 Let $(P^t)_{t \in \N}$ be the sequence of policies generated by~\Cref{alg: lec 1 PI}. Then 
         $\| v^{\pi} - v^{\pi,P^t}\|_{\infty} \leq \gamma^t \|v^{\pi}-v^{\pi,P^0}\|_{\infty}.$
\end{proposition}
Similarly as in~\Cref{chap: lec 0 mdp}, the proofs of the above convergence rates follow from the more general versions of the algorithms analyzed in Appendix~\ref{app: algorithms}.
We can still obtain a mathematical programming formulation for the policy evaluation program. However we need to be careful here: the operator $T^{\pi}$ is componentwise {\em concave} (as the minima over affine forms) whereas the operator $T^{P}$ is componentwise {\em convex} (as the maxima over affine forms). Hence to obtain a convex formulation for computing $v^{\pi}$ we need to use the alternate version of the contraction lemma mentioned in Section~\ref{sec: lec 0 fdt}, which yields the following proposition.
\begin{proposition}\label{prop: lec 1 math prog}
Let $\pi \in \PiS$ and let $\cP$ be a compact s-rectangular uncertainty set.
    The smallest worst-case value function $v^{\pi}$  is the unique optimal solution to the following convex program:
    \begin{equation}\label{eq: lec 1 math prog}
        \max \left\{ \sum_{s \in \cS} v(s) \; \middle| \; v(s) \leq \min_{P \in \cP} \sum_{a \in \cA} \pi(s,a)P(s,a) \tr (r(s,a) + \gamma v), \forall \; s \in \cS, v \in \R^{\cS}\right\}
    \end{equation}
\end{proposition}
\paragraph{Complexity.} The main difference with the complexity given for nominal MDPs in~\Cref{chap: lec 0 mdp} pertains to how we can evaluate the operator $T^{\pi}$ at each iteration of these algorithms. This requires solving an optimization problem with linear objective over the compact set $\cP$. This can be done efficiently for s-rectangular uncertainty sets with additional structure (e.g. based on $\ell_{\infty}$-balls, $\ell_{p}$-balls for $p \in \{1,2\}$, based on Kullback-Leibler divergence, etc.), and we point to \Cref{sec: lec 2 ref} for some references.
\begin{remark}
    In this section on the adversarial MDP we do not need convexity of $\cP$ for our results. However the Bellman operator $T^{\pi}$ only involves optimizing linear functions over $\cP$, hence we can replace $\cP$ by $\conv(\cP)$ in the expression of the Bellman operator $T^{\pi}$ (and therefore in the definition of the worst-case value function $v^{\pi}$ as the fixed point of $T^{\pi}$). This will not be true in the next section on robust MDPs (i.e., on the maxmin formulation~\eqref{eq: lec 1 rmdp}), as we will see next.
\end{remark}
\section{Optimal policies for rectangular robust MDPs}\label{sec: lec 2 rmdp}
We now get back to solving the robust MDP problem~\eqref{eq: lec 1 rmdp}. Similarly to the case of nominal MDPs and adversarial MDPs, robust MDPs can be analyzed via the properties and fixed point of the following {\em robust Bellman operator} $T:\R^{\cS} \to \R^{\cS}$, defined as 
\begin{equation}\label{eq: lec 2 T}
    T v(s) = \max_{\pi \in \PiS} T^{\pi}v(s) =  \max_{\pi \in \PiS} \min_{P \in \cP} T^{\pi,P}v(s), \forall \; v \in \R^{\cS}, s \in \cS.
\end{equation}
The robust Bellman operator inherits the monotonicity and contraction properties of the Bellman operator $T^{\pi,P}$, as shown in the next proposition.
\begin{proposition}\label{prop: lec 2 robust Bellman T}
        Let $\cP \subseteq \DSSA$ be a compact uncertainty set. The Bellman operator $T$ defined in~\eqref{eq: lec 2 T} is a monotone contraction with unique fixed point $v\opt \in \R^{\cS}$, satisfying the following robust Bellman equation in the variable $v \in \R^{\cS}$:
        \begin{equation}\label{eq: lec 2 Bellman rmdp}
            v(s) = \max_{\pi \in \PiS} \min_{P \in \cP} \sum_{a \in \cA} \pi(s,a) P(s,a) \tr (r(s,a) + \gamma v), \quad \forall \; s \in \cS
        \end{equation}
\end{proposition}
\begin{proof}
    Since $T^{\pi,P}$ is monotone we directly get that $T$ is monotone. To show that $T$ is a contraction, we note that for $v,w \in \R^{\cS}$ and $s \in \cS$ we have
\begin{align*}
    Tv(s) & = \max_{\pi \in \PiS} \min_{P \in \cP} T^{\pi,P}v(s) \\
    & \leq \max_{\pi \in \PiS} \min_{P \in \cP} \left(T^{\pi,P}w(s) + \gamma \| v-w\|_{\infty}\right) \\
    & \leq \gamma \| v-w\|_{\infty} + \max_{\pi \in \PiS} \min_{P \in \cP} T^{\pi,P}w(s) \\
    & = \gamma \| v-w\|_{\infty} + Tw(s) 
\end{align*}
where the first inequality is from the contraction property of $T^{\pi,P}$ for $\pi \in \PiS,P \in \cP$. This shows $Tv(s) - Tw(s) \leq \gamma \| v-w\|_{\infty}$. Reversing the roles of $v$ and $w$ above yields $Tw(s) - Tv(s) \leq \gamma \| v-w\|_{\infty}$, so that we have proved $\|Tv - Tw\|_{\infty} \leq \gamma \| v-w\|_{\infty}$ as claimed. Therefore there exists a unique fixed point $v\opt \in \R^{\cS}$ of the operator $T$ and this fixed point satisfies the equation~\eqref{eq: lec 2 Bellman rmdp}.
\end{proof}
Note that we do not need s-rectangularity for the previous proposition; this assumption will be needed for interpreting the fixed point $v\opt$ as the value function of an optimal policy.
Let us first investigate the case where $\cP$ is s-rectangular but not necessarily convex. We have already seen that an optimal policy for the robust MDP problem \eqref{eq: lec 1 rmdp} may need to be history-dependent.
\begin{proposition}\label{prop: lec 2 hist dep opt}
    Let $\cP \subseteq \DSSA$ be a compact s-rectangular uncertainty set. There may not exist any stationary optimal policies, i.e., we may have
    \[ \sup_{\pi \in \PiH} \inf_{P \in \cP} R(\pi,P) > \max_{\pi \in \PiS} \min_{P \in \cP} R(\pi,P).\]
\end{proposition}
\begin{proof}
This was already stated in \Cref{prop: lec 1 nonstationary opt policies}, where the proof considers an s-rectangular non-convex uncertainty set.
\end{proof}
We may still wonder if there is any useful information that we can obtain from the robust Bellman operator $T$ and its fixed point $v\opt$. These objects can be used to characterize ``the best policy among stationary ones", as proved next.
\begin{proposition}\label{prop: lec 2 maxmin piopt}
        Let $\cP \subseteq \DSSA$ be a compact s-rectangular uncertainty set. Let $\pi\opt \in \PiS$. The two statements are equivalent:
        \begin{enumerate}
            \item $\pi\opt \in \PiS$ is a best stationary policy across all initial distributions:
            \begin{equation}\label{eq: lec 2 unf opt PiS}
            \forall \; p_0 \in \DS, \pi\opt \in \arg \max_{\pi \in \PiS} \min_{P \in \cP} R(\pi,P).
        \end{equation}
        \item $T^{\pi\opt}v\opt = Tv\opt$.
        \end{enumerate}
\end{proposition}
\begin{proof}
First, we note that the set $\cap_{s \in \cS} \arg \max_{\pi \in \PiS} T^{\pi}v\opt(s)$ coincides with the set of policies $\bar{\pi} \in \PiS$ such that $Tv\opt = T^{\bar{\pi}}v\opt$.

 We show that any policy $\bar{\pi}$ in $\displaystyle \cap_{s \in \cS} \arg \max_{\pi \in \PiS} T^{\pi}v\opt(s)$ is uniformly optimal (among stationary policies); this set is non-empty because for any $s \in \cS$, the scalar $ T^{\pi}v\opt(s)$ only depends on the $s$-th component of $\pi \in \PiS$, and $\PiS$ is a Cartesian product.  Since $Tv\opt = v\opt$ by definition of $v\opt$, we get that $T^{\bar{\pi}}v\opt = v\opt$ so that $v\opt$ coincides with $v^{\bar{\pi}} \in \R^{\cS}$, the smallest worst-case value function of $\bar{\pi}$. Additionally, for any $\pi \in \PiS$ we have $v^{\pi}=T^{\pi}v^{\pi} \leq Tv^{\pi}$, where the inequality is by definition of $T$ as in \eqref{eq: lec 2 T}. Hence $v^{\pi} \leq Tv^{\pi}$, and since $T$ is monotone we have $v^{\pi} \leq T^n v^{\pi}$ for any $n \in \N$. Since $T$ is a contraction, taking $n \to + \infty$ yields $v^{\pi} \leq v\opt = v^{\bar{\pi}}$. Since $p_0 \in \Delta(\cS)$, we conclude that the set $\displaystyle \cap_{s \in \cS} \arg \max_{\pi \in \PiS} T^{\pi}v\opt(s)$ contains only policies that are uniformly optimal stationary policies. 

 Now let $\bar{\pi} \in \PiS$ be any stationary uniformly optimal policy. Since we must also have $v^{\bar{\pi}} \geq v^{\pi}$ for any other policy $\pi \in \PiS$, we get that $v^{\bar{\pi}} = v\opt$. Therefore $T^{\bar{\pi}}v\opt = T^{\bar{\pi}}v^{\bar{\pi}} = v^{\bar{\pi}}=v\opt = T v\opt$ so that $\bar{\pi} \in \cap_{s \in \cS} \arg \max_{\pi \in \PiS} T^{\pi}v\opt(s)$, which shows that $\cap_{s \in \cS} \arg \max_{\pi \in \PiS} T^{\pi}v\opt(s)$ is the set of uniformly optimal stationary policies.
\end{proof}
Note that in the proof of the above proposition we use ``optimal stationary policies" as solutions to $\sup_{\pi \in \PiS} \min_{P \in \cP} R(\pi,P)$, whereas in these notes ``stationary optimal policies" means stationary solutions to $\sup_{\pi \in \PiH} \inf_{P \in \cP} R(\pi,P)$. From \Cref{prop: lec 2 hist dep opt}, stationary optimal policies need not exist for s-rectangular compact sets (without convexity). 
We show below that they exist as soon as we assume convexity. 
\begin{theorem}\label{th: lec 2 main rmdps}
    Let $\cP \subseteq \DSSA$ be a compact convex s-rectangular uncertainty set. The following statements hold:
        \begin{enumerate}
            \item {\em Strong duality.}
            The following strong duality result holds: 
    \begin{equation}\label{eq: lec 2 strong duality}
        \max_{\pi \in \PiS} \min_{P \in \cP} R(\pi,P) = \min_{P \in \cP} \max_{\pi \in \PiS}  R(\pi,P). 
    \end{equation}
    \item {\em Optimality of stationary policies.}
    We have
    \[ \sup_{\pi \in \PiH} \inf_{P \in \cP} R(\pi,P) = \max_{\pi \in \PiS} \min_{P \in \cP} R(\pi,P).\]
    \item {\em Optimality conditions.}  $\pi \in \PiS$ is a stationary uniformly optimal policy for \eqref{eq: lec 1 rmdp} if and only if $T^{\pi}v\opt = Tv\opt$.

    Additionally, let $P\opt \in \cP$ such that $T^{P\opt} v\opt = Tv\opt$. Then: $P\opt$ is a uniform worst-case kernel for any stationary uniformly optimal  policy $\pi\opt \in \PiS$, and any stationary uniformly optimal  policy $\pi\opt \in \PiS$ is actually a {\em nominal} uniformly optimal policy in the nominal MDP with transition probabilities $P\opt$.
    \end{enumerate}
\end{theorem}
We encourage the reader to compare \Cref{th: lec 2 main rmdps} with the main theorem for adversarial MDPs (\Cref{cor: lec 2 adversarial mdp - alg}) and the main theorem for nominal MDPs (\Cref{th: lec 0 mdp main}); we refer to \Cref{tab: lec 2 comparisons mdp adv} for a comparison of the main objects for these three optimization problems.
\begin{table}[htbp]
\centering
\small
\setlength{\tabcolsep}{4pt}
\renewcommand{\arraystretch}{1.35}
\begin{adjustbox}{center,max width=\textwidth}
\begin{tabular}{@{}p{4.2cm}ccc@{}}
\toprule
& \textbf{Nominal MDP} & \textbf{Adversarial MDP} & \textbf{Robust MDP} \\
\midrule

States
& $\cS$
& $\cS$
& $\cS$ \\

Actions at state $s$
& $\cA$
& $\cP_s$
& $\cA$ \\

Stationary deterministic policies
& $\displaystyle \PiSD = \times_{s \in \cS} \cA$
& $\displaystyle \cP = \times_{s \in \cS} \cP_s$
& $\PiSD$ \\

Stationary policies
& $\displaystyle \PiS = \times_{s \in \cS} \Delta(\cA)$
& $\displaystyle \conv(\cP)=\times_{s \in \cS}\conv(\cP_s)$
& $\PiS$ \\

Set of operators
& $\displaystyle \{T^{\pi,P} \mid \pi \in \PiSD\}$
& $\displaystyle \{T^{\pi,P} \mid P \in \cP\}$
& $\displaystyle \{T^{\pi} \mid \pi \in \PiS\}$ \\

Optimality operator
& $\displaystyle T^{P}v(s)=\max_{\pi \in \PiSD}T^{\pi,P}v(s)$
& $\displaystyle T^{\pi}v(s)=\min_{P \in \cP}T^{\pi,P}v(s)$
& $\displaystyle Tv(s)=\max_{\pi \in \PiS}\min_{P \in \cP}T^{\pi,P}v(s)$ \\

Optimal value function 
& $v^{P}$ solution of \eqref{eq: lec 0 Bellman mdp} 
& $v^{\pi}$ solution of \eqref{eq: lec 2 adv bellman} 
& $v\opt$ solution of \eqref{eq: lec 2 Bellman rmdp} \\
\bottomrule
\end{tabular}
\end{adjustbox}
\caption{Comparisons of the properties of nominal MDPs, adversarial MDPs, and robust MDPs under an s-rectangular uncertainty set $\cP$. The operator $T^{\pi,P}:\RS \to \RS$ is defined as $T^{\pi,P}v(s) = \sum_{a \in \cA} \pi(s,a)P(s,a)\tr (r(s,a) + \gamma v)$, for $s \in \cS$ and $v \in \RS$. The column for adversarial MDPs requires compactness of $\cP$ while the column for robust MDPs requires compactness and convexity of $\cP$. For adversarial MDPs we identify stationary randomized policies (which should be $\times_{s \in \cS} \Delta(\cP_s)$) with their effective transition probabilities equal to $\times_{s \in \cS} \conv(\cP_s)$.}
\label{tab: lec 2 comparisons mdp adv}
\end{table}
The remainder of this section is devoted to the proof of \Cref{th: lec 2 main rmdps}. Before we do so, it is worth discussing the intuition and implications of this theorem:
\begin{itemize}
    \item Contrary to the nominal case, we may need randomized policies at optimality. We already observed this in \Cref{ex: lec 1 irrational}, where we computed the unique optimal policy and noticed it was randomized. From optimization and game-theoretic standpoints, this is not surprising, since the robust MDP problem~\eqref{eq: lec 1 rmdp} has a ``two-player zero-sum game" flavor and it is a maxmin problem. Randomization can also be understood ``geometrically": from the third statement in \Cref{th: lec 2 main rmdps}, a stationary uniformly optimal policy solves the maxmin problem $\max_{\pi \in \PiS} \min_{P \in \cP} T^{\pi,P}v\opt(s)$ for each $s \in \cS$, with $(\pi,P) \mapsto T^{\pi,P}v\opt(s)$ bilinear. If we dualize the inner minimization we obtain that uniformly optimal solutions are the optimal solutions to a (larger) maximization problem over $\pi \in \PiS$ and some additional constraints (coming from the dualization of the inner minimization over $P \in \cP$), and these additional constraints may prevent us from choosing deterministic policies (which correspond to the extreme points of $\PiS$).

    Additionally, for $(\pi\opt,P\opt)$ as in the last statement of \Cref{th: lec 2 main rmdps}, we get that $P\opt$ is a worst-case transition kernel for $\pi\opt$, and that $\pi\opt$ is a nominal optimal policy in the MDP with transitions $P\opt$. Therefore in general $P\opt$ must induce a particular MDP instance where an optimal policy can be chosen randomized, whereas we saw in \Cref{chap: lec 0 mdp} that a nominal optimal policy can always be chosen deterministic; we refer to \Cref{prop: lec 2 NE SD} for more detail on this.
    \item We didn't need the convexity assumption on $\cP$ in \Cref{sec: lec 2 adv}: this is because we reduced the policy evaluation problem to an MDP, and MDPs are tractable even when the set of actions is not convex (e.g. we derived the results from \Cref{chap: lec 0 mdp} with finite action sets). In contrast, we crucially need the convexity assumption for $\cP$ in Theorem~\ref{th: lec 2 main rmdps}. This is because we will rely on minimax arguments and saddle-point properties that require convex domains. For sa-rectangular sets (see next section) convexity is not needed because we don't rely on minimax theorems.
    \item Keep in mind that optimal policies may require irrational entries, as in~\Cref{ex: lec 1 irrational}, even for very simple convex compact s-rectangular uncertainty sets. In practice the algorithms we'll run only compute approximately optimal policies (see \Cref{sec: lec 2 algo}) so this is not an issue.
    \item \Cref{th: lec 2 main rmdps} justifies the existence of stationary optimal policies in \Cref{ex: lec 1 irrational}, as claimed in \Cref{chap: lec 1 robust MDP}.
    \item If one works with non-closed, convex s-rectangular uncertainty sets $\cP$ then the same techniques as we used for proving \Cref{th: lec 2 main rmdps} would still show that optimal policies can be chosen stationary (although we wouldn't guarantee the existence of a worst-case kernel anymore). This situation is akin to the case of robust MDPs with mean-payoff objective, see \Cref{sec: lec 3 sgs}.
    \end{itemize}

It is also worth discussing the interpretation of strong duality from a few different angles:
\begin{itemize}
    \item From an optimization standpoint, it is important to emphasize that in principle we do {\em not} expect such a maxmin duality result to hold: while $\cP$ and $\PiS$ are convex compact sets, the function $(\pi,P) \mapsto R(\pi,P)$ need not be concave-convex (nor quasi-concave quasi-convex for that matter, see \Cref{sec: lec 2 fdt}), so the usual von Neumann or Sion's minimax theorems do not apply readily. 
    \item From a game-theoretic standpoint, the left-hand side is the best worst-case return over stationary policies, i.e., the left-hand side is the largest scalar value that the decision maker can be guaranteed to obtain if the decision maker first chooses $\pi \in \PiS$ and then the adversary chooses $P \in \cP$. In contrast, the right-hand side is the worst-case best return, i.e., if the adversary first chooses $P \in \cP$ and then the decision maker chooses $\pi \in \PiS$, the right-hand side is the smallest scalar value that the adversary can ``pay" to the decision maker. While it is clear that the right-hand side (``worst-case best return") is always greater than the left-hand side (``best worst-case return"), it is by no means obvious that these two quantities are equal.
\end{itemize}
It is also crucial to reckon that \eqref{eq: lec 2 strong duality} is equivalent to the existence of a ``Nash equilibrium", i.e., to a pair $(\pi\opt ,P\opt) \in \PiS \times \cP$ that are ``best responses" to each other (I use quotation marks because the concept of best response is for games and I did not rigorously define the game here). We make this formal in the next proposition.
\begin{proposition}\label{prop: lec 2 NE SD}
Let $\cP \subseteq \DSSA$ be a compact set. 
\begin{enumerate}
    \item The following statements are equivalent:
\begin{enumerate}
    \item Strong duality holds: $\max_{\pi \in \PiS} \min_{P \in \cP} R(\pi,P) = \min_{P \in \cP} \max_{\pi \in \PiS}  R(\pi,P).$
    \item There exists a pair $(\pi\opt ,P\opt) \in \PiS \times \cP$ such that 
    \begin{equation}\label{eq: lec 2 best response}
        \begin{aligned}
            \pi\opt & \in \arg \max_{\pi \in \PiS} R(\pi,P\opt) \\
            P\opt & \in \arg \min_{P \in \cP} R(\pi\opt ,P)
    \end{aligned}
    \end{equation}
\end{enumerate}
\item If strong duality holds, the pairs $(\pi\opt,P\opt) \in \PiS \times \cP$ satisfying~\eqref{eq: lec 2 best response} are exactly the pairs $(\pi\opt,P\opt) \in \PiS \times \cP$ satisfying
\begin{equation}\label{eq: lec 2 maxmin minmax}
    \begin{aligned}
    \pi\opt  \in \arg \max_{\pi \in \PiS} \min_{P \in \cP} R(\pi,P) \\
    P\opt \in \arg  \min_{P \in \cP} \max_{\pi \in \PiS}  R(\pi,P). 
\end{aligned}
\end{equation}
\end{enumerate}
\end{proposition}
\begin{proof}
\begin{enumerate}
    \item We first prove the first implication. Assume that strong duality holds and let $(\pi\opt,P\opt)$ be as in~\eqref{eq: lec 2 maxmin minmax}; we can always find such a pair of policy-kernel because $\PiS,\cP$ are compact and $(\pi,P) \mapsto R(\pi,P)$ is continuous.

Then we have
\begin{align*}
   \max_{\pi \in \PiS} \min_{P \in \cP} R(\pi,P) = \min_{P \in \cP} R(\pi\opt,P) \leq    R(\pi\opt,P\opt) \leq \max_{\pi \in \PiS} R(\pi,P\opt) = \min_{P \in \cP} \max_{\pi \in \PiS}  R(\pi,P)
\end{align*}
Since we assumed strong duality holds then all terms above are equal so that 
\[R(\pi\opt,P\opt)= \max_{\pi \in \PiS} R(\pi,P\opt)  =\min_{P \in \cP} R(\pi\opt,P)\]
and we have proved~\eqref{eq: lec 2 best response}.

Let us now assume the existence of $(\pi\opt,P\opt)$ satisfying~\eqref{eq: lec 2 best response}. Then
\begin{align*}
   \max_{\pi \in \PiS} \min_{P \in \cP} R(\pi,P) \geq \min_{P \in \cP} R(\pi\opt,P) =    R(\pi\opt,P\opt) = \max_{\pi \in \PiS} R(\pi,P\opt) \geq \min_{P \in \cP} \max_{\pi \in \PiS}  R(\pi,P)
\end{align*}
and from weak duality we get that all terms above are equal, which shows that strong duality holds.
\item This follows from the analysis of the first statement of this proposition.
\end{enumerate}
\end{proof}
Note that the second statement in the above proposition shows the property that the set of pairs $(\pi\opt ,P\opt) \in \PiS \times \cP$ that satisfy~\eqref{eq: lec 2 best response} is actually a Cartesian product; this is surprising since \eqref{eq: lec 2 best response} states that $\pi\opt,P\opt$ are best responses to each other.

Finally, note that $\arg \min_{P \in \cP} R(\pi\opt ,P)$, the set of worst-case kernels for $\pi\opt $ a solution to the maxmin problem, contains at least one $P\opt \in \cP$ such that $P\opt \in \arg  \min_{P \in \cP} \max_{\pi \in \PiS}  R(\pi,P)$, but this need not be the case for {\em all} elements in $\arg \min_{P \in \cP} R(\pi\opt ,P)$. This situation is analogous for Nash equilibria in finite zero-sum games: the set of best-responses to a Nash equilibrium contains a Nash equilibrium (by definition) but may contain other strategy profiles too.

\subsection{Proof of strong duality}
We give two different proofs of strong duality. The first proof is based on ideas from robust optimization. The second proof is based on Kakutani's fixed point theorem. A third proof based on Brouwer's fixed point theorem is suggested in \Cref{sec: lec 2 fdt}.

\begin{proof}[First proof of strong duality (based on robust optimization).]
    
We will make repeated use of the contraction lemma (\Cref{lem: lec 0 contraction lemma}).
\begin{align*}
    \max_{\pi \in \PiS} \min_{P \in \cP} R(\pi,P) & = \max_{\pi \in \PiS} \min_{P \in \cP} p_0  \tr v^{\pi,P} \\
    & = \max_{\pi \in \PiS} \min_{P \in \cP} \min \;  \{ p_0  \tr v \; | \; v(s) \geq T^{\pi,P}v(s) \; \forall \; s \in \cS, v \in \R^{\cS}\} \\
    & = \max_{\pi \in \PiS} \min \;  \{ p_0  \tr v \; | \; v(s) \geq  \min_{P \in \cP} T^{\pi,P}v(s) \; \forall \; s \in \cS, v \in \R^{\cS}\} \\
    & =  \min \; \{ p_0  \tr v \; | \; v(s) \geq \max_{\pi \in \PiS} \min_{P \in \cP} T^{\pi,P}v(s)  \; \forall \; s \in \cS, v \in \R^{\cS}\} \\
    & =  \min \; \{ p_0  \tr v \; | \; v(s) \geq  \min_{P \in \cP} \max_{\pi \in \PiS} T^{\pi,P}v(s)  \; \forall \; s \in \cS, v \in \R^{\cS}\} \\
    & = \min_{P \in \cP}  \min \; \{ p_0  \tr v \; | \; v(s) \geq   \max_{\pi \in \PiS} T^{\pi,P}v(s)  \; \forall \; s \in \cS, v \in \R^{\cS}\} \\
    & = \min_{P \in \cP} \max_{\pi \in \PiS} \min \; \{ p_0  \tr v \; | \; v(s) \geq   T^{\pi,P}v(s)  \; \forall \; s \in \cS, v \in \R^{\cS}\} \\
    & =  \min_{P \in \cP} \max_{\pi \in \PiS} p_0  \tr v^{\pi,P} \\
    & = \min_{P \in \cP} \max_{\pi \in \PiS}  R(\pi,P)
\end{align*}
where the first equality is by definition of the return $R(\pi,P)$ and the value function $v^{\pi,P} \in \R^{\cS}$, the second equality is by the contraction lemma applied to the fixed point $v^{\pi,P}$ of the monotone contraction $T^{\pi,P}$, the third and fourth equalities will be justified next, the fifth equality is from Sion's minimax theorem applied to the concave-convex (bilinear) functions $(\pi,P) \mapsto T^{\pi,P}v(s)$ for each $s \in \cS$ and fixed $v \in \R^{\cS}$ (see Appendix~\ref{app:sec analysis}), the sixth and seventh equalities will be justified next, and the last two equalities are again by the contraction lemma and the definition of the return and value functions.

We are left with justifying the third, fourth, sixth, and seventh equalities. 

The third and the sixth proceed from the same argument, which is that $\cP$ is s-rectangular so that the $s$-th constraint $v(s) \geq T^{\pi,P}v(s)$ only depends on the $s$-th component of $P$, and a minimization program over two variables (here $P \in \cP$ and $v \in \R^{\cS}$) can be written as two consecutive minimization programs.

The fourth and seventh equalities proceed from the following facts: (a) the rectangularity of the set $\PiS = \times_{s \in \cS} \Delta(\cA)$ (b) the $s$-th constraint $v(s) \geq \min_{P \in \cP} T^{\pi,P}v(s)$ only depends on the $s$-th component of $\pi \in \PiS$ (c) the outer maximum over $\pi \in \PiS$ is always attained in the fourth and seventh equalities, i.e., there exists $\pi' \in \PiS$ such that 
\[\max_{\pi \in \PiS} \min \; \{ p_0  \tr v \; | \; v(s) \geq   T^{\pi,P}v(s)  \; \forall \; s \in \cS, v \in \R^{\cS}\}\]
is equal to
\[\min \; \{ p_0  \tr v \; | \; v(s) \geq   T^{\pi',P}v(s)  \; \forall \; s \in \cS, v \in \R^{\cS}\}\]
and there exists $\pi'' \in \PiS$ such that
\[\max_{\pi \in \PiS} \min \; \{ p_0  \tr v \; | \; v(s) \geq   \min_{P \in \cP} T^{\pi,P}v(s)  \; \forall \; s \in \cS, v \in \R^{\cS}\}\]
is equal to
\[\min \; \{ p_0  \tr v \; | \; v(s) \geq   \min_{P \in \cP} T^{\pi'',P}v(s)  \; \forall \; s \in \cS, v \in \R^{\cS}\}.\]
Notice the similarity between these last two facts and the attainability property; this concludes the proof.
\end{proof}
In the proof above, we emphasize that we use Sion's minimax theorem to prove strong duality as in~\eqref{eq: lec 2 strong duality}, except that we can't apply it directly to $(\pi,P) \mapsto R(\pi,P)$ but we can apply it to $(\pi,P) \mapsto T^{\pi,P}v(s)$ for $s \in \cS$ and a fixed $v \in \R^{\cS}$.

The second proof is an extension of the original proof for the existence of Nash equilibria in finite games to the case of robust MDPs.
We write $2^{\cX}$ for the power set of a set $\cX$, i.e., the set of all subsets of $\cX$. We need the notion of closed graph for set-valued functions, defined next.
\begin{definition}\label{def: lec 1 closed graph}
    Let $\cX \subset \R^n$ be a compact set and $F:\cX \to 2^{\cX}$. We say that $F$ has a closed graph if for any sequences $\left(x_{n}\right)_{n \in \N}$ and $\left(y_{n}\right)_{n \in \N}$ such that $\forall \; n \in \N, y_{n} \in F(x_{n})$ and $\lim_{n \rightarrow +\infty} x_{n} = x \in \cX$ and $\lim_{n \rightarrow + \infty} y_{n} = y$, we have $y \in F(x)$.
\end{definition}
We can now state Kakutani's fixed point theorem.
\begin{theorem}[Kakutani's fixed point theorem]\label{th: lec 1 kakutani}
    Let $\cX \subset \R^{n}$ be a nonempty compact convex set. Let $F:\cX \rightarrow 2^{\cX}$ be a function with a closed graph and assume that for any $x \in \cX$, $F(x)$ is a nonempty closed convex set.
Then there exists $x \in \cX$ such that $x \in F(x)$.
\end{theorem}
\begin{proof}[Second proof of strong duality (based on Kakutani).]
    We want to use Kakutani's fixed point theorem to show the existence of a pair $(\pi\opt ,P\opt) \in \PiS \times \cP$ such that $\pi\opt$ is an optimal policy for $P\opt$ and $P\opt$ is a worst-case transition kernel for $\pi\opt $:
    \[\pi\opt  \in \arg \max_{\pi \in \PiS} R(\pi,P\opt), P\opt \in \arg \min_{P \in \cP} R(\pi\opt ,P).\]
    In fact we will show a stronger result: we can pick such $\pi\opt ,P\opt$ as {\em uniformly} optimal policies and worst-case transition probabilities.
    For $(\pi,P) \in \PiS \times \cP$ we let $\PiS\opt(P)$ be the set of uniformly optimal stationary policies when the transition kernel is $P$, and $\cP\opt(\pi)$ be the set of uniform worst-case transition probabilities for the policy $\pi$. Note that both sets are the sets of uniformly optimal policies for some MDP instances with finitely many states, compact action sets, and continuous rewards and transitions (for $\PiS\opt(P)$ this is direct, and for $\cP\opt(\pi)$ this is from the interpretation of the policy evaluation problem as an adversarial MDP). From \Cref{lem: lec 0 convex sets unf policy} we get that both $\PiS\opt(P)$ and $\cP\opt(\pi)$ are convex, closed, and always non-empty.

    Let $F:\PiS \times \cP \to 2^{\PiS \times \cP}$ be defined as $F(\bar{\pi},\bar{P}) = \PiS\opt(\bar{P}) \times \cP\opt(\bar{\pi})$.
    To apply \Cref{th: lec 1 kakutani}, it remains to show that $F$ has a closed graph.  Let $(\pi_n,P_{n}) \in \PiS \times \cP$ and $(\hat{\pi}_n,\hat{P}_n) \in \PiS\opt(P_n) \times \cP\opt(\pi_n)$ for $n \in \N$. Let $(\pi,P) = \lim_{n \rightarrow + \infty} (\pi_n,P_n)$ and $(\hat{\pi},\hat{P}) = \lim_{n \rightarrow + \infty} (\hat{\pi}_n,\hat{P}_n)$ (the limits are assumed to exist). We want to show that $(\hat{\pi},\hat{P}) \in \PiS\opt(P) \times \cP\opt(\pi)$. 
    From $(\hat{\pi}_n,\hat{P}_n) \in \PiS\opt (P_n) \times \cP\opt(\pi_n)$, we know that for each $n \in \N$ we have
    \begin{align*}
T^{\hat{\pi}_n,P_n}v^{\hat{\pi}_{n},P_n}(s) \geq T^{\pi',P_n}v^{\hat{\pi}_{n},P_n}(s), \forall \; s \in \cS, \forall \; \pi' \in \PiS, \\
T^{\pi_n,P'}v^{\pi_{n},\hat{P}_n}(s) \geq T^{\pi_n,\hat{P}_n}v^{\pi_{n},\hat{P}_n}(s), \forall \; s \in \cS, \forall \; P' \in \cP.
    \end{align*}    
    Since $(\pi,P) \mapsto T^{\pi,P}v^{\pi,P}(s)$ is a continuous function for each $s \in \cS$, taking the limit as $n \rightarrow + \infty$ in the two inequalities above gives 
        \begin{align*}
T^{\hat{\pi},P}v^{\hat{\pi},P}(s) \geq T^{\pi',P}v^{\hat{\pi},P}(s), \forall \; s \in \cS, \forall \; \pi' \in \PiS, \\
T^{\pi,P'}v^{\pi,\hat{P}}(s) \geq T^{\pi,\hat{P}}v^{\pi,\hat{P}}(s) , \forall \; s \in \cS, \forall \; P' \in \cP
    \end{align*} 
    which shows that $\hat{\pi} \in \PiS\opt(P)$ and $\hat{P} \in \cP\opt(\pi)$, so that $F$ has a closed graph. 

    We can therefore apply \Cref{th: lec 1 kakutani} to obtain the existence of a pair $(\pi\opt ,P\opt) \in \PiS \times \cP$ such that $\pi\opt  \in \PiS\opt(P\opt),P\opt \in \cP\opt(\pi\opt)$, from which we conclude that strong duality holds following~\Cref{prop: lec 2 NE SD}.
   \end{proof}

\begin{remark}
 Recall that from \Cref{chap: lec 0 mdp} we know that given a nominal MDP, the set of optimal policies may not be convex (see~\Cref{lem: lec 0 convex sets unf policy}). Hence we cannot apply Kakutani's fixed point theorem to the correspondence that maps $(\bar{\pi},\bar{P}) \in \PiS \times \cP$ to the Cartesian product of the set of optimal policies for $\bar{P}$ and the set of worst-case transition probabilities for $\bar{\pi}$, i.e., to 
    $\displaystyle (\bar{\pi},\bar{P}) \mapsto \arg \max_{\pi \in \PiS} R(\pi,\bar{P}) \times \arg \min_{P \in \cP} R(\bar{\pi},P)$. This is why we considered the sets of {\em uniformly} optimal policies and worst-case transition probabilities, which are convex sets from~\Cref{lem: lec 0 convex sets unf policy}.
\end{remark}
\subsection{Proof of uniform optimality and optimality conditions}\label{sec: lec 2 prf unf opt}
We first prove that stationary optimal (not necessarily uniformly optimal) policies exist. We have
\begin{align*}
 \sup_{\pi \in \PiH} \inf_{P  \in \cP} 
 R(\pi,P) & \leq  \inf_{P  \in \cP} \sup_{\pi \in \PiH}  R(\pi,P) \\
 & \leq \inf_{P  \in \cP} \max_{\pi \in \PiS} 
 R(\pi,P) \\
 & = \max_{\pi \in \PiS} \min_{P  \in \cP}  
 R(\pi,P) \\
 &\leq \sup_{\pi \in \PiH} \inf_{P  \in \cP} 
 R(\pi,P)
\end{align*}
where the first inequality is always true from weak duality, the second inequality comes from the inner problem being a nominal MDP and therefore $\sup_{\pi \in \PiH} 
 R(\pi,P) = \max_{\pi \in \PiS} 
 R(\pi,P)$ for any fixed $P  \in \cP$, the equality comes from the strong duality result, and the last inequality comes from $\PiS \subset \PiH$. We conclude that all terms in the equations above are equal, hence $\displaystyle \sup_{\pi \in \PiH} \inf_{P  \in \cP} 
 R(\pi,P) = \max_{\pi \in \PiS} \inf_{P  \in \cP} 
 R(\pi,P)$ so that there exists a stationary optimal policy. \Cref{prop: lec 2 maxmin piopt} shows that the problem $\max_{\pi \in \PiS} \min_{P  \in \cP} 
 R(\pi,P)$ admits stationary uniformly optimal policies that are characterized as solutions to $T^{\pi}v\opt = Tv\opt$ (as an equation in $\pi \in \PiS$). The statement about $P\opt \in \cP$ such that $T^{P\opt}v\opt = Tv\opt$ follows analogously by \eqref{eq: lec 2 best response} in \Cref{prop: lec 2 NE SD}.
\subsection{Computing optimal policies}\label{sec: lec 2 algo}
\Cref{th: lec 2 main rmdps} shows that we can compute optimal policies for~\eqref{eq: lec 1 rmdp} by computing the fixed point $v\opt \in \R^{\cS}$ of the robust Bellman operator $T$ defined in~\eqref{eq: lec 2 T}. For nominal MDPs, we saw in \Cref{sec: lec 0 algo} that this can be done efficiently with value iteration, policy iteration, and linear programming. We describe the corresponding algorithms for robust MDPs next.
\paragraph{Iterative algorithms.}
The proofs of convergence of value iteration and policy iteration only rely on the fact that we want to compute the fixed point of an optimality operator that is attainable by a set of monotone contractions. In particular, the family $\{T^{\pi} \; | \; \pi \in \PiS\}$ satisfies the attainability property with respect to $T$:
\begin{lemma}
    Let $\cP$ be a compact uncertainty set. For any $v \in \R^{\cS}$, there exists $\pi \in \PiS$ such that $Tv = T^{\pi} v$.
\end{lemma}
\begin{proof}
    The proof follows directly from the compactness of $\PiS$ 
, the fact that the $s$-th component of $T^{\pi}$ only depends on $\pi(s)$, the fact that the $s$-th component of $T$ optimizes over $\pi(s) \in \DA$, and the fact that $\PiS$ is a rectangular set: $\PiS = \times_{s \in \cS} \DA$. \end{proof}
Similarly to the case of nominal MDPs and adversarial MDPs, the convergence of value iteration and policy iteration for robust MDPs (\Cref{alg: lec 2 VI} and \Cref{alg: lec 2 PI}) follows from the general convergence theorems given in Appendix~\ref{app: algorithms}.

\begin{algorithm}
\caption{Value iteration for robust MDPs}\label{alg: lec 2 VI}
\begin{algorithmic}[1]
\State Initialize $v^{0} \in \R^{\cS}$.
\For{$t \in \N$}
    \State $v^{t+1} = T v^{t}$
\EndFor
\end{algorithmic}
\end{algorithm}

\begin{algorithm}
\caption{Policy iteration for robust MDPs}\label{alg: lec 2 PI}
\begin{algorithmic}[1]
\State Initialize $\pi^0 \in \PiS$, $v^0 = v^{\pi^0}$
\For{$t \in \N$}
    \State Choose $\pi^{t+1} \in \PiS$ such that $Tv^{t}(s) = T^{\pi^{t+1}}v^{t}(s)$ for all $s \in \cS$  (and equal to $\pi^t$ if we can)
    \State Choose $v^{t+1} = v^{\pi^{t+1}}$
    \If{$\pi^{t+1}=\pi^{t}$}{ return $\pi^{t+1}$}
    \EndIf
\EndFor
\end{algorithmic}
\end{algorithm}

\begin{proposition}\label{prop: lec 2 VI rate}
     Let $\cP$ be a compact convex s-rectangular uncertainty set. Let $(v^t)_{t \in \N}$ be the sequence of vectors generated by \Cref{alg: lec 2 VI}.
    \begin{enumerate}
        \item The sequence $(v^t)_{t \in \N}$ converges to $v\opt $ at a geometric rate of $\gamma$:
        \[\|v\opt  - v^{t}\|_{\infty} \leq \gamma^t \|v\opt  - v^{0}\|_{\infty}\]
    \item Let $\epsilon>0$. Assume that $\|v\opt  - v^{t}\|_{\infty} \leq \epsilon$ for some $t \in \N$ and let $\pi \in \PiS$ such that $T^{\pi}v^t = Tv^t$. Then $\pi$ is a uniform $2\epsilon \gamma (1-\gamma)^{-1}$-optimal policy: $\| v\opt - v^{\pi}\|_{\infty} \leq 2\epsilon \gamma (1-\gamma)^{-1}$.

    Additionally, we can check that $\|v\opt- v^{t}\|_{\infty} \leq \epsilon$ by checking $\|T v^{t}-v^{t}\|_{\infty} \leq \epsilon (1-\gamma)$, and that last inequality holds after at most $\lfloor \frac{1}{1-\gamma} \log\left(\frac{\|v^1-v^0\|_{\infty}}{\epsilon (1-\gamma)}\right) \rfloor + 1$ iterations of \Cref{alg: lec 2 VI}.

    \end{enumerate}
\end{proposition}

\begin{proposition}\label{prop: lec 2 PI rate}
 Let $\cP$ be a compact convex s-rectangular uncertainty set. Let $(\pi^t)_{t \in \N}$ be the sequence of policies generated by~\Cref{alg: lec 2 PI}. Then 
         $\| v\opt - v^{\pi^t} \|_{\infty} \leq \gamma^t \| v\opt-v^{\pi^0}\|_{\infty}.$
\end{proposition}
A more refined convergence analysis can be given for policy iteration (\Cref{alg: lec 2 PI}), see the references in~\Cref{sec: lec 2 ref}.
Note that \Cref{alg: lec 2 VI} and \Cref{alg: lec 2 PI} require estimating the robust Bellman operator $T$, which is a maxmin optimization problem with bilinear objective in $(\pi,P)$ over the compact convex sets $\PiS$ and $\cP$, as can be seen from its definition~\eqref{eq: lec 2 T}. This family of optimization problems can be solved efficiently via a variety of methods. Without exploiting any further the structure of the problem, saddle-point optimization problems can be solved via first-order methods (mirror descent, mirror prox, extragradient algorithm, etc.). Over the years there have been several papers providing more efficient algorithms in the case where $\cP$ has some additional structure, e.g., if $\cP$ is based on some $\ell_{p}$-distance from some nominal kernel, see \Cref{sec: lec 2 ref} for the references there. Note that \Cref{alg: lec 2 PI} also requires computing a worst-case value function of the current policy at every iteration, which is equivalent to solving an MDP (following \Cref{cor: lec 2 adversarial mdp - alg}); thus we note that \Cref{alg: lec 2 PI} is a priori much more computationally intensive than \Cref{alg: lec 2 VI}.
\paragraph{Convex programming formulation.}
The linear programming formulation for MDPs obtained in \Cref{prop: lec 0 LP MDP} was obtained in two steps. In the first step we used the contraction lemma (\Cref{lem: lec 0 contraction lemma}) to obtain a mathematical program whose set of optimal solutions was exactly the optimal value function. In the second step we used the componentwise convexity of the nominal Bellman operator to conclude that this mathematical program was a convex (linear) program. 

The first step still applies for robust MDPs: since $T$ is a monotone contraction, we still have that
\begin{align*}
    \{v\opt\} = \arg \min \{ \sum_{s \in \cS} v(s) \; | \; v \geq Tv \} = \arg \max \{ \sum_{s \in \cS} v(s) \; | \; v \leq Tv \}
\end{align*}
where the first equality is the contraction lemma and the second optimization program is given in~\Cref{sec: lec 0 fdt}. However $T$ need not be componentwise convex anymore, nor componentwise concave, because of the maxmin formulation in~\eqref{eq: lec 2 T}, so we cannot conclude that the mathematical programs above are convex optimization problems. Some convex formulations can compute $\epsilon$-optimal policies but this requires more work, see \Cref{sec: lec 2 ref}.
\subsection{A closer look at non-convex s-rectangular uncertainty sets*}\label{sec: lec 2 noncvx srec}
We get back to non-convex s-rectangular sets. In this case we can still try to compute a best policy among stationary ones, i.e., to solve
\[\max_{\pi \in \PiS} \min_{P \in \cP} R(\pi,P).\]
From \Cref{prop: lec 2 maxmin piopt} this can be done by (a) computing $v\opt \in \RS$ the fixed point of the robust Bellman operator $T$, and (b) computing $\pi\opt$ such that $T^{\pi\opt}v\opt = Tv\opt$. As long as we can evaluate $Tv$ for any $v \in \R^{\cS}$, we can run value iteration or policy iteration to compute $v\opt$, similarly as in \Cref{sec: lec 2 algo}. Whether or not we can evaluate $Tv$ depends on the tractability properties of $\max_{ \pi \in \PiS} \min_{P \in \cP} T^{\pi,P}v(s)$. Even for $\cP$ non-convex we get that
\begin{align*}
    \max_{ \pi \in \PiS} \min_{P \in \cP} T^{\pi,P}v(s) = \max_{ \pi \in \PiS} \min_{P \in \conv(\cP)} T^{\pi,P}v(s) = \min_{P \in \conv(\cP)} \max_{ \pi \in \PiS}  T^{\pi,P}v(s) 
\end{align*}
where the first equality is because the inner minimization has a linear objective over $\cP$, while the second equality is from Sion's minimax theorem. The first equality implies that the tractability properties of $Tv$ are the same for $\cP$, s-rectangular and potentially non-convex, and $\conv(\cP)$, which is s-rectangular and convex! 

We can now take a second look at strong duality for rectangular models. For compact s-rectangular uncertainty sets $\cP$ we have
    \begin{align*}
        \max_{\pi \in \PiS} \min_{P \in \cP} R(\pi,P) & = \max_{\pi \in \PiS} \min_{P \in \conv(\cP)} R(\pi,P) \\
        & = \min_{P \in \conv(\cP)} \max_{\pi \in \PiS}  R(\pi,P)\\
        & \leq \min_{P \in P} \max_{\pi \in \PiS}  R(\pi,P)
    \end{align*}
    where the first equality is because the inner minimization is the adversarial MDP, so we can replace stationary deterministic policies of the adversary (the set $\cP$) by stationary randomized policies of the adversary (which yields $\conv(\cP)$), the second equality is from strong duality for s-rectangular convex compact sets, and the inequality is from $\cP \subseteq \conv(\cP)$.
    This shows that strong duality holds (for potentially non-convex s-rectangular compact sets) if and only if
    $\displaystyle \min_{P \in \conv(\cP)} \max_{\pi \in \PiS}  R(\pi,P) = \min_{P \in P} \max_{\pi \in \PiS}  R(\pi,P)$. Of course, if $\cP$ is convex this last equality is true, but we'll see in \Cref{sec: lec 2 sarec} this last equality also holds for sa-rectangular compact uncertainty sets.
\subsection{Additional properties of optimal solutions*}
\paragraph{Worst-case transition probabilities.}
For s-rectangular compact convex uncertainty sets, we have seen that a uniform worst-case transition kernel $P^{\pi}$ of a stationary policy $\pi \in \PiS$ can be computed as the solution of an MDP with action set $\cP_s$ at each state $s \in \cS$. In particular, this implies that we can always compute the worst-case probabilities among the set $\ext(\cP)$ of extreme points of $\cP$. However, if we are interested in finding {\em uniformly} optimal worst-case transition probabilities, from \Cref{th: lec 2 main rmdps} we want to find $P\opt \in \cP$ such that $T^{P\opt}v\opt = Tv\opt$, i.e., we want $P\opt$ to be a solution of the following optimization problem:
\[\min_{P \in \cP} \max_{\pi \in \PiS} \sum_{a \in \cA} \pi(s,a) P(s,a) \tr (r(s,a) + \gamma v\opt).\]
In general the minimum of the above optimization problem need not be attained at an extreme point of $\cP$ (because we'll dualize the inner maximization program), so if we're interested in finding uniformly optimal worst-case transition probabilities we may have to consider non-extreme points of $\cP$.
\paragraph{Adjacent improving policies.} Recall two policies are {\em adjacent} if they differ at exactly one state (we defined this notion in \Cref{chap: lec 0 mdp}). The next lemma generalizes the statement in \Cref{sec: lec 0 additional prop}.
\begin{lemma}\label{lem: lec 2 adjacent policies}
    Let $\cP$ be compact convex s-rectangular. Let $\pi \in \PiS$ be non-uniformly optimal. Then there exists $\pi' \in \PiS$, that differs from $\pi$ in only one state, such that $v^{\pi'} \geq v^{\pi}$, with strict inequality in at least one component.
\end{lemma}
\begin{proof}
Note that we always have $v^{\pi} = T^{\pi} v^{\pi} \leq Tv^{\pi}$, where the equality follows from the definition of $v^{\pi} \in \R^{\cS}$ as the fixed point of $T^{\pi}$ and the componentwise inequality follows from the definition of $T$.
    Since $\pi$ is not uniformly optimal, we must have that there exists a state $\bar{s} \in \cS$ such that $v^{\pi}(\bar{s}) < Tv^{\pi}(\bar{s})$. Choose $\pi'$ the policy that coincides with $\pi$ everywhere except at $\bar{s}$, where $\pi'$ chooses a randomized action attaining the maximum in $Tv^{\pi}(\bar{s})$. Then we must have $T^{\pi'}v^{\pi} \geq v^{\pi}$, with strict inequality at $\bar{s}$. Monotonicity and contraction of $T^{\pi'}$ imply that $v^{\pi'} \geq v^{\pi}$. Additionally $v^{\pi'}(\bar{s}) = T^{\pi'}v^{\pi'}(\bar{s}) \geq T^{\pi'}v^{\pi}(\bar{s}) > v^{\pi}(\bar{s})$, where the equality is from the fixed point equation for $v^{\pi'}$, the inequality is from $v^{\pi'} \geq v^{\pi}$, and the strict inequality is from the choice of $\bar{s}$ and $\pi'$.
\end{proof}
\Cref{ex: lec 0 ex 3} shows that the same statement is false if we replace worst-case value functions by worst-case returns.

\paragraph{Support of a randomized policy.} We have seen in Theorem~\ref{th: lec 2 main rmdps} that optimal policies can be chosen in $\PiS$, so that we may need randomization between several actions to achieve optimality in s-rectangular robust MDPs. Recall that from a mathematical programming standpoint, randomization can be traced to the dualization of the inner minimum over $P$ in the maxmin expression of the robust Bellman operator $T$. One may wonder between how many actions must we randomize; we provide a partial answer in the next proposition and relate the size of the randomization support with the {\em affine dimension} of each $\cP_s$.
\begin{proposition}\label{prop: lec 2 support randomization}
Let $\cP$ be compact convex s-rectangular. Let $s \in \cS$ and $v \in \R^{\cS}$. Let $\dima(\cP_s) \in \N$ be the affine dimension of $\cP_s$, i.e., the dimension of the affine hull of $\cP_s \subset \DSA$ as an affine subspace of $\R^{|\cS| \times |\cA|}$. For any $v \in \RS$ and $s \in \cS$ there exists a maximizer of
\[\max_{\pi(s) \in \DA} \min_{P(s) \in \cP_s} \sum_{a \in \cA} \pi(s,a) P(s,a) \tr (r(s,a) + \gamma v)\]
with support of cardinality at most $\min \{|\cA|, \dima(\cP_s)+1\}$.
\end{proposition}
\begin{proof}
Here we overload notation to write the maxmin problem from above as \[\max_{\pi(s) \in \DA} \min_{P(s) \in \cP_s} T^{\pi(s),P(s)}v(s)\] (in principle writing $``T^{\pi,P}"$ requires defining a full policy $\pi$ and kernel $P$, but the $s$-th component of $T^{\pi,P}$ only depends on the $s$-th component of $\pi$ and $P$). By strong duality of minmax problems with convex compact domains and bilinear objectives we have $\max_{\pi(s) \in \DA} \min_{P(s) \in \cP_s} T^{\pi(s),P(s)}v(s)=\min_{P(s) \in \cP_s} \max_{\pi(s) \in \DA} T^{\pi(s),P(s)}v(s)$. Let $P\opt(s) \in \cP_s$ in $\arg \min_{P(s) \in \cP_s} \max_{\pi(s) \in \DA} T^{\pi(s),P(s)}v(s)$. The minimax theorem for
\[\min_{P(s) \in \cP_s} \max_{\pi(s) \in \DA} T^{\pi(s),P(s)}v(s)\] gives the existence of some $\lambda \in \DA$ such that $\lambda(a)=0$ if $a \notin \arg \max_{a' \in \cA} T^{a',P\opt(s)}v(s)$ (noting $T^{a',P\opt(s)}v(s) = P\opt(s,a')\tr (r(s,a') + \gamma v)$) and
\[\sum_{a \in \cA} \lambda(a) \left( T^{a,P(s)}v(s) - T^{a,P\opt(s)}v(s) \right) \geq 0, \forall \; P(s) \in \cP_s.\]
Now note that $T^{a,P(s)}v(s) - T^{a,P\opt(s)}v(s) = T^{a,P(s)-P\opt(s)}v(s)$ because the same action $a$ is chosen at state $s$ in both operators. Hence we can interpret the left-hand side of the above inequality as a convex combination of elements of dimensions at most $\dima(\cP_s)$. From Carath{\'e}odory's theorem we can choose such $\lambda \in \Delta(\cA)$ with at most $\dima(\cP_s)+1$ nonzero entries.

Now we choose $\pi\opt(s) = \lambda$. The first-order condition gives that for any $P(s) \in \cP_s$ we have
\[\sum_{a \in \cA} \pi\opt(s,a)  T^{a,P(s)}v(s) \geq \sum_{a \in \cA} \pi\opt(s,a) T^{a,P\opt(s)}v(s) = \max_{a \in \cA} T^{a,P\opt(s)}v(s)\]
where the equality is because $\pi\opt(s)$ is supported on the active actions at $P\opt(s)$. Conversely for any $\pi(s) \in \DA$ we have
\[\min_{P(s) \in \cP_s} \sum_{a \in \cA} \pi(s,a) T^{a,P(s)}v(s) \leq \sum_{a \in \cA} \pi(s,a) T^{a,P\opt(s)}v(s) \leq \max_{a \in \cA} T^{a,P\opt(s)}v(s)\]
hence we conclude that if we construct $\pi\opt$ this way for $v\opt$ the fixed point of $T$ and for each state $s \in \cS$ we get an optimal policy, with support at most $\min \{|\cA|,\dima(\cP_s)+1\}$.
\end{proof}
\section{The case of sa-rectangular models}\label{sec: lec 2 sarec}
The previous sections have highlighted the good tractability and game-theoretic properties of s-rectangular models. The following {\em sa-rectangular} model imposes some additional structural assumption that the adversary can choose their decisions (worst-case transition probabilities) independently across state-action pairs.
\begin{definition}\label{dec: lec 2 sarec}
    An uncertainty set $\cP \subset \DSSA$ is sa-rectangular if we can write $\cP = \times_{(s,a) \in \cS \times \cA} \cP_{sa}$ for some sets $\cP_{sa} \subseteq \DS$ for each $(s,a) \in \cS \times \cA$.
\end{definition}
A typical example of sa-rectangular uncertainty set is based on $\cP_{sa} \subseteq \DS$ as
\[\cP_{sa} = \{ p \in \DS \; | \; d(p,\hat{P}(s,a)) \leq \epsilon_{sa}\}\]
for some reference transitions $\hat{P} \in \DSSA$. We can check ``visually" the sa-rectangularity assumption (for robust MDP instances represented as in the examples from \Cref{chap: lec 0 mdp} and \Cref{chap: lec 1 robust MDP}) if the parameters present in each figure representing the transition for each state-action pair are independent.
Note that sa-rectangular sets are s-rectangular, hence all the favorable properties (strong duality, existence of stationary optimal policies) outlined in the previous section apply (under the compactness and convexity assumptions). However, we will show below that for sa-rectangular uncertainty sets  $\cP$ we do {\em not} need $\cP$ to be convex to conclude as in \Cref{th: lec 2 main rmdps}.
\paragraph{Additional properties.}
The sa-rectangularity assumption provides some structural results on the Bellman operators. Let $\pi \in \PiS$, then
\begin{align*}
        T^{\pi} v(s) = \min_{P \in \cP} T^{\pi,P}v(s) & = \min_{(P(s,a))_{s,a} \in \times_{(s,a)} \cP_{sa}} \sum_{a \in \cA} \pi(s,a) P(s,a)\tr ( r(s,a) + \gamma v) \\
        & = \sum_{a \in \cA} \pi(s,a) \min_{P(s,a) \in \cP_{sa}} P(s,a)\tr ( r(s,a) + \gamma v)
\end{align*}
that is, we can interchange the summation over $a \in \cA$ and the minimization program over $P \in \cP$. The first equality above holds by definition of $T^{\pi}$, the second equality just uses the definition of $T^{\pi,P}$, and the third equality is because each term in the summation over $a \in \cA$ only depends on $P(s,a)$ and $\cP$ is a Cartesian product over $(s,a)$ (by the sa-rectangularity assumption). In particular, this implies that
   \begin{align}
        T v(s) = \max_{\pi \in \PiS} \min_{P \in \cP} T^{\pi,P}v(s) 
        & = \max_{\pi \in \PiS} \sum_{a \in \cA} \pi(s,a) \min_{P(s,a) \in \cP_{sa}} P(s,a)\tr ( r(s,a) + \gamma v) \nonumber \\
        & = \max_{a \in \cA} \min_{P(s,a) \in \cP_{sa}} P(s,a)\tr ( r(s,a) + \gamma v) \label{eq: lec 2 T sarec easier}
\end{align} 
i.e., to evaluate $Tv(s)$ we do not need to solve a ``true" saddle-point problem but we only need to compute the maximum of $|\cA|$ real numbers given by $\min_{P(s,a) \in \cP_{sa}} P(s,a)\tr ( r(s,a) + \gamma v)$ for each $a \in \cA$. In passing, we've obtained that 
\begin{equation}\label{eq: lec 2 Tv sarec}
    T v(s) = \max_{\pi \in \PiSD} \min_{P \in \cP} T^{\pi,P}v(s) 
\end{equation}
 for sa-rectangular compact uncertainty sets (this will be useful in the next section).
Since \Cref{th: lec 2 main rmdps} gives that the set of stationary uniformly optimal policies is the set of policies $\bar{\pi}$ such that $T^{\bar{\pi}}v\opt = Tv\opt$, the next lemma follows directly from \eqref{eq: lec 2 Tv sarec}.
\begin{lemma}\label{lem: lec 2 sa rec unif opt deter}
    Let $\cP \subset \DSSA$ be a convex compact sa-rectangular uncertainty set. Then there exists a stationary {\em deterministic} uniformly optimal policy, which can be found as any $\pi \in \PiSD$ such that $T^{\pi}v\opt = T v\opt$.
\end{lemma}
We give above the results for sa-rectangular compact {\em convex} uncertainty sets, because we required convexity for s-rectangular sets. However the convexity assumption is unnecessary for sa-rectangular robust MDPs, as we show next.
\begin{lemma}\label{lem: lec 2 T sa rec strong duality}
    Let $\cP$ be sa-rectangular compact. Let $v \in \R^{\cS}$ and $s \in \cS$. Then
    \begin{equation}\label{eq: lec 2 T sa rec strong duality}
    T v(s) = \max_{\pi \in \PiSD} \min_{P \in \cP} T^{\pi,P}v(s) = \min_{P \in \cP} \max_{\pi \in \PiSD} T^{\pi,P}v(s).
\end{equation}
\end{lemma}
\begin{proof}
From~\eqref{eq: lec 2 T sarec easier} we know that 
\[\max_{\pi \in \PiSD} \min_{P \in \cP} T^{\pi,P}v(s) = \max_{a \in \cA} \min_{P(s,a) \in \cP_{sa}} P(s,a) \tr (r(s,a) + \gamma v).\]
From weak duality we always have
\[\max_{a \in \cA} \min_{P(s,a) \in \cP_{sa}} P(s,a) \tr (r(s,a) + \gamma v) \leq \min_{P \in \cP} \max_{a \in \cA}  P(s,a) \tr (r(s,a) + \gamma v).\]
To show the equality we just need to find $\tilde{P} \in \cP$ such that, for all $s \in \cS$, we have
\[\max_{a \in \cA}  \tilde{P}(s,a) \tr (r(s,a) + \gamma v) = \max_{a \in \cA} \min_{P(s,a) \in \cP_{sa}} P(s,a) \tr (r(s,a) + \gamma v).\]
Select $\tilde{P} \in \cP$ such that, for each $(s,a) \in \cS \times \cA$ we have 
\begin{equation}\label{eq: lec 2 tildeP argmin sarec}
    \tilde{P}(s,a) \in \arg \min_{P(s,a) \in \cP_{sa}} P(s,a) \tr (r(s,a) + \gamma v).
\end{equation}
Since $\cP$ is sa-rectangular we get that $\tilde{P} \in \cP$ and then by construction we have
\[\max_{a \in \cA}  \tilde{P}(s,a) \tr (r(s,a) + \gamma v) = \max_{a \in \cA} \min_{P(s,a) \in \cP_{sa}} P(s,a) \tr (r(s,a) + \gamma v)\]
which shows that \eqref{eq: lec 2 T sa rec strong duality} holds.
\end{proof}
The proof of Lemma~\ref{lem: lec 2 T sa rec strong duality} may seem ``magic" at first glance, and it may seem surprising that we can get a strong duality result without the convexity assumption. Here is a game-theoretic interpretation of the proof: since the adversary can independently minimize each optimization problem over all pairs $(s,a)$, the adversary does not need to ``adapt" to the policy chosen by the decision maker. The fact that this ``worst-case transition kernel" (as defined in~\eqref{eq: lec 2 tildeP argmin sarec}) can be chosen independently of the chosen policy ensures that the value of the game is the same when the adversary plays first (i.e., when solving the minmax) and when the adversary plays second (i.e., when solving the maxmin). Still, note that the solution $\tilde{P}$ defined in~\eqref{eq: lec 2 tildeP argmin sarec} depends on the vector $v\in \R^{\cS}$, so we are not claiming that there exists some transition kernel that attains the worst-case return simultaneously across all policies. \Cref{lem: lec 2 T sa rec strong duality} easily implies that strong duality holds for $\cP$ a compact sa-rectangular uncertainty set.
\begin{proposition}\label{prop: lec 2 sa-rec strong duality}
    Let $\cP \subset \DSSA$ be compact sa-rectangular. We have
    \begin{equation}\label{eq: lec 2 strong duality sarec}
    \max_{\pi \in \PiSD} \min_{P \in \cP} R(\pi,P) =  \min_{P \in \cP} \max_{\pi \in \PiSD}  R(\pi,P).
\end{equation}
\end{proposition}
\begin{proof}
Note that $\min_{P \in \conv(\cP)} \max_{\pi \in \PiS}  R(\pi,P)$ admits an optimal solution at some $P \in \conv(\cP)$ that satisfies, for all $s \in \cS$, $P  \in \arg \min_{P \in \conv(\cP)} \max_{\pi \in \PiS} T^{\pi,P}v\opt(s)$.
    But for $\cP$ compact sa-rectangular we have
    \begin{align*}
        \min_{P \in \conv(\cP)} \max_{\pi \in \PiS} T^{\pi,P}v\opt(s) & \leq \min_{P \in \cP} \max_{\pi \in \PiS} T^{\pi,P}v\opt(s) \\
        & = \min_{P \in \cP} \max_{\pi \in \PiSD} T^{\pi,P}v\opt(s) \\
        & = \max_{\pi \in \PiSD} \min_{P \in \cP} T^{\pi,P}v\opt(s) \\
        & = \max_{\pi \in \PiSD} \min_{P \in \conv(\cP)} T^{\pi,P}v\opt(s) \\
        & \leq \max_{\pi \in \PiS} \min_{P \in \conv(\cP)} T^{\pi,P}v\opt(s)
    \end{align*}
where the inequality is from $\cP \subset \conv(\cP)$, the first equality is from the inner maximization problem having a linear objective for each fixed $P \in \cP$, the second equality is from~\eqref{eq: lec 2 T sa rec strong duality}, the third equality is from the inner minimization being linear in $P$, and the last inequality is from $\PiSD \subset \PiS$. So all terms above are equal and $\displaystyle \min_{P \in \conv(\cP)} \max_{\pi \in \PiS}  R(\pi,P) = \min_{P \in P} \max_{\pi \in \PiS}  R(\pi,P)$ holds for sa-rectangular compact uncertainty sets, without requiring the convexity assumption.
\end{proof}
Therefore all the results from \Cref{th: lec 2 main rmdps} and \Cref{lem: lec 2 sa rec unif opt deter} are still true for sa-rectangular compact (not necessarily convex) uncertainty sets.
\subsection{Algorithms for sa-rectangular robust MDPs}\label{sec: lec 2 sa-rec alg}
\paragraph{Iterative algorithms.}

The value iteration and policy iteration for s-rectangular robust MDPs (\Cref{alg: lec 2 VI} and \Cref{alg: lec 2 PI}) still apply for sa-rectangular models, the main difference being that it is usually more efficient to evaluate the operator $T$ in the case of sa-rectangular models, as explained after \eqref{eq: lec 2 T sarec easier}. There have been many papers working on giving fast algorithms for evaluating $Tv$ and $T^{\pi}v$ for sa-rectangular models with specific structural assumptions, and we refer to \Cref{sec: lec 2 ref} for more details.
\paragraph{Convex formulation.}
In the case of sa-rectangular models, under the assumption that each $\cP_{sa}$ is a polytope, we can obtain a convex formulation for computing an optimal policy relatively easily. The contraction lemma (\Cref{lem: lec 0 contraction lemma}, as extended in \Cref{sec: lec 0 fdt}) states that 
\[ \{v\opt\} = \arg \max \left\{ \sum_{s \in \cS} v(s) \; | \; v \leq Tv, v \in \R^{\cS}\right\}.\]
Let $\cV \coloneqq \{v \in \R^{\cS} \; | \; v \leq Tv\}$. Since $\displaystyle T v(s) = \max_{\pi \in \PiSD} \min_{P \in \cP} T^{\pi,P}v(s)$ we get that $\cV = \cup_{\pi \in \PiSD} \cV^{\pi}$, with $\cV^{\pi} \coloneqq \{v \in \R^{\cS} \; | \; v \leq T^{\pi}v\}$. For $\cVb^{\pi} = \cV^{\pi} \cap [-M,M]^{\cS}$ with $M = (1-\gamma)^{-1} \max_{s,a,s'} |r(s,a,s')|$ we still get $v^{\pi} \in \cVb^{\pi}$ for each $\pi \in \PiSD$ so that
\[\{v\opt\} = \arg \max \left\{ \sum_{s \in \cS} v(s) \; | \; v \in \conv(\cup_{\pi \in \PiSD} \cVb^{\pi}) \right\}\]
since optimizing a linear function over a set or over its convex hull yields the same value.

Now a theorem from Balas (see \Cref{sec: lec 2 ref}) provides the following polytope reformulation for the convex hull of a union of polytopes: for $Q_i = \{x \; | \; A_i x \leq b_i\}$ with $i =1,...,N$ we have that
\[\conv\left( \cup_{i=1}^{N} Q_i \right) = \{z \; | \; z = \sum_{i =1}^{N} z^i, A_i z^{i} \leq \lambda_i b_i, \forall \; i, \lambda \in \Delta(\{1,...,N\})\}.\]
Applying this formula for $\conv(\cup_{\pi \in \PiSD} \cVb^{\pi})$ when each $\cP_{sa}$ is a polytope yields a linear programming formulation for computing $v\opt$:
\begin{itemize}
    \item If $\cP_{sa}$ is represented as $\cP_{sa} = \conv\left(P^1(s,a), ..., P^{k_{sa}}(s,a)\right)$ for $k_{sa} \in \N$ for each state-action pair $(s,a)$ we get the following optimization program:
    \begin{align*}
        \max \; & \sum_{s \in \cS} v(s) \\
        & v = \sum_{\pi \in \PiSD} v^{\pi} \\
        & v^{\pi}(s) \leq \lambda_\pi P^k(s,\pi(s))\tr r(s,\pi(s)) + \gamma P^k(s,\pi(s))\tr v^{\pi}, \forall \; \pi \in \PiSD, s \in \cS, k \in \{1,..., k_{s,\pi(s)}\}, \\ 
        & - M \lambda_{\pi} \leq v^{\pi}(s) \leq M \lambda_{\pi}, \forall \; \pi \in \PiSD, s \in \cS, \\
        & \lambda \in \Delta(\PiSD), v^{\pi} \in \RS,\forall \; \pi \in \PiSD.
    \end{align*}
    Note that in the optimization program above $v^{\pi}$ is an auxiliary variable and need not be equal to a worst-case value function of $\pi \in \PiSD$.
    \item If $\cP_{sa}$ is represented as $\cP_{sa} = \{p \; | \; A_{sa}\tr p = b_{sa}, p \geq 0 \}$ then from strong duality we'd get $\min\{c_{sa}\tr p \; | \; A_{sa}\tr p = b_{sa}, p \geq 0 \} = \max \{ b_{sa} \tr y \; | \; A_{sa} y \leq c_{sa} \}$ hence we would arrive at the following linear program with unique optimal solution $v\opt$ (for the $v$-component of the solutions):
   \begin{align*}
        \max \; & \sum_{s \in \cS} v(s) \\
        & v = \sum_{\pi \in \PiSD} v^{\pi} \\
        & A_{s,\pi(s)}y^{\pi}_{s} \leq \lambda_\pi r(s,\pi(s)) + \gamma v^{\pi}, \forall \; \pi \in \PiSD, s \in \cS, \\
        & v^{\pi}(s) \leq b_{s \pi(s)}\tr  y^{\pi}_{s}, \forall \; \pi \in \PiSD, s \in \cS,\\
        & - M \lambda_{\pi} \leq v^{\pi}(s) \leq M \lambda_{\pi}, \forall \; \pi \in \PiSD, s \in \cS, \\
        & \lambda \in \Delta(\PiSD), v^{\pi} \in \RS, y^{\pi}_{s} \in \R^{m},\forall \; \pi \in \PiSD, s \in \cS.
    \end{align*}
\end{itemize}
Note that in both cases since $|\PiSD| = |\cA|^{|\cS|}$ the number of variables in the linear programming formulation is exponential in the number of states, so these reformulations are most likely impractical.
\section{Minimality of the rectangularity assumption*}\label{sec: lec 2 minimality rec}
The above sections have highlighted several of the favorable properties of rectangular (s-rectangular and sa-rectangular) uncertainty sets. It is natural to ask if these properties would still hold for other non-rectangular models of uncertainty. To formalize this question, we note that ideally we would like our model of uncertainty to satisfy the desirable properties (e.g., existence of stationary optimal policies, or existence of a robust Bellman equation for the optimal value functions) ``across all realistic situations we may face", i.e., for any instances we may have to solve. Hence in this section we characterize the class of uncertainty sets that have good properties {\em across all choices of the instantaneous rewards $r \in \R^{\cS \times \cA \times \cS}$, initial distribution $p_0 \in \DS$, and discount factor $\gamma \in [0,1)$}. In \Cref{sec: lec 2 ssp rec} we characterize the class of uncertainty sets for which the robust MDP problem~\eqref{eq: lec 1 rmdp} can be solved via a Bellman equation (also known as a {\em dynamic programming principle}), and in \Cref{sec: lec 2 unf opt rec} we characterize the class of uncertainty sets for which the robust MDP problem~\eqref{eq: lec 1 rmdp} always admits a stationary uniformly optimal policy. We will see that each of these properties separately characterizes the rectangularity assumptions.
\subsection{Minimality for the dynamic programming principle*}\label{sec: lec 2 ssp rec}
Let $\cP \subset \DSSA$ be a compact uncertainty set.
Note that we can always define the robust Bellman operator $T$ as in~\eqref{eq: lec 2 T} (even when $\cP$ is not s-rectangular), and $T$ is always a contraction, so that its fixed point $v\opt \in \R^{\cS}$ always exists. The main objective of this section is to characterize the set $\cP$ such that $v\opt$ {\em always} encodes the optimal return of the robust MDP problem, e.g., for which the following equality holds across all $r \in \R^{\cS \times \cA \times \cS}$, all $p_0 \in \DS$ and all $\gamma \in [0,1)$:
\begin{equation}\label{eq: lec 2 dp principle}
\sup_{\pi \in \PiH} \min_{P \in \cP} R(\pi,P) = p_0 \tr v\opt
\end{equation}
Equality~\eqref{eq: lec 2 dp principle} is referred to as a {\em dynamic programming principle}, as it relates the optimal objective of the robust MDP problem to the fixed point of the robust Bellman operator $T$. The next theorem uniquely characterizes the uncertainty set satisfying the dynamic programming principle across all instances. We give a slightly more general statement that allows us to characterize this property beyond convex sets.
\begin{theorem}\label{th: lec 2 dpp srec}
Let $\cP$ be a compact uncertainty set. The following statements are equivalent:
\begin{enumerate}
    \item {\em Dynamic programming principle.} For any $r \in \R^{\cS \times \cA \times \cS}$, any $p_0 \in \DS$ and any $\gamma \in [0,1)$, we have that $\displaystyle \sup_{\pi \in \PiH} \min_{P \in \conv(\cP)} R(\pi,P) = p_0 \tr v\opt$, with $v\opt \in \R^{\cS}$ the unique fixed point of the robust Bellman operator as defined in \eqref{eq: lec 2 T}.
    \item {\em Simultaneous solvability property.} The following simultaneous solvability property (SSP) holds:
    \begin{equation}\label{eq: lec 2 ssp srec}
        \cap_{s \in \cS} \arg \min_{P \in \cP} \sum_{a \in \cA} P(s,a)\tr U(s,a) \neq \emptyset, \quad \forall \; U \in \R^{\cS \times \cA \times \cS}.
    \end{equation}
    \item {\em s-rectangularity.} $\conv(\cP)$ is s-rectangular.
\end{enumerate}
\end{theorem}
\begin{proof}[Proof of \Cref{th: lec 2 dpp srec}]
{\em Dynamic programming implies simultaneous solvability.}

    Let $U \in \R^{\cS \times \cA \times \cS}$ and $s \in \cS$. Let $f_s:\DA\to \R$ be defined as, for $x \in \DA$,
    \[f_s(x) = \min_{P \in \cP} \sum_{a \in \cA} x(a) P(s,a)\tr U(s,a).\]
    Note that $f_s$ is continuous and concave (as the pointwise minimum of some linear functions). Let $\bar{x} \in \DA$ be the uniform distribution over $\cA$. The supporting hyperplane theorem (see Appendix~\ref{app:sec analysis}) gives that there exists a supergradient $g_s \in \R^{\cA}$ such that
    \begin{equation}\label{eq: lec 2 supergradient ineq}
        f_s(x) \leq f_s(\bar{x}) + g_s \tr (x-\bar{x}), \forall \; x \in \DA.
    \end{equation}
    Choose $\gamma=0$ and define the rewards as $r(s,a,s') = U(s,a,s') - g_s(a)$, so that 
    \[\min_{P \in \cP} \sum_{a \in \cA} x(a) P(s,a)\tr r(s,a) = f_s(x) - g_s \tr x.\]
    From \eqref{eq: lec 2 supergradient ineq} we get that $\bar{x}$ is an optimal policy so that $v\opt(s) = f_s (\bar{x}) - g_s \tr \bar{x}.$ For $p_0 \in \Delta(\cS)$ the uniform distribution over $\cS$, the dynamic programming principle gives $\displaystyle \sup_{\pi \in \PiH} \min_{P \in \conv(\cP)} R(\pi,P) = p_0 \tr v\opt$, and we conclude that the policy $\bar{\pi}$ choosing the uniform distribution $\bar{x}$ at every state is optimal, since we then get
    \[p_0 \tr v\opt = \sum_{s \in \cS} p_0(s) \min_{P \in \cP} v^{\bar{\pi},P}(s) \leq \min_{P \in \conv(\cP)} R(\bar{\pi},P) \leq \sup_{\pi \in \PiH} \min_{P \in \conv(\cP)} R(\pi,P) = p_0 \tr v\opt\]
    where the last equality is from the dynamic programming principle. Note that with $\gamma=0$ we have that $P \mapsto R(\bar{\pi},P)$ is linear so that $\min_{P \in \cP} R(\bar{\pi},P)=\min_{P \in \conv(\cP)} R(\bar{\pi},P)$. Hence for $\bar{P} \in \arg \min_{P \in \cP} R(\bar{\pi},P)$ a worst-case transition kernel for $\bar{\pi}$ we get that 
    \[0 = \sum_{s \in \cS} p_0 (s) \left( \sum_{a \in \cA} \bar{\pi}(s,a) \bar{P}(s,a)\tr U(s,a) - f_s (\bar{x})\right).\]
Each term of the sum is nonnegative, hence all terms are $0$. Since $\bar{\pi}(s,a)=1/|\cA|$ is independent of $s$, we get that $\bar{P} \in \cap_{s \in \cS} \arg \min_{P \in \cP} \sum_{a \in \cA} P(s,a) \tr U(s,a)$, so the simultaneous solvability property \eqref{eq: lec 2 ssp srec} holds.
\vspace{2mm}

{\em Simultaneous solvability implies s-rectangularity of the convex hull.}
To show that $\conv(\cP)$ is s-rectangular, we need to show that it is equal to the Cartesian product over its $s$-th marginals.
Let $\cPC_{s} = \{ P(s) \; | \; P \in \conv(\cP)\}$ be the set of $s$-th marginals of $\conv(\cP)$ (we use a superscript here to distinguish this set from the $s$-th marginal of $\cP$ defined for an s-rectangular set in \Cref{def: lec 2 srec}). We want to show that $\conv(\cP) = \times_{s \in \cS} \cPC_s$. Of course we always have $\conv(\cP) \subset \times_{s \in \cS} \cPC_s$ so we need to prove the reverse inclusion. 

We first note that the simultaneous solvability property implies that, for each $U \in \R^{\cS \times \cA \times \cS}$, we have
\begin{equation}\label{eq: lec 2 implication ssp srec}
    \min_{P \in \conv(\cP)} \sum_{s \in \cS} \sum_{a \in \cA} P(s,a)\tr U(s,a) = \min_{P \in \times_{s \in \cS} \cPC_s} \sum_{s \in \cS} \sum_{a \in \cA} P(s,a)\tr U(s,a)
\end{equation}
This is because (a) we are optimizing linear forms so we can replace minimization over $\cP$ by minimization over $\conv(\cP)$ (b) the simultaneous solvability property states that an optimal solution of the left-hand side above can be chosen independently for each state $s \in \cS$.

Now for the sake of contradiction, suppose that $\conv(\cP) \neq \times_{s \in \cS} \cPC_s$. Then we can find $P' \in \times_{s \in \cS} \cPC_s$ with $P' \notin \conv(\cP)$. Let $P'' \in \conv(\cP)$ be the orthogonal projection of $P'$ onto $\conv(\cP)$. 
The first-order optimality condition for the orthogonal projection (see Appendix~\ref{app:sec analysis}) implies that $\langle P'' - P', P \rangle \geq \langle P'' - P', P'' \rangle$ for any $P \in \conv(\cP)$.
However we get
\[\langle P'' - P', P' \rangle = \langle P'' - P', P'' \rangle - \| P'-P''\|_{2}^{2} < \langle P'' - P', P'' \rangle.\]
Therefore the linear form $P \mapsto \langle P'' - P', P \rangle$ has a strictly smaller minimum over $\times_{s \in \cS} \cPC_s$ than over $\conv(\cP)$, which contradicts \eqref{eq: lec 2 implication ssp srec}. Hence we must have $\times_{s \in \cS} \cP_{s}^{\sf c}=\conv(\cP)$, and $\conv(\cP)$ is s-rectangular.
\vspace{2mm}

{\em s-rectangularity implies dynamic programming principle.} This is the result from \Cref{th: lec 2 main rmdps} applied to $\conv(\cP)$.
\end{proof}
Note that the first implication only requires compactness; convexity of $\conv(\cP)$ is necessary in the second and third implications. Recall that $T$, as defined in \eqref{eq: lec 2 T}, takes componentwise minimization over $P \in \cP$ for each state $s \in \cS$. Hence the operator $T$ only ``sees" an {\em s-rectangular extension} of $\cP$, and the reader should not be surprised that the dynamic programming principle implies that $\conv(\cP)$ is s-rectangular.
The same lines of argument yield the following equivalent results for sa-rectangular uncertainty sets, except that the Bellman operator now has more structure. In particular, let $\hat{T}:\R^{\cS} \to \R^{\cS}$ be defined as
\begin{equation}\label{eq: lec 2 hat T}
   \hat{T} v(s) = \max_{\pi \in \PiSD} \min_{P \in \cP} T^{\pi,P}v(s), \quad \forall \; v \in \R^{\cS}, s \in \cS.
\end{equation}
Note that we do not need any rectangularity assumption on $\cP$ to define $\hat{T}$. However, since $\hat{T}$ is defined componentwise and the outer maximization is over $\PiSD$, we get that $\hat{T}$ only depends on the $(s,a)$-marginals of $\cP$. Hence in some sense $\hat{T}$ only depends on a {\em sa-rectangular extension} of $\cP$.
\begin{theorem}\label{th: lec 2 dpp sa rec}
Let $\cP$ be a compact uncertainty set. The following statements are equivalent:
\begin{enumerate}
    \item {\em Dynamic programming principle.} For any $r \in \R^{\cS \times \cA \times \cS}$, any $p_0 \in \DS$ and any $\gamma \in [0,1)$, we have that $\displaystyle \sup_{\pi \in \PiH} \min_{P \in \cP} R(\pi,P) = p_0 \tr \hat{v}\opt$, with $\hat{v}\opt \in \R^{\cS}$ the unique fixed point of the robust Bellman operator $\hat{T}$ as defined in \eqref{eq: lec 2 hat T}.
    \item {\em Simultaneous solvability property.} The following simultaneous solvability property (SSP) holds:
    \begin{equation}\label{eq: lec 2 ssp sarec}
        \cap_{(s,a) \in \cS \times \cA} \arg \min_{P \in \cP} P(s,a)\tr U(s,a) \neq \emptyset, \quad \forall \; U \in \R^{\cS \times \cA \times \cS}.
    \end{equation}
    \item {\em sa-rectangularity.} $\conv(\cP)$ is sa-rectangular.
\end{enumerate}
\end{theorem}
\begin{proof}
    {\em Dynamic programming implies simultaneous solvability.} Let $U \in \RSAS$ and use $\gamma=0$. Define $r(s,a,s') = U(s,a,s') - \min_{P \in \cP} P(s,a) \tr U(s,a)$ so that $\min_{P \in \cP} P(s,a) \tr r(s,a) = 0$. Note that this implies that $\hat{v}\opt(s)=0$ for each $s \in \cS$ (the worst-case value obtained from each state is $0$). Take $p_0$ the uniform initial distribution and $\bar{\pi}$ the uniform policy, so that 
    \[R(\bar{\pi},P) = \sum_{s,a} p_0(s)\bar{\pi}(s,a) P(s,a)\tr r(s,a) \geq 0\]
    and $\min_{P \in \cP} R(\bar{\pi},P)=0$ by construction of the rewards and $\hat{v}\opt=0_{\cS}$, since $p_0 \tr \hat{v}\opt=0$ and $p_0 \tr \hat{v}\opt = \sup_{\pi \in \PiH} \min_{P \in \cP} R(\pi,P)$ from the dynamic programming principle. For $\bar{P} \in \arg \min_{P \in \cP} R(\bar{\pi},P)$ we must have that for all state-action pair $(s,a)$ we have $\bar{P} \in \arg \min_{P \in \cP} P(s,a) \tr r(s,a) = \arg \min_{P \in \cP} P(s,a) \tr U(s,a)$, which shows that the simultaneous solvability property holds.

\vspace{2mm}
    {\em Simultaneous solvability implies sa-rectangularity of the convex hull.} This proof is similar to the second part of the proof of \Cref{th: lec 2 dpp srec}.

    \vspace{2mm}

    {\em sa-rectangularity implies dynamic programming principle.} 

    We first show that \eqref{eq: lec 2 ssp sarec} holds. Indeed for $U \in \RSAS$ we have
    \begin{align*}
        \min_{P \in \cP} \sum_{s,a} P(s,a) \tr U(s,a)  & = \min_{P \in \conv(\cP)} \sum_{s,a} P(s,a) \tr U(s,a) \\
        & =  \sum_{s,a} \min_{P \in \conv(\cP)} P(s,a) \tr U(s,a) \\
        & = \sum_{s,a} \min_{P \in \cP} P(s,a) \tr U(s,a)
    \end{align*}
    where the first and third equalities are because the objective function is linear, and the second equality is because $\conv(\cP)$ is sa-rectangular by assumption.
    The rest of the proof follows from considering \eqref{eq: lec 2 ssp sarec} applied to $U(s,a) = r(s,a) + \gamma \hat{v}\opt$. In this case we can find $P\opt \in \cP$ such that 
    \[P\opt \in \cap_{(s,a) \in \cS \times \cA} \arg \min_{P \in \cP} P(s,a)\tr U(s,a)\]
    and $\pi\opt \in \PiSD$ such that
    \[\pi\opt(s) \in \arg \max_{a \in \cA} \min_{P \in \cP} P(s,a) \tr (r(s,a) + \gamma \hat{v}\opt) = \arg \max_{a \in \cA} P\opt(s,a) \tr (r(s,a) + \gamma \hat{v}\opt) , \forall \; s \in \cS.\]
    Then $(\pi\opt,P\opt)$ is a saddle point: by construction $P\opt$ is a worst-case kernel for $\pi\opt$, while $\pi\opt$ is an optimal policy for $P\opt$, and $\hat{v}\opt = v^{\pi\opt,P\opt}$. Therefore $\displaystyle \sup_{\pi \in \PiH} \min_{P \in \cP} R(\pi,P) = p_0 \tr \hat{v}\opt$.
\end{proof}
\paragraph{Discussion.}
The simultaneous solvability properties \eqref{eq: lec 2 ssp srec} and \eqref{eq: lec 2 ssp sarec} say that we can optimize linear forms independently across the different dimensions of $\cP \subset \DSSA$. Hence, the SSPs are inherently statements about the set of extreme points of $\cP$; \Cref{th: lec 2 dpp srec} and \Cref{th: lec 2 dpp sa rec} characterize the convex hull of $\cP$ and not necessarily $\cP$. In particular if we don't impose the convexity assumption we can have non-rectangular sets satisfying the SSPs, simply by considering finite s-rectangular sets and adding one point in the middle, as shown in the next example.

\begin{example}\label{ex: lec 2 nons-rec ssp}
Take $\cS = \{s_0,s_1\}$ and assume that there is only one action $a_1$. We can identify $\DSSA$ with $\{ ((p,1-p),(q,1-q)) \; | \; p,q \in [0,1] \}$ where $p$ is the probability to transition to $s_0$ out of $s_0$ and $q$ is the probability to transition to $s_0$ out of $s_1$. We can identify $P \in \DSSA$ with the pair $(p,q) \in [0,1]^{2}$. For instance $P = (1,1)$ corresponds to $s_0$ absorbing and $s_1 \to s_0$ deterministically. 
Consider the following finite uncertainty set: $\displaystyle \cP_0 = \{(0,0),(0,1),(1,0),(1,1)\}.$ We can see $\cP_0$ as the set of four extreme points of the square in $[0,1]^{2}$ centered at $(1/2,1/2)$, and $\cP_0$ is compact s-rectangular (not convex since it is finite) and satisfies \eqref{eq: lec 2 ssp srec}. Now add the center and consider
$\displaystyle \cP = \{(0,0),(0,1),(1,0),(1,1),(1/2,1/2)\}.$ Note that $\cP$ is still compact and not convex, but it is not s-rectangular anymore because of the point $(1/2,1/2)$, e.g., $(1/2,1/2) \in \cP$ and $(0,0) \in \cP$ but $(0,1/2) \notin \cP$.
However the SSP \eqref{eq: lec 2 ssp srec} still holds for $\cP$, because the point $(1/2,1/2)$ is irrelevant from an optimization standpoint - it is never the {\em unique} minimizer of a linear form over $\cP$, and $\cP_0 = \cP \setminus \{(1/2,1/2)\}$ is s-rectangular.
\end{example}
\subsection{Minimality for the existence of uniformly optimal policies*}\label{sec: lec 2 unf opt rec}
For s-rectangular convex compact uncertainty sets, the existence of a stationary uniformly optimal policy is stated and proved in \Cref{th: lec 2 main rmdps}; this is of course very convenient from an implementation standpoint (we do not need to know where we start to play optimally). We show that an even weaker property (just the existence of uniformly optimal policies) completely characterizes s-rectangularity among compact convex uncertainty sets in the next theorem.
\begin{theorem}\label{th: lec 2 s-rec unf opt}
Let $\cA,\cS$ be the sets of actions and states with $|\cS| \geq 2, |\cA| \geq 2$.
Let $\cP$ be a compact convex uncertainty set. The following statements are equivalent:
\begin{enumerate}
\item For any choice of $r \in \R^{\cS \times \cA \times \cS},\gamma \in [0,1)$, there exists a uniformly optimal policy $\pi\opt$:
\[\forall\, r \in \R^{\cS \times \cA \times \cS}, \gamma \in [0,1) ,\exists\,\pi\opt\in\PiH,\forall\, p_0\in\Delta(\cS),\quad  \min_{P \in \cP} R(\pi\opt,P)= \sup_{\pi \in \PiH} \min_{P \in \cP} R(\pi,P).\]
\item $\cP$ is s-rectangular.
\end{enumerate}
\end{theorem}
\begin{proof}
    Recall that we have proved that the second statement implies the first one (in \Cref{th: lec 2 main rmdps}), so we are left to prove that the first statement implies the second one.

    Let $|\cA| \geq 2, |\cS| \geq 2$, and $\cP$ be compact convex. We proceed by contradiction. Assume that the first statement holds but that $\cP$ is not s-rectangular. Choose $\gamma=0$, so that the return of any history-dependent policy only depends on its chosen first action, so that optimal policies can be chosen stationary anyway. Here we write $\cP(s) \subset \DSA$ for the $s$-th marginal of $\cP$: $\cP(s) \coloneqq \{P(s) \; | \; \exists P' \in \cP, P'(s) = P(s)\}$.
    Following \Cref{th: lec 2 dpp srec}, we can choose $U^0 \in \R^{\cS \times \cA \times \cS}$ such that, for $\recs(\cP) = \times_{s \in \cS} \cP(s)$ the s-rectangular extension of $\cP$ (see \Cref{sec: lec 3 time var P} for more detail) we have
\begin{equation}\label{eq:s-rec ssp U0}
    \min_{P \in \cP} \langle U^0,P\rangle = \min_{P \in \cP}  \sum_{s \in \cS} \langle P(s), U^0(s)\rangle > \min_{P \in \recs(\cP)} \langle U^0,P\rangle =  \sum_{s \in \cS} \min_{P \in \cP} \langle P(s), U^0(s)\rangle.
\end{equation}
Note that we can perturb $U^0$ to a certain $U \in \R^{\cS \times \cA \times \cS}$ to both preserve the strict inequality and ensure uniqueness of minimizers on both sides. We can then scale and add a positive constant to every component of $U$ to ensure $1 \leq U \leq 2$ (since $P(s,a) \in \Delta(\cS)$). Let $P\opt$ be in $\arg \min_{P \in \cP} \langle U, P \rangle$. There exists a state $\bs \in \cS$ such that
\[\langle P\opt(\bs), U(\bs)\rangle > \min_{P \in \cP} \langle P(\bs), U(\bs)\rangle\]
and by construction of $U$ the objective on the right-hand side, involving only the $s$-th component of $P$, is minimized at a unique component $P'(\bs)$. Note that both vectors $(P\opt(\bs,a)\tr U(\bs,a))_{a \in \cA}$ and $(P'(\bs,a)\tr U(\bs,a))_{a \in \cA}$ are strictly positive (by $U \geq 1$ and $P(s,a) \in \DS$) and the first vector has a strictly larger sum than the second. If these vectors are proportional, we can modify the vector $U$ to ensure that the new corresponding vectors are non-proportional by adding $1$ to every component of $U(\bs,a_1) \in \RS$ for some particular $a_1 \in \cA$; this again does not change the strict gap in \eqref{eq:s-rec ssp U0}, and we now have $1 \leq U(s,a,s') \leq 3$ for any $s,a,s'$.

We can now construct a uniformly optimal policy for a particular choice of rewards (since we assumed the first statement). We will do so using Kakutani's fixed point theorem (\Cref{th: lec 1 kakutani}). Consider the correspondence $\PiS \to 2^{\PiS}$ that maps each $\pi \in \PiS$ to the set of stationary uniformly optimal policies in the robust MDP with uncertainty set $\cP$ and with instantaneous rewards defined as $r(s,a,s') = U(s,a,s') / \max \{ \pi(s,a), (6|\cA|)^{-1}\}$ (the maximum term in the denominator is here to allow for policies without full support). This correspondence satisfies the assumptions for \Cref{th: lec 1 kakutani} (recall that $\gamma=0$ so that the return is a linear function of the policies and the worst-case return is continuous and concave in the policy). Hence we get that there exists $\pi\opt \in \PiS$ that is a fixed point of this correspondence. We want to show that $\pi\opt(s,a) > (6|\cA|)^{-1}$ for each $(s,a)$ so that $\pi\opt$ has full support. Fix $(s,a)$ and let $\hat{\pi}$ that chooses action $a$ with probability $1$ in state $s$. We have
\begin{align*}
    \frac{1}{\max \{ \pi\opt(s,a), (6|\cA|)^{-1}\}} & \leq \min_{P \in \cP} \frac{P(s,a) \tr U(s,a)}{\max \{ \pi\opt(s,a), (6|\cA|)^{-1}\}}  \\
    & = \min_{P \in \cP} v^{\hat{\pi},P}(s) \\
    & \leq \min_{P \in \cP} v^{\pi\opt,P}(s) \\
    & \leq 3 \sum_{a'}  \frac{\pi\opt(s,a')}{\max \{ \pi\opt(s,a'), (6|\cA|)^{-1}\}} \\
    & \leq 3 |\cA|
\end{align*}
where the first inequality is from $U(s,a,s') \geq 1, P(s,a) \in \DS$, the equality is by definition of the instantaneous rewards and $\hat{\pi}$, the second inequality is from the uniform optimality of $\pi\opt$, the third inequality is from $U(s,a,s') \leq 3$, and the fourth inequality is always true. Hence we must have that $\pi\opt(s,a) \geq (3|\cA|)^{-1} > (6|\cA|)^{-1}$, $\pi\opt$ has full support and $r(s,a,s') = U(s,a,s')/\pi\opt(s,a)$; in particular we get a simple expression for $v^{\pi\opt,P}(s)$ as $v^{\pi\opt,P}(s) = \sum_{a \in \cA} P(s,a) \tr U(s,a)$.

We can now arrive at a contradiction. 
First, let $p_0$ be the uniform initial distribution so that $R(\pi\opt,P) = \frac{1}{|\cS|}\langle P, U \rangle$. From \eqref{eq:s-rec ssp U0} the unique worst-case kernel is $P\opt$, and from the analysis in \Cref{prop: lec 2 NE SD} we must have that $\pi\opt$ maximizes the nominal return against $P\opt$; in particular $\pi\opt(\bs) \in \DA$ maximizes the return against $P\opt(\bs)$, i.e.,
\[\pi\opt(\bs) \in \arg \max_{\pi(\bs) \in \DA} \sum_{a \in \cA} \pi(\bs,a) P\opt (\bs,a) \tr r(\bs,a)\]

Second, take $p_0$ starting in $\bs$ with probability $1$. Then $P'(\bs) \in \DSA$ is the unique $\bs$-component of a worst-case kernel (by construction of $U$). Therefore $\pi\opt(\bs)$ also maximizes the return against $P'(\bs)$. Now consider the optimality of $\pi\opt(\bs)$ in 
\begin{equation}\label{eq: lec 2 optimality piopt}
    \max_{\pi(\bs) \in \DA} \min_{P(\bs) \in \cP(\bs)} \sum_{a \in \cA} \pi(\bs,a) P(\bs,a)\tr r(\bs,a)
\end{equation}
For each $\bar{P}(\bs) \in \{P'(\bs),P\opt(\bs)\}$, the nominal optimality of $\pi\opt(\bs)$ above for the corresponding initial distribution implies the complementary-slackness condition:
\[0 = \sum_{a \in \cA} \pi\opt(\bs,a) \left( \max_{a' \in \cA} \bar{P}(\bs,a') \tr r(\bs,a') - \bar{P}(\bs,a) \tr r(\bs,a) \right). \]
Since $\pi\opt$ has full support we must then have 
\[ \max_{a' \in \cA} \bar{P}(\bs,a') \tr r(\bs,a') = \bar{P}(\bs,a) \tr r(\bs,a), \forall \; a \in \cA. \]
In particular, the scalar $c \coloneqq \bar{P}(\bs,a) \tr r(\bs,a)$ is independent of $a \in \cA$. Using the definition of the instantaneous rewards we know that $\bar{P}(\bs,a) \tr U(\bs,a) = \pi\opt(\bs,a) \bar{P}(\bs,a) \tr r(\bs,a)$. We can determine the constant $c$ as follows: since $\pi\opt(\bs) \in \DA$ we have
\[c = \sum_{a \in \cA} \pi\opt(\bs,a) c = \sum_{a \in \cA} \pi\opt(\bs,a) \bar{P}(\bs,a) \tr r(\bs,a) = \sum_{a \in \cA} \bar{P}(\bs,a) \tr U(\bs,a).\]
Applying the algebra above to $\bar{P} \in \{P',P\opt\}$ we conclude that
\begin{equation}\label{eq: lec 2 Popt P' aligned}
    \begin{aligned}
    P\opt(\bs,a) \tr U(\bs,a) & = \pi\opt(\bs,a) c\opt \\
    P'(\bs,a) \tr U(\bs,a) & = \pi\opt(\bs,a) c'
\end{aligned}
\end{equation}
with $c\opt \coloneqq \sum_{a \in \cA} P\opt(\bs,a) \tr U(\bs,a), c' = \sum_{a \in \cA} P'(\bs,a) \tr U(\bs,a)$. But \eqref{eq: lec 2 Popt P' aligned} means that the vectors $(P\opt(\bs,a) \tr U(\bs,a))_{a \in \cA}$ and $(P'(\bs,a) \tr U(\bs,a))_{a \in \cA}$ are both proportional to $\pi\opt(\bs)$, which is a contradiction.
\end{proof}
Recall that sa-rectangularity buys two additional properties, compared to s-rectangularity: (a) a stationary uniformly optimal policy can be chosen {\em deterministic}, (b) we do not need the convexity assumption.  
The next theorem shows a converse statement as \Cref{th: lec 2 s-rec unf opt} but for sa-rectangular uncertainty sets. We write $\PiHD$ for the set of deterministic history-dependent policies.
\begin{theorem}\label{th: lec 2 sa-rec unf opt}
Let $\cA,\cS$ be the sets of actions and states with $|\cS| \geq 2, |\cA| \geq 2$.
Let $\cP$ be a compact uncertainty set. The following statements are equivalent:
\begin{enumerate}
    \item For any choice of $r \in \R^{\cS \times \cA \times \cS}$ and $\gamma \in [0,1)$, there exists a uniformly optimal deterministic policy $\pi\opt$:
      \[ \forall\, r \in \R^{\cS \times \cA \times \cS}, \gamma \in [0,1), \exists\,\pi\opt\in\PiHD,\forall\,p_0\in\Delta(\cS),\quad  \min_{P \in \cP}  R(\pi\opt,P)= \sup_{\pi \in \PiH} \min_{P \in \cP} R(\pi,P).\]
      \item $\conv(\cP)$ is sa-rectangular.
\end{enumerate}
\end{theorem}
\begin{proof}
The fact that the second statement implies the first one was proved in \Cref{th: lec 2 dpp sa rec}, so we just need to prove that the first statement implies the second one.
Now assume that the first statement is true. Let $\gamma=0$ and $\pi \in \PiS$, we have 
\[R(\pi,P) = \sum_{s \in \cS} p_0 (s) \sum_{a \in \cA} \pi(s,a) P(s,a) \tr r(s,a)\] so that $\min_{P \in \cP} R(\pi,P) = \min_{P \in \conv(\cP)} R(\pi,P)$.

The necessity proof of \Cref{th: lec 2 s-rec unf opt} only uses $\gamma=0$, hence from \Cref{th: lec 2 s-rec unf opt} we must have that $\conv(\cP)$ is s-rectangular. It remains to show that it is sa-rectangular.
For the sake of contradiction, assume that $\conv(\cP)$ is not sa-rectangular. From the Simultaneous Solvability Property \eqref{eq: lec 2 ssp sarec} there exists $U \in \R^{\cS \times \cA \times \cS}$ such that 
    \begin{equation} \label{eq: SSP sa proof sa}
\cap_{(s,a) \in \cS \times \cA} \arg \min_{P \in \conv(\cP)} P(s,a)\tr U(s,a) = \emptyset .
    \end{equation}
    Note that we can assume that \begin{equation}\label{eq:min=0 across s,a}
           \min_{P \in \conv(\cP)} \; P(s,a)\tr U(s,a)=0, \forall (s,a) \in \cS \times \cA  .
       \end{equation}
    This follows from $P(s,a) \in \Delta(\cS)$ so we can subtract any sa-dependent constant from the sa-component of $U$ and this yields the same set $\arg\min_{P \in \conv(\cP)} P(s,a)\tr U(s,a)$. 
    We now show that because $\conv(\cP)$ is s-rectangular (despite not being sa-rectangular), there exists a state $s_0 \in \cS$ such that 
    $\cap_{a \in \cA} \arg \min_{P \in \conv(\cP)} P(s_0,a)\tr U(s_0,a) = \emptyset$.
    It also satisfies the  Simultaneous Solvability Property \eqref{eq: lec 2 ssp srec} for s-rectangular sets:
    \[\forall U'\in\R^{\cS\times\cA\times\cS}, \cap_{s\in\cS} \arg \min_{P\in \conv(\cP)} \sum_{a\in\cA} P(s,a) \tr U'(s,a) \neq \emptyset ,\]
    so there exists $P'\in\conv(\cP)$ such that for any $s\in\cS$, we have
    \[\sum_{a\in\cA} P'(s,a)\tr U(s,a) = \min_{P\in\conv(\cP)} \sum_{a\in\cA} P(s,a)\tr U(s,a).\]
    If we had that $\cap_{a \in \cA} \arg \min_{P \in \conv(\cP)} P(s,a)\tr U(s,a) \neq \emptyset$ for each $s\in\cS$, then 
    \[\min_{P \in \conv(\cP)} \sum_{a\in\cA} P(s,a)\tr U(s,a) =0\] for each $s\in\cS$ and $P'\in\cap_{(s,a) \in \cS \times \cA} \arg \min_{P \in \conv(\cP)} P(s,a)\tr U(s,a)$, which contradicts \eqref{eq: SSP sa proof sa}.
    Hence, there exists a state $s_0 \in \cS$ such that 
    $\cap_{a \in \cA} \arg \min_{P \in \conv(\cP)} P(s_0,a)\tr U(s_0,a) = \emptyset$.

Therefore, because all the $(P(s_0,a)\tr U(s_0,a))_{a\in\cA}$ are nonnegative:
   \[ \inf_{P \in \conv(\cP)} \sum_{a \in \cA} P(s_0,a)\tr U(s_0,a) >0 .\]
   
   Consider the following reward function: $r(s_0,a,s') = U(s_0,a,s')$ for any $(a,s')$ in $\cA \times \cS$, and $r(s,a,s')=0$ otherwise.
   Now consider $p_0$ the initial distribution that puts all mass at $s_0$. 
   The objective function of a policy $\pi \in \PiS$ is
   \[\inf_{P \in \conv(\cP)} \sum_{a \in \cA} \pi(s_0,a) \sum_{s'\in\cS} P(s_0,a,s') r(s_0,a,s').\]
   Any deterministic policy obtains a worst-case objective of $0$, following \eqref{eq:min=0 across s,a}. Consider a policy that uniformly randomizes over actions in $\cA$, i.e., that chooses each action $a \in \cA$ with probability $1/|\cA|$. Then its objective function is
   \[ \frac{1}{|\cA|} \inf_{P \in \conv(\cP)} \sum_{a \in \cA} P(s_0,a)\tr U(s_0,a) >0.\]
   This shows that in this example, a randomized policy may achieve positive return while all deterministic policies achieve a return equal to $0$, which concludes the proof by contradiction.
\end{proof}
\paragraph{Discussion.}
The theorems in this section show that the existence of uniformly optimal policies {\em for all reward parameters} is equivalent to rectangularity. The fact that we consider all possible values for $r \in \R^{\cS \times \cA \times \cS}$ is essential in our proof of \Cref{th: lec 2 s-rec unf opt} and \Cref{th: lec 2 sa-rec unf opt}, and means that we want the existence of uniformly optimal policies for all MDP instances we may encounter. However, one may be interested in a more refined characterization where the set of feasible rewards is a strict subset of $\R^{\cS \times \cA \times \cS}$, i.e., where we only consider a particular subset of instances. The case of rewards independent of the next states (i.e., restricting instantaneous rewards to belong to $\{r \in \R^{\cS \times \cA \times \cS} \; | \; r(s,a,s') = r(s,a,s''), \forall \; s,s',s'' \in \cS, a \in \cA\}$) is of particular interest and several non-rectangular models of uncertainty retain favorable properties for this class of rewards, see \Cref{sec: lec 3 discussion rewards}.
\section{Main takeaways}\label{sec: lec 2 main takeaways}
Before giving the exercises in the next section we review the main results of this chapter below:
\begin{itemize}
    \item {\em Policy evaluation (s-rec.).} Under the s-rectangularity assumption, the policy evaluation problem can be reformulated as an {\em adversarial MDP}, where $\cP$ is interpreted as a set of (stationary deterministic) policies (\Cref{th: lec 2 adversarial mdp}), hence we can compute the worst-case return of a stationary policy via the fixed point machinery from \Cref{chap: lec 0 mdp} for nominal MDPs (\Cref{cor: lec 2 adversarial mdp - alg}).
    \item {\em Policy optimization (s-rec.).} Under the s-rectangularity assumption, optimal policies may still need to be history-dependent (\Cref{prop: lec 2 hist dep opt}). However if we additionally assume convexity and compactness of the s-rectangular uncertainty set, we can find stationary uniformly optimal policies, and they are characterized by a robust Bellman equation (\Cref{th: lec 2 main rmdps}). This can be proved via strong duality in convex optimization or via the existence of a Nash equilibrium (based on Kakutani's fixed point theorem for correspondences). Optimal policies can then be computed via robust variants of value iteration and policy iteration (\Cref{alg: lec 2 VI}-\Cref{alg: lec 2 PI}), but it is less direct to obtain convex formulations. If the uncertainty set is s-rectangular but not convex, we can still compute the best stationary policies.
    \item {\em The case of sa-rectangular models.} Under the sa-rectangularity assumption, we don't need the convexity assumption anymore and uniformly optimal policies can be chosen stationary, because strong duality holds (\Cref{prop: lec 2 sa-rec strong duality}). Additionally optimal policies can be chosen deterministic. We describe algorithms based on value iteration, policy iteration and convex programming in \Cref{sec: lec 2 sa-rec alg}.
    \item (*) {\em Necessary and sufficient conditions for tractability.} In some sense made rigorous in \Cref{sec: lec 2 minimality rec}, rectangularity is the minimal assumption enabling dynamic programming algorithms to compute optimal policies (\Cref{th: lec 2 dpp srec}-\Cref{th: lec 2 dpp sa rec}), and enabling the existence of stationary uniformly optimal policies (\Cref{th: lec 2 s-rec unf opt}-\Cref{th: lec 2 sa-rec unf opt}).
\end{itemize}
\section{Exercises}\label{sec: lec 2 fdt}
\begin{enumerate}
    \item {\em Alternative definitions of rectangularity.}
    \begin{itemize}
        \item     Prove that $\cP$ is s-rectangular (as defined in \Cref{def: lec 2 srec}) if and only if for any two $P^1,P^2 \in \cP$ and any function $f:\cS \to \{1,2\}$, the transition kernel $P^3 = (P^{f(s)}(s))_{s \in \cS}$ belongs to $\cP$.
        \item Prove that $\cP$ is s-rectangular if and only if it is equal to the smallest (in the sense of the inclusion) s-rectangular uncertainty set containing $\cP$. Can you create an analogous definition with a notion of ``largest s-rectangular uncertainty set contained in $\cP$"?
    \end{itemize}

    \item {\em Evaluating the return of general policies.}
    \Cref{sec: lec 2 adv} focuses on solving the policy evaluation problem for a {\em stationary} policy $\pi \in \PiS$. If you are interested in the policy evaluation problem for a history-dependent policy $\pi \in \PiH$, you can think for yourself how we can even evaluate the return $R(\pi,P)$ for a fixed $P \in \cP.$
    \item {\em Properties of the return.} 
    Show that $\pi \mapsto R(\pi,P)$ need not be quasi-concave and that $P \mapsto R(\pi,P)$ need not be quasi-convex.

    Hint: use the properties of the sets of optimal policies of nominal MDPs and adversarial MDPs.
    \item {\em Strong duality from Brouwer.} Brouwer's fixed point theorem states that a continuous map from a nonempty convex compact set of a finite-dimensional Euclidean space to itself admits a fixed point. Use this theorem to prove strong duality for compact convex s-rectangular uncertainty sets.

    Hint: define the corresponding ``Nash improvement function".
    \item {\em Convexification of policy evaluation.}
    \begin{enumerate}
        \item Let $n \in \N$ and $\cE = \{1,...,m\}$ be a finite set of indices. Consider any sets $C_1,...,C_m \subseteq \R^{n}$. Show that the convex hull operation and the Cartesian product commute: $\conv\left(\times_{i \in \cE} C_i \right) = \times_{i \in \cE} \conv\left(C_i\right)$.
        \item Let $\cP$ be s-rectangular compact and let $\pi \in \PiS$.
        \begin{enumerate}
            \item Show that the worst-case Bellman operator $T^{\pi}$ remains the same if we replace $\cP$ by $\conv(\cP)$.
            \item Conclude that $\min_{P \in \cP} R(\pi,P) = \min_{P \in \conv(\cP)} R(\pi,P)$ for $\pi \in \PiS$, i.e., why we can replace $\cP$ by $\conv(\cP)$ in the policy evaluation problem.
            \item Explain why we cannot (in general) replace $\cP$ by $\conv(\cP)$ in the robust MDP problem $\sup_{\pi \in \PiH} \inf_{P \in \cP} R(\pi,P)$.
        \end{enumerate}
    \end{enumerate}
    \item {\em Optimistic MDPs.} The goal of this exercise is to leverage the analysis and tools of this chapter to easily derive the properties for the following {\em optimistic MDP} problem:
    \begin{equation}\label{eq: lec 3 optimistic}
        \sup_{\pi \in \PiH} \sup_{P \in \cP} R(\pi,P)
    \end{equation}
    given some uncertainty set $\cP \subset \DSSA$.
    \begin{enumerate}
        \item Show that for any $\cP \neq \emptyset$ (not necessarily s-rectangular, nor convex, nor compact), an optimal policy in \eqref{eq: lec 3 optimistic} exists and can be chosen stationary and deterministic.

        Hint: use the fact that two suprema commute.

        For the rest of the exercise we assume that $\cP$ is s-rectangular compact.
        \item Explain how we can solve, for some $\pi \in \PiS$, the Optimistic Policy Evaluation problem $\sup_{P \in \cP} R(\pi,P)$.
        \item Prove that we can interpret \eqref{eq: lec 3 optimistic} as a nominal MDP with state set $\cS$ and action set $\cA \times \cP_s$ available at each state $s \in \cS$. Write the instantaneous rewards and transition probabilities in this nominal MDP with augmented sets of actions.
        \item Give an optimality equation for the uniformly optimal stationary policies in \eqref{eq: lec 3 optimistic}.

        Hint: use the ``optimistic Bellman operator" $v \mapsto \left( \max_{a \in \cA} \max_{P \in \cP} P(s,a) \tr (r(s,a) + \gamma v) \right)_{s \in \cS}$.
        \item Adapt the algorithms from \Cref{sec: lec 2 algo} to solve \eqref{eq: lec 3 optimistic}. Can you easily obtain a convex programming formulation for solving \eqref{eq: lec 3 optimistic}?
    \end{enumerate}
    \item {\em A non-attainable set of operators.} Let $\cP$ be compact convex non-s-rectangular and let $T$ be defined as in \eqref{eq: lec 2 T}. Show that we can choose some rewards such that $\{T^{P} \; | \; P \in \cP\}$ does not satisfy the attainability property for $T$.

    Hint: use $\gamma=0$ and $v=(0,...,0) \in \R^{\cS}$.

    \item {\em Convex formulation of s-rectangular robust MDPs.}

    The goal is to outline the steps to obtain a convex formulation for approximating the optimal value function. Let $\cP \subset \DSSA$ be a compact convex s-rectangular uncertainty set. Let $\bar{\pi} \in \PiS$ be the uniform policy.
    \begin{enumerate}
        \item Let $s \in \cS$. Explain why $v \mapsto Tv(s)$ may not be convex nor concave.
        \item Let $b>0$ and consider the following operator $\tT:\R^{\cS} \to \R^{\cS}$ defined as
        \begin{equation}\label{eq: lec 2 T tilde}
            \tT v(s) = \max_{\pi \in \PiS} \min_{P \in \cP} T^{\pi,P}v(s) - b^{-1} \kl(\pi(s),\bar{\pi}(s))
        \end{equation}
        with the Kullback-Leibler divergence defined as $\kl(x,\bar{\pi}(s)) = \sum_{a \in \cA} x(a)\log(x(a)/\bar{\pi}(s,a))$ for $x \in \DA$, with $0\log(0)=0$ by convention.
        \begin{enumerate}
            \item Show that $\tT$ is a monotone contraction. How does $Tv$ compare to $\tT v$ for $v \in \R^{\cS}$?
            \item Let $\tv\opt \in \R^{\cS}$ be the unique fixed point of $\tT$. Show that 
            \[\tv\opt(s) \leq v\opt(s) \leq \tv\opt(s) + b^{-1}(1-\gamma)^{-1} \log(|\cA|),\quad \forall \; s \in \cS.\]
            \item Use the Donsker-Varadhan variational formula to show that
        \[\tT v(s) = b^{-1} \log \left( \min_{P(s) \in \cP_s} \frac{1}{|\cA|} \sum_{a \in \cA} \exp \left(b \cdot P(s,a) \tr (r(s,a) + \gamma v )\right) \right), \quad \forall \; s \in \cS.\]
        \end{enumerate}
        \item Let $\log_b$ and $\exp_b$ be defined as:
        \[
  \log_{b}(z) = \log(z)/b, \qquad  \exp_{b}(z') = \exp(b \cdot z'), \quad  \forall \; z>0, \forall \; z' \in \R.
\]
 Show that, coordinate-wise, we have
\[
  \tilde{T}v = \log_{b} \,\circ\, \tilde{t} \,\circ\, \exp_{b}v,
\]
with the operator $\tilde{t} :\left(\R_{+}^{*}\right)^{\cS} \rightarrow \left(\R_{+}^{*}\right)^{\cS}$ defined  for each state $s \in \cS$ and $x \in \left(\R_{+}^{*}\right)^{\cS}$ as
\begin{equation}\label{eq:definition-f}
  \tilde{t}x(s)=  \min_{P(s) \in \cP_s} \sum_{a \in \cA} \frac{1}{|\cA|} \exp_{b}\left(P(s,a)\tr r(s,a)\right)   \prod_{s' \in \cS} x(s')^{\gamma P(s,a,s')}.
\end{equation}
Prove that $\tilde{t}$ is a componentwise concave function.
\item Combine the above steps to obtain a convex formulation for computing an $\epsilon$-approximation of $v\opt$.

Hint: use the contraction lemma on $\tT$ to compute $\tv\opt$. Use $\tilde{t}$ to make a change of variable based on $\exp_{b}$; conclude for $b$ large enough.
\item (*) The approach above can be seen as using an entropic regularizer then doing a corresponding change of variables. Could the same approach work with a different choice of regularizer, for instance the squared-$\ell_2$ norm, if we could adapt the change of variables accordingly?
    \end{enumerate}
    \item {\em Complexity of checking the product structure.}
    \begin{enumerate}
        \item Assume that you are given a polytope as $Q=\{ x \in \R^{n} \; | \; Ax \leq b\}$ for some $b \in \R^{m},A \in \R^{m \times n}$. How can you check that $Q$ is a Cartesian product over $i \in \{1,...,n\}$?
        \item Assume that you are given a polytope as $Q = \conv(x_1,...,x_p)$ for some finite set $\{x_1,...,x_p\} \subset \R^{n}$. How can you check that $Q$ is a Cartesian product over $i \in \{1,...,n\}$?
    \end{enumerate}
\end{enumerate}

\section{References}\label{sec: lec 2 ref}
The first paper on robust MDPs is \cite{satia1973markov} where the authors consider what is now called sa-rectangular uncertainty sets and note the similarity with Shapley's stochastic games~\cite{shapley1953stochastic} (in fact sa-rectangular robust MDPs are closer to perfect-information stochastic games than to Shapley's concurrent stochastic games where players choose their decisions at the same time, see \Cref{chap: lec 3 connection}).  \cite{givan1997bounded} use such an sa-rectangular model (for ``interval MDPs" where the uncertainty set is such that each transition probability $P(s,a,s')$ belongs to a prescribed interval). The model now called ``sa-rectangular" comes from~\cite{nilim2005robust,iyengar2005robust} and the name ``s-rectangular" comes from~\cite{wiesemann2013robust} (\cite{le2007robust} uses the term ``state-rectangular").

 The first proof of strong duality for convex compact s-rectangular models presented above is from~\cite{nilim2005robust,wiesemann2013robust} and the second proof (based on Kakutani's fixed point theorem) and third proof (based on Brouwer's fixed point theorem, in \Cref{sec: lec 2 fdt}) are new (as far as I know); \cite{le2007robust} proves strong duality for s-rectangular models by invoking the seminal results from Shapley's paper~\cite{shapley1953stochastic}. The connection between Nash equilibria and strong duality in two-player zero-sum games is standard.

For the adversarial MDP in the s-rectangular case I follow~\cite{ho2021partial}. This idea is present in both~\cite{nilim2005robust,iyengar2005robust} for the sa-rectangular case (this is implied by the connection with perfect-information stochastic games given by these authors), and it is present in~\cite{goh2018data,goyal2023robust} for other models of uncertainty (r-rectangular under some additional assumptions on the rewards); it is worth noting that in stochastic games the idea that the decision problem faced by one player (fixing the opponent's strategy) is an MDP problem appeared from the start~\cite{shapley1953stochastic}.

Policy iteration as in \Cref{alg: lec 2 PI} is akin to the Strategy Iteration proposed by Hoffman and Karp~\cite{hoffman1966nonterminating} and~\cite{rao1973algorithms} for stochastic games. A convergence rate of Strategy Iteration is given in~\cite{hansen2013strategy} (see~\cite{akian2013policy} for an improved convergence rate) for perfect-information stochastic games, a setup close to sa-rectangular robust MDPs.

There are many papers that give efficient algorithms for computing the robust Bellman operator $T$ and the worst-case Bellman operator $T^{\pi}$. \cite{nilim2005robust,iyengar2005robust} give algorithms for $T^{\pi}$ for sa-rectangular uncertainty sets based on $\ell_2$-balls, $\ell_1$-balls, $\ell_{\infty}$-balls, entropy, finite models, and the Kullback-Leibler divergence.  \cite{givan1997bounded} focuses on sa-rectangular interval models (similar to models based on $\ell_{\infty}$-balls). \cite{ho2021partial,ho2022robust,behzadian2021fast} give algorithms for fast Bellman updates for s-rectangular robust MDPs based on $\ell_1$-balls, $\ell_{\infty}$-balls, and $\phi$-divergences.
A convex formulation for approximating the optimal value functions of s-rectangular robust MDPs is given in~\cite{grand2025convex}. The convex formulation based on disjunctive programming for sa-rectangular uncertainty sets is from some unpublished research notes with Marek Petrik; the theorem from Balas for disjunctive programming can be found as Theorem 5.1 in \cite{conforti2010extended}. The reader interested in the case of non-convex uncertainty sets can find some discussion in~\cite{wang2023foundation}. 

The connection between optimistic MDPs and robust MDPs (related to the exercise outlined in \Cref{sec: lec 2 fdt}) is well known, for instance \cite{givan1997bounded,tewari2007bounded} deal with both optimistic and robust (pessimistic) MDPs at the same time while \cite{iyengar2005robust} argue (in their Remark 2.1) that their results for robust MDPs imply the corresponding results for optimistic MDPs. What we show in the exercise in \Cref{sec: lec 2 fdt} is slightly stronger (in particular stationary deterministic optimal policies always exist for optimistic MDPs, in contrast to robust MDPs). The idea to interpret optimistic MDPs (with rectangular uncertainty sets) as nominal MDPs with augmented sets of actions was suggested to me by Raimundo Saona in research discussion and is also not new, see e.g. section 3.1.1 in the UCRL2 paper~\cite{jaksch2010near}.

The minimality of rectangularity properties is studied in~\cite{grand2024tractable} (for the dynamic programming principle for the policy evaluation problem) and \cite{benyamine2026existence} (for the existence of uniformly optimal policies); the reader interested in these types of necessary and sufficient conditions for tractability of robust MDPs may also be interested in \cite{benyamine2026dynamic} (which also highlights the attainability property) and \cite{li2025rectangularity}.
\chapter{Connection with other models}~\label{chap: lec 3 connection}
This chapter discusses the connection of the robust MDP model studied in these notes with several other models for sequential decision-making in uncertain environments:
\begin{itemize}
    \item In \Cref{sec: lec 3 uncertain rewards} and \Cref{sec: lec 3 uncertain rewards and transitions} we analyze robust MDPs where the uncertain rewards are unknown, both without and with uncertainty in the transition probabilities.
    \item In \Cref{sec: lec 3 multi-objective MDPs} we analyze the connection between multi-objective MDPs, Pareto optimality and robust MDPs.
    \item In \Cref{sec: lec 3 time var P} we study three different variants of robust MDPs where the transition probabilities may change over time and in \Cref{sec: lec 3 sgs} we highlight the connection with stochastic games.
    \item In \Cref{sec: lec 3 drmdps} we show that distributionally robust MDPs reduce to robust MDPs (despite being a generalization) and illustrate this with ambiguity sets based on Wasserstein distances.
    \item In \Cref{sec: lec 3 non-rec} we study three classes of non-rectangular uncertainty sets leading to favorable tractability properties (under some additional assumptions on the model).
    \item In \Cref{sec: lec 3 risk averse} we analyze the reformulations of risk-averse MDPs as robust MDPs for some particular classes of risk measures (coherent risk measures).
    \item In \Cref{sec: lec 3 discussion rewards} we discuss the natural assumptions that we can make about the reward models and analyze which ones are without loss of generality.
\end{itemize}
A ``visual" summary of the main takeaways is given in \Cref{sec: lec 3 main takeaways}. Throughout we assume $\cP \neq \emptyset$.
\section{Robust MDPs with uncertain rewards}\label{sec: lec 3 uncertain rewards}
In the previous chapters we assumed that the reward parameters are known to the decision maker and that only the transition probabilities are uncertain. Here we first consider the corresponding setting, where the transition probabilities are known to the decision maker but all we know about the rewards is that they belong to a set $\cR \subset \RSAS$. Since we keep $P$ fixed here and vary $r$, we write $R(\pi,r)$ and $v^{\pi,r}$ for the return and value function of a policy $\pi \in \PiH$ given a realization of $r \in \cR$. The goal in this section is to analyze the following optimization problem:
\begin{equation}\label{eq: lec 3 rmdp r}
    \sup_{\pi \in \PiH} \inf_{r  \in \cR} R(\pi,r)
\end{equation}
We have seen in \Cref{sec: lec 0 model} that we can express the return $R(\pi,r)$ of a policy $\pi \in \PiH$ as
\begin{equation}\label{eq: lec 3 R pir frequency}
    R(\pi,r) = \sum_{s \in \cS} \sum_{a \in \cA} \mu^{\pi}(s,a) P(s,a) \tr r(s,a)
\end{equation}
where $\mu^{\pi} \in \RSA$ are the state-action frequencies defined in \eqref{eq: lec 0 state-action occupancy}. An amazing property of the above expression for the return is that it is linear in the reward function $r$, which will prove central in all the arguments in this section. 
\paragraph{Policy evaluation.} We start by considering the problem of solving $\inf_{r \in \cR} R(\pi,r)$ for some fixed $\pi \in \PiS$. From \eqref{eq: lec 3 R pir frequency}, we get that, for $c^{\pi} \in \RSAS$ defined as $c^{\pi}(s,a,s') = \mu^{\pi}(s,a)P(s,a,s')$,  this optimization problem can be reformulated as 
\begin{equation}\label{eq: lec 3 MP PE r}
    \inf \{ c^{\pi \; \top} r \; | \; r \in \cR\}
\end{equation}
where we still use the scalar product notation ($\tr$) for  objects in $\RSAS$: 
\[c^{\pi \; \top} r = \sum_{s,a,s'} c^{\pi}(s,a,s')r(s,a,s').\]
The mathematical program \eqref{eq: lec 3 MP PE r} minimizes a linear form over the set $\cR$; whether this is doable efficiently (e.g., in polynomial time in the input size of the problem) depends on the properties of the set $\cR \subset \RSAS$, but for instance efficient algorithms exist if $\cR$ is a polyhedron. Of course to solve \eqref{eq: lec 3 MP PE r} we need to compute $c^{\pi}$ and therefore we need to compute $\mu^{\pi}$ first, but this can be done efficiently by solving a linear system of equations, as argued at the end of \Cref{sec: lec 0 model}. 
\paragraph{Properties of optimal policies.}
We first note that optimal policies in \eqref{eq: lec 3 rmdp r} can always be chosen stationary for a compact uncertainty set $\cR$, with a proof much simpler than for the case of uncertain transitions.
\begin{proposition}\label{prop: lec 3 stat opt policies r - model blind}
    Let $\cR \subset \RSAS$ be a compact set. Then there exists a stationary optimal policy:
    \[\sup_{\pi \in \PiH} \min_{r  \in \cR} R(\pi,r) = \max_{\pi \in \PiS} \min_{r  \in \cR} R(\pi,r).\]
\end{proposition}
\begin{proof}
This follows from
\begin{align*}
    \sup_{\pi \in \PiH} \min_{r  \in \cR} R(\pi,r) & = \sup_{\pi \in \PiH} \min_{r  \in \cR} \sum_{s,a} \mu^{\pi}(s,a) P(s,a) \tr r(s,a) \\
    & = \sup_{\mu \in \cQ} \min_{r  \in \cR} \sum_{s,a} \mu(s,a) P(s,a) \tr r(s,a) \\
    & = \sup_{\pi \in \PiS} \min_{r  \in \cR} \sum_{s,a} \mu^{\pi}(s,a) P(s,a) \tr r(s,a) \\
    &= \max_{\pi \in \PiS} \min_{r  \in \cR} R(\pi,r)
\end{align*}
where the first equality is via the expression of the return $R(\pi,r)$ as a function of the state-action frequencies as shown in \Cref{sec: lec 0 model}, and the second and third equalities are because the polytope $\cQ$ of the state-action frequencies induced by all $\pi \in \PiH$ coincides with the state-action frequencies induced by all $\pi \in \PiS$.
\end{proof}
Note that we may need randomization at optimality, as shown in the next simple example.
\begin{example}\label{ex: lec 3 randomized}
    Consider an MDP instance with one state $s_0$ and two actions $a_1,a_2$. The rewards are either $r_1$ with $r_1(s_0,a_1,s_0)=1,r_1(s_0,a_2,s_0)=0$ or $r_2$ with $r_2(s_0,a_1,s_0)=0,r_2(s_0,a_2,s_0)=1$. We can parametrize a stationary policy by $x \in [0,1]$ the probability to play $a_1$ in state $s_0$, which gives a worst-case return of $v^{x}(s_0) = \frac{1}{1-\gamma} \min \{x,1-x\}$. Hence an optimal policy chooses the randomized decision $x\opt=1/2$.
\end{example}
It is also worth emphasizing that without a rectangularity assumption for $\cR$ (of the kind $\cR = \times_{s \in \cS} \cR_s$ with $\cR_s \subset \R^{\cA \times \cS}$), there may not exist {\em uniformly} optimal policies in general, i.e., an optimal policy may need to depend on where we start.
\begin{example}\label{ex: lec 3 r not convex}
    Consider the MDP instance with two states $s_0,s_1$, two actions $a_1,a_2$ and take $\gamma=0$. The rewards are independent of the next state and are given by $r_t(s_0,a_1)=1-t,r_t(s_0,a_2)=-1,r_t(s_1,a_1)=0,r_t(s_1,a_2) = 2t-1$, parametrized by $t \in [0,1]$. In $s_0$ the adversary wants to choose $t$ large while in $s_1$ the adversary wants to choose $t$ small; this will prevent the existence of a uniformly optimal policy here. Note that we do not specify the transition probabilities because they are irrelevant ($\gamma=0$ and the rewards are independent of the next states). Let us write $x=\pi(s_0,a_1),y=\pi(s_1,a_1)$ with $x,y \in [0,1]$. 
    If we start in $s_1$ with probability $1$ then the worst-case return is $\min_{t \in [0,1]}(1-y)(2t-1) = y-1$ attained at $t=0$ and maximized at $y=1$.
    If we start in $s_0$ and $s_1$ each with probability $1/2$ then the worst-case return is $\min \{x+\frac{y}{2}-1, \frac{x-y}{2}\}$ which is maximized at $x=1$ and $y=1/2$. Hence there are no optimal policies in common for $p_0$ starting in $s_1$ and $p_0$ the uniform distribution over $\{s_0,s_1\}$, and therefore there is no uniformly optimal policy.
\end{example}

When the set $\cR$ is additionally convex compact, we obtain a strong duality result akin to that of \Cref{chap: lec 2 rectangular models}; again the proof is much simpler than the corresponding results for uncertain transitions because of the bilinearity of $R(\pi,r)$ in $\mu^{\pi}$ and $r$.
\begin{proposition}\label{prop: lec 3 strong duality}
    Let $\cR \subset \RSAS$ be a convex compact set. Then the following strong duality result holds:
    \begin{equation}\label{eq: lec 3 strong duality r}
        \max_{\pi \in \PiS} \min_{r \in \cR} R(\pi,r) = \min_{r \in \cR}  \max_{\pi \in \PiS}  R(\pi,r) 
    \end{equation}
\end{proposition}
\begin{proof}
    We have
    \begin{align*}
        \max_{\pi \in \PiS} \min_{r \in \cR} R(\pi,r) & = \max_{\mu \in \cQ} \min_{r \in \cR} \sum_{(s,a)} \mu(s,a) P(s,a) \tr r(s,a)  \\
        & = \min_{r \in \cR} \max_{\mu \in \cQ}  \sum_{(s,a)} \mu(s,a) P(s,a) \tr r(s,a) \\
        & = \min_{r \in \cR} \max_{\pi \in \PiS}  R(\pi,r)  
    \end{align*}
    where the first equality is via the expression of the return $R(\pi,r)$ as a function of the state-action frequency, the second equality is via strong duality for maxmin optimization problems with bilinear objectives and convex compact domains, and the third equality uses again the expression of $R(\pi,r)$.
\end{proof}

\paragraph{Computation of optimal policies.}
To compute an optimal policy we want to solve $\max_{\pi \in \PiS} \min_{r  \in \cR} R(\pi,r)$. As we saw in the proof of \Cref{prop: lec 3 stat opt policies r - model blind}, this can be reformulated as
\begin{equation}\label{eq: lec 3 maxmin r}
\max_{\mu \in \cQ} \min_{r  \in \cR} \sum_{(s,a)} \mu(s,a) P(s,a) \tr r(s,a).
\end{equation}
This last optimization problem has a bilinear objective and its tractability property will depend on the set $\cR$. In general we can write
\begin{align*}
\max_{\mu \in \cQ} \min_{r  \in \cR} \sum_{(s,a)} \mu(s,a) P(s,a) \tr r(s,a) = \max \; & t \\
& t \leq \min_{r \in \cR} \sum_{(s,a)} \mu(s,a) P(s,a) \tr r(s,a), \\
& t \in \R, \mu \in \cQ.
\end{align*}
The tractability property of the right-hand side can be described in some simple settings:
\begin{itemize}
\item When $\cR$ consists of a finite set of scenarios: $\cR = \{r_1,...,r_K\}$ for some $K \in \N$ and $r_i \in \RSAS$ for $i \in \{1,...,K\}$ we get a linear program:
\begin{align*}
 \max& \; t \\
& t \leq \sum_{(s,a)} \mu(s,a) P(s,a) \tr r_i(s,a), \forall \; i \in \{1,...,K\}, \\
& t \in \R, \mu \in \cQ.
\end{align*}
\item When $\cR$ is a polytope, we can use strong duality to turn the minimization problem over $r \in \cR$ in the constraints into a maximization problem, from which we can get again a linear program.
\item When $\cR$ is described by finitely many convex constraints and strong duality holds for mathematical programs with linear objective over $\cR$, then again we can turn \eqref{eq: lec 3 maxmin r} into a convex program.
\item When $\cR$ is a convex compact set one can also use primal-dual algorithms,  as long as we can optimize over $\cR$ (e.g. as long as we can compute proximal maps or linear-optimization oracles over $\cR$).
\end{itemize}
\paragraph{Comparison with robust MDPs with uncertain transitions.}
Robust MDPs with uncertain rewards appear to have much more favorable properties than robust MDPs with uncertain transitions, in that we need fewer assumptions on the set $\cR$ (no rectangularity assumption is necessary to derive the results in this section) to get tractable problems (i.e., with stationary optimal solutions that we can compute) than on $\cP$. This discrepancy comes from the linearity of the return $r \mapsto R(\pi,r)$ in $r$, whereas in general $P \mapsto R(\pi,P)$ is {\em not} a linear function of the transition probabilities (since we need to invert a matrix in \eqref{eq: lec 0 value function formula}).
\paragraph{What is observed? (part 1)}
It is worth taking a second look at \Cref{ex: lec 3 randomized}. In this example it is a bit puzzling that  non-stationary policies can't outperform stationary policies (as proved in \Cref{prop: lec 3 stat opt policies r - model blind}). In particular, if we could observe the instantaneous rewards obtained at each period, then we could determine whether the realization of $r$ chosen by the adversary is $r_1$ or $r_2$ (simply by playing $a_1$ with probability $1/2$ during the first period and checking whether we obtained a reward of $0$ or $1$), and then play accordingly onward (that is, if the realized reward is $1$ we keep playing the same action as in period $0$, otherwise the realized reward is $0$ and we should keep playing the other action). This yields a return of $\frac{1}{2} + \frac{\gamma}{1-\gamma}$, beating the stationary optimum when $\gamma > 0$. The subtlety lies in the definition of history-dependent policies. Recall that history-dependent policies $\PiH$ are defined in \Cref{sec: lec 0 model} as maps from all finite histories $(s_0,a_0,s_1,a_1,...,s_{t-1},a_{t-1},s_t)$ to mixed actions in $\Delta(\cA)$; in particular, the realized rewards $r(s_{\tau},a_{\tau},s_{\tau+1})$ for $\tau \leq t-1$ are {\em not} part of the histories, and so in principle history-dependent policies are not allowed to observe the rewards gathered so far at each period and adjust their decisions based on them. This is somewhat similar to robust MDPs with uncertain transitions, where we did not assume that the decision maker can observe the value of $P(s_{\tau},a_{\tau},s_{\tau+1})$ for $\tau \leq t-1$ along the trajectories $(s_0,a_0,s_1,a_1,...,s_{t-1},a_{t-1},s_t)$.

We can still obtain some results for the case where the decision maker is allowed to observe the instantaneous rewards at each period (in addition to the sequence of state-action pairs visited along the trajectory); let's call this class of policies {\em reward-augmented} history-dependent policies $\PiHra$. Indeed, in this model the set of stationary policies $\PiSra$, which only depend on the current state, is still identifiable with $\PiS = \times_{s \in \cS} \Delta(\cA)$, i.e., the sets of stationary policies coincide in the classical setup (with history-dependent policies as defined in \Cref{chap: lec 0 mdp}) and in the setup where the decision maker can observe the realized rewards. Therefore, for $\cR$ convex compact, strong duality as proved in \Cref{prop: lec 3 strong duality} would still hold, and we have seen at the beginning of \Cref{sec: lec 2 prf unf opt} that strong duality implies the existence of stationary optimal policies. Hence we conclude that for $\cR$ convex compact, an optimal solution to $\sup_{\pi \in \PiHra} \inf_{r \in \cR} R(\pi,r)$, with $\PiHra$ the set of reward-augmented history-dependent policies, can be chosen stationary: 
\[\sup_{\pi \in \PiHra} \inf_{r \in \cR} R(\pi,r) = \max_{\pi \in \PiS} \min_{r \in \cR} R(\pi,r)\]
and a stationary optimal policy can be computed via the methods described right above. If $\cR$ is not convex, then we may need non-stationary reward-augmented policies, as argued in the previous paragraph and exemplified in \Cref{ex: lec 3 randomized}.
\section{Robust MDPs with uncertain rewards and uncertain transitions}\label{sec: lec 3 uncertain rewards and transitions}
In this section, we treat the case where both the rewards $r \in \RSAS$ and the transitions $P \in \DSSA$ are uncertain. In this case, we are given an uncertainty set $\cU \subset \RSAS \times \DSSA$, which models the set of realizable pairs of rewards-transitions $(r,P) \in \RSAS \times \DSSA$. Note that in principle there may be joint constraints between $r$ and $P$ in $\cU$, i.e., we do not require $\cU$ to be a Cartesian product $\cR \times \cP$ with $\cR \subset \RSAS$ and $\cP \subset \DSSA$. 

We extend the return to $R(\pi,r,P)$ with $(r,P) \in \cU$, and the goal in this section is to solve
\begin{equation}\label{eq: lec 3 rmdp r p}
    \sup_{\pi \in \PiH} \inf_{(r,P) \in \cU} R(\pi,r,P)
\end{equation}
The first thing to notice is that given $(r,P) \in \cU$ we can replace $r \in \RSAS$ by $\hat{r}_{P,r} \in \RSAS$ the expected one-step reward as follows.
\begin{lemma}\label{lem: lec 3 R hat rP}
    Let $\pi \in \PiH, (r,P) \in \cU$. Then
    $R(\pi,r,P) = R(\pi,\hat{r}_{P,r},P)$ with $\hat{r}_{P,r} \in \RSAS$ defined as 
    \[\hat{r}_{P,r}(s,a,s') = P(s,a) \tr r(s,a), \forall \; s,a,s'.\]
\end{lemma}
\begin{proof}
    We already did such a calculation in \Cref{sec: lec 0 model} when dealing with state-action frequencies. By definition
    \[R(\pi,r,P) = \E^{\pi,P}\left[ \sum_{t=0}^{+\infty} \gamma^t r(s_t,a_t,s_{t+1}) \right] =  \sum_{t=0}^{+\infty} \gamma^t \E^{\pi,P}\left[ r(s_t,a_t,s_{t+1}) \right]\]
    and $\E^{\pi,P}\left[ r(s_t,a_t,s_{t+1}) \right] = \E^{\pi,P}\left[ P(s_t,a_t) \tr r(s_t,a_t) \right] = \E^{\pi,P}\left[ \hat{r}(s_t,a_t,s_{t+1}) \right].$
\end{proof}
Note $\hat{r}(s,a,s')$ is independent of the next state $s'$, so we could have identified $\hat{r}$ with an element of $\RSA$ instead of $\RSAS$. We still write $\hat{r} \in \RSAS$ for consistency (given it is an input of the return function) but we will allow ourselves to write $\hat{r}(s,a)$ instead of $\hat{r}(s,a,s')$ later in this section.

Given an uncertainty set $\cU \subset \RSAS \times \DSSA$ we can therefore construct an {\em effective uncertainty set} $\cUf$ as follows:
\begin{equation}\label{eq: lec 3 Uf}
    \cUf = \{ (\hat{r},P) \; | \; \exists (r,P) \in \cU, \hat{r}(s,a,s') = P(s,a) \tr r(s,a), \forall \; s,a,s' \}.
\end{equation}
Note that the first component of $\cUf$ has a bilinear expression in the original component $(r,P)$ in $\cU$.
What properties of $\cU$ may the effective uncertainty set $\cUf$ inherit?
\begin{itemize}
    \item Clearly if $\cU$ is compact then $\cUf$ is compact. 
    \item If $\cU$ is s-rectangular: $\cU = \times_{s \in \cS} \cU_s$ with $\cU_s \subset \{(r(s),P(s)) \; | \; r(s) \in \RAS, P(s) \in \DSA \}$ then $\cUf$ is s-rectangular: $\cUf = \times_{s \in \cS} \cUf_{s}$ with
    $    \cUf_s = \{ (\hat{r}(s),P(s)) \; | \; \exists (r(s),P(s)) \in \cU_s, \hat{r}(s,a,s') = P(s,a) \tr r(s,a), \forall \; a,s' \}.$
    \item The main issue is regarding convexity. If $\cU$ is convex, it may be the case that $\cUf$ is {\em not} convex, because of the bilinearity of $\hat{r}$, the reward-component of $\cUf$, in $(r,P) \in \cU$. We give an example where $\cU$ is convex but $\cUf$ is not convex below.
 \begin{example}\label{ex: lec 3 non cvx rp} Assume we are given a state-action pair $(s_0,a_0)$ and the next possible states are $s_0,s_1$. Let $(r^0,P^0) \in \cU$ such that $r^0(s_0,a_0,s_0) = r^0(s_0,a_0,s_1) =0, P^0(s_0,a_0,s_0)=1$, and let $(r^1,P^1) \in \cU$ such that $r^1(s_0,a_0,s_0) =0, r^1(s_0,a_0,s_1) =1, P^1(s_0,a_0,s_1)=1$. Let $\lambda \in [0,1]$ and $(r^{\lambda},P^{\lambda}) = (1-\lambda) (r^0,P^0) + \lambda (r^1,P^1)$. The expected instantaneous reward $\hat{r}^{\lambda}(s_0,a_0,s')$ is $\hat{r}^{\lambda}(s_0,a_0,s') = P^{\lambda}(s_0,a_0) \tr r^{\lambda}(s_0,a_0) = \lambda^2$. Hence the $(s_0,a_0)$-th component of $\cUf$ is $\{ (\lambda^2, \lambda ) \; | \; \lambda \in [0,1]\}$, which is not a convex set (e.g. the midpoint between $(0,0)$ and $(1,1)$ is not feasible).
 \end{example}
\end{itemize}
Note that {\em if} $\cUf$ is convex compact and s-rectangular, everything follows as in \Cref{chap: lec 2 rectangular models}. For $(\hat{r},P) \in \cUf, \pi \in \PiS$ we define the operators $\Tf^{\pi,\hat{r},P}, \Tf^{\pi}, \Tf:\R^{\cS} \to \RS$ as usual:
\begin{align*}
\Tf^{\pi,\hat{r},P}v(s) & = \sum_{a \in \cA} \pi(s,a) \left( \hat{r}(s,a) + \gamma P(s,a) \tr v \right) \\
\Tf^{\pi} v(s) & = \min_{(\hat{r},P) \in \cUf}\Tf^{\pi,\hat{r},P}v(s) \\
\Tf v(s) & = \max_{\pi \in \PiS} \min_{(\hat{r},P) \in \cUf}\Tf^{\pi,\hat{r},P}v(s)
\end{align*}
We of course get that $\Tf^{\pi,\hat{r},P}, \Tf^{\pi}, \Tf$ are monotone contractions, and we write $v\opt \in \R^{\cS}$ for the fixed point of $\Tf$.
\begin{theorem}\label{th: lec 3 main rmdps r P}
Let $\cU \subset \RSAS \times \DSSA$ be an uncertainty set for rewards-transition pairs. Assume that $\cUf$ is compact, convex and s-rectangular, with $\cUf$ defined as in \eqref{eq: lec 3 Uf}. The following statements hold:
\begin{enumerate}
    \item {\em Strong duality}. The following strong duality holds:
    \begin{equation}
        \max_{\pi \in \PiS} \min_{(\hat{r},P) \in \cUf} R(\pi,\hat{r},P) = \min_{(\hat{r},P) \in \cUf} \max_{\pi \in \PiS}  R(\pi,\hat{r},P).
    \end{equation}
    \item {\em Optimality of stationary policies.} We have
    \[ \sup_{\pi \in \PiH} \min_{(\hat{r},P) \in \cUf} R(\pi,\hat{r},P) = \max_{\pi \in \PiS} \min_{(\hat{r},P) \in \cUf} R(\pi,\hat{r},P).\]
    \item {\em Optimality conditions.} $\pi \in \PiS$ is a stationary uniformly optimal policy for \eqref{eq: lec 3 rmdp r p} if and only if $\Tf^{\pi}v\opt = \Tf v\opt$\end{enumerate}
\end{theorem}
The proof follows exactly the same lines as the proof of \Cref{th: lec 2 main rmdps}; we can still apply Sion's minimax theorem to show that in the definition of $\Tf$ we can turn the $\max \min$ into a $\min \max$, and this is where we require convexity of $\cUf$. It is worth emphasizing that we do require s-rectangularity of $\cR$ (the reward-component induced by $\cUf$) for \Cref{th: lec 3 main rmdps r P}, whereas in \Cref{sec: lec 3 uncertain rewards}, where $P$ was fixed and the rewards were unknown, we required $\cR$ to be convex and compact. This comes from the difference in solution methods between robust MDPs with uncertain rewards (which we solved in \Cref{sec: lec 3 uncertain rewards} based on state-action frequencies without requiring rectangularity) and robust MDPs with uncertain transitions (which we solved in \Cref{chap: lec 2 rectangular models} based on value functions, which are naturally componentwise objects and therefore require rectangularity to define proper worst-case value functions). To reconcile the two models and solve robust MDPs with uncertain rewards {\em and} uncertain transitions, we need stronger assumptions and therefore we assume s-rectangularity of $\cU$.

When $\cU$ is a Cartesian product of sa-rectangular sets, we do not need convexity (similar to the case of sa-rectangular uncertainty set in \Cref{sec: lec 2 sarec}), that is, if $\cU = \cR \times \cP$ with $\cR = \times_{s,a} \cR_{sa}, \cP = \times_{s,a} \cP_{sa}$ with $\cR_{sa} \subset \RS, \cP_{sa} \subset \DS$ compact sets then the conclusions of \Cref{th: lec 3 main rmdps r P} hold and additionally we get that an optimal policy can be chosen stationary {\em and} deterministic.

It remains to answer a central question: under what conditions is the effective uncertainty set $\cUf$ convex?
\paragraph{When is $\cUf$ convex?}
We have observed above that s-rectangularity and compactness of $\cU$ naturally imply s-rectangularity and compactness of $\cUf$. However obtaining the convexity of $\cUf$ is trickier. An interesting case is when the rewards are independent of the next state: If indeed we have that for any $(r,P) \in \cU$ then $r(s,a,s') = r(s,a,s''), \forall \; s,s',s'' \in \cS, a \in \cA$, we then get that $P(s,a) \tr r(s,a) = r(s,a,s')$ (where we can choose any $s'$ in the expression above), so that $\cUf = \cU$, so that  $\cUf$ is convex when $\cU$ is convex. To the best of our knowledge, most papers that focus on uncertain rewards {\em and} transition probabilities make this assumption, yet this is not without loss of generality, as we have seen repeatedly in these notes. We conclude this section with the following hardness result, which shows that policy evaluation is difficult when both $(r,P)$ are uncertain.
\begin{proposition}\label{prop: lec 3 hardness uncertain r p}
Let $\pi \in \PiSD$.
For a general convex compact s-rectangular polytope $\cU \subset \RSAS \times \DSSA$, it is in general NP-hard to decide $\inf_{(r,P) \in \cU} R(\pi,r,P) \leq c$.
\end{proposition}
\begin{proof}
    The proof is based on a reduction from the Clique problem. Consider an undirected graph $G=(V,E)$ with $V=\{1,...,n\}$ and $A \in \{0,1\}^{n \times n}$ its adjacency matrix. A clique is a subset of vertices all adjacent to each other. The size of a clique is its number of vertices. For some $k \in \N$, the Clique problem asks to decide whether there is a clique of size at least $k$ or not.

    We construct a reduction from the Clique problem to deciding $\inf_{(r,P) \in \cU} R(\pi,r,P) \leq c$ as follows. We consider an MDP instance where $\cS = V \cup \{s_0\}$ and there is one action $a$ (irrelevant here). The initial state is $s_0$, every state $j \in V$ is absorbing with reward $0$. The only uncertain probabilities and rewards are out of state $s_0$, and the set of feasible $(r(s_0),P(s_0)) \in \R^{\cA \times \cS} \times \DSA$ is given by
    \[\cU_{s_0} = \{ (r(s_0),P(s_0)) \; |\; x \in \Delta(V), (P(s_0,a,j))_{j \in V} = x, (r(s_0,a,j))_{j \in V} = - Ax, r(s_0,a,s_0)=0\}.\]
    Clearly $\cU$ is s-rectangular since the only state with some uncertainty is $s_0$. Moreover $\cU_{s_0}$ is an affine image of the polytope $\Delta(V)$ hence it is compact, convex and it is a polytope. We can take $\gamma=0$ (it does not play any role since the MDP goes to an absorbing state right after period $t=0$). The return $R(\pi,r,P)$ is equal to $R(\pi,r,P)=-x \tr Ax$ for some $x \in \Delta(V).$
    Now Motzkin-Straus theorem states that given a graph $G$ with adjacency matrix $A \in \{0,1\}^{n \times n}$ then $\max_{x \in \Delta(V)} x \tr Ax = 1 - \frac{1}{\omega(G)}$ with $\omega(G)$ the clique number of $G$, i.e., the maximum size of a clique in $G$. 
    Hence 
    \[\inf_{(r,P) \in \cU} R(\pi,r,P) = - 1 + \frac{1}{\omega(G)}\]
    and if we could compare $\inf_{(r,P) \in \cU} R(\pi,r,P)$ to $-1+1/k$ we could compare $\omega(G)$ to some $k$, which concludes our reduction.
\end{proof}
Note that in the reduction for the proof of the above proposition we get that the coordinate $s_0$ of $\cUf$ can be identified with $\{(-x \tr A x, x) \; | \; x \in \Delta(V)\}$ which need not be convex, similarly as in \Cref{ex: lec 3 non cvx rp}. This hardness result shows that non-convexity of $\cUf$ can be the main difficulty when studying robust MDPs with uncertain rewards and uncertain transitions, and general reward models.

\section{Multi-objective MDPs and Pareto optimality}\label{sec: lec 3 multi-objective MDPs}
\begin{figure}[h!]
    \begin{center}
    \begin{subfigure}{0.49\textwidth}
    \begin{center}
        \includegraphics[width=\linewidth]{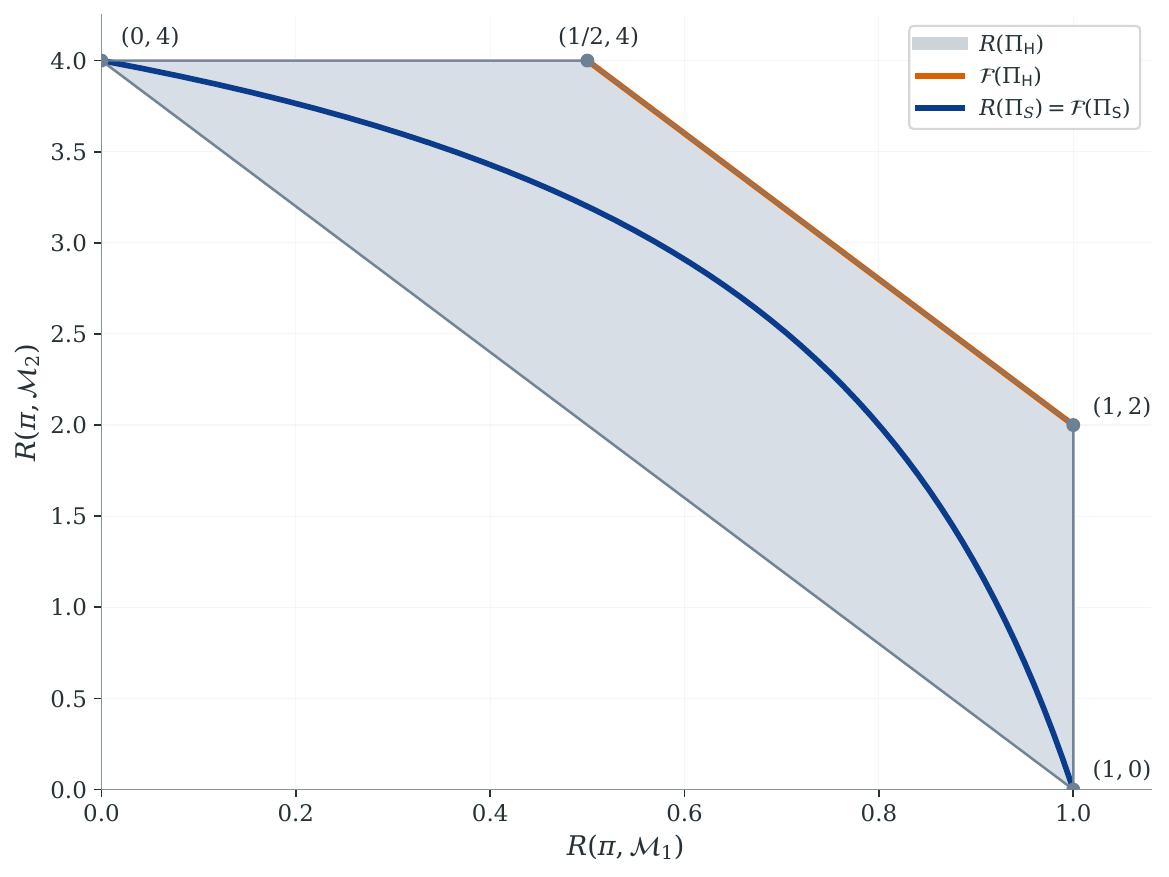}
    \end{center}
    \caption{Pareto front and achievable set}
    \label{fig: lec 3 momdp fig 1}
    \end{subfigure}
    \begin{subfigure}{0.49\textwidth}
    \begin{center}
        \includegraphics[width=\linewidth]{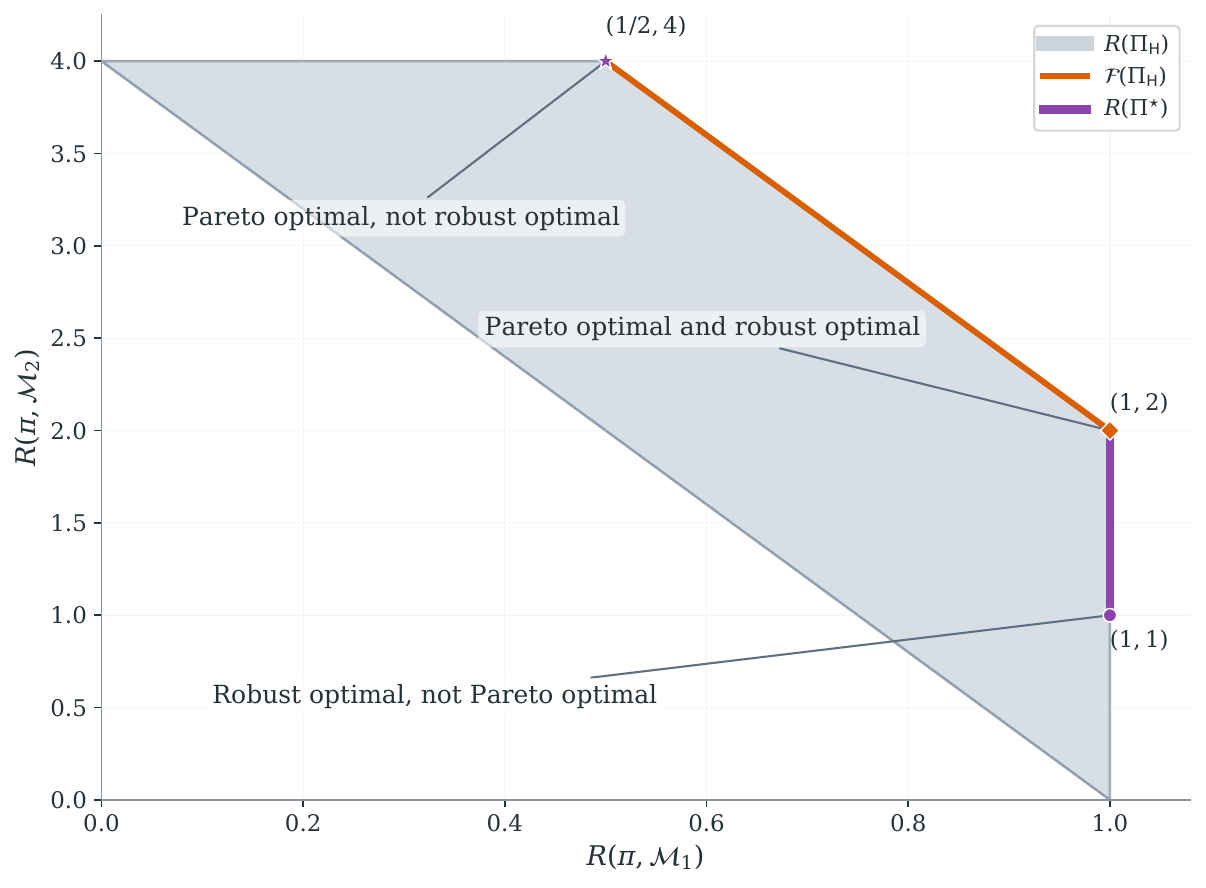}
    \end{center}
    \caption{Pareto front and solutions to \eqref{eq: lec 3 inf scalarized momdp}}
    \label{fig: lec 3 momdp fig 2}
    \end{subfigure}
    \end{center}
        \caption{Illustration of the achievable sets $R(\PiS),R(\PiH)$ and the Pareto fronts $\cF(\PiS),\cF(\PiH)$ (\Cref{fig: lec 3 momdp fig 1}), and the optimal solution of the scalarized problem with the infimum scalarization function as in \eqref{eq: lec 3 inf scalarized momdp}: $\Pi\opt = \arg \max_{\pi \in \PiH} \min \{ R(\pi,\cM_1),R(\pi,\cM_2)\}$. (\Cref{fig: lec 3 momdp fig 2}).}
        \label{fig: lec 3 pareto optimal}
    \end{figure}

This section focuses on {\em multi-objective MDPs}, where we are given several different MDP models $\cM_1,...,\cM_N$ with common state and action sets but (potentially) different instantaneous rewards and transition probabilities $\cM_i = (\cS,\cA,r^i,P^i,p_0)$ with $r^i \in \RSAS,P^i \in \DSSA$ for $i \in \{1,...,N\}$. We have some initial distribution $p_{0}$, common for each $i \in \{1,...,N\}$. Since the sets of states and actions are common across all MDP models here, it makes sense to talk about a policy $\pi \in \PiH$; however, it is unclear which objective we should pursue now. In particular a policy may perform well for some MDP model $\cM_i$ but poorly under some other model $\cM_j$, for some $j \neq i$. We will write $R(\pi,\cM_i) \in \R$ for the discounted return of $\pi$ in the $i$-th MDP model:
\[R(\pi,\cM_i) = \E^{\pi,P^i}\left[ \sum_{t=0}^{+\infty} \gamma^t r^i(s_t,a_t,s_{t+1}) \right].\]
We will write $R(\pi) \in \R^N$ for the vector of returns of a policy $\pi \in \PiH$ in each model $\cM_i$:
\[R(\pi) = (R(\pi,\cM_i))_{i \leq N}.\]
Naturally one arrives at the following definition of {\em Pareto optimality}.
\begin{definition}[Pareto domination and Pareto optimality]\label{def: lec 3 Pareto}
Let $\pi,\pi' \in \PiH$ and let $\Pi \subset \PiH$ be a class of policies (e.g., $\PiS$ or $\PiM$).
\begin{itemize}
    \item We say that $\pi$ {\em Pareto dominates} $\pi'$ if $R(\pi) \geq R(\pi')$ (componentwise) and $R(\pi,\cM) > R(\pi',\cM)$ for at least one model $\cM \in \{\cM_1,...,\cM_N\}$.
    \item We say that $\pi$ {\em strictly Pareto dominates} $\pi'$ if $R(\pi,\cM) > R(\pi',\cM)$ for any MDP model $\cM \in \{\cM_1,...,\cM_N\}$.
    \item When $\pi \in \Pi$, we say that $\pi$ is {\em Pareto optimal in $\Pi$} if no policy in $\Pi$ Pareto dominates it, and {\em weakly Pareto optimal} if no policy in $\Pi$ strictly Pareto dominates it.
\end{itemize}    
\end{definition}
With a slight abuse of terminology we will also say that ``$R$ Pareto dominates $R'$" for $R,R' \in \R^{N}$ when $R \geq R'$ componentwise and $R \neq R'$.
We also consider two sets of interest induced by a particular subclass of policies (think of the set of Markovian, or stationary, or stationary deterministic policies).
\begin{definition}[Pareto front and achievable set]
Let $\Pi \subset \PiH$ be a class of policies.
\begin{itemize}
    \item The {\em achievable set} $R(\Pi)$ is the set of vectors $R(\pi) \in \R^{N}$ induced by some $\pi \in \Pi$.
    \item The {\em Pareto front} $\cF(\Pi)$ is the set of all vectors $R(\pi) \in \R^{N}$ induced by policies $\pi \in \Pi$ that are Pareto optimal in $\Pi$.
\end{itemize}
\end{definition}
Note that we always have, by definition, $R(\Pi) \subseteq R(\PiH)$ when $\Pi \subset \PiH$.
One of the main questions in the literature on multi-objective MDPs is to compute Pareto optimal policies, and even better, when it is feasible, to enumerate the set of all Pareto optimal policies. Of course the definitions above are very general and can be extended to the case of MDPs with finite horizon, to MDPs with average (mean-payoff) return (see \Cref{sec: lec 3 sgs}), or to MDPs with risk-averse returns (see \Cref{sec: lec 3 risk averse}). Here, we focus on the discounted return, as this simplifies the exposition of the connection to the robust MDP model given above. We analyze a simple example below.
\begin{example}\label{ex: lec 3 momdp}
    Consider the MDP instance from \Cref{fig: lec 1 ex 3} with two MDP models that differ only in their transitions. We choose $\gamma=1/2$, we start in $s_0$ with probability $1$, and $p=1$ in $\cM_1$ and $p=0$ in $\cM_2$ (i.e., the effects of the actions $a_1$ and $a_2$ are reversed in $\cM_1$ and $\cM_2$). We only change the reward $r(s_0,a_2,s_1)$ from $1$ (in \Cref{fig: lec 1 ex 3}) to $r(s_0,a_2,s_1)=3$. Let $\pi \in \PiS$ be a stationary policy and $x = \pi(s_0,a_1)$. We write $v_{i}^{\pi}(s_0)$ for the value function derived by policy $\pi$ starting in $s_0$ in the MDP model $\cM_i$.  The Bellman equations give $v^{\pi}_1(s_0) = x + \frac{1}{2}(1-x)v^{\pi}_1(s_0), v^{\pi}_2(s_0) = 4(1-x) + \frac{x}{2}v^{\pi}_2(s_0)$. Therefore in this example
    \[R(\PiS) = \left\{ \left( \frac{2x}{1+x},\frac{8(1-x)}{2-x}\right)  \middle| \; x \in [0,1]\right\}.\]
    Note that the coordinates are respectively increasing and decreasing, hence $R(\PiS) = \cF(\PiS)$. It is easy to check that this set is not convex, see \Cref{fig: lec 3 pareto optimal}.

    We now focus on describing $R(\PiH)$. Given $\pi \in \PiH$ we write $x \in [0,1]$ for the probability that $\pi$ chooses $a_1$ at period $t=0$, and we now bound the return $R(\pi,\cM)$ for $\cM \in \{\cM_1,\cM_2\}$. In the MDP model $\cM_1$, with probability $x$ we obtain a return of $1$ (i.e. choosing $a_1$ at period $0$) and with probability $1-x$ we stay in $s_0$, and from then on we can obtain at most a return of $1/2$ (by choosing $a_1$). Hence $x \leq R(\pi,\cM_1) \leq x + \frac{1-x}{2}=\frac{1+x}{2}$. Similarly, in the MDP model $\cM_2$, if we choose $a_2$ at period $0$ then the return is $4$, otherwise we stay in $s_0$ and from there we can obtain at most $2$, by choosing $a_2$ (at period $t=1$). Hence $4(1-x) \leq R(\pi,\cM_2) \leq 4(1-x)+2x = 4-2x$. Hence $(R(\pi,\cM_1),R(\pi,\cM_2)) \in [x,(1+x)/2] \times [4(1-x),4-2x]$ for $x \in [0,1]$, and every point in the rectangle in right-hand side is attainable, by adjusting the probability to play $a_1$ or $a_2$ at period $t=1$, respectively after choosing $a_2$ or $a_1$ at period $t=0$. Therefore $R(\PiH) = \cup_{x \in [0,1]} [x,(1+x)/2] \times [4(1-x),4-2x]$. What does this look like? This is a quadrilateral with vertices $(0,4),(1/2,4),(1,2)$ and $(1,0)$, corresponding to the returns achievable by the four policies ``always play $a_2$", ``play $a_2$ then $a_1$", ``play $a_1$ then $a_2$", ``always play $a_1$", see \Cref{fig: lec 3 pareto optimal}. The Pareto front is the upper-right boundary of $R(\PiH)$, corresponding to $\cF(\PiH) = \left\{ \left( \frac{1+x}{2},4-2x \right) \; \middle| \; x \in [0,1] \right\}$.
\end{example}

We next investigate the properties of the sets $R(\Pi)$ and $\cF(\Pi)$ under various classes of policies $\Pi \subset \PiH$ and various assumptions on the unknown transitions and rewards. This will play an important role in the tractability (or intractability) of computing Pareto optimal policies, as we will see later in this section. 
\begin{theorem}\label{th: lec 3 R(Pi) properties}
Consider an instance $\cM_1,...,\cM_N$ of multi-model MDPs.
\begin{enumerate}
    \item Assume that the transition probabilities are the same across the MDP models $\cM_1,...,\cM_N$. Then:
    \begin{enumerate}
        \item The set $R(\PiH)$ is a polytope, i.e., it is the convex hull of finitely many points:
    \[R(\PiH) =  \conv\left( \{ \left(R\left(\pi,\cM_i\right)\right)_{i=1,...,N} \; | \; \pi \in \PiSD \} \right)\]
    \item $R(\PiH) = R(\PiS)$.
    \item $\cF(\PiH) = \cF(\PiS)$.
    \end{enumerate}
     \item Assume that the transition probabilities may differ across the MDP models $\cM_1,...,\cM_N$. Then:
    \begin{enumerate}
        \item The set $R(\PiH)$ is always compact and convex.
        \item We may have that $R(\PiS)$ is not convex. 
        \item We may have that $\cF(\PiH) \neq \cF(\PiS)$.
    \end{enumerate}
\end{enumerate}
\end{theorem}
\begin{proof}
    \begin{enumerate}
        \item We prove the three statements 1(a), 1(b), and 1(c) at the same time. Recall that by linearity of the expectation, the return is a linear function of the instantaneous rewards. Additionally, the state-action frequencies across all different models $\cM_1,...,\cM_N$ can be described by a single variable $\mu$ given that the transition probabilities are the same across all MDP models (we write this common kernel $P$ here). This means in particular that
\begin{align*}
    R(\PiH) = \{\left(R\left(\pi,\cM_i\right)\right)_{i=1,...,N} \; | \; \pi \in \PiH \} & = \{\left(R\left(\pi,\cM_i\right)\right)_{i=1,...,N} \; | \; \pi \in \PiS \} \\
    & = \left\{\left( \sum_{s,a} \mu(s,a)P(s,a)\tr r^i(s,a)  \right)_{i=1,...,N} \; | \; \mu \in \cQ \right\} \\
    & = \conv\left( \{\left(R\left(\pi,\cM_i\right)\right)_{i=1,...,N} \; | \; \pi \in \PiSD \} \right)
\end{align*}
where we recall that $\cQ$ is the set of state-action frequencies defined in \eqref{eq: lec 0 frequency polytope}, whose extreme points are the state-action frequencies induced by stationary deterministic policies. Hence $R(\PiH) = R(\PiS)$ is a polytope (it is the convex hull of finitely many points in $\R^{N}$) and every attainable return vector in $R(\PiH)$ can be achieved by some stationary policy $\pi \in \PiS$. From $R(\PiH)= R(\PiS)$ we directly get $\cF(\PiH) = \cF(\PiS)$. 
\item 
\begin{enumerate}
    \item 
We have seen in \Cref{prop: lec 1 continuity} that $\pi \mapsto R(\pi,\cM)$ is continuous for each model $\cM \in \{\cM_1,...,\cM_N\}$, and in \Cref{rmk: lec 0 topology} that $\PiH$ is compact (for the product topology). Hence the image $R(\PiH)$ is compact. We now show that it is a convex set. Let $\pi,\pi' \in \PiH, \lambda \in [0,1]$. We will construct a policy $\pi'' \in \PiH$ such that $R(\pi'',\cM) = \lambda R(\pi,\cM) +  (1-\lambda) R(\pi',\cM)$ simultaneously for all MDP models $\cM$. To do so, note that we just need to construct $\pi''$ such that, for any $t \in \N$ and any finite history $h_t = (s_0,a_0,...,s_{t-1},a_{t-1},s_{t})$, we have
\[\prod_{t'<t}\pi''(h_{t'})(a_{t'}) = \lambda \prod_{t'<t}\pi(h_{t'})(a_{t'}) + (1-\lambda) \prod_{t'<t}\pi'(h_{t'})(a_{t'}).\]
This can be done by choosing $\pi''(h_t)(a)$ proportional to $\lambda \left(\prod_{t'<t}\pi(h_{t'})(a_{t'}) \right) \pi(h_t)(a) + (1-\lambda) \left(\prod_{t'<t}\pi'(h_{t'})(a_{t'}) \right) \pi'(h_t)(a)$ for each $h_t$ and candidate action $a \in \cA$ when $\lambda \left(\prod_{t'<t}\pi(h_{t'})(a_{t'}) \right) + (1-\lambda) \left(\prod_{t'<t}\pi'(h_{t'})(a_{t'}) \right)$ is not null, otherwise $\pi''(h_t)(a)=1/|\cA|$ for each action $a$. By induction, the distribution of the trajectories under $\pi''$ is indeed the mixture of the laws under $\pi$ and $\pi'$,  hence we get that $R(\pi'',\cM) = \lambda R(\pi,\cM) +  (1-\lambda) R(\pi',\cM)$.
\item See \Cref{fig: lec 3 pareto optimal} and the analysis in \Cref{ex: lec 3 momdp}.
\item See \Cref{fig: lec 3 pareto optimal} and the analysis in \Cref{ex: lec 3 momdp}.
\end{enumerate}
    \end{enumerate}
\end{proof}
\Cref{th: lec 3 R(Pi) properties} gives a sufficient condition for $\cF(\PiS) = \cF(\PiH)$, i.e., for policies in $\PiS$ to suffice to attain every point of the Pareto front $\cF(\PiH)$. In particular, this is greatly dependent on the assumptions on the unknown transitions and rewards: the case of a common transition kernel has much more favorable properties than the case where the kernels differ across models. We next investigate the computation of Pareto optimal policies based on so-called ``scalarization functions".

\paragraph{Scalarization functions.}
A natural approach to deal with multi-objective optimization (beyond MDPs) is to introduce a {\em scalarization function}, which is simply a map $f:\R^{N} \to \R$. We then seek to find a policy that optimizes $f\left( R(\pi)\right)$ over a chosen policy class $\Pi \subset \PiH$, i.e., we want to solve
\begin{equation}\label{eq: lec 3 scalarized momdp}
    \sup_{\pi \in \Pi} f(R(\pi)).
\end{equation}
This model is of course very general: 
\begin{itemize}
    \item When $f(R) = R_i$ for $R \in \R^{N}$ then \eqref{eq: lec 3 scalarized momdp} recovers the nominal MDP model from \Cref{chap: lec 0 mdp}. More generally, when the scalarization function is linear: $f(R) = \sum_{i=1}^{N} w_i R_i$ for some $w \in \Delta(\{1,...,N\})$, then we obtain the so-called {\em multi-model MDPs}. Note that the class of functions $f:R \mapsto w \tr R$ for some $w \in \R^{N}_{+}$ yields the same optimizers of \eqref{eq: lec 3 scalarized momdp} as when $w \in \Delta(\{1,...,N\})$ (since we can just renormalize the weight vector).
    \item When $f(R) = \min_{i=1,...,N} R_i$ for $R \in \R^{N}$  then \eqref{eq: lec 3 scalarized momdp} recovers the robust MDP model with uncertain rewards and transitions from \Cref{sec: lec 3 uncertain rewards and transitions}, when there are finitely many unknown scenarios.
\end{itemize}
We now discuss two important questions:
\begin{itemize}
    \item {\bf Question 1:} Given a policy class $\Pi$, when does solving \eqref{eq: lec 3 scalarized momdp} over $\Pi$ yield Pareto optimal policies in $\Pi$? Can every Pareto optimal policy in $\Pi$ be computed this way for an appropriate scalarization function?
    \item {\bf Question 2:} Given a policy class $\Pi$, are there efficient algorithms for solving \eqref{eq: lec 3 scalarized momdp} (so that we can leverage the answer to the previous question to compute Pareto optimal policies)?
\end{itemize}
The answers to the above questions depend on (A) the geometry of $R(\Pi)$, which in turn depends on the assumptions made on the rewards and transition probabilities in the models $\cM_1,...,\cM_N$, and (B) the choice of the scalarization function $f$. We will study each set of assumptions on the uncertain models separately, and we will focus on linear scalarizations (since these are by far the most popular scalarization functions for multi-objective MDPs) and on worst-case scalarization (which recovers robust MDPs). We start with linear scalarization.
\subsection{Linear scalarization and multi-model MDPs}
Here we focus on the case where we are given some $w \in \Delta(\{1,...,N\})$ and where $f:R \mapsto w \tr R$, so that \eqref{eq: lec 3 scalarized momdp} becomes, for a chosen class of policies $\Pi \subset \PiH$, the following optimization problem:
\begin{equation}\label{eq: lec 3 linear scalarized momdp}
    \sup_{\pi \in \Pi} \sum_{i=1}^{N} w_i R(\pi,\cM_i)
\end{equation}
We next give a general response to the two questions outlined above, for the case of linear scalarization. We start with the first question (when linear scalarization recovers Pareto optimal solutions). 
\begin{theorem}\label{th: lec 3 momdp linear}
Let $\Pi \subseteq \PiH$ be a class of policies (e.g. $\Pi \in \{\PiS,\PiM,\PiH\}$). The following statements hold:
\begin{enumerate}
    \item Every optimal solution to \eqref{eq: lec 3 linear scalarized momdp} for some $w\in\Delta(\{1,\ldots,N\})$ is weakly Pareto optimal in $\Pi$. 
    \item If $w_i>0$ for every $i$, every optimal solution to \eqref{eq: lec 3 linear scalarized momdp} is Pareto optimal in $\Pi$.
    \item If $R(\Pi)$ is a convex set, then every weakly Pareto optimal policy in $\Pi$ is an optimal  solution of the scalarized problem for some $w\in\Delta(\{1,\ldots,N\})$ and $f:R\mapsto w \tr R$.
    \item If $R(\Pi)$ is a polytope, then every Pareto optimal policy in $\Pi$ is an optimal solution of \eqref{eq: lec 3 linear scalarized momdp} for some $w\in\Delta(\{1,\ldots,N\})$
    with $w_i>0$ for every $i$, and $f:R\mapsto w \tr R$.
\end{enumerate}
\end{theorem}
\begin{proof}
        \begin{enumerate}
        \item This first statement simply follows from the non-negativity of the weight vector $w \in \Delta(\{1,...,N\})$. Indeed if $\pi$ is optimal in \eqref{eq: lec 3 linear scalarized momdp} then if there is a policy in $\Pi$ that strictly Pareto dominates $\pi$ it would yield a strictly better objective value in \eqref{eq: lec 3 linear scalarized momdp}. Hence $\pi$ optimal in \eqref{eq: lec 3 linear scalarized momdp} implies that $\pi$ is weakly Pareto optimal. 
        \item In the case where $w_i>0$ for every model $i \in \{1,...,N\}$ then the same argument as for the first statement shows that $\pi$ optimal in \eqref{eq: lec 3 linear scalarized momdp} implies that $\pi$ is Pareto optimal.
        \item This is an implication of the hyperplane separation theorem (see Appendix \ref{app:sec analysis}).
        Let $R \in R(\Pi)$ be weakly Pareto optimal. Then $R(\Pi) - R$ and $(0,+\infty)^{N}$ are two disjoint convex sets. Therefore we can find a direction $w \in \R^{N}_{+}, w \neq (0,...,0)$ such that $w \tr (R'-R) \leq 0, \forall \; R' \in R(\Pi)$. Renormalizing $w$ to an element in the simplex $\Delta(\{1,...,N\})$ yields the conclusion.
        \item This is an application of duality in linear programming. Let $R$ be the vector of returns of a Pareto optimal policy. Then $R$ is the only feasible vector in the optimization program $\max \left\{ \sum_{i=1}^{N} R'_i \; | \; R' \in R(\Pi),R' \geq R \right\}$. From strong duality we get some nonnegative Lagrange multipliers $\lambda_1,...,\lambda_N \geq 0$ such that $R \in \arg \max \left\{ \sum_{i=1}^{N} (1+\lambda_i) R'_i \; | \; R' \in R(\Pi)\right\}$. This yields the optimization problem $\sup_{\pi \in \Pi} \left\{ \sum_{i=1}^{N} (1+\lambda_i) R(\pi,\cM_i) \; | \; \pi \in \Pi\right\}$ of the form \eqref{eq: lec 3 linear scalarized momdp} with positive coefficients.
    \end{enumerate}
\end{proof}
Combining \Cref{th: lec 3 momdp linear} with \Cref{th: lec 3 R(Pi) properties}, we get the following implications:
\begin{itemize}
    \item When the transition probabilities are common across all models, then we can recover every Pareto optimal policy in $\PiH$ as a solution of \eqref{eq: lec 3 linear scalarized momdp} for some positive weight vector $w$, and since $R(\PiS) = R(\PiH)$ in this case, stationary (randomized) policies suffice to attain every point of the Pareto frontier $\cF(\PiH)$.
    \item When the transition probabilities may be different across the models, then we can recover every weakly Pareto optimal policy in $\PiH$ as a solution of \eqref{eq: lec 3 linear scalarized momdp} for some nonnegative weight vector $w$. We may need history-dependent policies to reach some Pareto optimal vector in $\cF(\PiH)$.
\end{itemize}
We now describe the complexity properties of the optimization program in \eqref{eq: lec 3 linear scalarized momdp}. 
\paragraph{Tractability of \eqref{eq: lec 3 linear scalarized momdp}: the case of uncertain rewards.}
We assume here that the transition probabilities are the same across all MDP models $\cM_1,...,\cM_N$ and are equal to some $P \in \DSSA$.
In this case we can reformulate $R(\pi)$ given some $\pi \in \PiH$. Indeed, by linearity of the expectation, the return is a linear function of the rewards. Hence we get, for some $f:R \mapsto w \tr R$ with $w \in \R^{N}_{+}$, that
\begin{align*}
    f(R(\pi))  = \sum_{i =1}^{N} w_i R(\pi,\cM_i) & = \sum_{i =1}^{N} w_i  \E^{\pi,P}\left[ \sum_{t=0}^{+\infty} \gamma^t r^i(s_t,a_t,s_{t+1}) \right] \\
    & =   \E^{\pi,P}\left[ \sum_{t=0}^{+\infty} \gamma^t \left( \sum_{i =1}^{N} w_i r^i(s_t,a_t,s_{t+1})\right) \right] \\
    & = R(\pi,\bar{\cM}) 
\end{align*}
for $\bar{\cM}$ the MDP model $(\cS,\cA,\sum_{i=1}^{N} w_i r^i, P,p_0)$.
Therefore in this case \eqref{eq: lec 3 linear scalarized momdp} boils down to a nominal MDP model as studied in \Cref{chap: lec 0 mdp}, where the instantaneous rewards are the weighted sum of the instantaneous rewards of each model, so we can solve \eqref{eq: lec 3 linear scalarized momdp} efficiently (see \Cref{sec: lec 0 algo} for the description of the algorithms studied in these notes) for $\Pi \in \{\PiH,\PiM,\PiS\}$ (and more generally any class of policies containing $\PiSD$).
\paragraph{Tractability of \eqref{eq: lec 3 linear scalarized momdp}: the case of uncertain transitions.} 
When the transition probabilities are also uncertain and may differ across MDP models, things can become very complicated. For instance, it is known that solving the following optimization problem is NP-hard:
\[\sup_{\pi \in \PiH} \frac{1}{2}R(\pi,\cM_1) + \frac{1}{2}R(\pi,\cM_2)\]
 for the discounted return considered in this section, but this remains the case also for MDPs with finite horizon. Hence solving \eqref{eq: lec 3 linear scalarized momdp} may not be feasible when the transition probabilities are not the same across all models.
\subsection{Worst-case scalarization and robust MDPs}
This is the case most akin to robust MDPs as we studied here, and we consider solving the following optimization problem to compute Pareto optimal policies:
\begin{equation}\label{eq: lec 3 inf scalarized momdp}
    \sup_{\pi \in \Pi} \min_{i=1,...,N} R(\pi,\cM_i)
\end{equation}
We write $\Pi\opt = \arg \max_{\pi \in \Pi} \min_{i=1,...,N} R(\pi,\cM_i)$.
The motivation for solving the above mathematical program also sometimes comes from fairness considerations (we want to perform well across all scenarios and not just the most important ones).
The following theorem characterizes the Pareto optimality properties of optimal solutions in \eqref{eq: lec 3 inf scalarized momdp}.
\begin{theorem}\label{th: lec 3 momdp min}
Let $\Pi \subset \PiH$ be a policy class such that $R(\Pi)$ is compact. The following statements hold:
\begin{enumerate}
    \item Every optimal solution of \eqref{eq: lec 3 inf scalarized momdp} is weakly Pareto optimal in $\Pi$.
    \item There exists at least one optimal solution of \eqref{eq: lec 3 inf scalarized momdp} that is Pareto optimal in $\Pi$.
    \item Not every optimal solution of \eqref{eq: lec 3 inf scalarized momdp} is necessarily Pareto optimal in $\Pi$.
    \item Not every Pareto optimal policy in $\Pi$ can necessarily be recovered as a solution of \eqref{eq: lec 3 inf scalarized momdp}.
\end{enumerate}
\end{theorem}
\begin{proof}
\begin{enumerate}
    \item This is a direct reformulation of strict Pareto domination (and therefore weak Pareto optimality): if $\pi \in \Pi$ is an optimal solution to \eqref{eq: lec 3 inf scalarized momdp} then if $\pi' \in \Pi$ strictly Pareto dominates $\pi$ we would have $\min_{i} R(\pi',\cM_i) > \min_{i} R(\pi,\cM_i)$, which contradicts the optimality of $\pi$ in \eqref{eq: lec 3 inf scalarized momdp}.
    \item Consider $\Pi\opt = \arg \max_{\pi \in \Pi} \min_{i} R(\pi,\cM_i)$ the set of optimal solutions in \eqref{eq: lec 3 inf scalarized momdp}. The set $R(\Pi\opt)$ is compact by compactness of $R(\Pi)$ and continuity of $R \mapsto \min_i R_i$. Now consider the subset $\Pi\opt_{{\sf sum}} = \arg \max_{\pi \in \Pi\opt} \sum_{i=1}^{N} R(\pi,\cM_i)$, which corresponds to the set of optimal solutions in \eqref{eq: lec 3 inf scalarized momdp} that also maximize the sum. Then policies in $\Pi\opt_{{\sf sum}}$ are both Pareto optimal and, by construction, optimal solutions in \eqref{eq: lec 3 inf scalarized momdp}.
    \item See \Cref{fig: lec 3 pareto optimal} and the analysis in \Cref{ex: lec 3 momdp} and below.
    \item See \Cref{fig: lec 3 pareto optimal} and the analysis in \Cref{ex: lec 3 momdp} and below.
\end{enumerate}
\end{proof}
\begin{continuance}{ex: lec 3 momdp}
    We continue with the multi-objective MDP instance from \Cref{ex: lec 3 momdp}, and we compute the set $\Pi\opt = \arg \max_{\pi \in \PiH} \min \{R(\pi,\cM_1),R(\pi,\cM_2)\}$. We have shown in \Cref{ex: lec 3 momdp} that $R(\pi,\cM_1) \leq 1$ for any $\pi \in \PiH$, and we have seen that we can achieve $(1,2)$ with the policy ``play $a_1$ then $a_2$", hence $\max_{\pi \in \PiH} \min \{R(\pi,\cM_1),R(\pi,\cM_2)\}=1$. The policies achieving this maximum must be such that $R(\pi,\cM_1)=1$ (since this is the best one can achieve in $\cM_1$). Our characterization of $R(\PiH)$ from \Cref{ex: lec 3 momdp} gives that for such policies we must have $R(\pi,\cM_2) \in [0,2]$, and for a policy to be optimal in \eqref{eq: lec 3 inf scalarized momdp} we must also have $R(\pi,\cM_2) \in [1,2]$. Hence $R(\Pi\opt) = \{1\} \times [1,2]$, attained by policies choosing $a_1$ at period $0$, then if still at $s_0$ at period $1$ choosing $a_2$ with some probability $y \in [0,1]$, otherwise choosing $a_2$ at period $2$ (if still in $s_0$ at period $2$). This yields $R(\pi) = (1,1+y)$ for $y \in [0,1]$. Note that $(1,1) \in R(\Pi\opt)$ but this vector is Pareto dominated by $(1,2)$, therefore optimality in \eqref{eq: lec 3 inf scalarized momdp} does not imply Pareto optimality; additionally, $(1/2,4) \in \cF(\PiH)$ but is not attained by an optimal solution in \eqref{eq: lec 3 inf scalarized momdp}. This provides the proofs for the third and fourth statements of \Cref{th: lec 3 momdp min}.
\end{continuance}
\paragraph{Tractability properties of \eqref{eq: lec 3 inf scalarized momdp}: the case of uncertain rewards.}
In this case with $\Pi = \PiH$ we are back to the case of robust MDPs with uncertain rewards as studied in \Cref{sec: lec 3 uncertain rewards}. In particular for $\cR = \{r^1,...,r^{N}\}$ then \eqref{eq: lec 3 inf scalarized momdp} is equivalent to \eqref{eq: lec 3 rmdp r}, i.e., to $\sup_{\pi \in \PiH} \min_{r \in \cR} R(\pi,r)$; we refer to \Cref{sec: lec 3 uncertain rewards} for the favorable properties of this optimization problem (among other things, we can solve this problem via linear programming).
\paragraph{Tractability properties of \eqref{eq: lec 3 inf scalarized momdp}: the case of uncertain transitions.}
Here for $\Pi = \PiH$ the properties of \eqref{eq: lec 3 inf scalarized momdp} are similar to those of \eqref{eq: lec 3 linear scalarized momdp}; it is known (see references at the end of this chapter) that solving
\[\sup_{\pi \in \PiH} \min \{R(\pi,P^1), R(\pi,P^2)\}\]
is NP-hard in all generality. When the set $\cP = \{P^1,...,P^N\}$ is s-rectangular (not necessarily convex) and the rewards are common across the models, then \Cref{sec: lec 2 noncvx srec} provides algorithms for finding a best stationary policy (not necessarily optimal among history-dependent policies). If $\cP$ is sa-rectangular then the results in \Cref{sec: lec 2 sarec} apply directly (we don't need convexity for solving sa-rectangular robust MDPs).
If the rewards also differ across models, then we refer to \Cref{sec: lec 3 uncertain rewards and transitions}.
\paragraph{Connection with linear scalarizations.}
We conclude this section with a short comparison between linear scalarization and worst-case scalarization. Indeed we always have
\[\sup_{\pi \in \PiH} \min_{i=1,...,N} R(\pi,\cM_i) = \sup_{\pi \in \PiH} \min_{w \in \Delta(\{1,...,N\})} \sum_{i=1}^{N} w_i R(\pi,\cM_i)\]
so we can interpret the worst-case scalarization as the infimum over all linear scalarizations. This explains why it has ``weaker" properties (in terms of recovering Pareto optimal policies) than linear scalarization (as evident from comparing \Cref{th: lec 3 momdp linear} with \Cref{th: lec 3 momdp min}).

\section{Robust MDPs with time-varying transition probabilities}\label{sec: lec 3 time var P}
In \Cref{chap: lec 1 robust MDP} and \Cref{chap: lec 2 rectangular models} we assumed that the transition probabilities were chosen by the adversary once and for all before the decision maker starts choosing some actions. However we can imagine some setups where the transition probabilities are non-stationary and may change over time. We then need to choose a model for how the transition kernel can evolve.
\subsection{Three models of dynamic transition uncertainties}
Let $\cP \subset \DSSA$ be any uncertainty set. We will assume that at every period $t$, the transition kernel is some $P^t \in \cP$ chosen by the adversary, that is, the transition probabilities may change but always belong to the same set (this simplifies greatly the arguments in this section). We now need to rigorously describe how the transition $P^t$ for period $t$ is chosen by the adversary. Similarly to when we defined the classes of stationary, Markovian, and history-dependent policies for the decision maker in \Cref{chap: lec 0 mdp}, we need to describe what kind of information the adversary can observe when choosing the transition probabilities $P^t$ at each time period. Of course, at each period $t$ the transition probabilities $P^t$ need to be determined before the next state $s_{t+1}$ is realized. This leaves basically three levels of information that the adversary could use:
\begin{itemize}
    \item Either $P^t$ is chosen after observing the current $(s_t,a_t)$ state-action realization,
    \item Or $P^t$ is chosen after observing the current state $s_t$ but before the decision maker chooses their action $a_t$,
    \item Or $P^t$ is chosen before observing the current state $s_t$.
\end{itemize}
These three cases have deep connections with sa-rectangular, s-rectangular, and non-rectangular models, as we investigate next.
\paragraph{What is observed? (part 2)} 
We will use the notation
\[\cG_{t-1} = (\cS \times \cA)^{t}, \quad t \in \N\]
to denote the set of all finite histories $(s_0,a_0,...,s_{t-1},a_{t-1})$ up to time $t-1$, completed with the last chosen action (recall that by definition from \Cref{chap: lec 0 mdp} histories do not include the last chosen action). We write $g_t$ for a generic element in $\cG_{t-1}$.
For each $t \in \N$, a decision rule of the adversary at period $t$ maps a realization $(g_t,s_t,a_t) \in (\cS \times \cA)^{t+1}$ to some transition probabilities in $\DS$, representing the distribution of the next state $s_{t+1} \in \cS$ given $s_t,a_t$. A history-dependent policy $\kappa$ of the adversary is given by a sequence of decision rules $\left(\kappa_t\right)_{t \in \N}$ with $\kappa_t: (\cS \times \cA)^{t+1} \to \DS$ for each $t \in \N$. 

We will also need the following notation. For a given uncertainty set $\cP \subset \DSSA$, we consider its sets of $sa$-marginals $\cP_{s,a} \coloneqq \{P(s,a) \; | \; P \in \cP\} \subset \DS$ for each pair $(s,a) \in \cS \times \cA$, and its sets of $s$-marginals $\cP_{s} \coloneqq \{P(s) \; | \; P \in \cP \} \subset \DSA$ for each $s \in \cS$. 

We consider the following models of time-varying transition probabilities:
\begin{itemize}
\item {\em Model 1 (``post-state, post-action"): $\cPSA_{\sf H}$.} The set $\cPSA_{\sf H}$ coincides with the set of adversary policies that choose their decision rules consistent with $\cP$, after observing the current state-action pair $(s_t,a_t)$:
\begin{equation}\label{eq: lec 3 def cPSA}
    \cPSA_{\sf H} = \{ (\kappa_t)_{t \in \N} \; | \; \kappa_t(g_t,s_t,a_t) \in \cP_{s_t,a_t}, \forall t \in \N, g_t \in \cG_{t-1}, s_t \in \cS, a_t \in \cA\}.
\end{equation}
\item {\em Model 2 (``post state, pre-action"): $\cPS_{\sf H}$.} The set $\cPS_{\sf H}$ coincides with the set of adversary policies that choose their decision rules consistent with $\cP$, after observing the current state $s_t$ :
\begin{equation}\label{eq: lec 3 def cPS}
    \cPS_{\sf H} = \{ (\kappa_t)_{t \in \N} \; | \; (\kappa_t(g_t,s_t,a))_{a \in \cA} \in \cP_{s_t}, \forall t \in \N, g_t \in \cG_{t-1}, s_t \in \cS\}.
\end{equation}

\item {\em Model 3 (``pre-state"): $\cPG_{\sf H}$.} The set $\cPG_{\sf H}$ coincides with the set of adversary policies that choose their decision rules consistent with $\cP$, before observing the current state $s_t$ :
\begin{equation}\label{eq: lec 3 def cPG}
    \cPG_{\sf H} = \{ (\kappa_t)_{t \in \N} \; | \; (\kappa_t(g_t,s,a))_{s,a} \in \cP, \forall t \in \N, g_t \in \cG_{t-1}\}.
\end{equation}
\end{itemize}
Given that Model $i+1$ imposes some additional constraints compared to Model $i$ above (for $i \in \{1,2\}$), we easily get that $\cPG_{\sf H} \subseteq \cPS_{\sf H} \subseteq \cPSA_{\sf H}$.

The first important thing to notice is that in the definitions of Model 1 and Model 2, the uncertainty set $\cP$ only appears via its marginals ($\cP_{s,a}$ for each state-action pair $(s,a)$ for Model 1 and $\cP_s$ for each state $s$ in Model 2). This suggests that only certain ``rectangularizations" of $\cP$ matter in Model 1 and Model 2. We make this notion rigorous next.
\begin{definition}\label{def: lec 3 rectangular extensions}
    Let $\cP \subset \DSSA$ be an uncertainty set.
    \begin{itemize}
        \item The {\em s-rectangular extension} of $\cP$ is the uncertainty set $\recs(\cP) \subset \DSSA$ defined 
        \[\recs(\cP) = \times_{s \in \cS} \{P(s) \; | \; P \in \cP\}.\]
        \item The {\em sa-rectangular extension} of $\cP$ is the uncertainty set $\recsa(\cP) \subset \DSSA$ defined 
        \[\recsa(\cP) = \times_{(s,a) \in \cS \times \cA} \{P(s,a) \; | \; P \in \cP\}.\]
        \end{itemize}
\end{definition}
Equivalently, the s-rectangular extension $\recs(\cP)$ is the smallest s-rectangular set containing $\cP$, and the sa-rectangular extension $\recsa(\cP)$ is the smallest sa-rectangular set containing $\cP$. Note that obviously the s-rectangular extension of an s-rectangular set is itself, and the sa-rectangular extension of an sa-rectangular set is itself. The extensions play a crucial role for studying time-varying transitions, as shown in the next lemma.
\begin{lemma}\label{lem: lec 3 Model 12 eq extensions}
    Let $\cP \subset \DSSA$ be an uncertainty set.
    \begin{enumerate}
        \item We have $\cP_{\sf H}^{\sf SA} = \recsa(\cP)^{\sf SA}_{\sf H}$ and $\cP^{\sf S}_{\sf H} = \recs(\cP)^{\sf S}_{\sf H}$.
        \item If additionally $\cP$ is s-rectangular, then $\cP_{\sf H}^{\sf S}=\cP_{\sf H}^{\sf G}$.
        \item If additionally $\cP$ is sa-rectangular, then $\cP_{\sf H}^{\sf SA}=\cP_{\sf H}^{\sf S}=\cP_{\sf H}^{\sf G}$.
    \end{enumerate}
\end{lemma}
\begin{proof}
\begin{enumerate}
    \item This follows from the observation that the $s$-th marginal of $\recs(\cP)$ is the same as the $s$-th marginal of $\cP$, by construction. Since $\cPS_{\sf H}$ only depends on $(\cP_{s})_{s \in \cS}$ we must have that $\cPS_{\sf H} = \recs(\cP)^{\sf S}_{\sf H}$. The proof for the set $\cPSA_{\sf H}$ is similar.
    \item Let $t \in \N, g_t \in \cG_{t-1}$. For $\cP$ s-rectangular we get that the constraint $``(\kappa_t(g_t,s,a))_{s,a} \in \cP"$ in the definition of $\cPG_{\sf H}$ as in \eqref{eq: lec 3 def cPG} is equivalent to $``(\kappa_t(g_t,s,a))_{a} \in \cP_s, \forall \; s \in \cS"$, which is the constraint appearing in the definition of $\cPS_{\sf H}$ as in \eqref{eq: lec 3 def cPS}.
    \item This proof is similar to that of the statement right above.
\end{enumerate}
\end{proof}
Therefore in the models of time-varying adversaries described by Model 1 and Model 2, it is sufficient to consider sa-rectangular and s-rectangular uncertainty sets.  Additionally, for any of the three cases above, we can then ask that the decisions over time are always the same, in which case we say that the adversary policy is {\em stationary} and we identify the set of stationary policies with $\recsa(\cP),\recs(\cP)$ and $\cP$ respectively.
\begin{remark}
    It is tempting to also define {\em Markovian} adversaries as time-varying transition probabilities that only depend on the current states, since we also defined this notion for the decision maker; e.g. for Model 3 we could identify this set as $\cP^{\N}$, see  \Cref{sec: lec 3 ref} for some references on these models.
\end{remark}
\subsection{Main results for Model 1 and Model 2}
Given a policy of the adversary $\kappa \in \cPSA_{\sf H}$ (which contains the other two classes of history-dependent adversary policies) and a policy $\pi \in \PiH$ of the decision maker, we still get a distribution over the set of infinite trajectories, and therefore we can extend the function $P \mapsto R(\pi,P)$ defined on $\cP$ to $\kappa \mapsto R(\pi,\kappa)$ defined on $\cPSA_{\sf H}$ naturally as
\[R(\pi,\kappa)  = \E^{\pi,\kappa}\left[ \sum_{t=0}^{+\infty} \gamma^t r(s_t,a_t,s_{t+1}) \right].\]
The main goal in this section is to understand the properties of the policy evaluation problems and robust MDP problems with worst-case over time-varying transition probabilities. 
\[\inf_{\kappa \in \cP^{\sf i}_{\sf H} }R(\pi,\kappa)\]
and robust MDP problem
\[\sup_{\pi \in \PiH} \inf_{\kappa \in \cP^{\sf i}_{\sf H} }R(\pi,\kappa)\]
for $i \in \{{\sf SA, S, G}\}$.
We choose to organize the main results by model of time-varying transitions.
\paragraph{The case of Model 1.} We start with Model 1. We focus on the properties of the policy evaluation problem
\[\inf_{\kappa \in \cPSA_{\sf H} }R(\pi,\kappa)\]
and robust MDP problem
\begin{equation}\label{eq: lec 3 rmdp cPSA}
    \sup_{\pi \in \PiH} \inf_{\kappa \in \cPSA_{\sf H} }R(\pi,\kappa).
\end{equation}
\begin{theorem}[Robust MDPs with time-varying transitions, Model 1]\label{th: lec 3 Model 1}
    Let $\cP \subset \DSSA$ be compact. The following statements hold.
\begin{enumerate}
    \item {\em Policy evaluation.} Let $\pi \in \PiS$. Then the worst-case return of $\pi$ against $\kappa \in \cPSA_{\sf H}$ is attained at a stationary transition kernel in the sa-rectangular extension of $\cP$:
    \[\inf_{\kappa \in \cPSA_{\sf H}} R(\pi,\kappa) = \min_{P \in \recsa(\cP)} R(\pi,P).\]
    \item {\em Optimal policies.} An optimal policy for the optimization problem \eqref{eq: lec 3 rmdp cPSA} can always be chosen stationary:
    \begin{align*}
        \sup_{\pi \in \PiH} \inf_{\kappa \in \cPSA_{\sf H}} R(\pi,\kappa) & = \sup_{\pi \in \PiS} \inf_{\kappa \in \cPSA_{\sf H}} R(\pi,\kappa)
    \end{align*}
    In fact, a stationary uniformly optimal policy always exists for \eqref{eq: lec 3 rmdp cPSA} and it can be chosen to be also uniformly optimal for the robust MDP problem over the stationary sa-rectangular extension of $\cP$, i.e., $\sup_{\pi \in \PiH} \inf_{\kappa \in \cPSA_{\sf H}} R(\pi,\kappa) = \max_{\pi \in \PiS} \min_{P \in \recsa(\cP)} R(\pi,P)$.
\end{enumerate}
\end{theorem}
The above theorem is an immediate consequence of \Cref{lem: lec 3 Model 12 eq extensions} and the following proposition, specialized for sa-rectangular uncertainty sets (note that we don't need convexity here, similarly as in \Cref{sec: lec 2 sarec}).
\begin{proposition}[sa-rectangular time-varying transitions]\label{prop: lec 3 Model 1}
    Let $\cP \subset \DSSA$ be sa-rectangular compact. The following statements hold.
\begin{enumerate}
    \item {\em Policy evaluation.} Let $\pi \in \PiS$. Then the worst-case return of $\pi$ against time-varying transitions in $\cPSA_{\sf H}$ is attained at a stationary transition kernel in the sa-rectangular extension of $\cP$:
    \[\inf_{\kappa \in \cPSA_{\sf H}} R(\pi,\kappa) = \min_{P \in \cP} R(\pi,P).\]
    \item {\em Optimal policies.} An optimal policy for the optimization problem \eqref{eq: lec 3 rmdp cPSA} can always be chosen stationary:
    \begin{align*}
        \sup_{\pi \in \PiH} \inf_{\kappa \in \cPSA_{\sf H}} R(\pi,\kappa) & = \max_{\pi \in \PiS} \inf_{\kappa \in \cPSA_{\sf H}} R(\pi,\kappa)
    \end{align*}
    In fact, a stationary {\em uniformly} optimal policy always exists for \eqref{eq: lec 3 rmdp cPSA} and it can be chosen to be also uniformly optimal for the robust MDP problem over the stationary transition probabilities, i.e., $\sup_{\pi \in \PiH} \inf_{\kappa \in \cPSA_{\sf H}} R(\pi,\kappa) = \max_{\pi \in \PiS} \min_{P \in \cP} R(\pi,P)$.
\end{enumerate}
\end{proposition}
\begin{proof}
\begin{enumerate}
    \item We can interpret $\inf_{\kappa \in \cPSA_{\sf H}} R(\pi,\kappa)$ as an adversarial MDP where the set of states is $\cS \times \cA$ and the compact set of actions available at state $(s,a)$ is $\cP_{sa}$. The instantaneous reward for choosing action $P(s,a) \in \cP_{sa}$ at $(s,a)$ is $P(s,a)\tr r(s,a)$ (independent of the next state $(s',a')$) and the transition probabilities from $(s,a)$ to a next state $(s',a')$ is $P(s,a,s')\pi(s',a') \in [0,1]$ (at this point the reader should notice the similarity with the adversarial MDP defined in \Cref{sec: lec 2 adv} for s-rectangular model). We conclude that a stationary deterministic optimal adversarial policy exists, i.e., $\inf_{\kappa \in \cPSA_{\sf H}} R(\pi,\kappa) = \min_{P \in \cP} R(\pi,P)$.
    \item The proof is very close to the proof for the existence of stationary optimal policies for stationary kernels (beginning of \Cref{sec: lec 2 prf unf opt}). We have
\begin{align*}
 \sup_{\pi \in \PiH} \inf_{P  \in \cP} 
 R(\pi,P) & \leq  \inf_{P  \in \cP} \sup_{\pi \in \PiH}  R(\pi,P) \\
 & \leq \inf_{P  \in \cP} \max_{\pi \in \PiS} 
 R(\pi,P) \\
 & = \max_{\pi \in \PiS} \min_{P  \in \cP}  
 R(\pi,P) \\
  & = \max_{\pi \in \PiS} \inf_{\kappa  \in \cPSA_{\sf H}}  
 R(\pi,\kappa) \\
 &\leq \sup_{\pi \in \PiH} \inf_{\kappa  \in \cPSA_{\sf H}} 
 R(\pi,\kappa)
\end{align*}
where the first inequality is weak duality, the second inequality is because the inner maximization is an MDP, the first equality is strong duality for sa-rectangular compact uncertainty sets as in \Cref{sec: lec 2 sarec}, the second equality is via the first statement of this proposition, and the third is by $\PiS \subset \PiH$. Since the last term above is always smaller than the first term above (identifying $\cP$ with constant maps in $\cPSA_{\sf H}$), we get that all terms above are equal, so that a stationary optimal policy exists for~\eqref{eq: lec 3 rmdp cPSA} and can be computed as a solution to $\max_{\pi \in \PiS} \inf_{P  \in \cP} 
 R(\pi,P)$. This additionally shows that a stationary uniformly optimal policy exists (again via the results for sa-rectangular robust MDPs from \Cref{sec: lec 2 sarec}).
\end{enumerate}
\end{proof}
\paragraph{The case of Model 2.} We now focus on Model 2. We focus on the properties of the policy evaluation problem
\[\inf_{\kappa \in \cPS_{\sf H} }R(\pi,\kappa)\]
and robust MDP problem
\begin{equation}\label{eq: lec 3 rmdp cPS}
    \sup_{\pi \in \PiH} \inf_{\kappa \in \cPS_{\sf H} }R(\pi,\kappa).
\end{equation}
The analogous statements for robust MDPs with time-varying transition probabilities, as in Model 2, are given next.
\begin{theorem}[Robust MDPs with time-varying transitions, Model 2]\label{th: lec 3 Model 2}
    Let $\cP \subset \DSSA$ be compact. The following statements hold.
\begin{enumerate}
    \item {\em Policy evaluation.} Let $\pi \in \PiS$. Then the worst-case return of $\pi$ against time-varying transitions in $\cPS_{\sf H}$ is attained at a stationary transition kernel in the s-rectangular extension of $\cP$:
    \[\inf_{\kappa \in \cPS_{\sf H}} R(\pi,\kappa) = \min_{P \in \recs(\cP)} R(\pi,P).\]
\item {\em Optimal policies.} An optimal policy for the optimization problem \eqref{eq: lec 3 rmdp cPS} can always be chosen stationary:
    \begin{align*}
        \sup_{\pi \in \PiH} \inf_{\kappa \in \cPS_{\sf H}} R(\pi,\kappa) & = \max_{\pi \in \PiS} \inf_{\kappa \in \cPS_{\sf H}} R(\pi,\kappa)
    \end{align*}
    In fact, a stationary uniformly optimal policy always exists for \eqref{eq: lec 3 rmdp cPS}, and it can be chosen as a best policy among stationary policies for the robust MDP problem over the stationary s-rectangular extension of $\cP$, i.e., 
    \[\sup_{\pi \in \PiH} \inf_{\kappa \in \cPS_{\sf H}} R(\pi,\kappa) = \max_{\pi \in \PiS} \min_{P \in \recs(\cP)} R(\pi,P).\]
    If additionally $\recs(\cP)$ is convex, we additionally get that 
    \[\max_{\pi \in \PiS} \min_{P \in \recs(\cP)} R(\pi,P) = \sup_{\pi \in \PiH} \min_{P \in \recs(\cP)} R(\pi,P)\]
    and we can conclude that a stationary uniformly optimal policy for \eqref{eq: lec 3 rmdp cPS} can be chosen as a stationary uniformly optimal policy for the robust MDP problem with stationary transitions over $\recs(\cP)$.
\end{enumerate}
\end{theorem}
\begin{proof}
    \begin{enumerate}
        \item \Cref{lem: lec 3 Model 12 eq extensions} shows that $\inf_{\kappa \in \cPS_{\sf H}} R(\pi,\kappa) = \inf_{\kappa \in \recs(\cP)_{\sf H}^{\sf S}} R(\pi,\kappa)$ and we can interpret the left-hand side as an adversarial MDP (see \Cref{sec: lec 2 adv}) with a finite set of states, compact action sets, continuous rewards and transitions. Therefore $\inf_{\kappa \in \recs(\cP)_{\sf H}^{\sf S}} R(\pi,\kappa) = \inf_{P \in \recs(\cP)} R(\pi,P)$.
        \item The proof for the existence of stationary optimal policies is direct. Define $T$ the robust Bellman operator as in \eqref{eq: lec 2 T}. Since $T$ is defined componentwise and each $s$-component only depends on the $s$-th marginal of $\cP$, we get that the set $\cP$ and the set $\recs(\cP)$ induce the same robust Bellman operator. Let $v\opt$ be its fixed point, and $\pi\opt \in \PiS$ such that $Tv\opt = T^{\pi\opt}v\opt$. For any $\kappa \in \cPS_{\sf H}, t\in \N$ and finite history $h_t \in (\cS \times \cA)^{t} \times \cS$, we get that $\E^{\pi\opt,\kappa}\left[ r(s_t,a_t,s_{t+1}) + \gamma v\opt(s_{t+1}) \; | \; h_t \right] \geq Tv\opt (s_t) = v\opt(s_t)$, so by induction and taking the limit (similarly to the proof of \Cref{th: lec 0 mdp main}), we get that $R(\pi\opt,\kappa) \geq p_0 \tr v\opt$, i.e. we have shown $\inf_{\kappa \in \cPS_{\sf H}} R(\pi\opt,\kappa) \geq p_0 \tr v\opt$. Conversely, for $\pi \in \PiS$, the same argument combined with $v\opt=Tv\opt = T^{\pi\opt}v\opt \geq T^{\pi}v\opt$ yields $\inf_{\kappa \in \cPS_{\sf H}} R(\pi,\kappa) \leq p_0 \tr v\opt$. Therefore $\pi\opt \in \PiS$ is an optimal policy for \eqref{eq: lec 3 rmdp cPS}, and it has been constructed independently of $p_0$ hence it is uniformly optimal for \eqref{eq: lec 3 rmdp cPS}. The fact that $\pi\opt$ is additionally optimal in $\max_{\pi \in \PiS} \min_{P \in \recs(\cP)} R(\pi,P)$ has been shown in \Cref{sec: lec 2 rmdp} and follows easily from $Tv\opt = T^{\pi\opt}v\opt$.

        Finally, we know that for s-rectangular convex compact set an optimal policy for \eqref{eq: lec 1 rmdp} can be chosen stationary, which concludes the proof. Note that convexity of $\cP$ implies convexity of $\recs(\cP)$.
    \end{enumerate}
\end{proof}
We emphasize that without the convexity assumption on $\cP$, optimal policies for \eqref{eq: lec 3 rmdp cPS} against time-varying transitions and optimal policies for \eqref{eq: lec 1 rmdp} against stationary transitions in $\recs(\cP)$ may differ (indeed the latter set of optimal policies may require history-dependent policies, whereas the former set of optimal policies always contains a stationary policy). The next example illustrates a decision problem with $\cPS_{\sf H}$.
\begin{example}
  Consider the robust MDP instance from \Cref{ex: lec 0 ex 1}, continued in \Cref{sec: lec 1 properties optim sol} (right after \Cref{prop: lec 1 nonstationary opt policies}). There $\cP = \{P^1,P^2\}$ is s-rectangular but not convex. We saw that in this instance $\sup_{\pi \in \PiH} \min_{P \in \cP} R(\pi,P)$ is in general solved by non-stationary policies, which can exploit the information revealed after period $t=0$ to correctly guess the transition probabilities ($P^1$ or $P^2$) chosen by the adversary. In contrast \Cref{th: lec 3 Model 2} tells us that an optimal policy can be chosen stationary in \eqref{eq: lec 3 rmdp cPS}. We analyze this specific instance in more detail. Rewards are only obtained at odd periods. For $\pi \in \PiH$ and $t \in \N$ let $x_t \in [0,1]$ be the probability that the decision maker chooses action $a_1$ at period $2t$ and for $\kappa \in \cPS_{\sf H}$ let $p_t \in \{0,1\}$ be the transition probabilities chosen by the adversary for period $2t$; in principle, $x_t$ and $p_t$ depend on $h_{2t}$ the history up to period $2t$. Then 
  \begin{align*}
      R(\pi,\kappa) & = \E^{\pi,\kappa}\left[ \sum_{t=0}^{+\infty} \gamma^t r(s_t,a_t,s_{t+1}) \right] \\
      & = \sum_{t=0}^{+\infty} \gamma^{2t+1} \E^{\pi,\kappa}\left[ r(s_{2t+1},a_{2t+1},s_{2t+2}) \right] \\
      & = \sum_{t=0}^{+\infty} \gamma^{2t+1} \E^{\pi,\kappa}\left[x_t p_t + (1-x_t)(1-p_t) \right].
  \end{align*}
 Minimizing with respect to $\kappa$ we get that $\min_{\kappa \in \cPS_{\sf H}} R(\pi,\kappa)$, where in principle the adversary can choose $p_t=0$ at some period $2t$ and $p_{t'}=1$ at some other period $2t'$, is maximized with $x_t=1/2$ for all $t \in \N$ so that an optimal policy can be chosen stationary in this example. 
\end{example}
Note that if optimal policies for \eqref{eq: lec 1 rmdp} against stationary transitions can be chosen stationary, then in fact these stationary optimal policies coincide with stationary optimal policies for \eqref{eq: lec 3 rmdp cPS} against time-varying transitions.
An intermediate conclusion regarding Model 1 and Model 2 is as follows:
\begin{itemize}
    \item For both Model 1 and Model 2, all uncertainty sets are equivalent to their rectangular counterparts, and optimal policies can be chosen stationary (there is no value in learning with history-dependent policies since the transition probabilities may be changing with time), which then yields stationary transition probabilities as worst-cases;
    \item For Model 1 we additionally get that optimal policies for time-varying transitions coincide with optimal policies for stationary transitions over their sa-rectangular extension (as defined in \eqref{eq: lec 1 rmdp}), whereas for Model 2 the corresponding statement requires a convexity assumption, a contrast that we already encountered in \Cref{chap: lec 2 rectangular models}.
\end{itemize}
Perhaps surprisingly, this shows that with time-varying transition probabilities as Model 1 or Model 2, we actually obtain more favorable properties than for the original robust MDP problem with stationary transitions as in~\eqref{eq: lec 1 rmdp}! 
\subsection{The case of Model 3 and connection with POMDPs}\label{sec: lec 3 model 3}
We are left with considering time-varying transitions as in Model 3. The conclusions differ depending on the rectangularity properties of $\cP$:
\begin{itemize}
    \item When $\cP$ is sa-rectangular compact, from \Cref{lem: lec 3 Model 12 eq extensions} we know that $\cPG_{\sf H} = \cPS_{\sf H} = \cPSA_{\sf H}$ hence the conclusions of \Cref{th: lec 3 Model 1} apply verbatim (we just need to replace $\cPSA_{\sf H}$ by $\cPG_{\sf H}$ in the statements).
    \item When $\cP$ is s-rectangular compact, from \Cref{lem: lec 3 Model 12 eq extensions} we know that $\cPG_{\sf H} = \cPS_{\sf H}$ so that the conclusions from \Cref{th: lec 3 Model 2} apply verbatim (we just need to replace $\cPS_{\sf H}$ by $\cPG_{\sf H}$ in the statements).
\end{itemize}
Interesting properties occur when $\cP$ is not s-rectangular.
We start with policy evaluation. We first show that for non-rectangular uncertainty sets the worst-case return of a stationary policy may be achieved at a non-stationary adversary, i.e. we may have $\displaystyle \inf_{\kappa \in \cPG_{\sf H} }R(\pi,\kappa) < \min_{P \in \cP }R(\pi,P)$.
\begin{example}\label{ex: lec 3 tv kernel}
    Consider an MDP instance with two states $s_0,s_1$ and a single action $a_1$, with $\gamma \in (0,1)$. Assume that we start in $s_0$ with probability $1$. The instantaneous rewards are $r(s,a,s') = 1$ if $s= s'$ and $0$ otherwise, i.e., we obtain a reward of $1$ if we don't change states. Consider the uncertainty set $\cP = \conv(\{P^1,P^2\})$ where $P^1$ sends both $s_0,s_1$ to $s_0$ (we can identify $P^1$ with the transition matrix $\begin{pmatrix}
        1 & 0 \\
        1 & 0 \\
    \end{pmatrix}$)
   while $P^2$ sends both $s_0,s_1$ to $s_1$ (we can identify $P^2$ with the transition matrix $\begin{pmatrix}
        0 & 1 \\
        0 & 1 \\
    \end{pmatrix}$). Note that $\cP$ is not s-rectangular but it is convex. Then the adversary can choose the following feasible policy in $\cP_{\sf H}^{G}$: choose $P^2$ at even periods and $P^1$ at odd periods. This way we ensure that the sequence of visited states alternates between $s_0$ and $s_1$, so that the return is $0$. In contrast if the adversary always chooses the same transition probabilities in $\cP$ then the return is positive. Hence in this example the worst-case return of the policy $\{s_0,s_1\} \to \{a_1\}$ over $\kappa \in \cP_{\sf H}^{G}$ is $0$ and attained at a non-stationary adversary policy. 
\end{example}
From the adversary's standpoint, the optimization problem resembles a {\em Partially Observable MDP} (POMDP), which is a class of sequential decision-making problems where the decision maker does not observe the current state $s_t$ fully but only a ``masked" version (the so-called {\em ``observation"} $o_t$). The decision maker then forms a certain belief $b_t \in \DS$ of what is the current state at period $t$ and can act accordingly. Here we are in the case where the information is delayed, i.e., the state $s_{t}$ is revealed to the adversary only right before period $t+1$ (after the decision maker has chosen $a_t$ and after the adversary has chosen $P^t$). The next proposition shows how to apply classical results on POMDPs to our setup. We denote $C(\DS,\R)$ the set of continuous maps from $\DS \to \R$. For $\pi \in \PiS$ and $P \in \cP$ we write $B^{\pi,P}: C(\DS,\R) \to C(\DS,\R)$ the operator defined as, for $w \in C(\DS,\R)$,
\[ \left(B^{\pi,P} w \right)(b) = \sum_{s \in \cS} b(s) \sum_{a \in \cA} \pi(s,a) \left( P(s,a) \tr r(s,a) + \gamma w(P(s,a))\right), \forall \; b \in \DS.\]
In the equation above we interpret $b \in \DS$ as a certain belief of what will be the current state at the period of decision. Hence the right-hand side is the sum of the reward obtained in one period with a continuation value of $w \in C(\DS,\R)$. Hence the operator $B^{\pi,P}$ plays the same role as the operator $T^{\pi,P}$ defined in \eqref{eq: lec 0 T pi P} in \Cref{chap: lec 0 mdp}, and $\DS$ plays the role of $\cS$ (since the adversary can only form beliefs about the current state). Note that we use $P(s,a)$ as the argument for the next belief because $P(s,a)$ is the distribution of $s_{t+1}$ given a current state-action pair $(s,a)$, and the adversary, when choosing $P^{t+1}$ for period $t+1$, has observed $s_t$ and $a_t$ (following the definition of Model 3).
\begin{proposition} 
Let $\cP \subset \DSSA$ be a compact uncertainty set. Let $\pi \in \PiS$.
Let $B^{\pi}: C(\DS,\R) \to C(\DS,\R)$ be the operator defined as, for $w \in C(\DS,\R)$,
\[ \left(B^{\pi} w \right)(b) = \min_{P \in \cP} \left(B^{\pi,P} w \right)(b), \forall \; b \in \DS.\]
The following statements hold:
\begin{enumerate}
\item The operator $B^{\pi}$ is a monotone contraction for $\| \cdot \|_{\infty}$ on $C(\DS,\R)$, with unique fixed point $w^{\pi} \in C(\DS,\R)$.
\item We have 
\[\inf_{\kappa \in \cPG_{\sf H}} R(\pi,\kappa) = w^{\pi}(p_0)\]
and an optimal policy for the adversary can be chosen stationary in the belief $b_t \coloneqq P^{t-1}(s_{t-1},a_{t-1}) \in \DS$ at each $t \geq 1$ and $b_0 = p_0$.
\end{enumerate} 
\end{proposition}
\begin{proof}
\begin{enumerate}
\item The proof that $B^{\pi}$ is a monotone contraction follows the same lines as for these properties for $T^{\pi}$. The fact that $B^{\pi}w$ is always a continuous function (given $w \in C(\DS,\R)$) follows from the continuity of $(b,P) \mapsto \left(B^{\pi,P}w\right)(b)$. Now  $C(\DS,\R)$ is complete for $\| \cdot \|_{\infty}$ hence Banach's fixed point theorem gives the existence and uniqueness of the fixed point $w^{\pi} \in C(\DS,\R)$.
\item From the fixed point equation on $w^{\pi} $ and the definition of $B^{\pi}$ we get that $w^{\pi}(b_t) \leq (B^{\pi,P^t}w^{\pi})(b_t)$ for any admissible $P^t \in \cP$. That is, for any $g_t \in \cG_{t-1}$ and $\kappa \in \cPG_{\sf H}$ we get 
\[w^{\pi}(b_t) \leq \E^{\pi,\kappa}\left[ r(s_t,a_t,s_{t+1}) + \gamma w^{\pi}(b_{t+1}) \; | \; g_t \right].\] 
By induction we get
\[w^{\pi}(p_0) \leq \E^{\pi,\kappa}\left[ \sum_{t'=0}^{t}\gamma^{t'} r(s_{t'},a_{t'},s_{t'+1}) + \gamma^{t+1} w^{\pi}(b_{t+1})\right].\]
Since $w^{\pi}$ is a bounded function, taking the limit as $t \to + \infty$ yields that for any $\kappa \in \cPG_{\sf H}$ we have $w^{\pi}(p_0) \leq R(\pi,\kappa)$, hence $w^{\pi}(p_0) \leq \inf_{\kappa \in \cPG_{\sf H}} R(\pi,\kappa)$. 

To show the converse inequality we just choose $P\opt(b) \in \arg \min_{P \in \cP} (B^{\pi,P}w^{\pi})(b)$ for any belief $b \in \DS$ and construct $\kappa\opt \in \cPG_{\sf H}$ accordingly such that for any $g_t = (s_0,a_0,...,s_{t-1},a_{t-1})$ we have, $(\kappa\opt_t(g_t,s,a))_{s,a} = P^{\star}(b_t) \in \cP$, which concludes the proof.
\end{enumerate}
\end{proof}
The next example continues the analysis of \Cref{ex: lec 3 tv kernel} in light of the previous proposition.
\begin{continuance}{ex: lec 3 tv kernel}
Consider the setting of \Cref{ex: lec 3 tv kernel}. Since we start at $s_0$ with probability $1$ the belief for period $0$ is $b_0 = (1,0)$, and for $b=(b(s_0),b(s_1)) \in \DS$ the belief of where the current state is, we have
\[(Bw)(b) = \min_{\lambda \in [0,1]} b(s_0)\lambda + b(s_1) (1-\lambda) + \gamma w((\lambda,1-\lambda)).\]
At the fixed point we have that for $b=(1,0)$ the optimal $\lambda$ above is $\lambda=0$ and the next belief is $(0,1)$; for $b=(0,1)$ we are in the reverse situation and the optimal $\lambda$ is $1$ so that the next belief is $(1,0)$, that is, $w(b) = \min \{b(s_0),b(s_1)\}$. This is exactly the alternating adversary policy that we highlighted in \Cref{ex: lec 3 tv kernel}, where $P^2$ is chosen at even periods and $P^1$ is chosen at odd periods.
\end{continuance}
Again, for non-rectangular robust MDPs we are left with a simpler picture when the transitions are time-varying (since we have a connection with the well-studied field of POMDPs) than for stationary transitions. For the sake of comparison, note that past works by Littman and coauthors have shown that finding good stationary policies in POMDPs is hard (see reference section below).
We conclude with an example showing that an optimal policy for the decision maker need not be stationary in this setup.
\begin{example}\label{ex: lec 4 ex 1}
There are four states $s_0,s_1,G,B$ where $G,B$ are absorbing, and two actions $a_1,a_2$. The uncertainty set $\cP$ contains two transition kernels $P^1,P^2$ defined as follows. Under $P^1$, action $a_1$ yields the deterministic transitions $s_0 \to G$ and $s_1 \to B$, while action $a_2$ yields the transitions $s_0 \to B, s_1 \to B$. Under $P^2$, action $a_1$ yields the deterministic transitions $s_0 \to B$ and $s_1 \to s_0$, while action $a_2$ yields the transitions $s_0 \to G, s_1 \to s_0$. The rewards are $0$ everywhere except $r(s_0,a,G)=1$ for $a \in \{a_1,a_2\}$. Thus the decision maker wants to perform the transition $s_0 \to G$, while the adversary tries to guess their action and react accordingly (i.e. if the decision maker chooses $a_1$ in $s_0$ the adversary wants to choose $P^2$). The initial distribution is uniform over $\{s_0,s_1\}$. Take $\gamma=1/2$ to simplify our calculations. Let $\pi \in \PiH$.
Assume that $s_0$ is reached at period $t=1$. The decision maker faces the typical ``$ \max_x \min_p xp + (1-x)(1-p)$" situation we have encountered in many previous examples (e.g. right after \Cref{prop: lec 1 nonstationary opt policies}), with $x \in [0,1]$ the probability that the decision maker plays action $a_1$ and $p \in \{0,1\}$ modeling the choice of the adversary ($p=1$ for choosing $P^1$, $p=0$ otherwise). Hence the best worst-case continuation value from $s_0$ at period $t=1$ is $1/2$, attained by the uniform randomization $x\opt=1/2$.
Now let $z \in [0,1]$ the probability that the decision maker chooses $a_1$ at period $t=0$ when starting in $s_0$. Then $R(\pi,P^1)=z \cdot \frac{1}{2}$ (with the factor $1/2$ coming from the initial distribution), and $R(\pi,P^2) =  \frac{1-z}{2}+ \frac{1}{2} \cdot \frac{1}{2} \cdot  \frac{1}{2} = \frac{1-z}{2} + \frac{1}{8}$. Therefore, \[\sup_{\pi \in \PiH} \inf_{\kappa \in \cPG_{\sf H}} R(\pi,\kappa) = \max_{z \in [0,1]} \min \left\{ \frac{z}{2}, \frac{1-z}{2} + \frac{1}{8}\right\} = \frac{5}{16},\]
which is attained at $z\opt = 5/8$, that is, at the Markovian policy choosing $a_1$ with probability $5/8$ at period $t=0$ then $a_1$ with probability $1/2$ for $t \geq 1$.
In contrast, a stationary policy $\pi \in \PiS$ choosing always $a_1$ with probability $x \in [0,1]$ obtains a worst-case return of 
\[\max_{\pi \in \PiS} \inf_{\kappa \in \cPG_{\sf H}} R(\pi,\kappa) = \max_{z \in [0,1]} \min \left\{\frac{z}{2}, \frac{1-z}{2} + \frac{1}{4} \min \{z,1-z\}\right\}.\]
One can check that the maximum of this quantity over $z \in [0,1]$ is $\frac{3}{10}$ attained at $z=3/5$. Since $3/10 < 5/16$ we get that stationary policies are strictly suboptimal in this example.
\end{example}
The previous example has a particular structure, in that the transitions are deterministic: for all $P \in \cP, P(s,a,s') \in \{0,1\}, \forall \; (s,a,s') \in \cS \times \cA \times \cS$. In this case the adversary can correctly guess the current state from period $1$ onward; only the state at period $t=0$, sampled from $p_0 \in \DS$, is uncertain. Thus from period $1$ onward we are in a similar setup as Model 2 where the adversary can observe the current state and the decision maker can play stationary from period $t =1$ onward. This yields a value function $\tilde{v}\opt \in \RS$ for the rewards obtained starting from period $1$, and then we can play at period $0$ by solving
\[\max_{\pi \in \PiS} \min_{P \in \cP} \sum_{s \in \cS} p_0(s) * \sum_{a \in \cA} \pi(s,a) P(s,a) \tr (r(s,a) + \gamma \tilde{v}\opt)\]
which maximizes the sum of the instantaneous reward obtained at period $t=0$ and then the expected continuation values starting from $t=1$, where we know we can play stationary and obtain $\tilde{v}\opt \in \R^{\cS}$ as continuation values.
\section{Stochastic games}\label{sec: lec 3 sgs}
Stochastic games (SGs) were introduced by Shapley in 1953. We focus here on two-player zero-sum stochastic games, where players have opposite interests and choose some actions at each period, which influence the dynamics of the observed states, and players optimize for their gains.
\paragraph{Stochastic game model.}
Let us introduce this model more formally. In SGs it is common to distinguish between an {\em arena}, which is the setup in which players make their decisions, and the objective that the players optimize. A stochastic game consists of an arena and an objective. An arena of a two-player zero-sum stochastic game $\Gamma$ is determined by a tuple $\left(\Omega,\cI,\cJ,g,q\right)$ with $\Omega$ a finite set of states with cardinality $d$, $\cI$ the set of actions of the first player, $\cJ$ the set of actions of the second player, $g:  \cI \times \cJ \times \Omega \rightarrow \R$ is the instantaneous reward associated with each state $\omega \in \Omega$ and pair of actions $\left(i,j\right) \in \cI \times \cJ$, and $q: \cI \times \cJ \times \Omega \rightarrow \Delta(\Omega)$ gives the transition probabilities. Sometimes this kind of SG is called {\em concurrent} because players choose their actions simultaneously. For the sake of conciseness, we will refer to {\em two-player zero-sum stochastic games} simply as {\em stochastic games} (SGs). We will assume that the sets of actions $\cI,\cJ$ are compact and that the payoff and transition functions are continuous.

The game proceeds as follows. There is an initial state $\omega_{0}$ known to both players at period $t=0$. Both players simultaneously choose their actions $(x_{0},y_{0})$ in $\cI \times \cJ$, the payoff is $g(x_{0},y_{0}, \omega_{0})$, and the system transitions to $\omega_{1} \in \Omega$ with probability $q(\omega_{1}|x_{0},y_{0}, \omega_{0})$ at the next period. This sequence of events is repeated infinitely, and we note $H_{t} = \left(\Omega \times \cI \times \cJ\right)^{t} \times \Omega$ the set of histories of length $t \in \N$, $H = \cup_{t \in \N} H_{t}$ the set of all finite histories, and $H_{\infty} = \left(\Omega \times \cI \times \cJ\right)^{\N}$ the set of infinite histories. A strategy for player 1 (resp. for player 2) is a map $\sigma :H \rightarrow \Delta(\cI)$ (resp. $\tau: H \rightarrow \Delta(\cJ))$, and we write $\PiH^{i}$ for the set of all history-dependent, randomized strategies of player $i$, i.e., the set $\Delta(\cI)^{H}$ for the first player and the set $\Delta(\cJ)^{H}$ for the second player. We also write $\PiS^{i}$ for the set of stationary strategies of player $i$, which we can identify with the sets $\Delta(\cI)^{\Omega}$ and $\Delta(\cJ)^{\Omega}$, and $\PiSD^{i}$ for the set of stationary deterministic strategies, which we can identify with the sets $\cI^{\Omega}$ and $\cJ^{\Omega}$. Note that {\em history-dependent, stationary, and deterministic} strategies are usually called {\em behavioral, positional, and pure} in the SG literature. 

A pair of strategies $\left(\sigma,\tau\right) \in \PiH^{1} \times \PiH^{2}$ induces a distribution over the set of infinite histories. There are several ways to aggregate the random infinite stream of payoffs $(g(x_{t},y_{t},\omega_{t}))_{t \geq 0}$ induced by two strategies $\left(\sigma,\tau\right) \in \PiH^{1} \times \PiH^{2}$. The mean-payoff starting from state $\omega \in \Omega$ is 
\begin{equation}\label{eq: lec 3 mean payoff}
    \gamma(\sigma,\tau,\omega) \coloneqq \E^{\sigma,\tau} \left[ \limsup_{T \to + \infty} \frac{1}{T+1}\sum_{t=0}^{T} g(x_{t},y_{t},\omega_{t}) \;\middle|\;\omega_0 = \omega \right]
\end{equation}
and the discounted payoff is, for $\lambda \in (0,1)$,
\begin{equation}
    \gamma_{\lambda}(\sigma,\tau,\omega) \coloneqq \lambda\E^{\sigma,\tau}\left[\sum_{t=0}^{+\infty} (1-\lambda)^t g(x_{t},y_{t},\omega_{t}) \;\middle|\;\omega_0 = \omega \right].
\end{equation}

Of course robust MDPs greatly resemble this model despite different terminologies and notation (e.g., the terms ``policies" and ``history-dependent/stationary" correspond to ``strategies", ``behavioral/positional", the Shapley operator in SGs corresponds to the robust Bellman operator $T$, etc.), and this has been observed repeatedly in the past, as early as in the first paper on robust MDPs by Satia and Lave in 1973. 
We describe this connection for s-rectangular and sa-rectangular robust MDPs below.
\subsection{Perfect-information SGs and sa-rectangular robust MDPs} Consider an sa-rectangular compact uncertainty set $\cP = \times_{s,a} \cP_{sa}$, with state set $\cS$, action set $\cA$, rewards $r \in \RSAS$. Then we can model the corresponding robust MDP problem \eqref{eq: lec 1 rmdp} as a {\em perfect-information} stochastic game, which is a particular class of SGs where each state is controlled by exactly one player. In fact we will show that we can reformulate \eqref{eq: lec 1 rmdp} as a {\em turn-based} SG, where players play in turn: the first player chooses actions at even periods and the second player chooses actions at odd periods.

\begin{figure}[htb]
\begin{center}
    \includegraphics[width=0.7\linewidth]{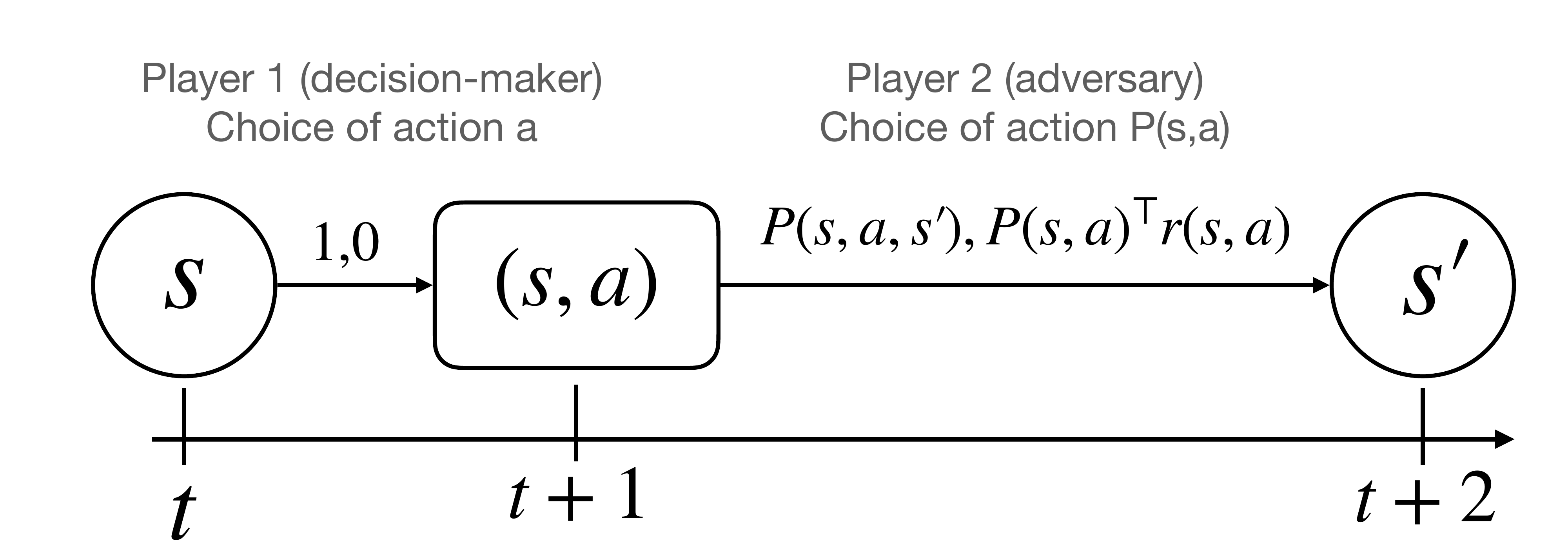}
    \end{center}
\caption{An sa-rectangular robust MDP as a perfect-information (turn-based) game. The decision-maker controls the states $s \in \cS$ and chooses an action $a \in \cA$. The adversary controls the states $(s,a) \in \cS \times \cA$ and chooses transition probabilities $P(s,a) \in \cP_{sa}$. The pairs above the arcs represent the transition probabilities and instantaneous rewards. }\label{fig: lec 3 sa-rec-as-sgs}
\end{figure}

Consider the following game instance. The first player is the decision maker and the second player is the adversary. The set of states is augmented to $\Omega = \cS \cup \cS \times \cA$. States $\omega \in \cS$ are visited at even periods (where the decision maker plays) and states $\omega \in \cS \times \cA$ are visited at odd periods (where the adversary plays). In each state $s \in \cS$ the available actions for the decision maker are $\cI = \cA$. If action $a \in \cA$ is chosen the dynamics transitions deterministically to the next state $\omega' = (s,a)$ and the reward for this transition is $0$. At a state $\omega'=(s,a)$ the set of actions for the adversary is $\cP_{sa}$ (note that we can allow for state-dependent action sets in SGs, it only complicates a bit the notation). If an ``action" $P(s,a) \in \cP_{sa}$ is chosen, the probability to transition to the next state $s'$ is $P(s,a,s')$ and the instantaneous reward is $P(s,a) \tr r(s,a)$ (which only depends on the current state $(s,a)$ and action $P(s,a)$ but not on the next state $s'$). The dynamics for this SG is represented in \Cref{fig: lec 3 sa-rec-as-sgs}.
To ensure that the discounted payoff for the SG is the discounted return for the robust MDP, we need to choose a discount factor $\gamma$ such that $1-\lambda = \sqrt{\gamma}$ (and divide the reward at the intermediate states by $(1-\sqrt{\gamma}) \sqrt{\gamma}$; this is not represented in \Cref{fig: lec 3 sa-rec-as-sgs} for the sake of readability) in this turn-based SG, because rewards are obtained only at odd periods (instead of every period in the original robust MDP instance).

There are two important caveats to the above constructions, pertaining to the sets of policies of each player:
\begin{itemize}
    \item In SGs each player should be allowed to play mixed strategies, i.e., to randomize over their sets of actions. For the decision maker this is akin to the situation in MDPs and robust MDPs where they can choose a randomized decision rule in $\DA$ at each period, instead of a fixed $a \in \cA$. For the adversary, in principle at each odd period we need to allow for a choice of randomized decision rule in $\Delta(\cP_{sa})$ (interpreted as the set of probability distributions over $\cP_{sa}$) instead of a choice of deterministic action in $\cP_{sa}$.
    \item In SGs each player should be allowed to play history-dependent strategies, whereas in classical robust MDPs as in \eqref{eq: lec 1 rmdp} the decision maker can choose history-dependent policies while the adversary must choose the transition probabilities once and for all. However even the definition of what constitutes a ``history" is different in SGs and in robust MDPs. Robust MDPs build upon nominal MDPs and define the set of all finite histories from the decision maker's perspective:
    \begin{align*}
        \cH &  = \cup_{t \in \N} \cH_t, \\
        \cH_t & = (\cS \times \cA)^{t} \times \cS
    \end{align*}
    where $\cH_t$ is the set of finite histories $(s_0,a_0,...,s_{t-1},a_{t-1},s_{t})$ defined in \Cref{sec: lec 0 model}: a history includes the state and the action sampled by the decision maker (first player) at each period. However in SGs a history must {\em also} include the past decisions of the second player! The first item above explains that the second player should choose a randomized decision rule in $\Delta(\cP_{sa})$ at each period, and sample a realization $P(s,a) \in \cP_{sa}$ (according to this distribution). Hence in principle a finite history for the first player in the SG representation of this sa-rectangular robust MDP is some sequence
    \[(s_0,a_0,(s_0,a_0),P^0(s_0,a_0),...,s_{t-1},a_{t-1},(s_{t-1},a_{t-1}), P^{t-1}(s_{t-1},a_{t-1}), s_{t})\]
    for some $t \in \N$ and $P^0(s_0,a_0) \in \cP_{s_0 a_0},...,P^{t-1}(s_{t-1},a_{t-1}) \in \cP_{s_{t-1} a_{t-1}}$ (note that we include the ``states" $(s_0,a_0),...,(s_{t-1},a_{t-1})$ in this history since these are states visited in the SG, where the adversary makes their decisions) and for the second player it is some sequence
    \[(s_0,a_0,(s_0,a_0),P^0(s_0,a_0),...,s_{t-1},a_{t-1},(s_{t-1},a_{t-1}), P^{t-1}(s_{t-1},a_{t-1}), s_{t},a_{t},(s_{t},a_{t})).\]
    This departs from the set of histories (see \Cref{sec: lec 0 model}) and time-varying transitions (see \Cref{sec: lec 3 time var P}) in robust MDPs, since in SGs both players have much more information (they can observe the realizations of the actions of the second player along each trajectory).
\end{itemize}
Fortunately, there are some special cases where the two caveats above do not matter for the interpretation of robust MDPs with sa-rectangular uncertainty sets as turn-based SGs. These are the cases of the objective functions for which we have that a pair of optimal strategies in the SG can be chosen stationary (only dependent on the current state) and deterministic. For instance the literature on SGs has shown that this is true in:
\begin{enumerate}
    \item Case 1: perfect-information SGs with the discounted objective when both players have compact action sets. For the sa-rectangular robust MDP it corresponds to compactness of $\cA$ and compactness of $\cP_{sa}$, which are usual assumptions.
    \item Case 2: perfect-information SGs with the mean-payoff objective when both players have finite action sets. For the sa-rectangular robust MDP it corresponds to finiteness of $\cA$ and finiteness of $\cP_{sa}$, and the average return $\Ravg(\pi,P)$ is defined as (compare with \eqref{eq: lec 3 mean payoff})
\[  \Ravg(\pi,P) =  \E^{\pi,P} \left[ \limsup_{T \rightarrow + \infty}\frac{1}{T+1} \sum_{t=0}^{T} r(s_{t},a_{t},s_{t+1})\right].
\]
\end{enumerate}
Stationarity of the first player's strategy means that their decision rules only depend on the current state $s_t$, and stationarity of the second player's strategy means that their decision rules depend on the current pair $(s_t,a_t)$, and not on the past histories. Hence for Case 1 and Case 2 above it is irrelevant that: (A) turn-based SGs allow for history-dependent adversaries, since an optimal strategy for the adversary in the turn-based SG formulation can be chosen stationary anyway, (B) the sets of histories for SGs and for robust MDPs differ, since optimal policies are stationary and only depend on the last states, where histories of MDPs and histories of SGs coincide. In particular, the adversary and the decision maker do not need to see the realizations of the past chosen transition probabilities to play optimally. Note Case 1 above corresponds to the setting we have studied in \Cref{sec: lec 2 sarec}, and to some extent the results we obtained there are relatively direct corollaries of existing results for stochastic games. Case 2 above yields the following proposition.
\begin{proposition}\label{prop: lec 3 mean payoff sa-rec}
Consider an sa-rectangular robust MDP with a finite set of transition probabilities. Then an optimal policy for the average return objective can be chosen stationary and deterministic:
\[\sup_{\pi \in \PiH} \inf_{P \in \cP} \Ravg(\pi,P) = \max_{\pi \in \PiSD} \min_{P \in \cP} \Ravg(\pi,P).\]
\end{proposition}
It is worth noting that even though Case 2 above requires finiteness of $\cP_{sa}$, it is not hard to extend the results to the case of compact $\cP_{sa}$
, although in this case worst-case kernels need not exist anymore (see references at the end of this section).
\subsection{Concurrent SGs and s-rectangular robust MDPs} Consider an s-rectangular convex compact uncertainty set $\cP$, with state set $\cS$, action set $\cA$, rewards $r \in \RSAS$. We can reformulate the corresponding robust MDP problem as a {\em concurrent} stochastic game, where both players choose their actions at the same time, as follows. The first player is the decision maker and the second player is the adversary. The set of states for the SG is now $\Omega=\cS$. In each state $s \in \cS$, the decision maker samples an action $a \in \cA$ from some (possibly randomized) decision rule in $\DA$, while the adversary chooses an action $P(s) \in \cP_s$. Since $\cP_s$ is convex we can interpret choosing an element $P(s) \in \cP_s$ as choosing the realization of a randomized decision over $\Delta(\cP_s)$ from some (possibly randomized) decision rule in $\Delta(\cP_s)$.

The instantaneous payoff is $P(s,a) \tr r(s,a)$ and the transition probability to the next state is $P(s,a,s')$. We can now use a discount factor of $\gamma$ (the same as for the robust MDP instance).
Note that the second caveat still applies here: a finite history for this SG consists of some sequence 
\[(s_0,(a_0,P^0(s_0)),s_1,(a_1,P^1(s_1)),...,s_{t-1},(a_{t-1},P^{t-1}(s_{t-1})),s_t)\]
for some $t \in \N$ and $P^0,...,P^{t-1} \in \cP$. However, for concurrent SGs with discounted return and a finite set of actions for the first player and a convex compact set of actions for the second player, it is known that stationary optimal strategies exist for both players. In principle this yields a mixed strategy $\mu \in \times_{s \in \cS} \Delta(\cP_s)$ for the adversary, but since $\mu$ only plays a role via its expectation (the robust Bellman operator is linear in $P$), convexity of $\cP$ enables us to replace the mixed strategy $\mu$ by its expected kernel $\bar{P} = \E_{P \sim \mu} \left[ P \right]$, interpreted as a stationary pure strategy for the adversary. Hence again the distinction between the histories for SGs and histories for robust MDPs does not matter in this setting, and again the adversary and the decision maker do not need to see the realizations of the past chosen transition probabilities to play optimally.

It is worth emphasizing that the rectangularity assumption, which plays a central role for robust MDPs, corresponds in SGs to the fact that each player can choose their actions independently at each state.
\paragraph{Applications.} Stochastic games date back to the 1950s, while the vast majority of the robust MDP literature has been published after the 2000s. Hence, the connection between robust MDPs and SGs is particularly useful when it enables us to claim new results for robust MDPs as more or less ``direct corollaries" of existing results for SGs (using the reformulations of robust MDPs as SGs with particular properties), or at least to use for robust MDPs some techniques developed for SGs. We did so above in \Cref{prop: lec 3 mean payoff sa-rec}. Some setups where results for SGs are more advanced than for robust MDPs include the mean-payoff objective (defined above for SGs), which corresponds in robust MDPs to the average return $\Ravg(\pi,P)$ defined above.
This kind of return has been studied by the SG community since the late 1950s, while this question has only recently been studied for robust MDPs. The same is true, for instance, for the Blackwell optimality criterion (close to the {\em uniform value} in SGs), SGs with incomplete information, and SGs that are not zero-sum. 

We conclude by emphasizing that sometimes it is not direct to ``translate" results for SGs to results for robust MDPs. For instance, if one is interested in the average return $\Ravg(\pi,P)$ for s-rectangular robust MDPs, naturally one would want to leverage existing results for SGs with mean-payoff objective. Perhaps the strongest result in this direction is Mertens-Neyman theorem, which states that concurrent SGs with finite sets of states and actions have a value {\em among history-dependent strategies}, that is, we have for each state $\omega \in \Omega$, that
\[\sup_{\sigma \in \PiH^{1}} \inf_{\tau \in \PiH^{2}} \gamma(\sigma,\tau, \omega) =  \inf_{\tau \in \PiH^{2}} \sup_{\sigma \in \PiH^{1}} \gamma(\sigma,\tau,\omega).\]
Hence the $\sup$ and $\inf$ need not be attained and in general it is required that the players can choose history-dependent (non-stationary, non-Markovian) strategies. Hence this result does not translate directly to studying $\sup_{\pi \in \PiH} \inf_{P \in \cP} \Ravg(\pi,P)$ because now optimal strategies in SGs make full use of the histories (in particular, using the other players' past decisions) whereas this is not allowed for the adversary in robust MDPs.
\section{Distributionally robust MDPs}\label{sec: lec 3 drmdps}
Distributionally robust MDPs (DRMDPs) assume that the transition probabilities are random and follow some uncertain distribution $\nu \in \cV$, where $\cV$ is a set of distributions over $\DSSA$, the set of all possible transition probabilities. The set $\cV$ is sometimes called the {\em ambiguity set}, not to be confused with the uncertainty set $\cP$ (note that these two sets don't live in the same space: $\cP$ is a set of transition probabilities whereas $\cV$ is a set of probability measures over transition probabilities). There are two main assumptions in the DRMDPs literature:
\begin{itemize}
    \item {\em Product distributions.} First it is assumed that $\cV$ is a Cartesian product: $\cV = \times_{s \in \cS} \cV_s$, where $\cV_s$ is a set of distributions over $\DSA$. Therefore the measure $\nu$ for the random transitions $P \in \DSSA$ is a product measure over $s \in \cS$.
 Each $\nu_s \in \cV_s$ represents a possible distribution for the $s$-th component $P(s)$ of the transition probabilities. Typically we will assume that the set of expectations $\cP_s = \{ \E_{P(s) \sim \nu_s} [P(s)] \; | \; \nu_s \in \cV_s\}$ is convex and compact.

    \item {\em Resampled kernels.} Second, the transition probabilities are resampled independently at every iteration from some $\nu^t \in \cV$ (chosen before the dynamics start), and the realizations are not observed by the decision maker.
\end{itemize}
Using the terminology from the seminal paper by Xu and Mannor we can then consider two classes of time-varying random transition probabilities:
\begin{itemize}
    \item {\em Stationary model.} In this case we assume that the adversary must choose the same distribution $\nu \in \cV$ for the transition probabilities across all periods. This set of distributions is written $\cVS$ with
    \[\cVS \coloneqq \{ (\nu)_{t \in \N} \; | \; \nu \in \cV\}.\]
    \item {\em Non-stationary model.} 
    Otherwise we assume that the adversary can choose a new distribution $\nu^t \in \cV$ at each period $t \in \N$.
    This set of distributions is written $\cVNS$ with
    \[\cVNS \coloneqq \{ (\nu^t)_{t \in \N} \; | \; \nu^t \in \cV, \forall \; t \in \N\} = \times_{t \in \N} \cV.\]
    Note the similarity with the definition of Model 2 from \Cref{sec: lec 3 time var P}.
    We could define more generally history-dependent policies for the adversary but we'll see that all these models collapse to stationary adversaries anyway. Note that $\cVS \subset \cVNS$.
\end{itemize}
\begin{remark}[Connection with multi-model MDPs]
Another natural formulation for MDPs with random transitions is to consider that the transition probabilities are sampled once and for all. In this case we want to solve $\displaystyle  \sup_{\pi \in \PiH} \inf_{\nu \in \cV} \E_{P \sim \nu} \left[R(\pi,P)\right]$. The special case where $\cV$ is the singleton $\{\nu\}$ yields the optimization problem $\displaystyle \sup_{\pi \in \PiH} \E_{P \sim \nu} \left[R(\pi,P)\right]$ which is known as {\em multi-model MDPs} and is in general intractable (even when $\nu$ has a support of cardinality two), see \Cref{sec: lec 3 multi-objective MDPs} and the reference at the end of this chapter.
\end{remark}
Recall that we have defined the return $R(\pi,P)$ with a fixed kernel $P \in \DSSA$. Given a sequence $(P^t)_{t \in \N} \in (\DSSA)^{\N}$ we can extend the return function to $R(\pi,(P^t)_{t})$ since the pair $(\pi,(P^t)_t)$ still induces a distribution over the set of trajectories. Therefore for $\mu \in \cVNS$ we can define 
\[R(\pi,\mu) \coloneqq \E_{(P^t)_t \sim \mu} \left[R(\pi,(P^t)_{t})\right].\]
The objective in distributionally robust MDPs is to solve 
\begin{equation}\label{eq: lec 3 drmdps stat}
\sup_{\pi \in \PiH} \inf_{\mu \in \cVS} R(\pi,\mu)
\end{equation}
and the same problem with worst-cases against time-varying distributions:
\begin{equation}\label{eq: lec 3 drmdps tv}
\sup_{\pi \in \PiH} \inf_{\mu \in \cVNS} R(\pi,\mu)
\end{equation}
We will see that both \eqref{eq: lec 3 drmdps stat} and \eqref{eq: lec 3 drmdps tv} are equivalent to robust MDPs, in a sense made precise below. We first show that given some $\mu \in \cVNS$, the return $R(\pi,\mu)$ only depends on the expected transition probabilities at each period, and not on their entire distribution.
\begin{lemma}\label{lem: lec 3 return reformulation}
Let $\pi \in \PiH$ and $\mu = (\nu^t)_{t \in \N} \in \cVNS$. Let $\bar{P}^t = \E_{P^t \sim \nu^t} \left[ P^t\right]$.
Then
\[R(\pi,\mu) = R(\pi, (\bar{P}^t)_{t \in \N}).\]
\end{lemma}
\begin{proof}
Recall that $R(\pi,\mu)$ is a function of $\bP^{\pi,\mu}$, the distribution over infinite histories induced by the pair $(\pi,\mu)$.
    Induction on the period $t$ shows that $(\bar{P}^t)_{t \in \N}$ and $\mu$ induce the same distribution over all finite histories $h_t$. Hence, they induce the same distribution over all infinite histories, and $R(\pi,\mu) = R(\pi, (\bar{P}^t)_{t \in \N})$.
\end{proof}
In particular, when $\mu \in \cVS$ we have $\mu= (\nu)_{t \in \N}$ for some $\nu \in \cV$ and $R(\pi,\mu) = R(\pi, \bar{P})$ with $\bar{P} = \E_{P \sim \nu} \left[ P\right]$. From the above lemma we easily see that distributionally robust MDPs can model robust MDPs as a special case: for $\cP \subset \DSSA$ a given s-rectangular uncertainty set, let $\cV = \times_{s \in \cS} \{ \delta_{\{P(s)\}} \; | \; P \in \cP\}$ with $\delta_{\{P(s)\} }$ the Dirac distribution which puts mass $1$ at $P(s) \in \cP_s$. Then $\bar{\cP} = \cP$ and from \Cref{lem: lec 3 return reformulation} we have $\inf_{\mu \in \cVS} R(\pi,\mu) = \inf_{P \in \cP} R(\pi,P)$. A somewhat surprising property is that the reverse is also true: distributionally robust MDPs will reduce to robust MDPs (this is because the distribution $\mu \in \cVNS$ only plays a role via its mean). We show that this is indeed the case in the next theorem, which is the fundamental result for distributionally robust MDPs.
\begin{theorem}\label{th: lec 3 drmdps}
Let $\cV = \times_{s \in \cS} \cV_s$ be a set of product distributions over $\DSSA$. 
Let $\bar{\cP} = \{\E_{P \sim \nu}\left[ P \right] \; | \; \nu \in \cV \}$ be the set of expected transition probabilities.
\begin{enumerate}
    \item The set $\bar{\cP}$ is always s-rectangular.
    \item Assume that $\bar{\cP}$ is compact and convex. Then
    
\begin{align*}
    \sup_{\pi \in \PiH} \inf_{\mu \in \cVS} R(\pi,\mu) & = \sup_{\pi \in \PiH} \inf_{\mu \in \cVNS} R(\pi,\mu) \\
    & = \sup_{\pi \in \PiH} \inf_{\bar{P} \in \bar{\cP}} R(\pi,\bar{\cP}) \\
    & = \max_{\pi \in \PiS} \min_{\bar{P} \in \bar{\cP}} R(\pi,\bar{\cP})
\end{align*}
and $\displaystyle \max_{\pi \in \PiS} \min_{\bar{P} \in \bar{\cP}} R(\pi,\bar{\cP})$ is a robust MDP problem with s-rectangular convex compact uncertainty set $\bar{\cP}$, which admits a stationary uniformly optimal policy.
\end{enumerate}
\end{theorem}
\begin{proof}
\begin{enumerate}
    \item The s-rectangularity of $\bar{\cP}$ is a direct consequence of the fact that each $\nu \in \cV$ is a product distribution with marginals that can be chosen independently in $\cV_s$ for each $s \in \cS$, so that each $s$-th component $P_s$ of the transition probabilities can be sampled independently from some $\nu_s \in \cV_s$.
    \item 
The proof can be done from scratch (using truncated returns as has been done several times in the rest of these notes) but it is faster to reuse the main results from \Cref{sec: lec 3 time var P}. We start by analyzing $\sup_{\pi \in \PiH} \inf_{\mu \in \cVNS} R(\pi,\mu)$. From \Cref{lem: lec 3 return reformulation} we can see this as a special case of robust MDPs with time-varying transition probabilities, where we have that at each period the transition probabilities must belong to $\bar{\cP}$ and are chosen dependent on time but not on the current state. Recall that Model 3 in \Cref{sec: lec 3 time var P} allows for transition probabilities at period $t$ dependent on $(s_0,a_0,...,s_{t-1},a_{t-1})$. Hence $\inf_{\mu \in \cVNS} R(\pi,\mu)$ can be interpreted as optimizing over a class of adversary policies $\bar{\cP}^{\sf G}_{\sf M'}$ that is a subclass of $\bar{\cP}^{\sf G}_{\sf H}$ (we use the subscript ``{\sf M'}" because it corresponds to a weaker class than Markovian adversary policies). Since $\bar{\cP}$ is convex compact and s-rectangular (by assumption), we get that $\bar{\cP}^{\sf G}_{\sf H} = \bar{\cP}^{\sf S}_{\sf H}$ (from Lemma~\ref{lem: lec 3 Model 12 eq extensions}), and from \Cref{th: lec 3 Model 2} we get that stationary policies are optimal against adversary policies in $\bar{\cP}^{\sf G}_{\sf H}$ and, in particular, are optimal against adversary policies in  $\bar{\cP}^{\sf G}_{\sf S}$. This shows 
\[\sup_{\pi \in \PiH} \inf_{\mu \in \cVNS} R(\pi,\mu) = \max_{\pi \in \PiS} \min_{\bar{P} \in \bar{\cP}} R(\pi,\bar{\cP}).\]
In particular this shows that the worst-case over $\cVNS$ is attained at a stationary distribution in $\cVS$ when $\pi \in \PiS$, so that $\displaystyle \sup_{\pi \in \PiH} \inf_{\mu \in \cVNS} R(\pi,\mu) = \sup_{\pi \in \PiH} \inf_{\mu \in \cVS} R(\pi,\mu).$  The rest of the theorem follows from \Cref{th: lec 2 main rmdps}.
\end{enumerate}
\end{proof}
We proved that distributionally robust MDPs, where transition probabilities are sampled at each period, are actually equivalent to robust MDPs, where the transitions are stationary! An optimal policy for a distributionally robust MDP can therefore be computed by the algorithms presented in \Cref{sec: lec 2 algo}. As evident from the proof of \Cref{th: lec 3 drmdps}, this equivalence between distributionally robust MDPs and robust MDPs remains true even if we allow for more complicated adversary policies (e.g. history-dependent policies that can observe the current state but not the current action chosen by the decision maker), because the analysis we did for Model 2 in \Cref{sec: lec 3 time var P} also works for history-dependent adversaries.
\paragraph{Example: ambiguity sets based on Wasserstein distances.} Wasserstein distances have become popular in distributionally robust optimization (beyond distributionally robust MDPs), and we illustrate the reductions outlined above in this case. To simplify the notation, we use $\cD$ for the set of all possible transition probabilities out of a state: $\cD = \Delta(\cS)^{\cA}$.
Assume we are given a metric $d$ induced by a norm over $\cD$ and an integer $p \in \N \setminus \{0\}$ (the ``type"). 
The Wasserstein distance $W_p(\nu,\nu')$ between two distributions $\nu,\nu'$ over $\cD$  is 
\begin{align*}
W_{p}(\nu,\nu') = \min & \; \left(\E_{(P,P') \sim \kappa}  \left[ d(P,P')^{p} \right]\right)^{1/p} \\
& \kappa \in \Delta(\cD \times \cD) \\
& \Pi_{1} \kappa = \nu, \Pi_{2} \kappa = \nu'
\end{align*}
where $\Pi_{1} \kappa$ and $\Pi_{2} \kappa$ are the first and second marginals for a coupling $\kappa$ on $\cD \times \cD$. We can also allow $p = + \infty$ in which case we define $W_{\infty}(\nu,\nu')$ by
\begin{align*}
W_{\infty}(\nu,\nu') = \min & \; \kappa\textrm{-ess.sup}(d) \\
& \kappa \in \Delta(\cD \times \cD) \\
& \Pi_{1} \kappa = \nu, \Pi_{2} \kappa = \nu'
\end{align*}
with $\kappa\textrm{-ess.sup}(d)$ the essential supremum of $d$ under $\kappa$, equal to
$\inf \{ c \in \R \; | \; \kappa \left( \{ (x,y) \; | \;  d( x,y) > c \} \right) = 0 \}.$
Suppose that we have a nominal estimate $\bnu(s) \in \Delta(\cD)$ for each $s \in \cS$, and we are interested in ambiguity sets in the form of balls for the Wasserstein distance: for $p \in \N \cup \{+\infty \}$, define
\[\cV = \times_{s \in \cS} \cV_s, \cV_s = \{\nu \in \Delta(\cD) \; | \; W_{p}(\nu,\bnu(s)) \leq \alpha\}\]
with $\bnu = \otimes_{s \in \cS} \bnu(s)$ (a product distribution).
Assume that $\bnu$ is of the form $\bnu(s) = \frac{1}{N}\sum_{i=1}^{N} \delta \{\bar{P}_{i}(s)\}$ for some $\bar{P}_{i}(s) \in \cD$. Then the set of expected kernels $\bar{\cP}$ induced by $\cV$ is $\cP_p = \times_{s \in \cS} \cP_p(s)$, with, when $p \in \N \setminus \{0\}$,
\[\cP_p(s) = \left\{ \frac{1}{N}\sum_{i=1}^{N} P_{i}(s) \; \middle| \; \frac{1}{N}\sum_{i=1}^{N} d\left(P_i(s),\bar{P}_{i}(s)\right)^{p} \leq \alpha^{p}, P_i(s) \in \cD, \forall \; i \in \{1,...,N\} \right\} \]
and with $p=+\infty$ we have that $\cP_{\infty}(s)$ is defined by
\[\cP_{\infty}(s) = \left\{ \frac{1}{N}\sum_{i=1}^{N} P_{i}(s) \; \middle| \; d\left(P_i(s),\bar{P}_{i}(s)\right) \leq \alpha, P_i(s) \in \cD, \forall \; i \in \{1,...,N \} \right\}.\]
\section{Special classes of non-rectangular uncertainty sets}\label{sec: lec 3 non-rec}
We have seen in \Cref{chap: lec 1 robust MDP} that for general uncertainty sets, several optimization problems related to robust MDPs can be difficult, e.g., computing the worst-case transition probabilities for some policy may be NP-hard (\Cref{prop: lec 1 np hard}) and we may need history-dependent policies to play optimally. In fact \Cref{sec: lec 2 minimality rec} shows that rectangular models of uncertainty are ``necessary" (in a certain sense made precise there) to obtain favorable properties, such as  the validity of Bellman (dynamic programming) equations and the existence of uniformly optimal policies. 
In this section we will focus on some special cases of non-rectangular uncertainty sets that retain these favorable properties under some additional assumptions.
\subsection{r-rectangularity and generalizations}\label{sec: lec 3 r-rec}
We start with the definition of r-rectangular uncertainty sets.
\begin{definition}\label{def: lec 3 r-rec}
    We say that $\cP \subset \DSSA$ is r-rectangular if there exist $r \in \N$, some sets $\cW^{1},...,\cW^{r} \subset \DS$, and some coefficients $u_{sa}\in \Dr$ for each state-action pair $(s,a) \in \cS \times \cA$ such that
    \begin{equation}\label{eq: lec 3 def r rec}
        \cP = \left\{ P \; | \; P(s,a) = \sum_{i \in [r]} u_{sa}(i) w^{i}, w^{i} \in \cW^{i}, \forall \; i \in [r], \forall \; (s,a) \in \cS \times \cA \right\}
    \end{equation}
\end{definition}
The vectors $w^{1},...,w^{r}$ are called the {\em factors}; note that the coefficients ($u$'s) are fixed and known, but the factors can vary. If the sets $\cW^1,...,\cW^r$ are compact then $\cP$ is compact; vice versa if $\cP$ is compact we can replace $\cW^{1},...,\cW^{r}$ by their closures. Hence when we assume that $\cP$ is r-rectangular and compact we will assume that $\cW^{1},...,\cW^{r}$ are compact. There are a few interpretations for the r-rectangularity assumption:
\begin{itemize}
    \item Equation \eqref{eq: lec 3 def r rec} states that for each state-action pair $(s,a)$, the transition probabilities $P(s,a) \in \DS$ toward the next state can be factorized as first sampling some $i \in [r]$ following the distribution $u_{sa} \in \Delta([r])$ and then sampling the next state according to the factor distribution $w^i \in \DS$. 
    \item We can rewrite $P(s,a) = \sum_{i \in [r]} u_{sa}(i) w^{i}$ as $P(s,a) = W u_{sa}$ with $W \in \R^{\cS \times [r]}$ the matrix with columns $w^{1},...,w^{r} \in \DS$, hence $P(s) = W u_s$ with $P(s) \in \DSA$ and $u_s \in (\Dr)^{\cA}$. Hence we can interpret r-rectangularity as some form of non-negative matrix factorization, which is the optimization problem of factoring a nonnegative matrix $A$ as a product of non-negative matrices $B,C$: $A=BC$; in our case one matrix is kept fixed (the coefficients $u$) while the other can vary (the factor matrix $W$).
    \end{itemize}

\paragraph{Relations with rectangular models.}
It should be clear that in general r-rectangular uncertainty sets are neither s-rectangular nor sa-rectangular. For instance if $r=1$ then we have a ``rank one" uncertainty set where for every state-action pair $(s,a)$ we have $P(s,a) = w^{1},w^{1} \in \cW^{1}$. Clearly this set is neither s-rectangular nor sa-rectangular (unless $\cW^{1}$ is a singleton), i.e. we have the following lemma.
\begin{lemma}
    There exists $\cP \subset \DSSA$ that is r-rectangular but neither s-rectangular nor sa-rectangular.
\end{lemma}
However we always have that $(w^1,...,w^r)$ belongs to the Cartesian product set $\cW^{1} \times ... \times \cW^{r}$, even though $\cP$ itself may not be a Cartesian product set under the r-rectangularity assumption. Hence we can interpret r-rectangularity as a form of ``underlying" rectangularity. In particular r-rectangular models contain sa-rectangular models as a ``limiting case" ($r=|\cS| \times |\cA|$), as we show next.
\begin{lemma}
    Let $\cP \subset \DSSA$ be sa-rectangular. Then $\cP$ is r-rectangular and we can choose $r=|\cS| \times |\cA|$.
\end{lemma}
\begin{proof}
Let $\cP = \times_{s,a} \cP_{sa}$. Then we can just use $r = |\cS| \times |\cA|$ and ``enumerate". We have a bijection $\phi:\cS \times \cA \to [r]$. Then we can write $\cP$ as an r-rectangular uncertainty set by constructing the coefficients $u_{sa}$ as $u_{sa}(i)=1$ if $\phi((s,a))=i$ and $u_{sa}(i)=0$ otherwise, and $w^{i}=P(\phi^{-1}(i)), \cW^{i} = \cP_{\phi^{-1}(i)}$.
\end{proof}
The fact that there exist s-rectangular uncertainty sets that are not r-rectangular will be shown after \Cref{cor: lec 3 properties rrec}.
It is worth emphasizing that sa-rectangularity and s-rectangularity are geometric properties (being equal to its s-rectangular or sa-rectangular extensions), while r-rectangularity is an ``existence property" (see \Cref{def: lec 3 r-rec}), which may be harder to disprove. So the easiest way to prove that a set is not r-rectangular is to prove that it does not have some properties implied by r-rectangularity.
\paragraph{r-rectangular models are intractable in general.}
\Cref{sec: lec 2 minimality rec} tells us that we should not expect r-rectangular models to have favorable tractability properties in all generality. Indeed we first show that we may need history-dependent policies even in rank one uncertainty sets.
\begin{proposition}\label{prop: lec 3 r rank non-stationary}
    There exists a compact convex r-rectangular uncertainty set $\cP$ with $|\cS|=|\cA|=2$ and rank $r=1$, for which for some choice of rewards $r \in \RSAS$ and initial distribution $p_0 \in \DS$, an optimal policy must be non-stationary for any choice of $\gamma \in (0,1)$.
\end{proposition}
\begin{proof}
\begin{figure}
    \centering
    \includegraphics[width=0.5\linewidth]{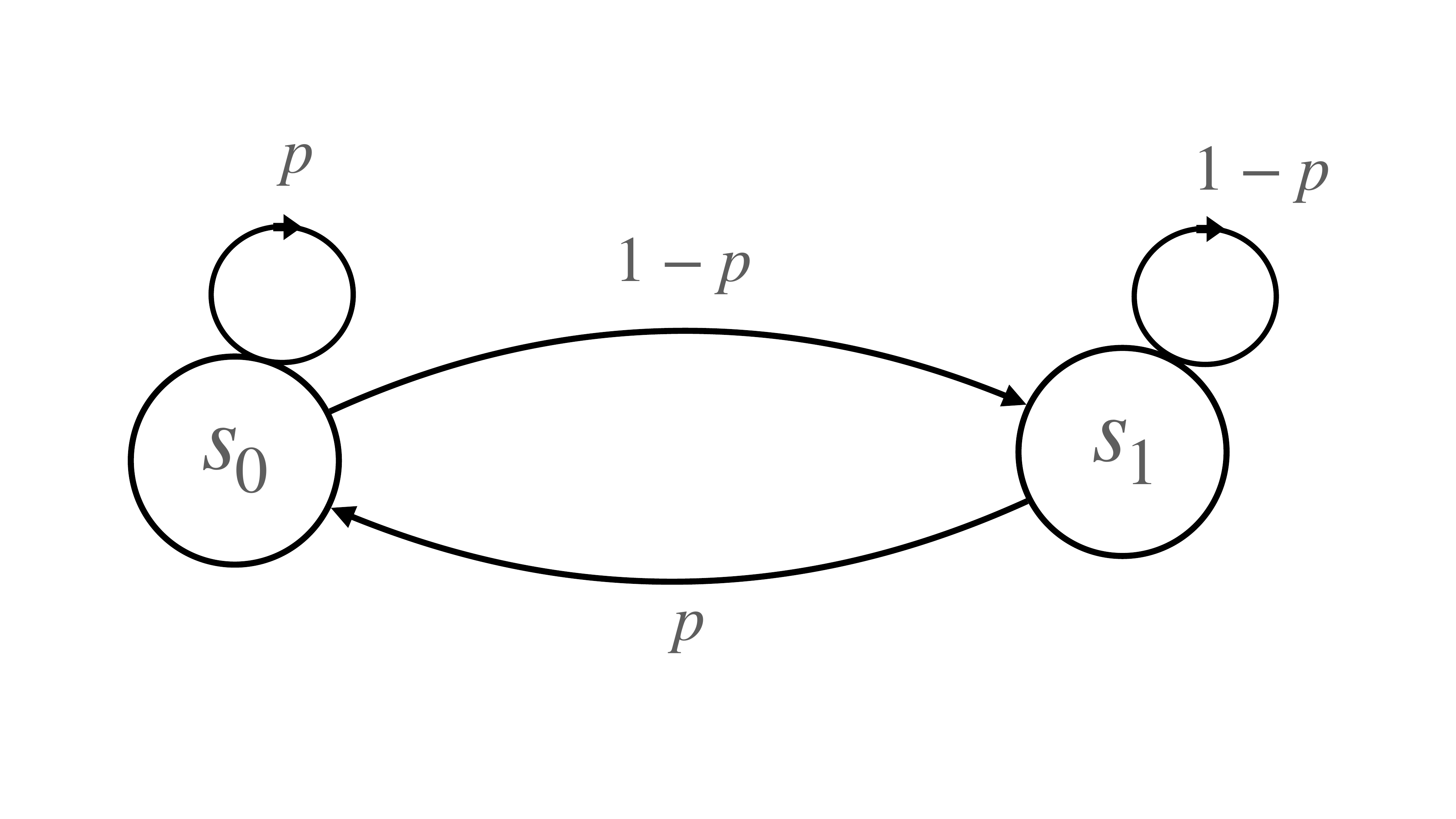}
    \caption{A rank one r-rectangular robust MDP.}
    \label{fig: lec 3 rrec}
\end{figure}

We consider an r-rectangular MDP with two states $\cS = \{s_0,s_1\}$ and two actions $\{a_1,a_2\}$. With a rank one r-rectangular robust MDP we get that starting from period $t\geq 1$ the next state is sampled following $w=(p,1-p) \in \Delta(\cS)$, with $p \in [0,1]$ the probability that the next state is $s_0$; the state at period $t=0$ is sampled from $p_0$ and we take here $p_0 = (1,0)$. This is represented in \Cref{fig: lec 3 rrec}.  We choose the rewards as follows:
  \begin{align*}
      & r(s_0,a_1,s_0)=0, r(s_0,a_1,s_1) = -1 \\
      & r(s_0,a_2,s_0) = -1, r(s_0,a_2,s_1) = 1 \\
      & r(s_1,a_1,s_0)=-1, r(s_1,a_1,s_1) = 0 \\
      & r(s_1,a_2,s_0) = 1, r(s_1,a_2,s_1) = -1.
      \end{align*}
The decision maker is facing the following tradeoff. Action $a_1$ is the best action to choose in either state if we think that we are going to remain in that same state at the next period; otherwise we better choose $a_2$. The adversary chooses the probability $p \in [0,1]$ that the next state is $s_0$ (otherwise it is $s_1$ with probability $1-p$). Direct calculations give that for $\pi \in \PiS$ we have
\begin{align*}
    (1-\gamma)R(\pi,P) = & (1-\gamma + \gamma p) \left( 1-2p + (3p-2)\pi(s_0,a_1)\right) \\
    & + \gamma (1-p) \left( 2p-1 + (1-3p)\pi(s_1,a_1) \right)
\end{align*}
and therefore the optimal policy among stationary ones is to choose $\pi(s_0,a_1) = \pi(s_1,a_1) = 2/3$  with a worst-case return of $-\frac{1}{3}\frac{1}{1-\gamma}$.

Now consider the following history-dependent policy $\pi'$: in the first two periods we match the best stationary policies: at $t=0$ play $\pi(s_0,a_1) = 2/3$, at $t=1$ play the optimal stationary policy again, then for $t\geq 2$: if the state visited at period $1$ is equal to $s_0$ then play forever $\pi(s_0,a_1) = 1, \pi(s_1,a_2)=1$ (that is, we bet that $p \geq 1/2$, and that we will often ``stay" if in $s_0$ and ``move away" if we're in $s_1$) otherwise play forever $\pi(s_0,a_2)=1,\pi(s_1,a_1)=1$ (that is, we bet that $p \leq 1/2$, and that we will often ``move away" if in $s_0$ and ``stay" if in $s_1$).

For this history-dependent policy, we have a return of $R(\pi',P) = - \frac{1}{3}-\frac{\gamma}{3} - \frac{\gamma^2 p (1-p)}{1-\gamma}$, with worst-case attained at $p=1/2$ and a worst-case return of $- \frac{1}{3} - \frac{\gamma}{3} - \frac{\gamma^2}{4 (1-\gamma)}$, which is strictly larger than the best worst-case return for a stationary policy computed above. This shows that non-stationary policies strictly outperform stationary ones in this small example.
\end{proof}
We next prove that solving the policy evaluation problem is NP-hard for r-rectangular uncertainty sets, even when $r=1$.
\begin{proposition}\label{prop: lec 3 hardness r rec}
Let $\pi \in \PiSD$.
For a general convex compact r-rectangular polytope $\cP \subset \RSAS $ with $r=1$, it is in general NP-hard to decide $\inf_{P \in \cP} R(\pi,P) \leq c$.
\end{proposition}
\begin{proof}
    The proof is similar to the hardness result from \Cref{prop: lec 3 hardness uncertain r p} and uses a reduction from the Clique problem. Let $G=(V,E)$ be an instance of a Clique decision problem with $V=\{1,...,n\}$ and $A \in \{0,1\}^{n \times n}$ the adjacency matrix. We take $\cS = \{s_0\} \cup V,\cA = \{a_0\},p_0 = \delta_{s_0},\gamma = 1/2$ (we could take any $\gamma \in (0,1)$ here). The uncertainty set is parametrized by $p \in \Delta(V)$ and is defined as 
    \[ \cP = \{ P \; | \; \exists \; p \in \Delta(V), P(s,a_0,s_0)=0,(P(s,a_0,s'))_{s' \in V} = p, \forall \; s \in \cS \}.\]
    Hence $\cP$ is r-rectangular of rank $r=1$ and it is a compact convex polytope. We choose the rewards as $r(s,a_0,s') = 1 - A_{ss'}$ for $s,s' \in V$ and $0$ everywhere else. With these definitions the expected one-step instantaneous reward is $\sum_{s,s' \in V} p_{s} p_{s'} (1-A_{ss'}) = 1 - p \tr A p$ from period $1$ onward, when the transition probabilities are parametrized by $p \in \Delta(V)$. 
    Hence $R(\pi,P) = \sum_{t=1}^{+\infty} 2^{-t} (1-p \tr Ap) = 1 - p \tr Ap$. Motzkin-Straus theorem gives again that $\max_{p \in \Delta(V)} p \tr A p = 1 - \omega(G)^{-1}$ with $\omega(G)$ the clique number of $G$. So again if we could decide $\min_{P \in \cP} R(\pi,P) \leq c$ with $c = 1/k$ for $k \in \N$ we could decide $\omega(G) \geq k$, which concludes the reduction.
\end{proof}
\paragraph{Tractability under some additional assumptions.}
The previous propositions suggest that we need some additional assumptions on the structure of the MDP instance to obtain favorable properties for r-rectangular uncertainty sets. Indeed, the original papers on robust MDPs with r-rectangular uncertainty sets claim tractability of the policy evaluation problem and of computing an optimal policy. This holds when we use the assumption that the rewards are independent of the next state, that is, that $r(s,a,s') = r(s,a,s'')$ for all $s,s',s'' \in \cS$ and all $a \in \cA$. In fact in this case, we can show that r-rectangular uncertainty sets are equivalent to sa-rectangular uncertainty sets, in two different ways:
\begin{itemize}
    \item We can augment the set of states to show that an r-rectangular model of uncertainty on the original robust MDP instance yields an sa-rectangular uncertainty set on the augmented robust MDP;
    \item Even without state-augmentation, we can show that an r-rectangular uncertainty set $\cP$ yields the same Bellman operator as its sa-rectangular extension $\recsa(\cP)$ (defined in \Cref{sec: lec 3 time var P}).
\end{itemize}
We make these two statements precise next. We start with the first statement below.
\begin{proposition}\label{prop: lec 3 rrec as sarec}
   Consider a robust MDP instance $ (\cS,\cA,\cP,r,\gamma)$ with $\cP \subset \DSSA$ r-rectangular compact. Assume that the rewards are independent of the next state and write $r(s,a) \in \R$ for the rewards at the state-action pair $(s,a)$. Consider the following robust MDP instance $(\hat{\cS},\hat{\cA},\hat{\cP},\hat{r},\hat{\gamma})$:
    \begin{itemize}
        \item The set of augmented states is $\hat{\cS} = \cS \cup [r]$
        \item The set of actions is $\cA$ at $\hat{s} \in \hat{\cS}$ such that $\hat{s} \in \cS$, and $\{1\}$ for $\hat{s} \in [r]$ (that is the decision maker has no choice at states in $[r]$). We can identify the set of stationary policies $\hat{\Pi}_{\sf S}$ in this augmented robust MDP instance with $\{\phi(\pi)  \; | \; \pi \in \PiS\}$ for $\phi: \PiS \to \R^{(|\cS|+r) \times (|\cA|+1)}$ defined as
        \[\phi(\pi)=\begin{pmatrix}
            \pi & 0_{\cS} \\
            0_{[r]\times \cA} & 1_{[r]} 
        \end{pmatrix}, \forall \; \pi \in \PiS.\]
        \item For $(s,a) \in \cS \times \cA$ we choose $\hat{r}(s,a) = r(s,a)$ and $\hat{r}(i,1)=0$ for $i \in [r]$.
        \item For $(s,a) \in \cS \times \cA$  we define $\hat{\cP}_{sa} = \{(0_{\cS},u_{sa})\}$ and for $i \in [r]$ we define $\hat{\cP}_{i 1} = \{ (w^{i},0_{[r]}) \; | \; w^i \in \cW^{i}\}$.
        \item The initial distribution $\hat{p}_0$ assigns probability $p_0(s)$ on state $s \in \cS$.
        \item The discount factor $\hat{\gamma}$ is $\hat{\gamma}= \sqrt{\gamma}$.
    \end{itemize}
Then $\hat{\cP}$ is sa-rectangular and $\phi$ is a bijection between stationary policies in the original robust MDP instance and in the lifted robust MDP instance. 
\end{proposition}
\begin{proof}
    The proof is a direct reformulation and is illustrated in \Cref{fig: lec 3 rrec lifted MDP}.

        \begin{figure}[h!]
    \begin{center}
    \begin{subfigure}{0.35\textwidth}
    \begin{center}
        \includegraphics[width=\linewidth]{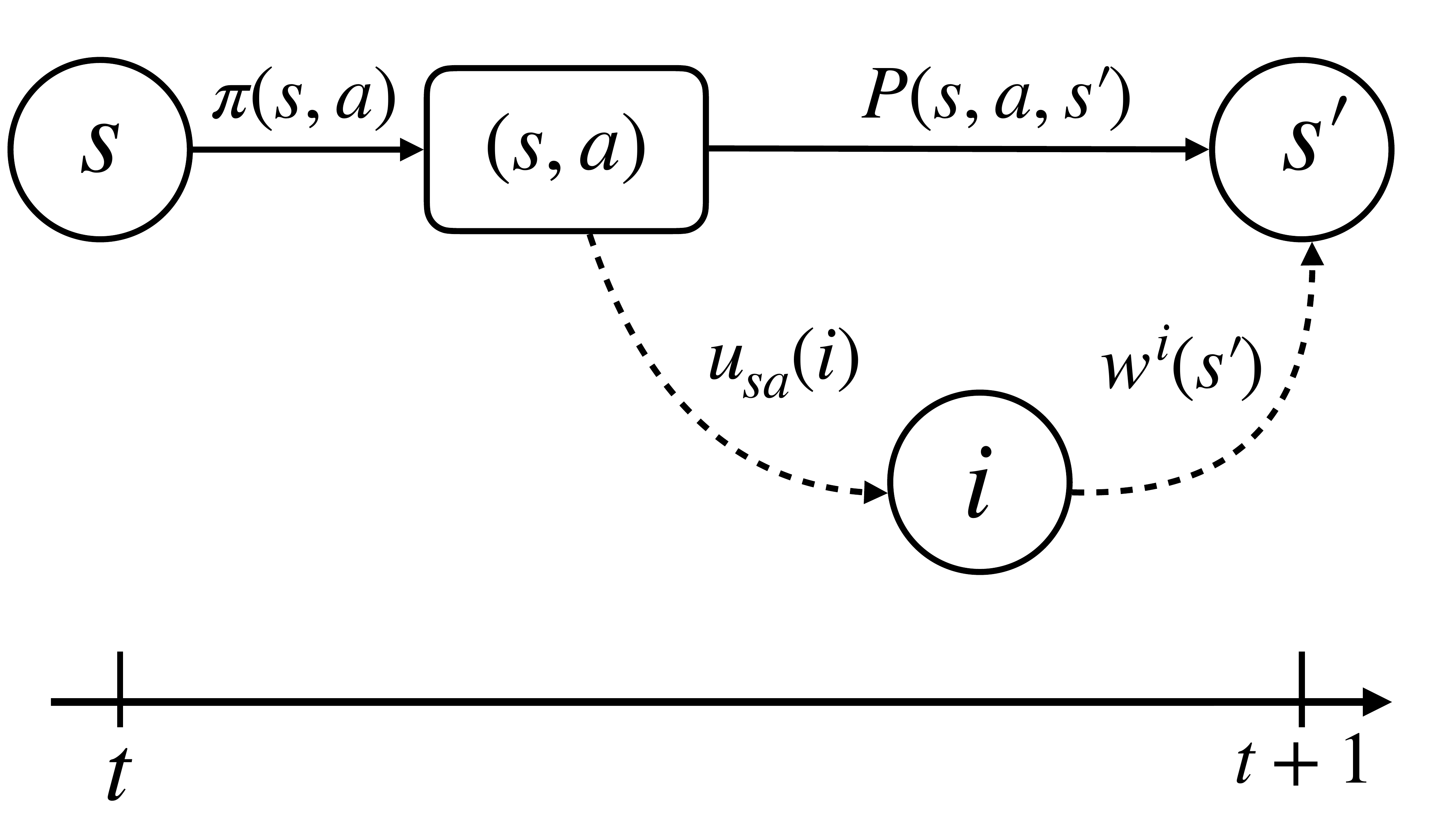}
    \end{center}
    \caption{Transitions and rewards for an r-rectangular uncertainty set}
    \label{fig: lec 3 factor MDP}
    \end{subfigure}
    \begin{subfigure}{0.55\textwidth}
    \begin{center}
        \includegraphics[width=\linewidth]{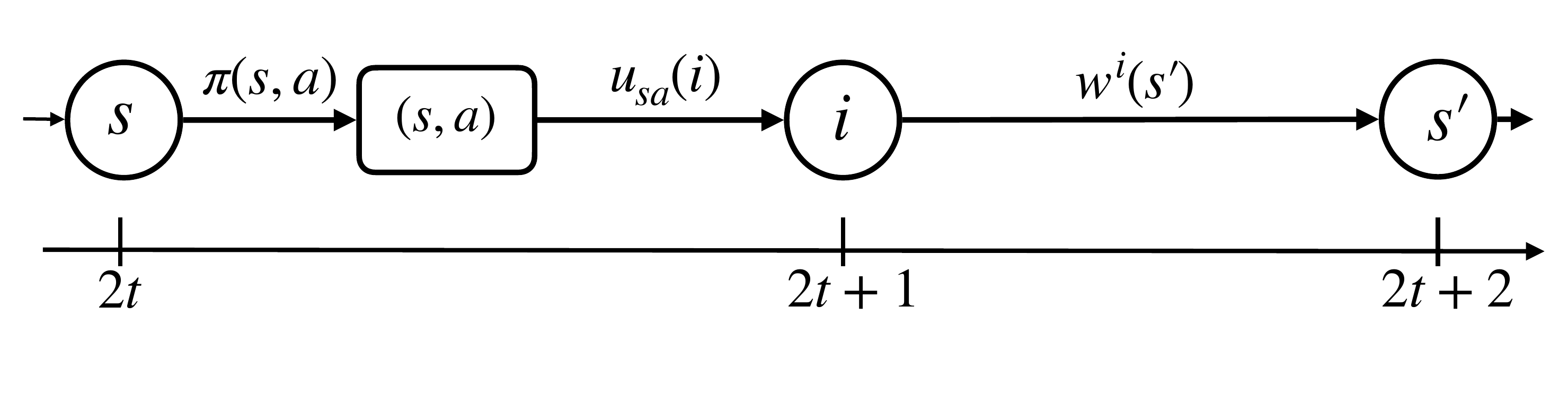}
    \end{center}
    \caption{Transitions and rewards in the augmented robust MDP instance}
    \label{fig: lec 3 lifted MDP}
    \end{subfigure}
    \end{center}
        \caption{Illustration for the state-augmentation reformulation from \Cref{prop: lec 3 rrec as sarec}.}
        \label{fig: lec 3 rrec lifted MDP}
    \end{figure}
\end{proof}
This proposition shows that r-rectangularity of $\cP$ is just sa-rectangularity on an augmented state space; that is, the non-rectangularity of $\cP$ on the original set of states comes from the ``projection" of some components from the augmented set of states. So already we can see that the favorable tractability properties of sa-rectangularity will hold for r-rectangular models (under a compactness assumption so that we can apply the results from \Cref{sec: lec 2 sarec}), that is, we directly obtain the existence of stationary uniformly optimal policies for $\cP$ r-rectangular compact uncertainty sets. However we can prove even stronger results: r-rectangular uncertainty sets are equivalent to sa-rectangular uncertainty sets. This is the purpose of the next proposition.
\begin{proposition}\label{prop: lec 3 rrec sa-rec extension}
 Let $\cP$ be a compact r-rectangular uncertainty set and assume that the rewards are independent of the next state. Then
 \[\forall \; v \in \RS, \exists P' \in \cP, \forall \; \pi \in \PiS, \forall s \in \cS,  T^{\pi,P'}v(s)=\min_{P \in \recsa(\cP)} T^{\pi,P}v(s).\]
 where $\recsa(\cP)$ is the sa-rectangular extension of $\cP$ defined in \Cref{sec: lec 3 time var P}. In particular this implies
 \[\forall \; v \in \RS, \forall \; \pi \in \PiS, \forall s \in \cS, \min_{P \in \cP} T^{\pi,P}v(s) = \min_{P \in \recsa(\cP)} T^{\pi,P}v(s).\]
\end{proposition}
\begin{proof}
We have by definition of $\recsa(\cP)$ that
\begin{align*}
    \min_{P \in \recsa(\cP)} T^{\pi,P}v(s) & = \sum_{a \in \cA} \pi(s,a) r(s,a) + \gamma \sum_{a \in \cA} \pi(s,a) \min_{P \in \cP} P(s,a) \tr v \\
    & = \sum_{a \in \cA} \pi(s,a) r(s,a) + \gamma \sum_{a \in \cA} \pi(s,a) \sum_{i \in [r]} u_{sa}(i) \min_{w^i \in \cW^i} w^{i \; \top} v
\end{align*}
By choosing factors in $\arg \min_{w^i \in \cW^i} w^{i \; \top} v$ for each $i \in [r]$ we recover $P'$ feasible in $\cP$, independent of $\pi,s$, such that $T^{\pi,P'}v(s)=\min_{P \in \recsa(\cP)} T^{\pi,P}v(s), \forall \; s \in \cS$ as claimed.
The second implication is because $\cP \subset \recsa(\cP)$.
\end{proof}
A rephrasing of the above proposition is that the set of operators $\{T^{\pi,P} \; | \; P \in \cP\}$ satisfies the attainability property with the contraction $v \mapsto (\min_{P \in \recsa(\cP)} T^{\pi,P}v(s))_{s \in \cS}$. Note that we are not claiming that $\cP$ and $\recsa(\cP)$ coincide (in particular $\cP$ is only r-rectangular so it may not be sa-rectangular), but rather that minimizing $P \mapsto T^{\pi,P}v(s)$ over $\cP$ and $\recsa(\cP)$ yields the same optimal objectives.
The above proposition is enough to show that $\cP$ can be replaced by its sa-rectangular extension in the following sense.
\begin{proposition}\label{prop: lec 3 rrec sarec extension}
     Let $\cP$ be a compact r-rectangular uncertainty set and assume that the rewards are independent of the next state. Let $\pi \in \PiS$. Then 
     \[ \min_{P \in \cP} R(\pi,P) = \min_{P \in \recsa(\cP)}  R(\pi,P).\]
\end{proposition}
\begin{proof}
    We have $\displaystyle \min_{P \in \recsa(\cP)}  R(\pi,P) = p_0 \tr \hat{v}^{\pi}$ with $\hat{v}^{\pi}$ the unique fixed point of 
    \[v \mapsto \left(\min_{P \in \recsa(\cP)} T^{\pi,P}v(s)\right)_{s \in \cS}.\] But we know that $\{T^{\pi,P} \; | \; P \in \cP\}$ satisfies the attainability property with this contraction. Hence $\hat{v}^{\pi} = v^{\pi,P^{\pi}}$ for some feasible $P^{\pi} \in \cP$. Hence $\min_{P \in \recsa(\cP)} R(\pi,P) = R(\pi,P^{\pi}) \geq \min_{P \in \cP} R(\pi,P)$. Since the reverse inequality is always true (given $\cP \subset \recsa(\cP)$) we get that $\min_{P \in \cP} R(\pi,P) = \min_{P \in \recsa(\cP)}  R(\pi,P)$.
\end{proof}
Combining the arguments above proves the following corollary.
\begin{corollary}\label{cor: lec 3 properties rrec}
    Let $\cP$ be an r-rectangular compact uncertainty set and assume that the rewards are independent of the next state. Then the following properties hold:
    \begin{enumerate}
       \item We have
       \begin{align*}
           \sup_{\pi \in \PiH} \inf_{P \in \cP} R(\pi,P) = \max_{\pi \in \PiS} \min_{P \in \cP} R(\pi,P) = \max_{\pi \in \PiS} \min_{P \in \recsa(\cP)} R(\pi,P).
       \end{align*}
        \item {\em Optimal policies.} 
        The set of stationary deterministic uniformly optimal policies for $\cP$ coincides with the set of stationary deterministic uniformly optimal policies for $\recsa(\cP)$. In particular the set of stationary deterministic uniformly optimal policies for $\cP$ is non-empty and can be computed with the algorithms outlined in \Cref{sec: lec 2 sarec}.
    \end{enumerate}
\end{corollary}
\begin{proof}
\begin{enumerate}
    \item 
    Let $\pi\opt$ be a stationary deterministic uniformly optimal policy for $\recsa(\cP)$, and $P\opt$ such that $P\opt \in \arg \min_{P \in \recsa(\cP)} R(\pi\opt,P),\pi\opt \in \arg \max_{\pi \in \PiSD} R(\pi,P\opt)$.

    We know that $R(\pi,P\opt) \leq R(\pi\opt,P\opt)$ for any $\pi \in \PiH$. From \Cref{prop: lec 3 rrec sa-rec extension} we can actually choose $P\opt \in \cP$, hence $ \sup_{\pi \in \PiH} \min_{P \in \cP} R(\pi,P) \leq R(\pi\opt,P\opt)$. But $R(\pi\opt,P\opt) = \max_{\pi \in \PiSD} \min_{P \in \recsa(\cP)} R(\pi,P) = \sup_{\pi \in \PiH} \min_{P \in \recsa(\cP)} R(\pi,P)$, hence we have proved
    \[\sup_{\pi \in \PiH} \min_{P \in \cP} R(\pi,P) \leq \max_{\pi \in \PiSD} \min_{P \in \recsa(\cP)} R(\pi,P) = \sup_{\pi \in \PiH} \min_{P \in \recsa(\cP)} R(\pi,P).\]
    Since $\cP \subset \recsa(\cP)$ we must have that all terms above are equal, which concludes the proof for the first statement.
    \item \Cref{prop: lec 3 rrec sarec extension} shows that every stationary policy has the same worst-case value over $\cP$ and over $\recsa(\cP)$, and that the two maxmin problems have the same optimal value. Hence a stationary deterministic policy is uniformly optimal for $\cP$ if and only if it is uniformly optimal for $\recsa(\cP)$.
\end{enumerate}
\end{proof}
The fact that optimal policies for r-rectangular compact uncertainty sets can be chosen {\em deterministic} implies that some s-rectangular sets are not r-rectangular (for instance we saw in the previous chapter that this is the case for \Cref{ex: lec 0 ex 1} with $p \in [0,1]$).
\paragraph{Generalizations.}
\Cref{prop: lec 3 rrec sarec extension} has highlighted an interesting property: some uncertainty sets $\cP$ are such that the induced set of operators $\{T^{\pi,P} \; | \; P \in \cP\}$ has the attainability property for $v \mapsto (\min_{P \in \recsa(\cP)} T^{\pi,P}v(s))_{s \in \cS}$ (in particular this is true for r-rectangular uncertainty sets when the rewards are independent of the next state). \Cref{cor: lec 3 properties rrec} shows that this attainability property yields that $\cP$ inherits the good properties of $\recsa(\cP)$ in terms of optimal policies and tractability (we don't use at all r-rectangularity for \Cref{cor: lec 3 properties rrec}, only the attainability outlined above). This yields the following two interesting questions:
\begin{enumerate}
    \item First question: are there other classes of uncertainty sets $\cP$ such that, when the rewards are independent of the next state, we have
     \[\forall \; v \in \RS, \exists P' \in \cP, \forall \; \pi \in \PiS, \forall s \in \cS,  T^{\pi,P'}v(s)=\min_{P \in \recsa(\cP)} T^{\pi,P}v(s).\]
    \item Second question: how about the same question with s-rectangular extensions, that is, are there other classes of uncertainty sets $\cP$ such that when rewards are independent of the next state, we have
    \[\forall \; v \in \RS, \forall \; \pi \in \PiS, \exists P' \in \cP, \forall s \in \cS,  T^{\pi,P'}v(s)=\min_{P \in \recs(\cP)} T^{\pi,P}v(s).\]
\end{enumerate}
    Note the slightly weaker formulation for the second question compared to the first one (I invert the positions of the quantifiers for the kernel $P' \in \cP$ and the policy $\pi \in \PiS$); this is because the second question involves the s-rectangular extension of $\cP$, so that the choice of policy may matter for this formulation, whereas the sa-rectangular extension is state-action wise.
The answers to both questions are positive, and several classes of non-rectangular uncertainty sets have been proposed that satisfy these attainability properties. I give here two very simple examples (left as an exercise):
\begin{enumerate}
    \item Consider $\cP \subset \DSSA$ such that
\begin{equation}\label{eq: lec 3 generalized rectangularity}
    \begin{aligned}
    \cP &  = \left\{ \left.\left(\sum_{i=1}^r u_{sa}(i)w_{i} \right)_{sa} \right|(w^{1},\ldots,w^{r}) \in \cW,(\bm{u}_{sa})_{s,a} \in \mathcal{U} \right\},\\
    \cW & = \times_{i \in [r]} \cW^{i}, \cW^{i} \subset \DS, \forall \; i \in [r], \\
    \cU & =  \times_{(s,a) \in \cS \times \cA} \cU^{sa}, \cU^{sa} \subset \Dr, \forall \; (s,a) \in \cS \times \cA.
    \end{aligned}
\end{equation}
Then $\cP$ provides a positive answer to the first question (when $\cW^{i}$ and $\cU^{sa}$ are compact sets for $i \in [r]$ and $(s,a) \in \cS \times \cA$).
    \item Let $\cP \subset \DSSA$ be an uncertainty set and assume that $\cS$ can be partitioned into $\cS_1$ and $\cS_2$ such that the transitions out of the states in $\cS_1$ are s-rectangular while the transitions out of $\cS_2$ are r-rectangular:
    
        \begin{equation}\label{eq: lec 3 (s+r)-rec}
    \begin{aligned}
    \cP &  = \cP_{1} \times \cP_{2}, \cP_{1} \subset \Delta(\cS)^{\cS_{1} \times \cA}, \cP_{2} \subset \Delta(\cS)^{\cS_{2}\times \cA}, \\
    \cP_{1} & = \times_{s \in \cS_{1}} \cP_{s}, \cP_{s} \subset \Delta(\cS)^{\cA}, \\
    \cP_{2} & = \left\{ \left({\sum_{i=1}^{r}}
    u_{sa}(i)w^{i}\right)_{s \in \cS_{2},a \in \cA} \; | \; (w^{1},...,w^{r}) \in \cW^{1} \times ... \times \cW^{r} \right\}
\end{aligned}
\end{equation}
    Then $\cP$ gives a positive answer to the second question (when $\cP_1$ and $\cP_2$ are compact sets).  
\end{enumerate}
It is not hard to come up with other classes of uncertainty sets that give a positive answer to the first or second question above (see references of this chapter). However in some sense (more precisely: regarding policy evaluation), all of these uncertainty sets are ``equivalent" to their sa-rectangular extensions (for those with a positive answer to the first question above) or their s-rectangular extensions (for those with a positive answer to the second question above), so perhaps it is easier to work with s-rectangular or sa-rectangular uncertainty sets directly.
\subsection{Lightning-Does-not-Strike-Twice uncertainty set}\label{sec: lec 3 k-rec}
\begin{figure}
    \centering
    \includegraphics[width=0.85\linewidth]{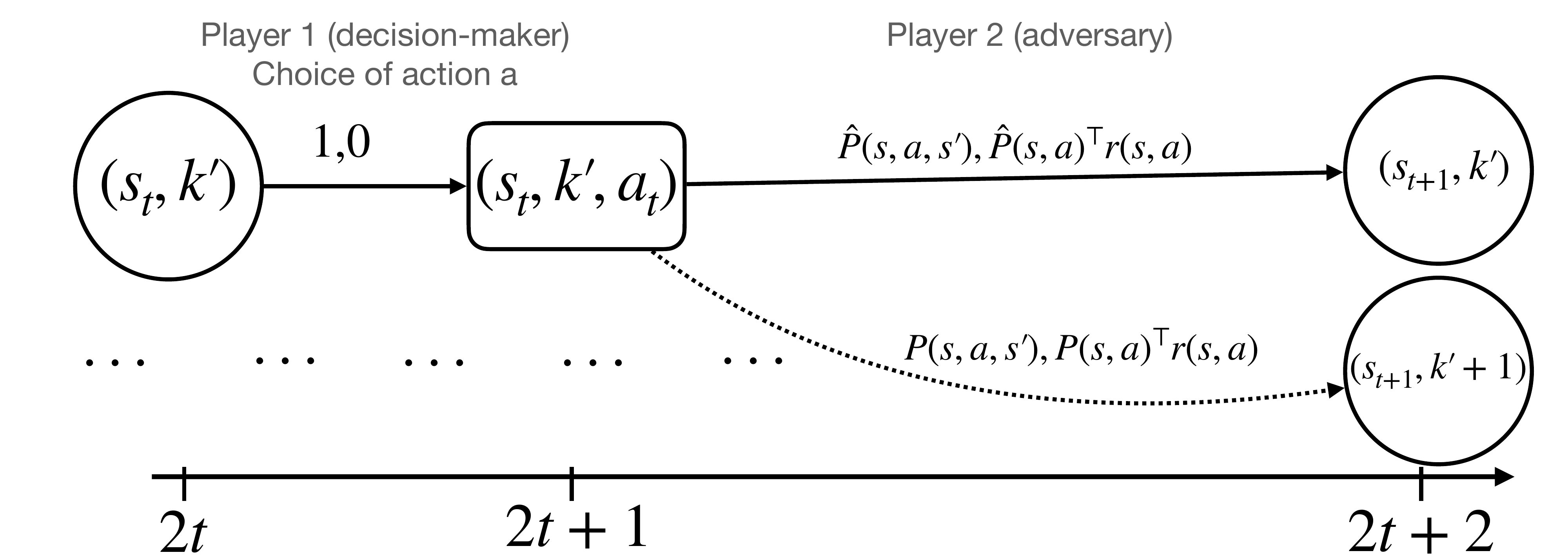}
    \caption{The perfect-information stochastic game reformulation of robust MDPs with Lightning-Does-not-Strike-Twice (LDST) model of uncertainty.}
    \label{fig: lec 3 LDST}
\end{figure}

Another popular model is the so-called {\em Lightning-Does-not-Strike-Twice} (LDST) model of uncertainty, and its generalization, {\em k-rectangularity.}
In the (adaptive) LDST model we are given a nominal kernel $\hP \in \DSSA$, and some s-rectangular uncertainty set $\cP = \times_{s} \cP_{s}$, assumed to be compact here, such that $\hP \in \cP$. The transition probabilities are allowed to vary over time and the dynamics unfold as follows: at period $t$, the decision maker observes the current state $s_t$ and all the past choices $P^{t'}(s_{t'}) \in \DSA$ chosen by the adversary at the past periods $t' < t$ out of the past states $s_0,...,s_{t-1}$. Based on this information, the decision maker chooses $a_t$, the adversary observes $s_t,a_t$ and chooses the transition probabilities $P^{t}(s_t)$ out of $s_t$ for period $t$. The next state $s_{t+1} \in \cS$ is sampled from $P^t(s_t,a_t) \in \DS$. The constraint for the adversary is that along each realized infinite trajectory, the transition probabilities can deviate at most $k$ times from the nominal kernel $\hP$, where the integer $k \in \N$ is given and known to the decision maker and the adversary. Note the difference with the previous robust MDP model considered in these notes: in the LDST model, at any period $t \in \N$ the decision maker is supposed to be able to know the cardinality of $\{t' \; | \; P^{t'}(s_{t'})\neq \hP(s_{t'}), t' < t\}$, i.e., the decision maker knows how many times the kernel has deviated from $\hP$ and can calculate how much of its ``budget of deviation" the adversary has used until period $t$. We emphasize that this is a significant departure from the robust MDP model studied in the rest of these notes, where the decision maker makes their decision at period $t \in \N$ based on $(s_0,a_0,s_1,a_{1},...,s_{t-1},a_{t-1},s_{t})$ the current history up to period $t$, whereas in the LDST model the decision maker makes its decision at period $t$ based on $(s_0,a_0,k_0,s_1,a_1,...,s_{t-1},a_{t-1},s_{t},k_t)$ where $k_t=\sum_{t'=0}^{t-1} 1_{\{P^{t'}(s_{t'})\neq \hP(s_{t'})\}}$, and the adversary bases its decision on $(s_0,a_0,k_0,s_1,a_1,...,s_{t-1},a_{t-1},s_{t},k_t,a_t)$. In particular one may wonder how can the decision maker, which in principle only observing the sequence of states, could infer that the transitions $(s_t,a_t) \to s_{t+1}$ happened based on some particular distribution $P(s_t,a_t) \in \DS$ over the next states. That could be the case if $s_{t+1}$ does not belong to the support of $\hP(s_t,a_t)$. Here it is interesting to compare the LDST model for the sets of information available to the decision maker and adversary, and the connection we made between sa-rectangular robust MDPs and perfect-information stochastic games in \Cref{sec: lec 3 sgs}. In particular, it is straightforward to reformulate robust MDPs with the LDST model of uncertainty as a perfect-information stochastic game as follows:
\begin{itemize}
    \item The set of states of the perfect-information stochastic game is $\cS \times \{0,...,k\} \cup (\cS \times \{0,...,k\} \times \cA)$, which we can represent visually as two stacks: a first stack of $k+1$ copies of $\cS$ and a second stack of $k+1$ copies of $\cS \times \cA$. The decision maker controls the states in $\cS \times \{0,...,k\}$ and the adversary controls the states in $\cS \times \{0,...,k\} \times \cA$.
    \item Let $(s,k') \in \cS \times \{0,...,k-1\}$. When the decision maker chooses action $a \in \cA$ in the state $(s,k')$ the transition is deterministic to $(s,k',a)$. When in $(s,k',a)$ the adversary chooses the current transition probabilities $P(s,a) \in \cP_{sa}$. If $\hP(s,a) \neq P(s,a)$ then the transition is to $(s',k'+1)$ with probability $P(s,a,s')$, otherwise the transition is to $(s',k')$ with probability $P(s,a,s')$. In both cases, the rewards are obtained at the transition $(s,k',a)$ to $(s',k'')$ with $k'' \in \{k',k'+1\}$ (its value depending on whether $\hP$ was chosen at period $t$) and are equal to $P(s,a)\tr r(s,a) \in \R$.

    For the states in the form $(s,k)$ the situation is the same for the decision maker but the adversary is constrained to choosing $\hP(s,a)$, i.e., the MDP becomes ``stuck" in the last level of the stacked copies of $\cS$ and $\cS \times \cA$ (since the budget of deviations is exhausted).
    \item Because the rewards are obtained every other period we choose a discount factor of $\sqrt{\gamma}.$
\end{itemize}
The resulting perfect-information SG is represented in \Cref{fig: lec 3 LDST}. For the sake of conciseness, for a given $v \in \R^{\cS \times \{0,...,k\}}$ we write $v(:,k') \in \R^{\cS}$ for $(v(s,k'))_{s \in \cS}$, and similarly we write $v(:,k',a) \in \R^{\cS}$ for $(v(s,k',a))_{s \in \cS}$ given $v \in \R^{\cS \times \{0,...,k\} \times \cA}$. With this notation, the optimality (fixed point) equations for $(v_1,v_2) \in \R^{\cS \times \{0,...,k\}} \times \R^{\cS \times \{0,...,k\} \times \cA}$ are
\begin{align*}
    v_1(s,k') & = \max_{a \in \cA} \sqrt{\gamma} \cdot v_2(s,k',a), \forall \; s \in \cS, k' \in \{0,...,k\},\\
    v_2(s,k',a) & = \min \{ \hP(s,a) \tr (r(s,a) + \sqrt{\gamma} v_1(:,k')), \\
    &\min_{P(s,a) \in \cP_{sa}} P(s,a) \tr (r(s,a) + \sqrt{\gamma} v_1(:,k'+1))\}, \quad \forall \; (s,a) \in \cS \times \cA, k' < k, \\
    v_2(s,k,a) & = \hat{P}(s,a) \tr (r(s,a) + \sqrt{\gamma} v_1(:,k)), \quad \forall \; (s,a) \in \cS \times \cA
\end{align*}
The third line above is the case where the adversary has exhausted all the budget, in which case $k'=k$ and the adversary must choose $\hat{P}$. The second line above is the case where already $k'<k$ deviations have occurred in the past and the adversary can choose to ``consume" one more unit of deviation (choose some $P(s,a)$, eventually with $P(s,a) \neq \hat{P}(s,a)$) or not (stick to the nominal kernel). Note that a pair $(v\opt_1,v\opt_2)$ of solutions to the fixed point equations above exists because the relevant operators are contractions, and in fact on the augmented state set we have a problem similar to Model 1 from \Cref{sec: lec 3 time var P}, so that we can invoke \Cref{th: lec 3 Model 1} and \Cref{prop: lec 3 Model 1} to conclude that there exists a uniformly optimal stationary policy (stationary in the augmented states). Note that this does not require convexity of $\cP$ in the first place, only compactness to invoke results for sa-rectangular robust MDPs. This is because the decisions of the adversary happen {\em after} observing the action $a_t$ picked by the decision maker. In particular, given a pair $(v\opt_1,v\opt_2)$ of solutions to the above equations we get that an optimal policy for the decision maker can be chosen as a function of the current state $s$ and the current number of ``kernel deviations" $k'$ in the $\arg \max$ in the first equation above, and similarly for the adversary. It is worth emphasizing that the resulting policy of the decision maker is {\em not} stationary in the sense defined in \Cref{chap: lec 0 mdp}, since it may in general depend on $k'$, so that when visiting the same state $s$ over time an optimal decision maker may choose different decision rules. However, an optimal policy is stationary in the augmented state $(s,k')$, so in some sense in the LDST model we have a ``summary statistic" that we can use to make our decisions: even though an optimal policy may be non-stationary over $\cS$ we still need only a small amount of information to make our decisions. 

\begin{remark}
    In principle, we should divide the rewards in the stochastic game (SG) by $\sqrt{\gamma}$ so that the return in the LDST robust MDP and the SG are the same; in the construction above there is an additional $\sqrt{\gamma}$ multiplicative factor for the return in the SG compared to that of the robust MDP; of course this does not change the set of optimal policies, so we avoid this rescaling above for the sake of readability.
\end{remark}
\begin{remark}
    The dynamic programming equations count the number of deviations of the vector $P(s,a) \neq \hat{P}(s,a)$ instead of the ``block" $P(s) \neq \hat{P}(s)$. However for LDST the adversary chooses $P(s_t)$ after observing $a_t$ so these two models are the same.
\end{remark}
\subsection{k-rectangular uncertainty set}
The model of {\em k-rectangularity} is a natural extension of the LDST model of uncertainty for the case of {\em finite} horizon robust MDPs. Assume that there are $H$ periods of decisions $\{0,...,H-1\}$ and the sets of states $\cS_0,...,\cS_{H}$ visited at each period are all pairwise disjoint (i.e., a state is visited at most once along a trajectory); the set $\cS_H$ is the set of terminal states. We write $\cS = \cS_0 \cup ... \cup \cS_H$ for the set of states. We are given an uncertainty set $\cP$ for the transition probabilities across states. The dynamics unfold as follows: at each period $t$, in the current state $s_t \in \cS_t$,
\begin{itemize}
    \item the decision maker will pick an action $a_t$ and the adversary will pick a block $P(s_t) \in \DSA$ for the transitions out of $s_t$.
    \item  The decision maker can observe the trajectory $(s_0,a_0,...,s_t)$ up to period $t$ as well as the past choices $P(s_0),...,P(s_{t-1}) \in \DSA$ of the adversary.
    \item  The decision maker first chooses their action $a_t$ based on this information, then the adversary observes $a_t$ and chooses $P(s_t) \in \DSA$. The next state $s_{t+1}$ is sampled from $P(s_t,a_t) \in \DS$ and $P(s_t)$ is revealed to the decision maker.
\end{itemize}
The objective here is the finite-horizon (total) return
\[R(\pi,P) = \E^{\pi,P}\left[ \sum_{t=0}^{H-1} r(s_t,a_t,s_{t+1}) \right] \]
Note that this is close to Model 1 in \Cref{sec: lec 3 time var P} except that the decision maker can observe the realizations of the past transition probabilities (instead of the realized next state).
Note that we can write $P(s_t)$ instead of $P^t(s_t)$ because each state is visited at most once. There are two main departures from the robust MDP model studied in these notes (apart from the finite-horizon assumption, which we can always model as an infinite-horizon problem anyway by adding some terminal absorbing states with rewards $0$):
\begin{itemize}
    \item Along each realized trajectory and at each period $t$ the choice of $P^t(s_t)$ by the adversary must be ``compatible" with at least one element of the original set $\cP$ of admissible transition probabilities.
    \item Along each realized trajectory, the decision maker observes the realizations of the successive adversary choices of the entire $P(s_t) \in \DSA$ at each period $t$ after taking their action $a_t$, instead of solely the realization of $s_{t+1}$ as in the robust MDP problem studied so far.
\end{itemize}
Given the first condition, it makes sense to assume that we visit each state-action pair at most once (this is the assumption that $\cS_0,...,\cS_H$ are pairwise disjoint; otherwise the adversary would be constrained to take the same choice again, to keep compatibility with a kernel in $\cP$).
This first condition naturally leads to the definition of the {\em projection} of the uncertainty set $\cP$ onto a subset of states. Given $\cS' \subset \cS$ and $P \in \cP$ we write $P_{\cS'} = (P(s))_{s \in \cS'}$ for the transition probabilities out of states in $\cS'$. The projection $\cP_{\cS'} \subset (\Delta(\cS))^{\cS' \times \cA}$ of $\cP$ onto $\cS'$ is the set of feasible transition probabilities out of $\cS'$, generated by transition probabilities in $\cP$:
\begin{equation*}
    \cP_{\cS'} = \{P_{\cS'} \; | \; \exists P_{\cS \setminus \cS'} \text{ such that } (P_{\cS'},P_{\cS \setminus \cS'}) \in \cP\}.
\end{equation*}
Now in the model of time-varying adversary for k-rectangular sets, we are interested in counting, given the observed realizations of the transition probabilities out of states in some $\cS_0 \subset \cS$ (think of all visited states so far in the MDP up to the current period), the number of feasible realizations for the transition probabilities out of states in some other $\cS\opt \subset \cS \setminus \cS_0$ (think of the future states visited during the later dynamics). This leads to the definition of {\em conditional} projection onto a set of states, given a realization of the transition probabilities outside of this set of states. In particular, let $P_{\cS_0} \in \cP_{\cS_0}$ be a realization of the parameters in $\cS_0$. Then the conditional projection of $\cP$ onto $\cS\opt$ given $P_{\cS_0}$ is
\begin{equation*}
    \cP_{\cS\opt}(P_{\cS_0}) = \{ P_{\cS\opt} \; | \; \exists P_{\cS \setminus \left( \cS_{0} \cup \cS\opt \right)}, (P_{\cS_0}, P_{\cS\opt},P_{\cS \setminus \left( \cS_{0} \cup \cS\opt \right)}) \in \cP\}
\end{equation*}
This is the set of feasible (future) realizations for the transitions out of $\cS\opt$ given the observed (past) realizations $P_{\cS_0}$ out of $\cS_0$. A natural ``consistency identity", which plays a crucial role in the later state-augmentation technique for solving k-rectangular robust MDPs, is that given $\cS_1,\cS_2,\cS_3 \subset \cS$ pairwise disjoint and $P_1 \in \cP_{\cS_1},P_2 \in \cP_{\cS_2}$ jointly feasible, we have
\begin{equation}\label{eq: lec 3 continuation krec}
    \cP_{\cS_3}((P_1,P_2)) = \left(\cP_{\cS_2 \cup \cS_3}(P_1)\right)_{\cS_3}(P_2)
\end{equation}
The left-hand side $\cP_{\cS_3}((P_1,P_2))$ is just the transition probabilities out of $\cS_3$ compatible with the realization $P_1,P_2$ out of $\cS_1$ and $\cS_2$, while the right-hand side states that we decompose this set as first the conditional projection of $\cP$ onto $\cS_2 \cup \cS_3$ given $P_1$, which we project onto $\cS_3$ given $P_2$ afterwards. We can already see here some form of dynamic programming principle: as information (realization of the current transition probabilities) is revealed over time, we can decompose the set of transition probabilities out of the remaining non-visited states sequentially.

We are now ready to define the k-rectangularity assumption. Given $\cS\opt \subset \cS$ (representing the set of states that we have not visited yet), we define the following collection of sets:
\[\cE_{\cS\opt} = \{\cP_{\cS\opt}(P_{\cS_0}) \; | \; \cS_0 \subset \cS \setminus \cS\opt, P_{\cS_0} \in \cP_{\cS_0}\}\]
The set $\cE_{\cS\opt}$ contains all the conditional projections of $\cP$ induced by the realizations of some transition probabilities outside of $\cS\opt$.
\begin{definition}
    An uncertainty set $\cP$ is {\em k-rectangular} if $|\cE_{\cS\opt}| \leq k$ for every $\cS\opt \subset \cS$.
\end{definition}
This means that for any subset of states $\cS\opt \subset \cS$ for which we do not know the realizations of the transition probabilities, the realizations of the transitions on any combination of states outside of $\cS\opt$ yield at most $k$ distinct sets of transitions out of $\cS\opt$. Note that here for each realization of transitions outside of $\cS\opt$ we may have an infinite set of possible kernels out of $\cS\opt$; it is the number of such possible sets that is finite, not the number of possible kernels. It is worth investigating some examples.
\begin{example}[s-rectangular uncertainty sets]
    Assume that $\cP$ is s-rectangular. For $\cS_0,\cS\opt$ disjoint we have $\cP_{\cS\opt}(P_{\cS_0}) = \times_{s \in \cS\opt} \cP_s$ (independently of the realization of $P_{\cS_0} \in \cP_{\cS_0}$). Hence $|\cE_{\cS\opt}|=1$ and s-rectangular uncertainty sets are k-rectangular with $k=1$.
\end{example}
\begin{example}[LDST uncertainty sets]
Consider an LDST uncertainty set where we have at most $K \in \N$ deviations along any trajectory. Then given some realizations of the transition probabilities at $\cS_0$, we can compute $K_0$ the number of deviations from the nominal $\hat{P}$ for the states in $\cS_0$. The projection is
\[\cP_{\cS\opt}(P_{\cS_0}) = \{P_{\cS\opt} \; | \; \sum_{s \in \cS\opt} 1_{ \{ P(s) \neq \hat{P}(s) \}} \leq K - K_0 ,P(s) \in \cP_s, \forall \;s \in \cS\opt\}\]
hence there are at most $K+1$ possible realizations of the projections. Here LDST uncertainty sets are k-rectangular with $k=K+1$.
\end{example}
There is an interesting algebra for the set of k-rectangular sets for finite values of $k$, which is stable by union, intersection, Cartesian product, and Minkowski averaging; we refer to Section 4 in the original paper on k-rectangularity (see \Cref{sec: lec 3 ref}) for more details on this.
\paragraph{State-augmentation technique.} We now describe how we can compute an optimal policy for the decision maker for a k-rectangular uncertainty set. Assume that we are at a given period $t \in \N$ and that we have visited the states $s_0,...,s_{t-1}$ so far, and that the decision maker knows the corresponding realizations of the transition probabilities $P(s_0),...,P(s_{t-1})$. The set of all possible future states is $\cS_{t:H} \coloneqq \cS_t \cup ... \cup \cS_H$, and the uncertainty set of all feasible future transition probabilities compatible with what we observed so far is $\cP_{\cS_{t:H}}(P_{s_0},...,P_{s_{t-1}}).$ By definition $\cP_{\cS_{t:H}}(P_{s_0},...,P_{s_{t-1}}) \in \cE_{\cS_{t:H}}$ and from the k-rectangularity assumption we know that $|\cE_{\cS_{t:H}}| \leq k$. We will write $x_t$ for a general element of $\cE_{\cS_{t:H}}$.
Once the decision maker has chosen $a_t$ and the adversary has chosen $P_{s_t} \in \cP_{\cS_{t:H}}(P_{s_0},...,P_{s_{t-1}})_{\{s_t\}}$, the set of feasible future transition probabilities compatible with what we observed so far becomes 
$(\cP_{\cS_{t:H}}(P_{s_0},...,P_{s_{t-1}}))_{\cS_{t+1:H}}(P(s_t))$, which, given \eqref{eq: lec 3 continuation krec}, is equal to $\cP_{\cS_{t+1:H}}(P_{s_0},...,P_{s_t})$. Therefore, once $x_t$ is known, the distribution of the next state only depends on the current state-action $(s_t,a_t)$, the current $x_t$, and the choice of transition probabilities at period $t$; we do not need to know the rest of the history up to period $t$. 
This means that we can reformulate the robust MDP problem over a compact k-rectangular uncertainty set with state set $\cS$ as a robust MDP problem with a compact {\em sa-rectangular} uncertainty set over the augmented set of states $\hat{\cS} = \hat{\cS}_{0} \cup ... \cup \hat{\cS}_{H}$ with $\hat{\cS}_t = \cS_t \times \cE_{\cS_{t:H}}$ for each period $t$, with time-varying adversaries, for which we can apply the dynamic programming equations and reductions outlined in \Cref{sec: lec 3 time var P} for $\cPSA_{\sf H}$ (which can be adapted to deal with finite-horizon problems by adding absorbing states with rewards $0$ and adjusting the instantaneous rewards at each period to account for discounting). The k-rectangularity assumption implies that the augmented set of states has finite cardinality, since $|\hat{\cS}_{t}| = |\cS_t| \times |\cE_{\cS_{t:H}}| \leq k |\cS_t|$ so that $\hat{\cS} \leq k | \cS|$. More precisely, the lifted set of states at period $t$ is $\cS_t \times \cE_{\cS_{t:H}}$, the set of actions available at this augmented state is $\cA$, the transition probability $(s,x) \to (s',x')$ when choosing action $a$ at period $t$ is $P(s,a,s')1_{\{ x' = x_{\cS_{t+1:H}}(P(s))\}}$ and the lifted instantaneous rewards are $r((s,x),a,(s',x')) = r(s,a,s')$. Since the adversary chooses their decisions after observing the action chosen by the decision maker, the dynamic programming equation for this model is, for $t \in \{0,...,H-1\}$ and, for $s \in \cS_t$ and $x \in \cE_{\cS_{t:H}}$,
\begin{equation*}
    v(s,x) = \max_{a \in \cA} \inf_{P(s) \in x_{\{s\}}} \sum_{s' \in \cS_{t+1}} P(s,a,s')\left( r(s,a,s') + v\left(s',x_{S_{t+1:H}}(P(s))\right) \right)
\end{equation*}
and $v(s,x) = 0, \forall \; s \in \cS_{H}$.
For instance in the LDST problem we have seen that $x_t$ can be mapped to the number of past deviations from the nominal kernels.
\paragraph{Remarks on k-rectangular models.}
It is important to recall again that the information model of k-rectangularity assumes that the decision maker can observe the realizations of the transition probabilities along the trajectories, which is radically different from the other models of robust MDPs studied in these notes. In particular if the decision maker only observes the current state $s_t$ then they can't build the variable $x_t$ (label of the current conditional projection of the feasible remaining kernels) and cannot play optimally. This is an important limitation of the k-rectangular model of uncertainty, along with the fact that some uncertainty sets are k-rectangular only for very large $k$ (e.g. a union of $J \in \N$ s-rectangular sets is k-rectangular with $k$ at most $2^J -1$ which is very large), so that the resulting number of augmented states may be very large. Additionally, it is unclear how one can extend this modeling assumption to the case of infinite-horizon problems, where we cannot assume that each state is visited once (otherwise we have an infinite number of states, so the methods developed in these notes do not apply).
\section{Risk-averse MDPs}\label{sec: lec 3 risk averse}
The study of risk-averse MDPs is a large and active field of research; as such we are by no means trying to be exhaustive in this section and our only goal is to highlight the well-known connection between risk-averse MDPs and robust MDPs. We will use a tilde symbol to denote random variables in this section.
\paragraph{Preliminaries on risk measures.} We introduce here the definition of risk measures and coherent risk measures. Risk measures are real functions of random variables, which we can interpret as assigning a scalar ``score" to a random variable.
\begin{definition}
    Let $\left(\Omega,\cF,\bP\right)$ be a probability space. We write $\cL(\Omega,\cF,\bP;\R)$ for the set of essentially bounded random variables, defined as the set of $\cF$-measurable functions $\Omega \to \R$ that are bounded outside a set of measure $0$ (for $\bP$). A {\em risk measure} is a map from $\cL(\Omega,\cF,\bP;\R)$ to $\R$.
\end{definition}
The reader should already be familiar with the most classical examples such as the expectation: $\tX \mapsto \E_{\bP}\left[ \tX \right]$ with $\E_{\bP}\left[ \tX \right] = \int_{\omega \in \Omega} \tX(\omega)d\bP(\omega)$, the variance $\tX \mapsto 
\E_{\bP} \left[ \left( \tX - \E_{\bP} \left( \tX\right) \right)^2  \right] $, the Value-at-Risk, the Conditional-Value-at-Risk, etc.
A risk measure $\rho$ is said to be:
\begin{itemize}
    \item {\em monotone} if for any random variables $\tX,\tY$, $\tX \leq \tY$ a.s. implies $\rho(\tX) \leq \rho(\tY)$;
    \item {\em translation-invariant} if $\rho(\tX + c) = \rho(\tX)+c$ for any $c \in \R$ and for any random variable $\tX$;
    \item {\em concave} if it is monotone, translation invariant, and additionally for any random variables $\tX,\tY$ and scalar $\lambda \in [0,1]$ we have $\rho(\lambda \tX + (1-\lambda)\tY) \geq \lambda \rho(\tX) + (1-\lambda)\rho(\tY)$;
    \item {\em coherent} if it is monotone, translation-invariant, concave, and additionally it is positively homogeneous: $\rho(\alpha \tX) = \alpha \rho(\tX)$ for any $\alpha \geq 0$ and any random variable $\tX$.
\end{itemize}
The reader can check that the expectation is coherent but not the variance; the essential infimum $\essinf(\tX) = \sup \{ b \in \R \; | \; \bP(\tX < b) = 0\}$ is coherent. It is worth emphasizing that we are focusing on random variables that model the return (that we want to maximize), whereas a large part of the literature on risk measures focuses on random variables modeling losses. So, for instance, we focus on concave risk measures instead of convex risk measures. Similarly, our definition of Conditional-Value-at-Risk will focus on the smallest (bottom) $\alpha$-fraction of the return instead of the largest (top) $\alpha$-fraction of the return (we would use the later if we were focusing on losses or costs, i.e., on quantities we want to minimize).

The following theorem, known as the {\em dual representation of coherent and concave risk measures}, provides the fundamental tools toward connecting risk-averse MDPs and robust MDPs.
\begin{theorem}\label{th: lec 3 dual representation risk averse}
    Let $(\Omega,\cF,\bP)$ be a probability space with $\Omega$ finite: $\Omega = \{1,...,m\}$, $\cF = 2^{\Omega}$, and let $\rho$ be a risk measure.
    \begin{enumerate}
        \item If $\rho$ is a concave risk measure, there exists $\alpha:\Delta(\Omega) \to \R\cup \{+\infty\}$ such that        
        \begin{equation}\label{eq: lec 3 dual rep concave}
            \rho(\tX) = \min_{\bP' \in \Delta(\Omega)} \E_{\bP'}\left[ \tX \right] + \alpha(\bP')
        \end{equation}
        \item If $\rho$ is a coherent risk measure, there exists a nonempty compact convex subset $\cQ^{\rho} \subseteq \Delta(\Omega)$ of measures over $\Omega$, called the {\em risk envelope}, such that
        \begin{equation}\label{eq: lec 3 dual rep coherent}
            \rho(\tX) = \min_{\bP' \in \cQ^{\rho}} \E_{\bP'} \left[ \tX \right]
        \end{equation}
    \end{enumerate}
\end{theorem}
For two distributions $\bP,\bP'$ in $\Delta(\Omega)$ we note $\bP' \ll \bP$ if $\bP'$ is absolutely continuous with respect to $\bP$, i.e., if for any $\omega \in \Omega$ we have $\bP(\omega)=0 \Rightarrow \bP'(\omega)=0$.
It is worth studying in more detail the dual representations of some classical risk measures:
\begin{itemize}
    \item For a given $\lambda \in (0,1]$, $\cvar_{\lambda}(\tX)$ is defined as the expectation of the worst $\lambda$-fraction of the random variable. The risk envelope theorem gives
    \[\cvar_{\lambda}(\tX) = \inf \{ \E_{\bP'}\left[\tX\right] \; | \; \bP' \in \Delta(\Omega), 0 \leq \bP'(\omega) \leq \lambda^{-1} \bP(\omega), \forall \; \omega \in \Omega\}.\]
    \item For a given $\lambda>0$, the {\em entropic risk measure} (\erm) is defined as $\erm_{\lambda}(\tX) = - \frac{1}{\lambda}\log\left( \E_{\bP}\left[\exp(-\lambda \tX)\right] \right).$ We already encountered in the exercises from \Cref{sec: lec 2 fdt} the Donsker-Varadhan variational formula:
    \[\erm_{\lambda}(\tX) = \inf \{\E_{\bP'}\left[\tX\right] + \lambda^{-1} \kl(\bP',\bP) \; | \; \bP' \in \Delta(\Omega), \bP' \ll \bP \}\]
    with the Kullback-Leibler divergence defined as
    \[\kl(\bP',\bP) = \sum_{\omega \in \Omega} \bP'(\omega) \log (\bP'(\omega)/\bP(\omega))\]
    under the convention that $0 \log(0) = 0$.
    \item For $\lambda \in (0,1]$, we define the {\em entropic value at risk} (\evar) as $\evar_{\lambda}(\tX) = \sup_{\theta>0} \{\erm_{\theta}(\tX) + \theta^{-1}\log(\lambda)\}.$ The dual representation of \evar{} is
    \[\evar_{\lambda}(\tX) = \inf \{\E_{\bP'}\left[\tX \right] \; | \; \bP' \in \Delta(\Omega), \bP' \ll \bP, \kl(\bP',\bP) \leq - \log(\lambda) \}.\]
\end{itemize}
At this point \eqref{eq: lec 3 dual rep coherent} already suggests some connections with robust MDPs, where $\cQ^{\rho}$ could play the role of an uncertainty set. Indeed the worst-case Bellman operator can be understood as a worst-case expectation for the continuation values (see \Cref{sec: lec 2 adv}). However, $\cQ^{\rho}$ is not a set of transition probabilities (we didn't define any notions of dynamic processes anyway so far in this section), it is a set of measures over the universe $\Omega$. 
Later in this section we will see that we can convert $\cQ^{\rho}$ into an uncertainty set in some specific cases. We also note that \eqref{eq: lec 3 dual rep concave} suggests an alternative model to robust MDPs, where instead of having hard constraints on the transition probabilities the adversary instead incurs a ``cost" $\alpha(\bP')$ when choosing $\bP'$. This model is sometimes called ``soft-robust MDPs" (although there are multiple meanings of this term, which is also used to denote multi-model MDPs).
\paragraph{Finite-horizon robust MDPs: static risk measure.} Consider a finite-horizon MDP where the dynamics end at period $H$ and the discount factor is $\gamma=1$. The set of states is $\cS$. We write $\PiMD$ for the set of Markovian, deterministic policies. Consider a fixed Markovian {\em deterministic} policy $\pi=(\pi_0,...,\pi_{H-1}) \in \PiMD$ and a {\em deterministic} initial distribution $p_0$ (both of the assumptions can be removed but they simplify the later calculations). 
Assume that the nominal transition kernel is $P^0 \in \DSSA$. Given $\pi,P^0$, if we are in state $s_t$ at period $t$ the next state $\ts_{t+1}$ is random with a probability distribution $P^0(s_t,\pi_t(s_t))$. Hence the total return $\tR(\pi,P^0)$ is a random variable, with
\begin{equation}\label{eq: lec 3 random return}
    \tR(\pi,P^0) = \sum_{t=0}^{H-1} r(\ts_t,\pi_t(\ts_t),\ts_{t+1})
\end{equation}
where $(\ts_t)_{t \leq H}$ is a random sequence of states. This random variable is defined for the probability space with universe $\Omega = \cS^{H+1}$ corresponding to all possible realizations of the random sequence of states visited at period $t=0,t=1,...,t=H$, with the natural filtration 
\[\cF_t = \sigma(\ts_0,...,\ts_t) \subseteq 2^{\Omega}\]
for $t=0,...,H$, and the nominal probability measure $\bP^{\pi,P^0}$ on $\Omega$ is 
\[\bP^{\pi,P^0}(s_0,...,s_H) = p_0(s_0)\prod_{t=0}^{H-1} P^0(s_t,\pi_t(s_t),s_{t+1}).\]
In risk-averse MDPs we are interested in evaluating $\rho(\tilde{R}(\pi,P^0))$ for some risk measure $\rho$. Note that for $\rho=\E_{\bP^{\pi,P^0}}$ this is just a nominal MDP where the objective function $\E_{\bP^{\pi,P^0}}\left[\tilde{R}(\pi,P^0)\right]$ coincides with our original definition of the expected return as in \Cref{sec: lec 0 model}. We also want to optimize $\pi \mapsto \rho(\tilde{R}(\pi,P^0))$ but this requires more formal definitions (since we need to define $\tilde{R}(\pi,P^0)$ and the probability space for history-dependent randomized policies) that are not relevant for the connection with robust MDPs.
Given a coherent risk measure on $(\Omega,\cF,\bP^{\pi,P^0})$, the risk envelope $\cQ^{\rho}$ is a set of distributions over $\Omega$ (which we can identify with the set of trajectories), and $\rho(\tilde{R}(\pi,P^0))$ is such that
\begin{equation}\label{eq: lec 3 raMDP 0}
    \rho(\tilde{R}(\pi,P^0)) = \inf_{\bP' \in \cQ^{\rho}} \E_{\bP'}\left[ \sum_{t=0}^{H-1} r(\ts_t,\pi_t(\ts_t),\ts_{t+1})\right].
\end{equation}
Note that $\cQ^{\rho}$ is a set of probability measures over trajectories, so we cannot interpret the right-hand side in the above equation as a policy evaluation problem in a robust MDP. To connect risk-averse MDPs and robust MDPs we want to reformulate the right-hand side of the above equation as a minimization over {\em transition probabilities} (possibly history-dependent). This can be done with the chain rule, as explained in the next proposition.
\begin{proposition}\label{prop: lec 3 disintegration}
    \begin{enumerate}
        \item Let $\bP'$ be a probability law on $\Omega = \cS^{H+1}$ with initial marginal $p_0$. For any $h_t=(s_0,...,s_t) \in \cS^{t+1}$ with $\bP'(h_t)>0$ (in particular we must have $s_0$ equal to the state chosen deterministically by the initial distribution $p_0$), define $P(h_t) \in \Delta(\cS)$ as $P(h_t,s') = \bP'(h_t,s')/\bP'(h_t)$ for any $s' \in \cS$. Then for any $(s_0,...,s_H) \in \Omega$ we have
    \begin{equation}\label{eq: lec 3 disintegration}
            \bP'(s_0,...,s_H) = p_0(s_0)\prod_{t=0}^{H-1} P(h_t,s_{t+1}).
        \end{equation}
        \item Conversely, for any collection of $P(h_t) \in \Delta(\cS)$ defined for each $h_t \in \cS^{t+1}$ for each $t < H$, we can define a probability measure $\bP'$ on $\Omega = \cS^{H+1}$ by the formula \eqref{eq: lec 3 disintegration}.
    \end{enumerate}
\end{proposition}
\begin{proof}
    Both statements are just plain rewriting of the definitions of $P(h_t)$ as $P(h_t,s') = \bP'(h_t,s')/\bP'(h_t)$.
\end{proof}
Note that in general we need history-dependent transition probabilities, i.e., $P(h_t) \in \Delta(\cS)$ for the distribution over the state at period $t+1$ after observing history $h_t = (s_0,...,s_t)$, and not necessarily stationary transition probabilities. We show how this works for \cvar{} and \evar{} below:
\begin{enumerate}
    \item For \cvar, using the dual representation theorem (\Cref{th: lec 3 dual representation risk averse}) to obtain \eqref{eq: lec 3 raMDP 0} and then \Cref{prop: lec 3 disintegration} to convert a minimization problem over measures in $\Delta(\cS^{H+1})$ into transition probabilities, we get that $\cvar_{\lambda}(\tilde{R}(\pi,P^0))$ is equal to
    \begin{align*}
         \inf \; & \E^{\pi,P}\left[\sum_{t=0}^{H-1} r(\ts_t,\pi_t(\ts_t),\ts_{t+1})\right] \\
        & \prod_{t=0}^{H-1} \frac{P(s_0,...,s_t,s_{t+1})}{P^0(s_t,\pi_t(s_t),s_{t+1})} \leq  \lambda^{-1}, \forall \; (s_t)_{t\leq H} \in \cS^{H+1}, \\
        & P(s_0,...,s_t) \in \Delta(\cS), P(s_0,...,s_t) \ll P^0(s_t,\pi_t(s_t)), \forall \; (s_t)_{t\leq H} \in \cS^{H+1}, t \leq H
    \end{align*}
    where $\E^{\pi,P}$ is the expectation with respect to the probability measure over $\cS^{H+1}$ defined by $\pi \in \PiMD$ and $P$ satisfying the constraints in the above optimization program.
    This is a policy evaluation problem (in the finite-horizon setting) with an uncertainty set characterized by the set of constraints in the optimization problem above. The resulting uncertainty set has history-dependent transition probabilities and is non-rectangular given we have constraints across the transition probabilities out of different states (linked by the product $\prod_{t=0}^{H-1}$) and because the constraints in this uncertainty set are trajectory-wise (instead of state-wise). Note that we already encountered robust MDPs with time-varying transition probabilities in \Cref{sec: lec 3 time var P}, but there we assumed that the transition probabilities chosen at period $t$ out of $s_t$ were always constrained to belong to the same set $\cP_{s_t}$, we have constraints over the transition probabilities at different periods and when visiting the same state over time the sets of constraints for the transition probabilities out of this state change, hence the results from \Cref{sec: lec 3 time var P} do not apply directly.
    \item For \evar, using the dual representation theorem (\Cref{th: lec 3 dual representation risk averse}) we get that $\evar_{\lambda}(\tilde{R}(\pi,P^0))$ is equal to 
    \begin{align*}
         \inf \; & \E^{\pi,P}\left[\sum_{t=0}^{H-1} r(\ts_t,\pi_t(\ts_t),\ts_{t+1})\right] \\
        & \sum_{t=0}^{H-1} \E^{\pi,P} \left[\kl\left(P(\ts_0,...,\ts_t),P^0(\ts_t,\pi_t(\ts_t)) \right) \right] \leq  -\log(\lambda), \\
        & P(s_0,...,s_t) \in \Delta(\cS), P(s_0,...,s_t) \ll P^0(s_t,\pi_t(s_t)), \forall \; (s_{t'})_{t' \leq H} \in  \cS^{H+1}, t \leq H
    \end{align*}
    This again yields an uncertainty set with history-dependent, non-rectangular transition probabilities.
\end{enumerate}
So we have seen that in general risk-averse MDPs with a {\em coherent} risk measure will induce a particular type of robust MDPs where the transition probabilities are allowed to be history-dependent and constrained in a ``nasty" way: there may be constraints across different periods and constraints across the transitions out of different states, which typically yields difficult optimization problems, as we have seen in the previous chapters.
For a {\em concave} risk measure $\rho$ we obtain an optimization problem sometimes called {\em soft-robust MDPs}, where there is a ``soft" penalization term $\alpha(\bP) \in \R \cup \{+\infty\}$ in the objective:
\[\rho(\tilde{R}(\pi,P^0)) = \inf_{\bP' \in \Delta(\Omega)} \E_{\bP'}\left[ \sum_{t=0}^{H-1} r(\ts_t,\pi_t(\ts_t),\ts_{t+1})\right] + \alpha\left( \bP' \right).\]
In principle the term $\alpha(\bP')$ need not decompose period-wise, so that evaluating $\rho(\tilde{R}(\pi,P^0))$ may be difficult, let alone optimizing it over $\pi$; \erm\; is one particular example where the penalty term $\alpha$ decomposes across periods (it is a KL divergence). Note that hard constraints can be encoded this way by taking $\alpha(\bP') = + \infty$ for some forbidden measures $\bP'$ outside of some prescribed set $\cQ^{\rho}$.
\paragraph{Finite-horizon robust MDPs: conditional and Markov risk measures.} We have seen in the examples above that coherent risk measures can be related to robust MDPs with non-rectangular uncertainty sets, which are in general difficult to optimize over. A natural idea to recover robust MDPs with rectangular uncertainty sets from risk-averse MDPs is the following: instead of applying the risk measure $\rho$ to the return of the entire trajectory, we apply $\rho$ inductively, e.g., informally we consider the following objective (for $s_0 \in \cS$ the state chosen deterministically by $p_0$):
\begin{equation}\label{eq: lec 3 dynamic return informal}
    \rho\left(r(s_0,\ta_0,\ts_1) + \rho(r(\ts_1,\ta_1,\ts_2) + \rho( ... + \rho(r(\ts_{H-1},\ta_{H-1},\ts_{H}))))\right)
\end{equation}
This construction is known as nested {\em conditional} risk measures (this choice of terminology comes from the fact that the risk measure $\rho$ is applied at every period as information about the current state is revealed, rather than to the entire trajectory).
More formally, in the inner evaluation above each risk measure is actually a {\em conditional} risk measure $\rho(\cdot | \cF_t): \cL(\Omega,\cF_{t+1},\bP^{\pi,P^0}; \R) \to \cL(\Omega,\cF_{t},\bP^{\pi,P^0}; \R)$ in the sense of Ruszczy{\'n}ski–Shapiro (see references at the end of this chapter): for an $\cF_{t+1}$-measurable random variable $\tX$, we get that $\rho(\tX | \cF_t)$ is an $\cF_t$-measurable random variable, i.e. its universe is the set of finite histories up to period $t$. Hence in \eqref{eq: lec 3 dynamic return informal} the outer $\rho$ is actually $\rho(\cdot | \cF_0)$, the first inner $\rho$ is actually $\rho(\cdot | \cF_1)$, etc.
To construct a family of {\em coherent} conditional risk measures in our context, following the dual representation theorem (\Cref{th: lec 3 dual representation risk averse}) we can simply choose, for each period $t$ and history $h_t=(s_0,a_0,...,s_{t})$ a history-dependent risk envelope 
\[\cQ^{\rho}(h_t) \subseteq \{ p \in \DS \; | \; p \ll P^0(s_t,\pi_t(s_t))\}\]
and then define, for a given $\cF_{t+1}$-measurable random variable $\tX$, 
\begin{equation}\label{eq: def conditional risk measures}
    \left(\rho(\tX|\cF_t)\right)(h_t) \coloneqq \inf \left\{ \E_{\ts \sim P(h_t)}\left[ \tX(h_t,\ts)\right] \; \middle| \; P(h_t) \in \cQ^{\rho}(h_t) \right\}.
\end{equation}
With this construction, the objective \eqref{eq: lec 3 dynamic return informal} is well-defined, and the monotonicity of the conditional risk measures ensures that {\em time-consistency}, i.e., dynamic programming, holds. We can additionally assume that the risk measures satisfy the following Markov property, that the risk envelopes $\cQ^{\rho}(h_t)$ (with $h_t = (s_0,a_0,...,s_t)$) only depend on the last state $s_t$ and the deterministic choice of action $\pi_t(s_t)$, i.e., that there exists some $\cP^{\rho}_{s_t \pi_t(s_t)} \subseteq \DS$ (say compact) such that $\cQ^{\rho}(h_t) = \cP^{\rho}_{s_t \pi_t(s_t)}$. 
For instance, we get the following constraints on the corresponding uncertainty $\cP^{\rho}$ for conditional \cvar{} and conditional \evar:
\begin{align*}
    \cP^{\cvar}_{sa} & = \{ p \in \DS \; | \; 0 \leq p(s') \leq \lambda^{-1} P^0(s,a,s'), \forall \; s' \in \cS\} \\
    \cP^{\evar}_{sa} & = \{ p \in \DS \; | \; \kl(p,P^0(s,a)) \leq - \log(\lambda)\}
\end{align*}
Applying this in \eqref{eq: def conditional risk measures}, the reader should now recognize something very close to an sa-rectangular robust MDP formulation. Indeed, we now get that
\[\left(\rho(\tX|\cF_t\right)(h_t) \coloneqq \min \left\{ \E_{\ts \sim p}\left[ \tX(h_t,\ts)\right] \; \middle| \; p \in \cP^{\rho}_{s_t \pi_t(s_t)} \right\}.\]
In particular for $v \in \R^{\cS}$ we have
\[\rho\left(r(\ts_t,\pi_t(\ts_t),\ts_{t+1}) + v(\ts_{t+1}) \; | \; \cF_t\right) = T_t^{\pi,\rho}v(\ts_t)\]
for the operator $T_t^{\pi,\rho}:\R^{\cS} \to \R^{\cS}$ as
\begin{equation*}
    T_t^{\pi,\rho}v(s) \coloneqq \min_{p \in \cP^{\rho}_{s \pi_t(s)}} p \tr \left( r(s,\pi_t(s)) + v \right), \forall \; v \in \RS, s \in \cS.
\end{equation*}
With these definitions, we finally obtain
\begin{equation*}
    \left(T_0^{\pi,\rho} \circ T_1^{\pi,\rho} \circ ... \circ T_{H-1}^{\pi,\rho}\right)0_{\cS}(s_0) = \rho\left(r(s_0,\ta_0,\ts_1) + \rho(r(\ts_1,\ta_1,\ts_2) + \rho( ... + \rho(r(\ts_{H-1},\ta_{H-1},\ts_{H}))))\right)
\end{equation*}
with $0_{\cS} = (0,...,0) \in \R^{\cS}$; this can be compared with the same kind of expressions we obtained for nominal objectives right before \Cref{def: lec 0 monotone contraction}. In particular, the left-hand side above can be interpreted as the (total, undiscounted) return of a policy $\pi \in \PiMD$ obtained after $H$ periods of decisions in a finite-horizon robust MDP where the transitions are chosen adversarially at each period in an sa-rectangular uncertainty set; this is akin to Model 1 from \Cref{sec: lec 3 time var P}, and we know that, in the case of {\em finite} horizon, an optimal policy for this model can be chosen Markovian deterministic (see \Cref{th: lec 3 Model 1}).

One of the main issues of the above approach is that the nested objective \eqref{eq: lec 3 dynamic return informal} need not be law-invariant, in that the returned scalar value may depend on the timing at which the rewards are obtained (and not solely on the distribution of the random return). In particular \eqref{eq: lec 3 dynamic return informal} may differ from $\rho\left(\sum_{t=0}^{H-1} r(\ts_t,\ta_t,\ts_{t+1})\right)$. 
It is known that, under fairly general assumptions (see references at the end of this chapter), the only conditional Markov law-invariant risk measure is the entropic risk measure (extended to $\lambda = + \infty$ and $\lambda=0$ by continuity), and if additionally we assume that the corresponding law-invariant risk measure must be coherent, only the expectation and the essential infimum satisfy this time-consistency property. So while conditional risk measures ``restore" the dynamic programming property for risk-averse MDPs, in some sense the fact that almost all of the conditional risk measures are not law-invariant induces some strange estimation properties for the return.
\paragraph{Infinite-horizon problems.} 
We now explain how to extend these formulations for the case of infinite-horizon problems. We  assume $\gamma<1$ and we still have $\pi \in \PiMD$.
In this case the random return is $\tR(\pi,P^0)=\sum_{t=0}^{\infty} \gamma^t r(\ts_t,\pi_t(\ts_t),\ts_{t+1})$ and the universe is the set of all infinite trajectories $\cS^{\N}$, which is of course not finite anymore. Static coherent risk measures still admit some dual representations (if we additionally assume the Fatou property) but the resulting uncertainty set reformulations are harder to describe:
\begin{enumerate}
    \item For \cvar{} we get the following uncertainty set of history-dependent transition probabilities:
    \begin{align*}
        \cP^{\cvar}_{\sf H} = \{ P \; | \; & P(s_0,...,s_t) \in \DS, P(s_0,...,s_t) \ll P^0(s_t,\pi_t(s_t)),  \forall (s_0,...,s_t) \in \cS^{t+1}, \forall\; t \in \N, \\
        & \prod_{t=0}^{H-1} \frac{P(s_0,...,s_t,s_{t+1})}{P^0(s_t,\pi_t(s_t),s_{t+1})} \leq \lambda^{-1}, \forall \; H \geq 1, \forall h_H \in \cS^{H+1}\}. 
    \end{align*}
    \item For \evar{} we obtain
    \begin{align*}
        \cP^{\evar}_{\sf H} = \{ P \; | \; & P(s_0,...,s_t) \in \DS, P(s_0,...,s_t) \ll P^0(s_t,\pi_t(s_t)),  \forall (s_0,...,s_t) \in \cS^{t+1}, \forall\; t \in \N, \\
        & \sum_{t=0}^{+\infty} \E^{\pi,P} \left[  \kl\left(P(\ts_0,...,\ts_t),P^0(\ts_t,\pi_t(\ts_t))\right) \right] \leq -\log(\lambda)\}.
    \end{align*}
\end{enumerate}
Hence for static \cvar{} and static \evar{} we obtain reformulations of risk-averse MDPs as robust MDPs with history-dependent transition probabilities that have constraints across both periods and states; arguably this should constitute difficult optimization problems.
In contrast in the case of conditional Markov risk measures we obtain a very simple reformulation: for $\pi \in \PiSD$ the corresponding Bellman operator becomes $T^{\pi,\rho}:\RS \to \RS$ defined as
\[T^{\pi,\rho}v(s) = \rho(r(s,\pi(s),\ts') + \gamma v(\ts'))\]
where $\ts' \sim P^0(s,\pi(s))$; note that in the equations above we slightly abuse notation since we note that the random variable $\ts'$ has different laws $P^0(s,a) \in \DS$ for different state-action pairs $s,a$, hence in principle we should define one risk measure for each state-action pair; this is not necessary if we use law-invariant risk measures.
The operator $T^{\pi,\rho}$ has several important properties:
\begin{itemize}
    \item $T^{\pi,\rho}$ is a monotone contraction on $\RS$ as long as $\rho$ is monotone and translation-invariant (this is related to Crandall-Tartar theorem as in \Cref{sec: lec 0 fdt}), hence it admits a unique fixed point $v^{\pi,\rho} \in \RS$.
    \item If $\rho$ is a coherent risk measure then we additionally have
    \[T^{\pi,\rho}v(s) = \min_{P(s,\pi(s)) \in \cP_{s\pi(s)}} P(s,\pi(s))\tr (r(s,\pi(s)) + \gamma v) \]
    for some nonempty convex compact set $\cP_{s\pi(s)}$, hence we again obtain that risk-averse MDPs with Markovian coherent risk measures are equivalent to sa-rectangular robust MDPs.
\end{itemize}
As a short conclusion we note that all the analysis above can be made for history-dependent randomized policies $\pi \in \PiH$ and randomized initial distribution $p_0 \in \Delta(\cS)$ (instead of $\pi \in \PiMD$ and deterministic initial distributions). The main changes are that the trajectory law of the next state also depends on the random variables $\ta_t$ for the actions chosen following $\pi$ and depends on $p_0$, and therefore the corresponding uncertainty sets become policy-dependent, initial-distribution dependent, and the action probabilities also change, which departs even more from our model of robust MDPs with unknown stationary transition probabilities. Also we only focused on the ``policy evaluation problem" for risk-averse MDPs (i.e., estimating $\rho(\tilde{R}(\pi,P^0))$) and not on the policy optimization problem (i.e., optimizing $\pi \mapsto \rho(\tilde{R}(\pi,P^0))$); we point the reader to \Cref{sec: lec 3 ref} for references on these topics.
\section{Assumptions on the instantaneous rewards*}\label{sec: lec 3 discussion rewards}
I suggest reading \Cref{sec: lec 0 additional prop} and \Cref{sec: lec 3 uncertain rewards and transitions} before diving into this topic. We have encountered a few cases where things are considerably harder for general rewards than for rewards that are independent of the next state:
\begin{itemize}
    \item In \Cref{sec: lec 2 minimality rec}, the minimality of s-rectangularity and sa-rectangularity for dynamic programming and the existence of uniformly optimal policies crucially relies on general reward models. If we instead restrict ourselves to next state-independent rewards then the so-called {\em r-rectangular} model of uncertainty and its generalization also satisfy these desirable properties (see \Cref{sec: lec 3 r-rec}). In turn, r-rectangular models are NP-hard to solve for general models of rewards (\Cref{prop: lec 3 hardness r rec}).
    \item For uncertain rewards {\em and} uncertain transitions, general rewards make the policy evaluation problem NP-hard even for convex compact s-rectangular uncertainty sets (\Cref{prop: lec 3 hardness uncertain r p}) while under the same assumptions on the uncertainty set, next-state independent rewards make both policy evaluation and computing an optimal policy tractable (\Cref{th: lec 3 main rmdps r P}).
\end{itemize}
This may seem a bit puzzling, since we saw in \Cref{sec: lec 0 additional prop} that for {\em nominal MDPs} we could have assumed that the rewards were next-state independent, without loss of generality. This fact was proved with two different techniques: state-augmentation (preserving the expression of the discounted return) or replacing rewards by expected rewards (and noticing the Bellman operator $T^P$ didn't change). We explain here in more detail when these techniques may work or fail:
\begin{itemize}
    \item Assume first that the rewards are fixed.

    For $\cP$ s-rectangular convex the state-augmentation technique yields an s-rectangular convex uncertainty set over the augmented state set, where rewards are dependent only on the current state, leading to a tractable problem (it falls under the umbrella of models studied in \Cref{chap: lec 2 rectangular models}). The ``expected reward" technique yields a robust MDP with s-rectangular uncertainty set $\cUf$ over the uncertain rewards and uncertain transitions (expected rewards are uncertain because they involve $P \in \cP$), and $\cUf$ is convex; this falls under the umbrella of tractable problems described in \Cref{sec: lec 3 uncertain rewards and transitions}. Therefore when the rewards are fixed and $\cP$ is s-rectangular convex we can safely assume that the rewards only depend on the current state.

    \item Now assume that the rewards are not fixed but belong to some $\cR \subset \RSAS$.

    For the ``expected reward" technique we've seen that if $\cP$ and $\cR$ are s-rectangular convex the corresponding $\cUf$ is s-rectangular but it may not be convex, so this method may fail in this case (it may not lead to tractable reformulations). Note that from \Cref{sec: lec 3 uncertain rewards and transitions} this method works when $\cP$ is sa-rectangular compact and $\cR$ is sa-rectangular compact (since we don't need convexity of the resulting $\cUf$ in this case).

    In contrast the state-augmentation may fail even for sa-rectangular sets $\cP$ and $\cR$. This is because the rewards out of the different fictitious states $(s,a,s')$ can be ``correlated" in $\cR_{sa}$ so that the corresponding set of uncertain rewards in the augmented state robust MDPs is not sa-rectangular. For this technique to work a sufficient condition is $\cR = \times_{s,a,s'} \cR_{sas'}$ with $\cR_{sas'} \subset \R$, which appears to be a very strong assumption (the adversary can choose the rewards independently across all $(s,a,s')$).
\end{itemize}
This is summarized in \Cref{tab:reward-reformulations}.
\begin{table}[htp]
    \centering
    \caption{Tractability of the expected-reward and state-augmentation
    reformulations under different rectangularity assumptions. A checkmark ($\checkmark$) means that the resulting reformulation leads to a convex compact s-rectangular robust MDP or a compact sa-rectangular robust MDP, and a cross ($\times$) means the contrary (the reformulation need not lead to a convex compact s-rectangular robust MDP). This table
    assumes $\cU = \cR \times \cP$ (using the notation of
    \Cref{sec: lec 3 uncertain rewards and transitions}), with $\cR,\cP$ compact; $\cP$ must be convex for the first column, whereas $\cP$ need not be convex for the second column.}
    \label{tab:reward-reformulations}

    \footnotesize
    \renewcommand{\arraystretch}{1.32}
    \setlength{\tabcolsep}{4pt}

    \begin{tabularx}{0.78\textwidth}{
        @{}
        >{\raggedright\arraybackslash}p{0.19\textwidth}
        c c c
        @{}
    }
        \toprule
        &
        \multicolumn{1}{c}{Fixed rewards}
        &
        \multicolumn{2}{c}{Uncertain rewards $r \in \cR$}
        \\
        \cmidrule(lr){2-2}
        \cmidrule(lr){3-4}

        Reformulation
        &
        \multicolumn{1}{c}{$\cP$ s-rectangular}
        &
        \multicolumn{1}{c}{$\cP$ sa-rectangular}
        &
        \multicolumn{1}{c}{$\cP$ s-rectangular}
        \\
        \midrule

        Expected-reward reduction
        &
        $\checkmark$
        &
        $\checkmark$

        (if $\cR$ sa-rec.)
        &
        $\times$ (in general)
        \\

        \addlinespace[3pt]

        State augmentation
        &
        $\checkmark$
        &
        $\times$ (in general)
        &
        $\times$ (in general)
        \\

        \bottomrule
    \end{tabularx}
\end{table}

\section{Main takeaways}\label{sec: lec 3 main takeaways}
This chapter has highlighted the connections between robust MDPs and several models of sequential decision-making under uncertainty. We summarize the main reductions in \Cref{fig: lec 3 main takeaways} and below:
\begin{itemize}
    \item \Cref{sec: lec 3 uncertain rewards}: in the case where the transition probabilities are fixed and the rewards are unknown, optimal policies can be chosen stationary robust MDPs can be solved via mathematical programming, even when the uncertainty set (for the rewards) is non-rectangular. This approach is based on the fact that the return is linear in the state-action frequencies and in the instantaneous rewards.
    \item \Cref{sec: lec 3 uncertain rewards and transitions} : when both the rewards and the transitions are unknown, things are much more complicated. In particular the Bellman operators now involve a term $P(s,a)\tr r(s,a)$ which is bilinear in the variables $P$ and $r$. In general this leads to difficult (NP-hard) optimization problems even when the uncertainty sets for $P$ and $r$ are s-rectangular, except when the ``effective uncertainty set" (see \eqref{eq: lec 3 Uf}) is convex (and s-rectangular).
    \item \Cref{sec: lec 3 multi-objective MDPs}: for multi-objective MDPs, the infimum as a scalarization function yields robust MDPs, which can be interpreted as ``worst-case" multi-model MDPs (with linear scalarization function). We can recover Pareto optimal policies by solving robust MDPs, but in many cases it is NP-hard to do so (except when only the rewards are uncertain and the transitions are common across models).
    \item \Cref{sec: lec 3 time var P}: robust MDPs with time-varying transition probabilities reduce to robust MDPs with stationary probabilities (studied in \Cref{chap: lec 1 robust MDP} and \Cref{chap: lec 2 rectangular models}) by replacing the uncertainty set by its rectangular counterpart when the adversary makes its decision after observing the current state. Otherwise the decision problem of the adversary is close to a POMDP.
    \item \Cref{sec: lec 3 sgs}: for the discounted return and the mean-payoff, sa-rectangular robust MDPs can be formulated as perfect-information stochastic games (SGs). For the discounted return, s-rectangular robust MDPs can be formulated as concurrent stochastic games. The ``information" available to the decision maker (in robust MDPs) and the players (in stochastic games) differ: in SGs, histories include the past decisions chosen by all players, whereas in robust MDPs we cannot observe the values of the past transition probabilities. This difference is irrelevant when stationary optimal strategies exist in SGs, as is the case for the objectives described above.
    \item \Cref{sec: lec 3 drmdps}:  distributionally robust MDPs is a natural extension to robust MDPs that allow for distributional information over the set of feasible transition probabilities. Perhaps surprisingly, distributionally robust MDPs (where the transition probabilities are sampled independently across each period) can be reformulated as robust MDPs with stationary transition probabilities.
    \item \Cref{sec: lec 3 non-rec}: under an additional assumption that the instantaneous rewards are independent of the next state, several non-rectangular uncertainty sets (including r-rectangular and some extensions of this model) can be solved by dynamic programming over the original set of states and in fact can be replaced by their sa-rectangular extension. The model of Lightning-Does-not-Strike-Twice and its generalization of k-rectangular uncertainty sets require a state augmentation to recover dynamic programming equations.
    \item \Cref{sec: lec 3 risk averse}: risk-averse MDPs with coherent risk measures can be reformulated as robust MDPs with non-rectangular uncertainty; to recover rectangular uncertainty sets one can resort to nested coherent risk measures. These results are only derived for the policy evaluation problem (i.e., for a fixed policy of the decision maker) to simplify the exposition.
    \item (*) \Cref{sec: lec 3 discussion rewards}: the assumption that the instantaneous reward is independent of the next state can be made without loss of generality in some setups (for instance when the rewards are fixed and the uncertainty set is s-rectangular) and not without loss of generality in some other setups (for instance when the rewards are also unknown).
\end{itemize}
\begin{figure}[htb!]
    \centering
    \includegraphics[width=1\linewidth]{figures/lecture_3_main_takeaways_v5.pdf}
    \caption{Main reductions proved in this chapter. An arrow indicates that an optimization problem reduces to another one; next to an arrow is either an assumption under which this reduction works (e.g., ``$\cP$ cvx" for $\cP$ a convex set) or a label of the kind ``$\cP \gets \recsa(\cP)$", which means that this reduction always works by replacing $\cP$ by $\recsa(\cP)$ in the decision problem at the end of the arrow (and similarly for ``$\cP \gets \recs(\cP)$" and ``$\cP \gets \bar{\cP}$"). The s-rectangular and sa-rectangular extensions $\recs(\cP)$ and $\recsa(\cP)$ have been defined in \Cref{sec: lec 3 time var P}. The set $\bar{\cP}$ is the set of expected kernels induced by a set of measures $\cV \subset \Delta(\DSSA)$ as defined in \Cref{sec: lec 3 drmdps}. The reduction from risk-averse MDPs to robust MDPs for nested coherent risk measures has only been sketched in \Cref{sec: lec 3 risk averse} and is added here for the sake of completeness.}
    \label{fig: lec 3 main takeaways}
\end{figure}
\section{References}\label{sec: lec 3 ref}

For \Cref{sec: lec 3 uncertain rewards} I follow the model description in \cite{mcmahan2003planning} (robust MDPs with uncertain rewards), \Cref{sec: lec 3 uncertain rewards and transitions} extends this model to joint uncertainty in rewards and transitions. The optimality of stationary policies for the case of uncertain rewards and fixed transitions is in \cite{gadot2024solving}, and interestingly uncertainty in rewards is related to state entropy regularization~\cite{ashlag2026state}. It is worth mentioning that some seminal robust MDP papers allow (such as \cite{xu2010distributionally}) both $r$ and $P$ to be unknown but do assume that the rewards do not depend on the next state.
To the best of my knowledge, the hardness result for the case of convex rectangular robust MDPs with uncertain rewards and transitions (\Cref{prop: lec 3 hardness uncertain r p}) is new; it is worth mentioning a closely related result (Proposition 10 in \cite{kitamura2026revisiting}) for the case of finite uncertainty sets.

For \Cref{sec: lec 3 multi-objective MDPs}, I refer to \cite{roijers2013survey} for a survey, to \cite{randour2015percentile} for a very complete paper studying several different objectives (mean-payoff, limsup and liminf of the payoff, discounted returns, etc), and to \cite{ropke2024divide,lu2023multi,delgrange2020simple} for some more recent works. The vast majority of the papers on multi-objective MDPs (including the ones cited in the previous sentence) focus on the case where the transitions are common across all models, with the notable exception of \cite{scheftelowitsch2017multi,hahn2017multi} who focus on interval MDPs. The relation between the convexity properties of the achievable set $R(\Pi)$ and the possibility to recover Pareto optimal policies from linear scalarization is well-known, see e.g. \cite{feinberg1995constrained,hernandez2004scalarization}. The case of {\em multi-model MDPs}, who study the scalarized problem \eqref{eq: lec 3 linear scalarized momdp} has been studied under different names over the years: \cite{derman2018soft,lobo2020soft} call it soft-robust MDPs, \cite{buchholz2020concurrent,buchholz2019computation} call it concurrent MDPs, \cite{raskin2014multiple} call it multiple-environment MDPs. There exist many papers showing that the scalarized problem (with a linear function) is difficult, see e.g. \cite{le2007robust}, \cite{buchholz2019computation}, as well as \cite{li2026computational} (for the case of the scalarized problem with the infimum function).

For \Cref{sec: lec 3 time var P} (time-varying transitions) I follow Section 2 in \cite{wang2023foundation} for the description of the sets of time-varying uncertainty, which provides a rigorous analysis of the several possible corner cases that occur when inspecting various restrictions on the sets of policies (e.g., the decision maker can only use deterministic policies potentially history-dependent) and uncertainty sets (s-rectangular, sa-rectangular, non-rectangular, nonconvex ...). \cite{li2025rectangularity} also focuses on the relation between time-varying transition probabilities, rectangularity, and so-called ``game formulations" of robust sequential decision problems.
The connection between the policy evaluation problem for time-varying non-rectangular uncertainty sets and POMDPs is from conversations with Marek Petrik and with St{\'e}phane Gaubert and Marianne Akian; I believe that it is known more widely but I am not aware of a more formal argument. The two papers~\cite{littman1994memoryless,vlassis2012computational} show hardness results for stationary decision makers in POMDPs. 

For \Cref{sec: lec 3 sgs} (connection with stochastic games) I refer to \cite{neyman2003stochastic,mertens2015repeated} for classical textbooks on SGs and to \cite{laraki2015advances,renault2019tutorial} for shorter recent surveys. The connection between stochastic games and robust MDPs can be found in several places, including \cite{satia1973markov}, section~5 of \cite{iyengar2005robust}, the last paragraph of page 7 of \cite{xu2010distributionally}, the last paragraph of section 2 in \cite{wiesemann2013robust}, and the introduction of \cite{grand2025convex}. \Cref{prop: lec 3 mean payoff sa-rec} is in \cite{grand2026beyond}. The Mertens-Neyman theorem is proved in~\cite{mertens1981stochastic}.

For \Cref{sec: lec 3 drmdps} (distributionally robust MDPs) I adapted~\cite{xu2010distributionally} to the case of uncertain transitions (\cite{xu2010distributionally} focuses on uncertain rewards and uncertain transitions, but assumes that rewards are independent of the next state, which is not without loss of generality, see \Cref{sec: lec 3 uncertain rewards and transitions}).
I refer to~\cite{le2007robust,buchholz2020concurrent} for hardness results for multi-model MDPs with discounted return. 

For \Cref{sec: lec 3 non-rec} (non-rectangular uncertainty sets), I follow \cite{grand2024tractable} for the analysis of r-rectangular models. The model of r-rectangularity was introduced in \cite{goh2018data} for policy evaluation, and the analysis of optimal policies was done in \cite{goyal2023robust}; the equivalence with sa-rectangularity over an extended state set was suggested to the author by Erick Delage in 2019. For the connections between non-rectangular models and their rectangular extensions when rewards are independent of the next state I follow \cite{grand2024tractable} (see also \cite{li2025rectangularity} for the finite-horizon case and connection with multistage stochastic optimization). The Lightning-Does-not-Strike-Twice model is from~\cite{mannor2012lightning} and the k-rectangular model of uncertainty is by the same authors in~\cite{mannor2016robust}.

For \Cref{sec: lec 3 risk averse} (connection with risk-averse MDPs) I follow \cite{chow2015risk} for the formulation of \cvar{} risk-averse MDPs as a robust MDP with non-rectangular uncertainty sets; recent such formulations can be found for instance in~\cite{shapiro2025distributionally,shapiro2025risk}. Conditional risk measures are introduced in~\cite{ruszczynski2006conditional}. The reader interested in Markovian risk measures for MDPs can refer to \cite{ruszczynski2010risk}; the main results about time-consistent, law-invariant risk measures are in \cite{kupper2009representation,shapiro2012time,marthe2023beyond}. An excellent textbook on risk measures is \cite{follmer2016stochastic}. 
\appendix
\part*{Appendices}
\chapter{Background in topology and convex analysis}\label{app:sec analysis}
\paragraph{Compactness and continuity.} The following are used in \Cref{chap: lec 1 robust MDP}.  We want to study the compactness of $\PiH$. Note that $\PiH$, the set of all history-dependent policies, is contained in an infinite-dimensional space since $\PiH = \times_{h \in \cH} \Delta(\cA)$ with $\cH = \cup_{t \in \N} (\cS \times \cA)^t \times \cS$. Therefore we need to be careful when describing its ``compactness" (a notion that in all generality depends on the topology chosen over $\PiH$). We use the product topology: convergence in this topology means convergence of the action distributions at each fixed finite history. The product topology has the following important property.
\begin{theorem}[Tykhonoff's theorem]\label{app:th tykhonoff}
    The product of any collection of compact topological spaces is compact for the product topology.
\end{theorem}
We recall the definition of upper semicontinuity here.
\begin{definition}[Upper semicontinuity]
 Let $f:\cX \to \R \cup \{+\infty,-\infty\}$. The function $f$ is {\em upper semicontinuous} on $\cX$ if for every $y \in \R$ the set $\{x \in \cX \; | \; f(x) < y\}$ is open in $\cX$.
\end{definition}
In the definition above $\cX$ is a topological space, so the upper semicontinuity property of $f$ depends on the topology chosen on $\cX$. The product topology makes $\PiH$ compact (from Tykhonoff's theorem above), and then we can use Weierstrass' theorem in the following general version.
\begin{theorem}[Weierstrass theorem]
Let $\cX$ be a nonempty compact topological space and $f:\cX \to \R \cup \{+\infty,-\infty\}$ be upper semicontinuous. Then $f$ attains its supremum over $\cX$.
\end{theorem}
\paragraph{Minimax theorems.}
In \Cref{chap: lec 2 rectangular models} and \Cref{chap: lec 3 connection} we use strong duality in minimax problems a few times. The first such result was obtained by von Neumann in 1928 for continuous quasi-convex quasi-concave functions over simplexes (although people generally refer to this result only for bilinear objectives); for our purpose we use the following version of {\em Sion's theorem}, which allows for more general domains. I only give the version we need (i.e., the domains belong to Euclidean spaces and the objective is continuous in each variable separately), but Sion's theorem applies more broadly.
\begin{theorem}[Sion's theorem (for our purpose)]
    Let $\cX,\cY$ be nonempty convex compact subsets of Euclidean spaces and $f: \cX \times \cY \to \R$. Assume that:
    \begin{itemize}
        \item For any $y \in \cY$, $x\mapsto f(x,y)$ is continuous and quasi-concave on $\cX$,
        \item For any $x \in \cX$, $y\mapsto f(x,y)$ is continuous and quasi-convex on $\cY$.
    \end{itemize}
    Then we have that $\sup_{x \in \cX} \inf_{y \in \cY} f(x,y) = \inf_{y \in \cY}  \sup_{x \in \cX} f(x,y)$ and all the suprema and infima are attained.
\end{theorem}
\paragraph{Convex geometry and analysis.}
In \Cref{chap: lec 2 rectangular models} and \Cref{chap: lec 3 connection} we use the following results from convex analysis and convex geometry. We start with Banach's fixed point theorem.
\begin{theorem}[Banach's fixed point theorem]
    Let $\cX$ be a complete metric space. A contraction $T:\cX \to \cX$ has a unique fixed point $v_T \in \cX$, and for any $v_0 \in \cX$, $\displaystyle \lim_{t \to + \infty} T^t v_0$ exists and $\displaystyle \lim_{t \to + \infty} T^t v_0 = v_T$.
\end{theorem}
\begin{theorem}[Supporting hyperplane theorem for concave functions]
    Let $\cX $ be a nonempty convex compact subset of $\R^{n}$ and let $f:\cX \to \R$ be a continuous concave function. For every $x_0$ in the relative interior of $\cX$ there exists $g \in \R^{n}$ such that $f(x) \leq f(x_0) + g \tr (x-x_0), \forall \; x \in \cX$.
\end{theorem}
The vector $g$ in the above theorem is a supergradient of $f$ at $x_0$; we refer to \cite{rockafellar2025convex}, Theorem 23.4 (which is stated in the more familiar case of subgradients of proper continuous convex functions).

The next theorem is known as the {\em separating hyperplane} theorem; we give it here its statement for Euclidean space. We refer to Section 11 in ~\cite{rockafellar2025convex} for more general statements and to Section 2.5 in \cite{vandenberghe2004convex} for the statement below. 
\begin{theorem}\label{th: app separating hyperlane}
    Let $\cX,\cY \subset \R^n$ be two nonempty disjoint convex sets. Then there exists $w \in \R^n, w \neq 0$ and $\alpha \in \R$ such that $w \tr x \leq \alpha \leq w \tr y, \forall \; (x,y) \in \cX \times \cY$.
\end{theorem}
The next result provides the first-order optimality conditions for orthogonal projections (here it's useful to draw a picture, for instance of the orthogonal projection onto the unit disk in $\R^2$), see pages 270-271 in \cite{rockafellar2025convex}.
\begin{theorem}
    Let $\cX $ be a nonempty convex closed subset of $\R^{n}$ and let $x \in \R^{n}$. The orthogonal projection $\bar{x}$ of $x$ onto $\cX$ is the unique minimizer of $y \mapsto \|x-y\|_{2}$ over $y \in \cX$ and is characterized by $\bar{x} \in \cX$ and
    \[\langle \bar{x}-x,y-\bar{x} \rangle \geq 0, \forall \; y \in \cX.\]
\end{theorem}
\chapter{Algorithms for fixed point computations}\label{app: algorithms}
In this appendix we present two general versions of value iteration and policy iteration and prove their convergence guarantees. Doing so enables us to analyze simultaneously the convergence rates of value iteration and policy iteration in the three distinct setups of interest in these notes: (a) nominal MDPs (\Cref{chap: lec 0 mdp}), (b) adversarial MDPs under s-rectangularity (\Cref{sec: lec 2 adv}), (c) robust MDPs  under s-rectangularity (\Cref{sec: lec 2 rmdp}).

We will use the following notation: we study families of monotone contractions $\{T^{\theta} \; | \; \theta \in \Theta\}$ indexed by a compact set $\Theta$, with $T^{\theta}:\R^{\cS} \to \R^{\cS}$ (for $\cS$ a finite set). The rate of contraction for $\| \cdot \|_{\infty}$ is $\gamma \in [0,1)$ for each contraction $T^{\theta}$, and we write $v^{\theta} \in \R^{\cS}$ for its fixed point. We will consider the ``componentwise supremum" of the family, defined as the operator $\bT:\R^{\cS} \to \R^{\cS}$ with
\[\bT v(s) = \sup_{\theta \in \Theta} T^{\theta} v(s), \forall \; v \in \R^{\cS}, s \in \cS.\]
We assume that this supremum is attained on each component here (this is true for our applications for MDPs and robust MDPs).
We know that $\bT$ is a monotone contraction (the arguments follow similarly as in the proof of \Cref{prop: lec 0 operators contraction}), and we define $\bv \in \R^{\cS}$ as its unique fixed point. Finally, the attainability property is defined as follows.
\begin{definition}
    The family $\{T^{\theta} \; | \; \theta \in \Theta\}$ satisfies the attainability property with respect to the contraction $\bT$ if for any $v \in \R^{\cS}$, there exists $\btheta \in \Theta$ such that $\bT v = T^{\btheta}v$.
\end{definition}
 The value iteration and policy iteration algorithms are two methods to compute $\bv$ for families of monotone contractions satisfying the attainability property. 

\paragraph{Value iteration.} Here we rederive the convergence rate of Picard's iteration.
\begin{algorithm}
\caption{Value iteration (general algorithm)}\label{alg: app VI (general)}
\begin{algorithmic}[1]
\State {\bf Input:} A compact set $\Theta$, a family of monotone contractions $\{T^{\theta} \; | \; \theta \in \Theta\}$ with constant $\gamma \in [0,1)$, satisfying the attainability property with respect to $\bT$
\State Initialize $v^{0} \in \R^{\cS}$.
\For{$t \in \N$}
    \State $v^{t+1} = \bT v^{t}$
\EndFor
\end{algorithmic}
\end{algorithm}

\begin{proposition}\label{prop: app VI rate}
    Let $(v^t)_{t \in \N}$ be the sequence of vectors generated by \Cref{alg: app VI (general)}.
    \begin{enumerate}
        \item The sequence $(v^t)_{t \in \N}$ converges to $\bv$ at a geometric rate of $\gamma$:
        \[\|\bv - v^{t}\|_{\infty} \leq \gamma^t \|\bv - v^{0}\|_{\infty}\]
    \item Let $\epsilon>0$. Assume that $\|\bv - v^{t}\|_{\infty} \leq \epsilon$ for some $t \in \N$ and let $\theta \in \Theta$ such that $T^{\theta}v^t = \bT v^t$. Then the fixed point of $T^{\theta}$ satisfies $\|\bv-v^{\theta}\|_{\infty} \leq 2\epsilon \gamma (1-\gamma)^{-1}$.

    Additionally, we can check that $\|\bv - v^{t}\|_{\infty} \leq \epsilon$ by checking $\|\bar{T}v^{t}-v^{t}\|_{\infty} \leq \epsilon (1-\gamma)$, and that last inequality holds after at most $\lfloor \frac{1}{1-\gamma} \log\left(\frac{\|v^1-v^0\|_{\infty}}{\epsilon (1-\gamma)}\right) \rfloor + 1$ iterations of \Cref{alg: app VI (general)}.
    \end{enumerate}
\end{proposition}
\begin{proof}
\begin{enumerate}
    \item We have
        \[\|\bv - v^{t+1}\|_{\infty} = \|\bT \bv - \bT v^t\|_{\infty} \leq \gamma \|\bv - v^{t}\|_{\infty}\]
        where the inequality uses the contraction property of $\bT $. We get the desired inequality by induction.
        \item For $\theta$ constructed as in the statement of the proposition, we have
\begin{align*}
    \|v^{\theta} - \bv\|_{\infty} = \|v^{\theta} - T^{\theta}v^t + T^{\theta}v^t - \bv\|_{\infty} & \leq \|v^{\theta} - T^{\theta}v^t \|_{\infty}  + \|T^{\theta}v^t - \bv\|_{\infty} \\
    & = \|T^{\theta}v^{\theta} - T^{\theta}v^t \|_{\infty}  + \|\bT v^t - \bT \bv\|_{\infty} \\
    & \leq \gamma \|v^{\theta} - v^t \|_{\infty}  + \gamma \|v^t - \bv\|_{\infty} \\
    & \leq \gamma \|v^{\theta} - \bv \|_{\infty}  + 2 \gamma \|v^t - \bv\|_{\infty}
\end{align*}
and rearranging the terms yields the first desired inequality.

Now note that
\[\|\bv - v^{t}\|_{\infty} = \|\bv - \bar{T}v^t +\bar{T}v^t- v^{t}\|_{\infty} \leq \gamma \|\bv - v^{t}\|_{\infty} + \|\bar{T}v^t- v^{t}\|_{\infty} \]
so that rearranging we obtained 
\[\|\bv - v^{t}\|_{\infty} \leq \frac{\|\bar{T}v^t- v^{t}\|_{\infty}}{1-\gamma}.\]
Now $\|\bar{T}v^t- v^{t}\|_{\infty} = \|v^{t+1} - v^{t} \|_{\infty} \leq \gamma  \|v^{t} - v^{t-1} \|_{\infty}$ hence $\|v^{t+1} - v^{t} \|_{\infty} \leq \gamma^t  \|v^{1} - v^{0} \|_{\infty}$.
 To ensure $\|\bv - v^{t}\|_{\infty} \leq \epsilon$ it is sufficient that $\gamma^{\bar{t}}  \|v^{1} - v^{0} \|_{\infty} \leq \epsilon (1-\gamma)$, for which we can take $\bar{t} = \lfloor \frac{1}{1-\gamma} \log\left(\frac{\|v^1-v^0\|_{\infty}}{\epsilon (1-\gamma)}\right) \rfloor + 1$.
\end{enumerate}
\end{proof}

\paragraph{Policy iteration.}
\begin{algorithm}
\caption{Policy iteration (general algorithm)}\label{alg: app PI (general)}
\begin{algorithmic}[1]
\State {\bf Input:} A compact set $\Theta$, a family of monotone contractions $\{T^{\theta} \; | \; \theta \in \Theta\}$ with constant $\gamma \in [0,1)$, satisfying the attainability property with respect to $\bT$
\State Initialize $\theta^0 \in \Theta$, $v^{0}$ the fixed point of $T^{\theta^0}$
\For{$t \in \N$}
    \State Choose $\theta^{t+1} \in \Theta$ such that $\max_{\theta \in \Theta} T^{\theta}v^{t}(s) = T^{\theta^{t+1}}v^{t}(s)$ for all $s \in \cS$ (and equal to $\theta^t$ if we can)
    \State Choose $v^{t+1} \in \R^{\cS}$ as the fixed point of $T^{\theta^{t+1}}$
\EndFor
\end{algorithmic}
\end{algorithm}
The convergence of \Cref{alg: app PI (general)} is a little longer to prove than for \Cref{alg: app VI (general)} but uses the same ideas; again we only need the monotonicity, contractivity, and attainability properties.
\begin{proposition}\label{prop:app PI rate}
 Let $(\theta^t)_{t \in \N}$ be the sequence generated by \Cref{alg: app PI (general)}. Then 
         $\|\bv-v^{\theta^t}\|_{\infty} \leq \gamma^t \|\bv-v^{\theta^0}\|_{\infty}.$
\end{proposition}
\begin{proof}
    Let $t \in \N$. Note that $v^{\theta^{t}} \leq v^{\theta^{t+1}}$. Indeed
\[T^{\theta^{t+1}}v^{\theta^{t}} = \bT v^{\theta^{t}} \geq T^{\theta^{t}}v^{\theta^{t}} = v^{\theta^{t}}.\]
Using then the monotonicity and contraction properties of the operator $T^{\theta^{t+1}}$, we can apply $T^{\theta^{t+1}}$ on both sides of the inequality 
$T^{\theta^{t+1}}v^{\theta^{t}} \geq v^{\theta^{t}}$. This yields $v^{\theta^{t}} \leq v^{\theta^{t+1}}$.

We now prove the desired convergence rate. We first note that 
\begin{align*}
    \bv - v^{\theta^{t+1}} = \bv  - T^{\theta^{t+1}}v^{\theta^{t+1}}= \bv  - T^{\theta^{t+1}}v^{\theta^{t}} + T^{\theta^{t+1}}v^{\theta^{t}} -  T^{\theta^{t+1}}v^{\theta^{t+1}}.
\end{align*}
Hence we conclude that $0 \leq \bv  - v^{\theta^{t+1}} \leq \bv  - T^{\theta^{t+1}}v^{\theta^{t}}$ since $T^{\theta^{t+1}}v^{\theta^{t}}\leq  T^{\theta^{t+1}}v^{\theta^{t+1}}$ (from monotonicity of $T^{\theta^{t+1}}$ and $v^{\theta^{t}} \leq v^{\theta^{t+1}}$).
Now we get
\[ \|\bv  - T^{\theta^{t+1}}v^{\theta^{t}}\|_{\infty}  = \|\bT \bv  - \bT v^{\theta^{t}}\|_{\infty}\leq \gamma \| \bv  - v^{\theta^{t}}\|_{\infty}\]
By induction $\|\bv  - T^{\theta^{t+1}}v^{\theta^{t}}\|_{\infty} \leq \gamma^{t+1} \| \bv  - v^{\theta^{0}}\|_{\infty}$ and since $0 \leq \bv  - v^{\theta^{t+1}} \leq \bv  - T^{\theta^{t+1}}v^{\theta^{t}}$ we get the desired inequality on $\|\bv  - v^{\theta^{t}}\|_{\infty}$.
\end{proof}

\newpage

\bibliographystyle{alpha}
\bibliography{ref}

\end{document}